\documentclass{memo-l} 
\usepackage{epsf,amscd,amssymb} 
 
\newtheorem{thm}{Theorem}[chapter] 
\newtheorem{lem}{Lemma}[chapter] 
\newtheorem{cor}{Corollary}[chapter] 
\newtheorem{prop}{Proposition}[chapter] 
\newtheorem{claim}{Claim}[chapter] 
\newtheorem{decthm}{Theorem A} 
 
\newtheorem{crscor}{Corollary A} 
 
\newtheorem{topcor}{Corollary B}

\theoremstyle{definition} 
 
\newtheorem{defn}{Definition}[chapter] 
 
\theoremstyle{remark}  
 
\newtheorem{rem}{Remark}[chapter] 
\newtheorem{hyp}{Hypothesis} 
\newtheorem{exmp}{Example}[chapter] 

\numberwithin{figure}{chapter}
\numberwithin{equation}{chapter}
 
\newcommand\bR{\mathbf{R}}

\newcommand\bZ{\mathbf{Z}}
\newcommand\bV{\mathbf{V}}

\newcommand\rp{\bR P}
\newcommand\GL{\mathrm{GL}}
\newcommand\SL{\mathrm{SL}}
\newcommand\PGL{\mathrm{PGL}}
\newcommand\PSO{\mathrm{PSO}}

\newcommand\PSL{\mathrm{PSL}}
\newcommand\SLt{\mathrm{SL}}

\newcommand\dev{\mathbf{dev}}
\newcommand\SI{\mathbf{S}}

\newcommand\inte{\mathrm{int}}
\newcommand\Bd{\mathrm{bd}}
\newcommand\clo{\mathrm{Cl}}
\newcommand\bdd{\mathbf{d}}
\newcommand\ideal[1]{\tilde #1_\infty}
\newcommand\ra{\rightarrow}

\newcommand\che{\check}
\newcommand\emp{\emptyset}
\newcommand\eps{\epsilon}

\newcommand\vth{\vartheta}
\newcommand\vpi{\varphi}
\newcommand\Aff{\mathrm{Aff}}
\newcommand\ovl{\overline}
\newcommand\Aut{\mathrm{Aut}}
\newcommand\Idd{\mathrm{I}}

\newcommand\awlg{may assume without loss of generality }
\newcommand\cO{{\rm ,}\ }
\newcommand\hll{\hfill\break}

\newcommand{\tri}{\triangle} 
\newcommand{\brno}{\bR^n - \{O\}} 
\newcommand{\brto}{\bR^3 - \{O\}} 
\newcommand{\SII}{{{\SI}^{n-1}_\infty}} 
\newcommand{\SIT}{{{\SI}^2_\infty}} 
 
\newcommand{\hideal}[1]{{#1}_{h \infty}} 
\newcommand{\Iideal}[1]{\ideal{#1}^{\mathsf{i}}} 
\newcommand{\hIideal}[1]{{#1}^{\mathsf{i}}_{h \infty}} 
\newcommand{\fideal}[1]{\ideal{#1}^{\mathsf{f}}} 
\newcommand{\hfideal}[1]{{#1}^{\mathsf{f}}_{h \infty}} 
\newcommand{\help}{\psi} 
\newcommand{\ini}{\xi} 
\newcommand{\fin}{\eta} 
\newcommand{\upp}{u} 
\newcommand{\dow}{d} 
\newcommand{\deck}{\vth} 
\newcommand{\Lamp}{\Lambda^{\mathsf{p}}} 
\newcommand{\Lambc}{\Lambda^{\mathsf{c}}} 
\newcommand{\alin}{\alpha^{\mathsf{i}}} 
\newcommand{\alfi}{\alpha^{\mathsf{f}}} 
\newcommand{\alfiin}{\alpha^{\mathsf{fi}}}  
\newcommand{\dfLambc}{\delta_\infty^{\mathsf{f}} \Lambc} 
\newcommand{\dinLambc}{\delta_\infty^{\mathsf{i}} \Lambc} 
\newcommand{\dfiLambc}{\delta_\infty^{\mathsf{fi}} \Lambc} 
\newcommand{\dfLamp}{\delta_\infty^{\mathsf{f}} \Lamp} 
\newcommand{\dinLamp}{\delta_\infty^{\mathsf{i}} \Lamp} 
\newcommand{\dfiLamp}{\delta_\infty^{\mathsf{fi}} \Lamp} 
\newcommand{\lora}{\longrightarrow}

\begin{document}

\frontmatter

\title[Radiant affine manifolds]  
{The decomposition  and classification of  radiant affine 
$3$-manifolds\footnote{\Large With Appendix C: 
``Radiant affine $3$-manifolds with boundary, and certain   
radiant affine $3$-manifolds,'' \hfill\break
by Thierry Barbot (UMPA, \'Ecole Normale Sup\'erieure de Lyon, 
46, all\'ee d'Italie, LYON, France)
and Suhyoung Choi. 
}
}     
\author{Suhyoung  Choi}    
\address{Department  of Mathematics \\   
College   of   Natural  Sciences \\     
Seoul  National University \\   
151--742 Seoul,   Korea}  
\email{shchoi@math.snu.ac.kr} 
\date{\today}     
\thanks{Version 5}     
\thanks{Research partially supported by GARC-KOSEF, 
the Ministry of Education 1997-001-D00036, 
and the BK21 program of the Ministry of Education}   
\thanks{A part of this result was announced at the conference   
``Combinatorial problems arising in   knots  and 
$3$-manifolds'' January 21-24, 1997 MSRI, Berkeley, CA}  
 
\begin{abstract} 
An affine manifold is a manifold with torsion-free flat affine connection. 
A geometric topologist's definition  of an affine manifold   
is a  manifold with an atlas of  charts to the affine space   
with affine transition functions;  
a radiant affine manifold  is an affine manifold with a
holonomy group consisting of affine  transformations  fixing a common fixed 
point. We decompose  a closed  radiant affine $3$-manifold 
into radiant $2$-convex  affine manifolds  and radiant concave  affine 
$3$-manifolds along mutually  disjoint totally geodesic  tori or Klein 
bottles  using   the   convex and concave decomposition   of   
real    projective $n$-manifolds developed   earlier.    
Then we  decompose  a  $2$-convex 
radiant   affine manifold into convex  radiant affine manifolds and 
concave-cone affine manifolds.  To do this, we will obtain 
certain  nice geometric objects  in  the  Kuiper completion  of 
a holonomy cover.  The equivariance and local finiteness property of the 
collection of such  objects will show that  their union covers a  
compact submanifold of codimension zero, the complement of  
which is convex. Finally, using the results of Barbot, we will show  
that a closed radiant affine $3$-manifold admits a total cross-section,  
confirming a conjecture of Carri\`ere, and hence every closed 
radiant affine $3$-manifold is homeomorphic to 
a Seifert fibered space with trivial Euler  
number, or a virtual bundle over a circle with fiber homeomorphic to  
a Euler characteristic zero surface.  
In Appendix C, Thierry Barbot and the author show 
the nonexistence of certain radiant affine $3$-manifolds and that 
compact radiant affine $3$-manifolds with nonempty totally geodesic 
boundary admit total cross-sections, which are key results 
for the main part of the paper. 
\end{abstract}   
 
\maketitle  
 
\tableofcontents 
 
\mainmatter 

\setcounter{chapter}{-1} 
\chapter{Introduction} 
The  classical $(X, G)$-geometry is  the study of invariant properties 
on a space $X$ under  the  action of a Lie group $G$ as Felix  Klein  
proposed.  An  $(X, G)$-structure on a manifold   
$M$ prescribes a local  identification of 
$M$ with $X$ with transition  functions in $G$.
An $(X,G)$-structure is also given by  an immersion, 
so-called developing map, from the universal cover of  
the manifold $M$ to $X$ equivariant with respect to a homomorphism,  
so-called holonomy homomorphism, from the deck-transformation group  
$\pi_1(M)$ to $G$. Differential geometers can define $(X, G)$-structure  
using higher order curvatures. 
   
Two of the fundamental general questions  
of geometric  topologists are to find 
a topological  obstruction  to the   existence of structures  locally 
modeled  on $(X,  G)$ and  classify  such  structures on a  given 
manifold up to $(X, G)$-self-diffeomorphisms. 

There are many important $(X, G)$-structures. In particular, 
when $X$ has a $G$-invariant complete metric, then 
the study is equivalent to the study of discrete cocompact subgroups 
of Lie groups which is well-established and gave us many 
fruitful topological implications. 
 
We will be working in real projective and affine structures on 
manifolds, a very much open area presently.
Recall that real projective geometry is given by a pair  
$(\bR P^n, \PGL(n+1, \bR))$ where $\PGL(n+1, \bR)$ is  
the projectivized general linear group acting on $\bR P^n$.  
An affine geometry is given by a pair $(\bR^n, \Aff(\bR^n))$  
where $\Aff(\bR^n)$ is the group of affine transformations of 
$\bR^n$; i.e., transformations of form $x \mapsto Ax + b$ for 
$A$ an element of the general linear group $\GL(n, \bR)$ and $b$  
an $n$-vector.  
A {\em real projective $n$-manifold}\/ is an $n$-manifold  with  
a geometric structure modeled on real projective geometry;   
an {\em affine $n$-manifold}\/ one  with that modeled  
on affine geometry.   (In differential geometry, an affine 
manifold is defined as a manifold with a flat torsion-free  affine 
connection. A real projective $n$-manifold is equivalent to 
a manifold with a projectively flat torsion-free affine connection.)     

An important distinguishing feature of these geometric manifolds is 
that geodesics are defined by connections. Basically, they correspond 
to usual straight lines under charts.

A  sphere  with real projective  structure  induced from  the standard
double covering map is said to be  a {\em real  projective sphere}.  A
compact real projective disk  with geodesic boundary real projectively
diffeomorphic to a closed hemisphere  in $\SI^2$ is  said to be a {\em
real projective hemisphere}. These are trivially unique structures.
  
Affine and real projective structures on annuli or tori  
were classified by Nagano-Yagi \cite{NY} and Goldman \cite{Gsthesis}.  
Benz\'ecri \cite{Ben} showed that any affine surface must 
be homeomorphic to tori or annuli.
 
A real projective manifold is {\em convex} if its universal cover 
is real projectively homeomorphic to a convex domain in  
an affine patch of $\bR P^n$. They form the most important class 
of real projective manifolds. 
As hyperbolic geometry is a geometry included in projective  
geometry by the Klein model, we see that every hyperbolic  
$n$-manifold carry a canonical real projective structure. 
Also, there are nontrivial deformations of these real projective 
manifolds to convex real projective manifolds as shown by  
Kac-Vinberg \cite{KV}, Koszul \cite{Kz} and Goldman \cite{Gconv}. 
 
Many other examples of affine    manifolds     were   constructed   by 
Sullivan-Thurston \cite{ST}, where they construct many examples of 
complex projective or M\"obius structures on surfaces and manifolds. 

The {\em holonomy group} is the image of the holonomy homomorphism. 
The purpose of this paper is to show that   
the above two questions can be answered for  radiant  affine  
structures  on $3$-manifolds. A so-called radiant affine manifold  is   
an affine manifold whose holonomy group fix a common point of   
$\bR^n$; i.e., an atlas can be chosen so that  
transition functions are in the subgroup $\GL(n, \bR)$ of 
$\Aff(\bR^n)$. One can think of a radiant affine  
manifold as a manifold with $(\bR^n, \GL(n, \bR))$-structure,  
where $\GL(n, \bR)$ acts on $\bR^n$ in the standard manner. 

We recall that  
examples of radiant affine $(n+1)$-manifold can be obtained  
by so-called generalized affine suspensions over real projective  
$n$-manifolds: We realize a real projective manifold 
times a real line as an affine manifold of dimension $n+1$  
by adding the radial direction to in a natural manner.  
Then we take a quotient by an infinite cyclic group acting  
affinely and mostly in the radial direction to make a compact  
affine quotient manifold. A generalized affine suspension is said  
to be a {\em Benz\'ecri suspension} if the group acts  
purely in the radial direction up to finite order 
(see Chapter \ref{ch:racom} and Appendix \ref{app:radd} 
for more details).  

A submanifold of a real projective or affine manifold is 
{\em totally geodesic}\/ if each point of it has 
a neighborhood with a chart mapping the manifold into 
a subspace of the real projective space or 
the affine space respectively.
Given a real projective  or affine manifold  $M$, let $N$ be a totally 
geodesic $(n-1)$-dimensional manifold properly imbedded in 
the interior $M^o$ of $M$. Then 
we  define  the splitting  of $M$  along $N$ as  taking the disjoint 
unions of the completions of all components of $M - N$ 
(see Chapter 1 or \cite{psconv} for  more details).   
We say that a 
real projective or  affine manifold $M$ {\em  decomposes}\/ along $N$ to 
another real  projective or affine manifold  $M'$  if $M'$ is obtained 
from $M$ by splitting along $N$.  
In this paper, we will prove the following result:  
\begin{decthm}
Let  $M$ be a compact radiant affine $3$-manifold with empty or totally 
geodesic  boundary.  Then $M$  decomposes along the  union of finitely 
many  disjoint totally geodesic tori or  Klein bottles,  
tangent to the radial flow, into  
\begin{enumerate} 
\item convex radiant affine $3$-manifolds  
\item generalized affine  suspensions   of  real projective  spheres,  real 
projective planes, real projective hemispheres, or $\pi$-annuli {\rm 
(}or M\"obius  bands\/{\rm )}  of type  C\/{\rm ;}  or affine tori,  
affine Klein bottles, or affine annuli  
{\rm (} or  M\"obius   bands\/{\rm  )} with geodesic boundary.  
\end{enumerate} 
\end{decthm}

Note this  decomposition is not canonical; however, the statements can be 
made into ones about canonical decomposition with additional requirements.  
A {\em $\pi$-annulus}\/ of  type  C is the third  type  of a $\pi$-annulus in 
Section 3 of \cite{cdcr2}, i.e., a $\pi$-annulus  that is  the sum of 
two elementary annuli of type IIb and with nondiagonalizable  
holonomy for the generator of the fundamental group. 
A {\em $\pi$-M\"obius band of type C}\/  is a  real  projective  
M\"obius  band  that is doubly covered by a $\pi$-annulus of type C.  

The radiant vector field on $\bR^n$ as given by 
$\sum_{i=1}^n x_i \partial/\partial x_i$ is $\GL(n, \bR)$-invariant, 
and hence induces a vector field and a flow on $M$ by $\dev$ and 
the covering map. The vector field is said to be 
a {\em radiant vector field}, and the flow the {\em radial flow}.
A {\em total cross-section} is a closed transverse submanifold that meets 
every flow line. 
Carri\`ere \cite{Carlet} asked  whether  every compact radiant  affine 
$3$-manifold admits a total cross-section to 
the radial  flow. In the paper, 
Carri\`ere showed that if the radial flow preserves a volume form,  
then this is true. In  dimension 6, Fried produced  
a  counter-example \cite{fried}.   

Barbot has shown that the conjecture is true if the holonomy group  
is virtually solvable (and more generally, if some finite-index subgroup  
preserves a plane of ${\bR}^{3}$), or if the manifold is  
homeomorphic to a Seifert $3$-manifold 
(see \cite{barbot2} and \cite{resolu}). 
Barbot confirms the conjecture if the $3$-manifold has  
a totally geodesic surface tangent to the radial flow, 
if there exists a closed orbit of non-saddle type, or if  
the $3$-manifold is convex (see Theorem \ref{thm:barbot}  
and \cite{barbot1}). 
Barbot and Choi in Appendix \ref{app:radd}, i.e.,  
Theorem \ref{thm:barbot2II}, show that if the boundary 
is totally geodesic and nonempty, then the $3$-manifold has  
to be an affine suspension. 
These results
enable us to answer the Carri\`ere conjecture in the positive: 
\begin{crscor}
\label{cor:crscor}
Let $M$ be a compact radiant affine $3$-manifold with empty or totally 
geodesic boundary.  Then $M$ admits a total cross-section to  
the radial flow. As a consequence $M$ is affinely diffeomorphic  
to one of the following affine manifolds\/{\rm :} 
\begin{itemize} 
\item a Benz\'ecri suspension over a  real projective surface of 
negative Euler characteristic with   empty or geodesic boundary,  
\item a generalized affine  suspension over a  real  projective sphere,  
a real  projective plane, or a hemisphere,   
\item a generalized affine suspension over a $\pi$-annulus   
\/{\rm  (} or  M\"obius band\/{\rm )} of type C\/{\rm ;}  
or an affine torus, an affine Klein bottle,  
an affine annulus \/{\rm (} or M\"obius band\/{\rm )}  
with geodesic boundary. 
\end{itemize} 
\end{crscor}  
 
The classification of Benz\'ecri suspensions and generalized 
affine  suspensions  over real projective  surfaces of  zero Euler 
characteristic reduces  to   the  classification  of  real  projective 
orbifolds of negative Euler characteristic and the classification of  
real  projective automorphisms of real projective tori or Klein bottles 
(see \cite{barbot2}). 
 
We obtain the following topological characterization of radiant affine 
$3$-manifolds.  
(See Scott \cite{Sc} for definition of a Euler number of  
a Seifert fibration.) 
\begin{topcor}
Let $M$ be a compact radiant affine $3$-manifold with empty or
totally geodesic boundary. Then $M$ is homeomorphic  
to a Seifert space of Euler number $0$ or a virtual bundle over  
a circle with fiber homeomorphic to a Euler characteristic zero 
compact surface.  
\end{topcor}

The purpose  of  this  paper is   to prove Theorem  A, approaching  by 
geometric   techniques   developed initially   for two-dimensions  in 
\cite{cdcr1} and \cite{cdcr2} and  later generalized to $n$-dimensions 
in the monograph \cite{psconv}.  
 
In Chapters  \ref{ch:prel}, \ref{ch:psconv}, \ref{ch:racom} 
and  Appendices \ref{app:dipping} and   
\ref{app:seqballs},  we will discuss  general 
$n$-dimensional real projective and affine  manifolds while there  are 
no losses from studying in the general dimensions.  The main arguments 
for Theorem A are  carried out in  Chapters  
\ref{ch:tdim}-\ref{ch:obtain}: First, using the 
decomposition result obtained in \cite{psconv}, we decompose a 
radiant  affine  manifold  into  $2$-convex  submanifolds  and radiant 
concave affine submanifolds (see Definition \ref{defn:cafm} in Chapter 
\ref{ch:psconv}). A real projective $3$-manifold is {\em $2$-convex} 
if a nondegenerate real projective map from  
$T^o \cup F_2 \cup F_3 \cup F_4$ for an affine $3$-simplex $T$  
in $\bR^3$ with sides $F_1, \dots, F_4$ always extends to one from $T$. 
A {\em concave affine} manifold is a real projective  
manifold with a very special affine structure so that any point  
of its universal cover is covered by domains affinely 
homeomorphic to affine half-spaces. 
Next,  we decompose   a $2$-convex radiant  affine 
manifold   into  convex  affine   manifolds  and  concave-cone  affine 
manifolds (see Definition  \ref{defn:ccaff}  
in Chapter \ref{ch:obtain}) in  Chapters \ref{ch:2conv}-\ref{ch:obtain}. 
A concave cone affine $3$-manifold is a real projective manifold 
with a special affine structure so that its universal cover is  
covered by domains affinely homeomorphic to affine quarter spaces. 
We  will prove Theorem A in  the final part   of  
Chapter \ref{ch:obtain} using results  from  Chapters \ref{ch:caff} 
and \ref{ch:pcrc}.   
In  Chapter  \ref{ch:caff}, we  show that  radiant 
concave affine manifolds  and concave-cone affine manifolds are  
generalized affine suspensions.  In Chapter \ref{ch:pcrc}, we discuss  
a step needed in Chapters \ref{ch:2conv} and \ref{ch:pcrc}.  
In Appendix C, we show
the nonexistence of affine manifolds whose developing maps
are universal covering maps of $\bR^3$ with a line removed.
More importantly, we show that the Carri\`ere conjecture holds 
if the affine manifold has a nonempty totally geodesic boundary 
of sufficiently general type; i.e., they admit total cross-sections.
 
In Chapter \ref{ch:prel}, we give definitions of  
real  projective   structures, 
developing maps, holonomy covers, and so on. We   discuss lifts  of  the   
developing map $\tilde M \ra \rp^n$ to $\SI^n$, and the holonomy  
homomorphism to the group $\Aut(\SI^n)$ of projective automorphisms 
of $\SI^n$. A {\em holonomy covering}\/ $M_h$ 
of $M$ is the cover of  $M$ corresponding to  the kernel of a holonomy 
homomorphism $h$.  $M_h$ is well-defined since $h$  can change only by 
a conjugation by a  real projective automorphism.  Then $\dev$ induces 
a metric on $M_h$ from the standard Riemannian metric on $\SI^n$,  
and the completing $M_h$ for the metric, we obtain the  Kuiper completion 
$\che{M_h}$ of $M_h$. We discuss the ideal set or frontier  of $M_h$ in 
$\che{M_h}$, convex subsets in the Kuiper completions, and 
extending maps  from one ball to  another one. This type of completion  
was first introduced by Kuiper \cite{Kp}  and is the most useful  
technical device in studying the incompleteness in geometric  
structures.  (Since this is the  most important technical  machinery 
needed in this paper, we  advise  our readers to become completely  
comfortable with involved notions.)  The  ideal set   
can be  thought of as points infinitely far away from points of  
$M_h$ if $M_h$ is endowed with a Riemannian metric $d_M$ induced from  
that of the compact manifold  $M$. We also define various  
polyhedra in $\SI^n$ and hence ones on $\che M_h$ as well.  
In particular, a {\em bihedron} is an $n$-ball in $\SI^n$  
that is the closure of a component of the complement of  
two distinct great $(n-1)$-spheres in $\SI^n$. A {\em trihedron}  
is one that is the closure of a component of the complement of 
three great $(n-1)$-spheres in general position.  
A {\em tetrahedron} is one for four great $(n-1)$-spheres in  
general position. Corresponding objects in $\che M_h$ are also defined.  
(This section can be browsed through if the reader is familiar 
with the field.)
 
In  Chapter \ref{ch:psconv}, we discuss $(n-1)$-convexity,   
and the decomposition of real projective $n$-manifolds into  
$(n-1)$-convex submanifolds and concave affine $n$-manifolds.      
However, the   submanifold  along  which  we 
decompose may not be totally geodesic.   
 
In Chapter \ref{ch:racom}, radiant affine manifolds, regarded as  
real projective manifolds, are discussed. We define the term 
generalized affine suspensions of a real projective surface, and show that  
a generalized affine suspension over a negative Euler characteristic real  
projective surface is a Benz\'ecri suspension over a real projective  
surface. We discuss how the radial flow  extends to  the projective  
completions, and define the origin, radiant sets, finitely ideal  sets,      
infinitely  ideal sets, and discuss the relationship between   
radiant  $n$-bihedra  and  $n$-crescents.  A so-called  
$n$-crescent is essentially  a convex $n$-bihedron one  of 
whose sides is  in the ideal set. 
 
In  Chapter \ref{ch:tdim},  we   apply  materials of  
Chapters  \ref{ch:prel}, \ref{ch:psconv}, and \ref{ch:racom} to 
three-dimensional  radiant affine  manifolds.   We show  that a  radiant 
affine   $3$-manifold   decomposes  into   $2$-convex  radiant  affine 
$3$-manifolds and radiant  concave affine $3$-manifolds along mutually 
disjoint union of totally geodesic tori or Klein bottles. It is  
proved in Chapter \ref{ch:totgeo} that the decomposing surface  
is totally geodesic unlike the ones in the general real projective 
manifolds. 
 
In Chapters \ref{ch:2conv}-\ref{ch:obtain},  
we study exclusively $2$-convex radiant affine $3$-manifolds 
with empty or totally geodesic boundary. We assume  that  $M$ is  
not convex in these chapters.  The so-called crescent-cone  is  
a  radiant trihedron with  three sides (which are lunes)   
so that two of its sides are  in  the  ideal   set respectively   
(see   Definition    \ref{defn:crscone}   and     Figure 
\ref{fig:crescent-cone}).   We will prove  Theorem  A and Corollary A, 
following the basic strategy analogous to that of papers \cite{cdcr1} 
and \cite{psconv}:  
\begin{enumerate} 
\item[(1)] We find a radiant tetrahedron detecting nonconvexity in  
$\che M_h$ for a compact radiant affine manifold $M$.  
\item[(2)]  By taking a suitable  subsequence of deck transformations, 
we  replace   it by a   crescent-cone in $\che M_h$ (by a blowing up   
argument as in \cite{cdcr1}).  
\item[(3)]  Using the equivariance and  local finiteness properties of 
the  collection  of  crescent-cones,  we  obtain a disjoint collection  
of concave-cone  affine submanifolds in $M$ of codimension $0$ with  
totally geodesic boundary.  
\item[(4)] We obtain the decomposition: the closures of components  
of the complement in $M$ of the union consist of convex radiant  
affine $3$-manifolds.  
\end{enumerate} 
More details are given below:  
\begin{enumerate}  
\item[(a)] Since $M$  is  not convex,  we  obtain a triangle in  $M_h$ 
detecting the nonconvexity of $M$.  
\item[(b)] Obtain  a   radiant  tetrahedron $F$ by radially extending the  
triangle.  The  tetrahedron has four sides $F_1, F_2, F_3,$ and $F_4$,  
where $F_4$  is a subset of infinitely ideal set, and $F_1^o, F_2^o$    
those of $M_h$, and $F_3$ meets the finitely ideal set and $M_h$---(1).  
\item[(c)] We choose an appropriate sequence of points  
$p_i \in F_3 \cap M_h$ leaving every compact subset of $M_h$.  
Choosing a fixed fundamental domain, a sequence of equivalent  
points $q_i$ in it, and the deck transformation $\deck^i$  
so that $\deck^i(p_i) = q_i$, we obtain a sequence of objects  
in $F$ pulled along with $p_i$ by $\deck^i$. We show  
that $\dev(\deck^i(F))$  
converges to a radiant $3$-ball in $\SI^3$. Hence, there exists  
a radiant convex $3$-ball $F^u$ in $\che M_h$ to which $F^i$ ``converges'' 
(see Appendix \ref{app:seqballs}). The convex $3$-ball $F^u$ is  either  
a radiant tetrahedron, or a radiant trihedron as $F^i$ are radiant  
tetrahedra. Using the following claim, we show that  there exists  
a so-called pseudo-crescent-cone or crescent-cone, which is a radiant  
tetrahedron or a radiant trihedron with exactly one side  
intersecting $M_h$ and the remaining sides in $\hideal{M}$.  
By Chapter \ref{ch:pcrc}, we rule out pseudo-crescent-cones---(2). 
\item[(d)] {\bf  Claim}\/: {\em We can  choose an appropriate sequence 
of points  $p_i \in F_3 \cap M_h$ leaving every compact subset   
of $M_h$ such that  the sequence  of  the $d_M$-distances  from $p_i$  
to  $(F_1 \cup F_2)\cap M_h$ goes to  infinity}\/  
(Proposition \ref{prop:MainEqu}):\hll  
To prove  the claim, we  suppose  that for every  sequence of points $p_i$,   
the $d_M$-distance   is bounded  above, and choose    $p_i$  in a  certain 
manner.  
We show this implies contradictions in Chapters \ref{ch:tric} and 
\ref{ch:rtrh}
when $F^u$ is a radiant tetrahedron and a radiant trihedron respectively. 
\item[(e)] In Chapter \ref{ch:obtain}, 
using the equivariant and locally finite collection of crescent-cones,   
we obtain a radiant set which covers a concave-cone affine  $3$-manifold 
in $M$ ---(3).  From  the  classification of  radiant  concave  affine 
manifolds and concave-cone affine manifolds in   
Chapter \ref{ch:caff}, we complete the proof of Theorem A.  ---(4).  
\end{enumerate} 
 
In Chapter \ref{ch:caff}, we classify radiant  concave affine manifolds  
and concave-cone   affine manifolds.    We   show  that  they  are   
generalized affine suspensions  of   Euler  characteristic   
nonnegative   real projective surfaces with geodesic boundary.   
The method uses the work by Tischler and heavily that by Barbot
and Choi in Appendix C.  
 
In Chapter \ref{ch:pcrc}, we show that pseudo-crescent-cones do   
not occur if $M$ is not convex.  
 
The  materials in  Appendices \ref{app:dipping}  and \ref{app:seqballs}  
are    adopted from the  paper \cite{cdcr1} to general $n$-dimension.   
The proofs are identical.  (We include these for  reader's convenience.)   
In  Appendix  \ref{app:dipping}, we discuss dipping intersection of  
two balls, and transversal intersection of two $n$-crescents.   
In Appendix \ref{app:seqballs},  we discuss sequences of compact convex 
$n$-ball in the Kuiper completion, and discuss their  convergence 
property when they include a common open ball.  
(See \cite{psconv} for details.)  

In Appendix C, written by Barbot and Choi, we show the
nonexistence of a radiant affine $3$-manifold whose developing map
is a universal covering map to $\bR^3$ with a line removed. 
The proof follows from the fact that such a manifold has a total
cross-section, which uses the work of Barbot by showing that 
the holonomy group of such a manifold is either solvable or 
has a hyperbolic element. 
The second part of Appendix C is devoted to showing 
that a compact radiant affine manifold with nonempty boundary 
has a total cross-section. The boundary here is assumed to be totally
geodesic and convex or must be a two-dimensional Hopf manifold.
We do this by showing that a certain connected abelian Lie group
of rank $2$ acts on the manifold giving us a foliation by tori or
a decomposition of the manifold into pieces which admit 
total cross-sections. Again, the theory of Barbot developed 
for the case when the radiant affine manifold contains 
a totally geodesic surface tangent
to the radial flow is used. These are key steps for 
the main part of the paper.

\section*{Acknowledgement} 
The author would like to  thank  Boris Apanasov, Thierry Barbot,  Yves 
Carri\`ere, William Goldman,  Yoshinobu Kamishima, Hyuk Kim, Sadayoshi 
Kojima, Shigenori Matsumoto,  William Thurston,  and Abdelghani Zeghib 
for    many  helpful   discussions. The    author  particularly enjoyed 
the cooperation with Thierry  Barbot generously explaining and supplying   
him with short arguments,  
in many occasions enabling the author  
to complete  this work. The author also thank him for writing 
an appendix with him, providing key results in the boundary case 
using his approach essentially.
Furthermore, the author likes to thank  
Yves Carri\`ere  and David Fried for posing many extremely  interesting   
conjectures in  this    field.   As usual, we need    more    
well-posed problems  in  this field. With  their  guidance, the field   
has  become  more  mature and fertile. Also, the   author thanks  
Sang-Eun   Lee for much help on  my computer  running Linux.   
The figures  were drawn  by xfig  and  Maple. 
Finally, a final
portion of this work was completed during the author's visit 
to the IHES and the ENS-Lyon. The author appreciates the hospitality 
given by both institutes with great tradition of excellence 
in geometry of all kinds. The author also thanks GARC for generous 
financial support. Finally, the author is very thankful for 
the support of Haeyon while the author spent numerous hours 
writing and revising this paper.
 
\chapter{Preliminary} 
\label{ch:prel} 
 
In this chapter, we define developing maps and holonomy and discuss   
the  natural relation  between  affine  manifolds   and  real projective  
manifolds, some history on the subject of affine and real projective
manifolds, developing maps and holonomy homomorphisms lifted to  
$\SI^n$ and  the group $\Aut(\SI^n)$ of projective  automorphisms,  
the completion $\che M$ of the universal cover $\tilde M$ 
of a real projective manifold $M$ due to Kuiper, and the holonomy  
cover $M_h$ of $M$. We define convex sets in $\SI^n$ and $\che{M}$  
and the Kuiper completion $\che M_h$ of $M_h$. We introduce an important   
lemma that how two balls in $\che{M_h}$ meet may be completely  
read from their images under the developing map $\dev$ of $M$ under 
very natural circumstances.   
 
An $(X, G)$-structure on a manifold $M$ is given by a maximal atlas of 
charts to    $X$ where the transition  functions   are in $G$.   
An $(X, G)$-map is a local diffeomorphism preserving 
$(X,  G)$-structures locally. (Excellent treatments can be found in 
Ratcliffe \cite{RT} and Thurston \cite{ThB}.)   (This will be our  view 
point of $(X,    G)$-structures in this   paper.  However,  there  
are differential-geometry version of $(X, G)$-structures which we will not 
use.)  
 
Given an   $(X,  G)$-structure on  $M$,  we  can  associate  to  it an 
immersion $\dev: \tilde M \ra X$, called {\em a developing map}, and a 
homomorphism $h$ from the group $\pi_1(M)$  of deck transformations to 
$G$, called  a {\em holonomy  homomorphism}\/ satisfying $h(\gamma)\circ 
\dev = \dev \circ  \gamma$ for each $\gamma  \in \pi_1(M)$.  $\dev$ is 
obtained  by analytically extending coordinate   charts of $\tilde  M$ 
altered by post-composition with an element of  $G$ so as to extend in 
sequence.  
Since given any deck transformation  $\gamma$, $\dev \circ \gamma$  is 
again    a developing  map, it    equals   $h(\gamma) \circ \dev$  for 
$h(\gamma) \in  G$.   One can see  easily that  $h$ is  a homomorphism 
$\pi_1(M) \ra G$.  $h$ is  said to be  the {\em holonomy homomorphism}\/ 
and  its image the {\em holonomy  group}.  
Given an $(X,  G)$-structure on $M$, 
$(\dev, h)$ is  determined  up  to  the equivalence relation   $(\dev, 
h(\cdot)) \sim (\vth \circ \dev, \vth \circ h(\cdot) \circ \vth^{-1})$ 
where $\vth$ is an element of $G$. (Fundamental class of examples  
are given by Sullivan-Thurston \cite{ST}.)

A  real projective  space $\rp^n$ is  given  as the  quotient space of 
$\bR^{n+1} - \{O\}$ for the origin $O$ by relation generated by  
scalar multiplications. The group of general  linear transformations  
$\GL(n+1, \bR)$ descends to the group $\PGL(n+1, \bR)$ acting on  
$\rp^n$ as projective transformations. An affine space $\bR^n$ is a  
Euclidean space  $\bR^n$ with the group of affine transformations  
$\Aff(\bR^n)$ acting on it.  
 
A {\em real projective $n$-manifold}\/ is an  $n$-manifold  with  
an $(\rp^n, \PGL(n+1, \bR))$-structure;  an  {\em affine $n$-manifold}\/   
is  one with  an $(\bR^n,    \Aff(\bR^n))$-structure.     

A prominent feature of real projective or affine manifolds 
is that geodesics are defined by connections. 
A curve in a real projective manifold is {\em geodesic}\/ if 
it corresponds to a straight line in the projective space 
under charts, and similarly for curves in affine manifolds.
{\em Segments} or {\em lines} are imbedded images of geodesics.

For affine manifolds, the completeness
of geodesics makes sense since an affine parameter of a geodesic 
is defined. By a {\em complete} real line, we mean a geodesic
in an affine manifold which has affine parameter from $-\infty$
to $\infty$. An affine manifold is often incomplete even if 
it is closed (see Kobayashi \cite{Kob:84}, Kobayashi-Nagano
\cite{KN}, and Sullivan-Thurston \cite{ST}.)

Let    $\rp^{n-1}$     be   an 
$(n-1)$-dimensional subspace of $\rp^n$. Then $\rp^n - \rp^{n-1}$ has a 
natural structure of an  affine space $\bR^n$ so  that the subgroup of 
projective transformations acting   on the space equals  $\Aff(\bR^n)$ 
and   projective    and   affine  geodesics   agree     on  $\rp^n   - 
\rp^{n-1}$ up to parameterization. Obviously, an  affine  structure on   
a manifold  induces a unique real  projective structure since charts  
to $\bR^n$ are naturally  charts to  $\rp^n$ and  affine  
transition maps are real projective ones by the above  identification.   
We will always regard an affine  manifold as a real projective   
manifold obtained in this natural manner.  
 
We assume in  this paper that the  manifold-boundary $\delta M$ of a 
real projective manifold $M$  is empty or  totally geodesic; i.e., for 
each point of $\delta M$, there exists  an open neighborhood $U$ and 
a  chart  $\phi: U  \ra \rp^n$  so that  $\phi(U)$  does  not meet the 
subspace $\rp^{n-1}$ and  is a convex open  domain  intersected with a 
closed affine half-space in $\rp^n  - \rp^{n-1}$.  We also assume that 
affine  manifolds have empty or  totally geodesic boundary; i.e., they 
have empty or totally geodesic boundary as real projective manifolds.  
 
We give some simple examples and a short history on real projective
and affine manifolds relevant to this paper.
\begin{exmp}\label{exmp:hopf} 
Let $\vth$ be a matrix with the absolute value of  
all its eigenvalues greater than one. 
The quotient of $\bR^n - \{O\}$ by $\langle \vth \rangle$ is  
homeomorphic to $\SI^{n-1} \times \SI^1$, and has  
a natural radiant affine structure. The manifold is said to be  
a {\em Hopf manifold}. 
 
If $\vth$ preserves an $(n-1)$-dimensional subspace,  
then the quotient of $H-\{O\}$ where $H$ is a closed affine  
half-space bounded by the subspace is homeomorphic to  
an $(n-1)$-ball times $\SI^1$. This is said to be 
a {\em half-Hopf}\/ manifold. This has totally geodesic boundary.  
 
Let $x_1, \dots, x_n$ be coordinates of $\bR^n$. 
Let $H$ be an open upper half-space with the positive  
$x_n$-axis removed. Let $\vth$ be given by  
\[(x_1, \dots, x_{n-1}, x_n) \mapsto (2x_1, \dots, 2x_{n-1}, x_n)\]  
and $\vpi$ by  
\[(x_1, \dots, x_{n-1}, x_n) \mapsto (2x_1,\dots, 2x_{n-1}, 2x_n).\]  
Then $H/\langle \vth, \vpi \rangle$ is homeomorphic to 
the $\SI^{n-2} \times \SI^1 \times \SI^1$, and has  
a radiant affine structure. Denote this radiant affine  
manifold by $\mathcal{E}_2$. 
\end{exmp} 

One of the simplest example  of an annulus with 
a  real projective structure is
given as the quotient  of an invariant triangle  by an infinite cyclic
group of  real  projective transformations  with representing matrices
conjugate to diagonal matrices  with positive distinct  eigenvalues; a
compact annulus with  interior real projectively homeomorphic  to such
an annulus is said to be an {\em elementary}\/ annulus (see \cite{cdcr2}
and \cite{CG}). 

The  Euler characteristic  zero closed  surfaces  with real projective
structures   were classified by    Goldman \cite{Gconv}   and Nagano-Yagi
\cite{NY}; such a  surface either admits  a natural compatible  affine
structure or  is  built  up from   elementary  annuli by  gluing  (see
\cite{Gconv} and \cite{CG}).  If a closed surface $\Sigma$ with a real
projective structure has a negative  Euler characteristic, then it  is
shown  that  $\Sigma$  is obtained   from  surfaces  with convex  real
projective structures and elementary annuli  by projective
gluing (see \cite{CG}, \cite{cdcr1},  \cite{cdcr2},  
and  \cite{cdcr3}).   Since      Goldman
\cite{Gconv} classified  and  constructed all convex  real  projective
structures on surfaces  up to isotopy,  we can classify and  construct
all real projective   structures on closed  surfaces.  The deformation
space of real   projective structures on $\Sigma$, i.e.,   the space of
equivalence classes of real projective   structures on $\Sigma$ up  to
isotopy  with appropriate topology,  is shown to  be homeomorphic to a
countably infinite  disjoint  union of cells  of dimension  $16g - 16$
where $g$ is the genus of $\Sigma$ \cite{CG}. 

A source of examples of real projective  manifolds is obtained from the
Klein  model  of  hyperbolic space:   Given  a  Lorentzian  metric  on
$\bR^{n+1}$, the   positive part  of the hyperboloid  $-  x_0^2 +  x_1^2 +
\cdots + x_n^2 = -1$ is the  hyperbolic $n$-space $H^n$, and the group
$\PSO(1,  n)$ acts on  $H^n$  as a   group  of isometries.  Under  the
projection $\bR^{n+1} - \{O\}  \ra \rp^n$, $H^n$ becomes identified with
an open ball  $B^n$ in $\rp^n$ and $\PSO(1,  n)$ with the natural copy
of   $\PSO(1, n)$ in  $\PGL(n+1, \bR)$.  Therefore,  it follows that a
hyperbolic structure on a    manifold $M$ induces a   real  projective
structure    on $M$. When   $n  =  2$,  Koszul \cite{Kz},  Kac-Vinberg
\cite{KV}, and Goldman   \cite{Gconv} found parameters  of  nontrivial
deformations  of  hyperbolic    real   projective  structures  to   
non-hyperbolic  but  convex   ones.   It  is  thought   that many such
deformations exist in  dimensions $\geq 3$ 
giving  us a plethora  of  real projective structures. 
Koszul found some higher dimensional deformations, and Goldman
\cite{Glett} parameters of deformations of real projective   
structures on the complement  of the figure eight knot. 

Many examples of affine manifolds were   constructed   by
Sullivan-Thurston \cite{ST} from real projective manifolds and 
complex projective manifolds by suspension-like constructions. 
If a  closed  surface admits  an affine
structure, then its   Euler  characteristic  is   zero as   shown   by
Benz\'ecri   \cite{Ben}.   The    outstanding  conjecture   on  affine
$n$-manifolds is Chern's conjecture that  the Euler characteristic of
a  closed affine manifold is  zero.  Kostant-Sullivan \cite{KS} showed
this to be  true if the  affine manifold  is complete,  i.e.,  it is a
quotient of $\bR^n$  by a properly discontinuous and  free action of a
subgroup   of  $\Aff(\bR^n)$.  If the fundamental    group of a closed
affine manifold is amenable (or isomorphic to a free product
of  groups of polynomial growth), then  Chern's conjecture  is true by
Hirsch-Thurston    \cite{HT}      (see  also    Kim-Lee    \cite{KL1},
\cite{KL2}). However, this conjecture is yet to  be proved for general
case, and strangely there seem  to be really no  tools to study  this
question. 

Other examples of affine and projective manifolds were constructed 
by Smillie \cite{Smdiag}, in particular with diagonal holonomy, 
and later more general class were found by Benoist \cite{Beno:94}
and \cite{Beno2}. 

Another conjecture due to L. Markus states that  an affine manifold is
complete if and only if it has a  volume form parallel with respect to
the affine connection. This  conjecture was settled if the fundamental
group  is   virtually abelian, nilpotent,  or solvable   of rank 
$\leq n$ by Smillie \cite{smillie} and 
Fried, Goldman, and Hirsch \cite{FGH}. 
The conjecture is verified if the holonomy group lies in the group
of Lorentzian transformations of $\bR^n$ with a flat Lorentz metric by
Carri\`ere \cite{Card}. 

One  question is to  find an example of a three-manifold  not admitting a
real projective or affine structure, as asked by Goldman. Smillie  
\cite{Sm2} showed that a connected sum of lens spaces does not admit 
an affine structure.  Choi \cite{uaf} showed that an affine $3$-manifold 
with parallel volume form must be   irreducible and has a   $3$-cell 
as a  universal cover, answering Carri\`ere's question  \cite{Card}. 
(See \cite{uaf2} for a summary.) 

Going to the main part of this section,
we now define  the splitting of a real  projective manifold  $M$ along 
a submanifold $N$ of codimension one  in   $M^o$. We take a   regular 
neighborhood  $K_i$ of   each component  $S_i$   of $N$.  If $S_i$  is 
two-sided,   then  we take  the closure in $K_i$ of each component 
of $K_i - S_i$.  
If $S_i$ is one-sided, then we take  a double  cover  $(K'_i, p_i)$   
of  $K_i$ corresponding to the index-two-subgroup of $\pi_1(K_i)$ given by   
the homomorphism $\pi_1(\delta K_i) \ra \pi_1(K_i)$ induced by the  
inclusion map and remove  $S'_i$ the part corresponding to $S_i$ in   
$K_i$.  Take  a component $J$ of $K'_i - S'_i$.  Then the closure  
$J'$  of $J$ has  a  natural  real  projective or affine structure   
with totally geodesic boundary equivalent to $S'_i$. 
$J$ is obviously identical  
with $K_i - S_i$. For each component of $K_i - S_i$, we identify 
it with the appropriate subset of the closures as above to obtain the split  
manifold. Our construction gives us a real projective or affine  
manifold with totally geodesic boundary; the resulting manifold is  
said to be obtained from {\em splitting along $N$}.  
 
There exists  a standard  double cover  $p: \SI^n  \ra \rp^n$ from the 
standard  unit sphere $\SI^n$  in  $\bR^{n+1}$ with a standard  Riemannian 
metric  $\mu$  which is  projectively  flat and has  the same geodesic 
structure as one induced by $p$ from $\rp^n$, i.e., $\mu$-geodesics in 
$\SI^n$ agree   with     projective  geodesics in  $\SI^n$     up   to 
parameterization.  
The metric  $\mu$  also induces one on  $\rp^n$, 
and  the geodesic structure of this  Riemannian metric agree with that 
of  the  given    real  projective  structure    on  $\rp^n$   up   to 
parameterization. Denote  by $\bdd$  the  distance metric  on  $\SI^n$ 
induced from $\mu$.  $\SI^n$ can also  be identified with the quotient 
of $\bR^{n+1} -\{O\}$ by the  equivalence relation  generated by  
positive scalar multiplications. The general linear group  
$\GL(n+1, \bR)$ acts on $\bR^{n+1} - \{O\}$ and  hence on its quotient     
space $\SI^n$. The induced self-diffeomorphisms form a group  
$\Aut(\SI^n)$ of all projective self-diffeomorphisms of $\SI^n$, which   
obviously is isomorphic to the subgroup $\SL_\pm(n+1, \bR)$ of  
$\GL(n+1, \bR)$ consisting of linear maps with determinant $\pm 1$.  
 
A homeomorphism $f: X \ra Y$ for metric  spaces $(X, \bdd_X)$ and $(Y, 
\bdd_Y)$ are said to be {\em quasi-isometric}\/ if  
\[C^{-1}\bdd_X(x, y) \leq \bdd_Y(f(x), f(y)) \leq C\bdd_X(x, y), x, y \in X\] 
hold for a uniform   positive constant $C$ independent of  $x$ and  $y$.    
We introduce a Kuiper  completion  of   $\tilde  M$ (see \cite{Kp}).   
Given a development pair  $(\dev, h)$ of a real projective 
manifold $M$, $\dev$ lifts to an immersion $\dev': \tilde M \ra \SI^n$ 
and  $h$  to  $h':\pi_1(M)   \ra  \Aut(\SI^n)$ so  that   they satisfy 
$h'(\gamma) \circ \dev'  = \dev'  \circ  \gamma$ for each $\gamma  \in 
\pi_1(M)$.  We induce the Riemannian  metric $\mu$ on $\tilde M$  from 
$\mu$ on  $\SI^n$, complete the distance  metric $\bdd$ on $\tilde M$ induced 
from $\mu$,  and denote by  $\che M$ the Cauchy completion  of $\tilde M$ and 
$\bdd$ the completed metric.  By an abuse  of notation, we will simply 
denote $(\dev',    h')$ by  $(\dev, h)$.    Since   $\dev$ is distance 
non-increasing,   $\dev$ extends to  the  completion $\che M$ and since 
each deck transformation $\vth$ is a quasi-isometry with respect to $\bdd$, 
it extends to a self-homeomorphism of  $\che M$.   We will denote the 
extended maps by same symbols $\dev$ and $\vth$ respectively (see \cite{KT} 
and \cite{psconv}).  
 
Given a real projective  structure on $M$, any other development 
pair   $(\dev'':\tilde M \ra  \SI^n, h''(\cdot))$  equals $(\vth \circ 
\dev, \vth \circ h(\cdot) \circ \vth^{-1})$ for an element $\vth \in 
\Aut(\SI^n)$. The    completion  with     respect  to  $\dev''$    is 
quasi-isometric to $\che M$ with $\bdd$ above since the metric induced 
from $\dev''$ is quasi-isometric with respect to the  original $\mu$.  
We will  say that $(\che M, \bdd)$ is a {\em Kuiper completion}\/ of  
$\tilde M$.  
 
\begin{figure}[t] 
\centerline{\epsfxsize=3.3in \epsfbox{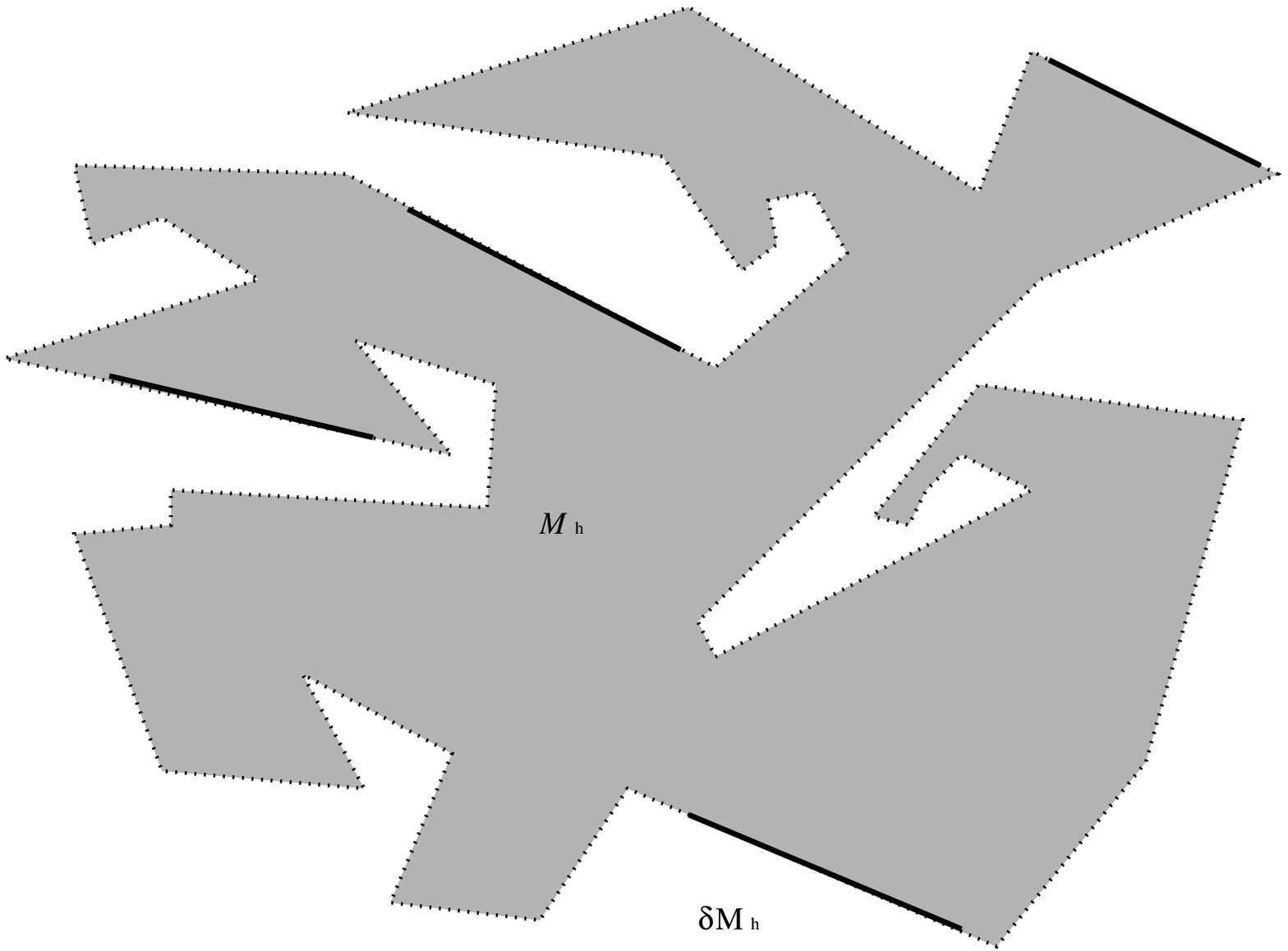}}  
\caption{\label{fig:figm}   A  figure of   $\che M$.    The dark lines 
indicate   $\delta \tilde M$ and  the  dotted  lines  the  ideal  boundary 
$\ideal{M}$.}  
\end{figure} 
\typeout{<<rfig1.eps>>}
 
Recall  that the complement in  $\rp^n$ of a codimension-one  
subspace $\rp^{n-1}_\infty$ can be  identified with  $\bR^n$,   and    
the group of   real   projective transformations    of    $\rp^n$   
acting    on  $\rp^{n-1}_\infty$  to 
$\Aff(\bR^n)$ (see Berger  \cite{B}).  The subspace $\rp^{n-1}_\infty$ 
corresponds to the great sphere  $\SII$ under the double covering  map 
$\SI^n  \ra  \rp^n$,  and $\bR^n$  to an   open  hemisphere bounded by 
$\SII$.  Let $\mathcal{H}$ denote the closed hemisphere including this, and 
$\Aut(\mathcal{H})$  the    group   of  orientation-preserving  projective 
transformations acting on $\mathcal{H}$.  Under the double covering map the 
interior $\mathcal{H}^o$ of $\mathcal{H}$ corresponds  to  the affine  space 
$\bR^n$ and $\Aut(\mathcal{H})$ to  $\Aff(\bR^n)$ in a one-to-one manner, 
preserving the geodesic structures  up to parameterization.    
(We  will  fix our choice  of  the hemisphere $\mathcal{H}$ and the 
identification of $\mathcal{H}^o$ with $\bR^n$, which determines
of $\SII$,  and we use $\bR^n$ and $\mathcal{H}^o$ interchangeably, and
so use $\Aff(\bR^n)$ and $\Aut(\mathcal{H})$.)  

The correspondence between $\bR^n$ and the open $n$-hemisphere
is realized by moving the subspace $\bR^n$   
in $\bR^{n+1}$ in the orthogonal direction by a unit and   
stereographically project from the origin onto a hemisphere $\SI^n$
(see \cite{psconv}).

Let $M$ be an affine $n$-manifold. By above consideration, when  
$\mathcal{H}^o$  is identified with $\bR^n$ and the affine transition  
functions with the real projective ones, $M$ has a naturally induced real  
projective structure.   Thus,  there  exists a development  pair 
$(\dev, h)$  for $M$ considered as  a real projective manifold so that 
$\dev: \tilde  M  \ra  \mathcal{H}$  and 
$h(\gamma)  \in  \Aut(\mathcal{H})$ for 
$\gamma \in \pi_1(M)$, and clearly $\dev$ maps $\che M$ into $\mathcal{H}$. 
 
We introduce the   holonomy   cover  $M_h$   of a   real    projective 
$n$-manifold $M$: Since the holonomy homomorphism may change only by a 
conjugation by  an element of $\Aut(\SI^n)$,  the kernel $K_M$  of the 
holonomy homomorphism    $h$    is independent   of  the  choice    of 
$h$. Considering  $\pi_1(M)$ as the  group of deck transformations, we 
define   the {\em   holonomy cover}\/ $M_h$  as the   cover  of   $M$ 
corresponding to $K_M$; that is, we define  $M_h = M/K_M$ and identify 
$K_M$   with  $\pi_1(M_h)$.   While  we   have  $\dev  \circ   \vth  = 
h(\vth)\circ  \dev = \dev$ for  any deck transformation $\vth$ in the 
kernel,  it   follows that   $\dev$   induces  a   well-defined  local 
diffeomorphism $\dev_h: M_h \ra \SI^n$.  Given $\dev_h$, we can define 
a holonomy homomorphism   $H_h:  \pi_1(M)/\pi_1(M_h) \ra  \Aut(\SI^n)$ 
from the deck-transformation group $\pi_1(M)/\pi_1(M_h)$  of $M_h$ as 
we did for $(\dev, h)$; the resulting pair  $(\dev_h, H_h)$ is said to 
be a {\em  development pair}.  Clearly $H_h$ is induced homomorphism  
from $h$. Similarly  to above, other  development 
pair $(\dev'_h, H'_h(\cdot))$  equals $(\vth \circ \dev_h,  \vth \circ 
H_h(\cdot) \circ \vth^{-1})$  for  an element $\vth  \in \Aut(\SI^n)$. 
(For convenience,   we  drop  the subscripts   $h$   from the notation 
$\dev_h$ from now on and write $h$ for $H_h$.)  
 
The immersion $\dev$ induces a Riemannian metric on $M_h$ from $\SI^n$,  
which we denote by $\mu$, and  $\mu$ induces a distance metric  $\bdd$ 
on $M_h$; a  completion  of $M_h$ is denoted  by  $\che M_h$ and   the 
completed metric by $\bdd$; $\hideal{M}$ denotes the ideal set $\che 
M_h - M_h$. As before, the developing map $\dev:M_h \ra \SI^n$ extends  
to a map $\che M_h \ra \SI^n$ and each deck transformation $\vth$  
extends to a self-diffeomorphism $\che M_h \ra \che M_h$. We will denote  
the extensions by $\dev$ and $\vth$ respectively. 
 
\begin{exmp}\label{exmp:compl} 
The universal cover of a Hopf manifold $M$ may be identified with  
$\brno$, and $\che M_h$ with the $n$-hemisphere that is  
the closure of $\brno$. The ideal set $\hideal{M}$ is the union 
of $\{O\}$ and the sphere $\SI^{n-1}$ which is the boundary.  
 
The holonomy cover of $\mathcal{E}_2$ may be identified with  
$U - l$ where $U$ is given by $x_n > 0$ and $l$ the $x_n$-axis.  
$\che \mathcal{E}_{2, h}$ equals the closure of $U$ in $\SI^n$,  
and the ideal set the union of the boundary of $U$ in $\SI^n$  
and $l$ intersected with the closure of $U$. 
\end{exmp}

A segment  or a line in  $\SI^n$  is said to   be {\em convex}\/  if its 
$\bdd$-length is $\leq \pi$.   
A segment  of $\bdd$-length $\leq \pi$ is always mapped to   
one of $\bdd$-length $\leq \pi$ by projective automorphisms of  
$\SI^n$. We define a {\em convex subset}  
of $\SI^n$ as a   subset such that  any two  points  of  the set can    
be connected by a segment in the set of $\bdd$-length $\leq \pi$.  
 
\begin{figure}[t] 
\centerline{\epsfxsize=3.0in \epsfbox{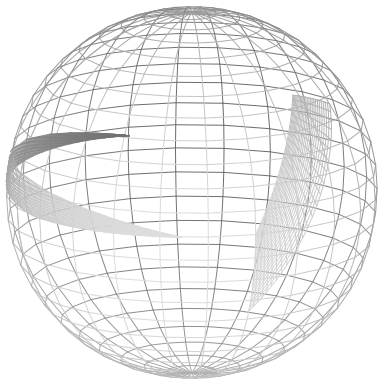}}  
\caption{\label{fig:cnv} Some  examples   of  convex   sets in 
$\SI^2$. The left one indicate a lune and the right one 
a simply convex set.}  
\end{figure}

\begin{figure}[t] 
\centerline{\epsfxsize=3in \epsfbox{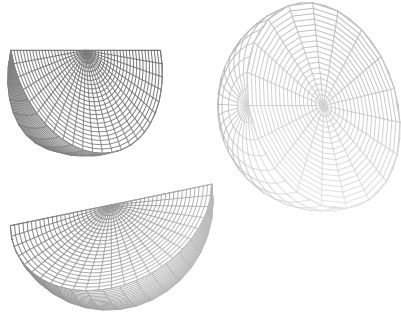}}  
\caption{\label{fig:cnvsets} Some examples  of convex sets in  $\SI^3$ 
stereographically projected  to $\bR^3$  and translated: Two  trihedra 
and one bihedron. Some of their sides are  drawn as totally geodesic  
sides for computer graphics convenience.}  
\end{figure}

A  {\em great $0$-sphere}\/ is  the set of  antipodal points, which is 
not  convex.  A  {\em  simply convex}\/ subset   of $\SI^n$ is  a convex 
subset of a $\bdd$-ball of radius $< \pi/2$,  i.e., it is a precompact 
subset of  an open   hemisphere,   and when  the open  hemisphere   is 
identified with an affine space, it is  a bounded affinely convex set. 
(See Figures \ref{fig:cnv} and \ref{fig:cnvsets}).  
 
It is easy  to categorize all compact convex subsets of $\SI^n$; they 
are homeomorphic to an $i$-ball   or an $i$-sphere always.   
A compact convex  subset of $\SI^n$ is  either a {\em  great $i$-sphere}\/  
$i \geq 1$, i.e.,  a  totally geodesic  $i$-sphere;   
an  {\em $i$-hemisphere}, i.e., the closure of a component  
of a great  $i$-sphere with a great $(i-1)$-sphere in it removed;  
or a  convex, proper compact subset of an $i$-hemisphere. 
Hence a compact convex subset is always homeomorphic to an $i$-ball 
or an $i$-sphere, and hence is a manifold, and we can define  
the boundary $\delta A$ of a compact convex subset $A$. 
A convex subset of $\SI^n$ is either a great $i$-sphere or  
a convex subset of a convex $i$-ball for some integer $i$. 
A convex,  proper compact subset of an $i$-hemisphere 
contains a unique great sphere of dimension $j$, $0 \leq j  < i$ or is 
a compact simply convex $i$-ball. (See \cite{psconv} for details.) 
 
See  Chapter  6  of  Ratcliffe \cite{RT} for   more  details on  convex 
sets. Note that our definition of convexity is slightly different from 
the book.  Assuming  that a subset  $A$ of $\SI^3$ is   not a pair  of 
antipodal points, a subset $A$ is convex if and  only it is ``convex'' 
in the  sense of  Ratcliffe \cite{RT}.  A  pair of  antipodal points is 
``convex''   according to the  book.   A minor  difficulty is that the 
intersection of two  convex sets may  not  be convex according to  our 
definition. But most of   the theory passes   to ours by mild  and 
obvious modifications, which we will not endeavor to list.  
 
A {\em side}\/ of a convex compact subset $C$ of $\SI^n$ is a  nonempty, 
maximal convex subset of $\delta C$.  A  {\em convex polyhedron} $P$ 
in $\SI^n$ is a nonempty, compact,  convex subset of $\SI^n$ such that 
the collection of its  sides is finite. A  side of a convex polyhedron 
of dimension $n$ is  again a convex polyhedron of  dimension $n-1$.  A 
convex polyhedron is always homeomorphic to a ball or equal to a great 
$i$-sphere, $i \geq 1$.  
 
When our dimension $n$  is $3$, we will  use the  terms as follows:  A 
{\em convex polygon}\/ in $\SI^3$ is a convex disk in a great $2$-sphere 
with finitely many sides, and a {\em convex  polyhedron}\/ in $\SI^3$ is a 
convex $3$-ball in $\SI^3$ with finitely many sides. 
A {\em lune} is the closure of a component of a great $2$-sphere in  
$\SI^3$ (or sometimes in $\SI^2$) with two distinct great  
circles removed. A {\em triangle} is the closure of a component of one  
with three distinct great circles in general position removed. 
 
A convex but not simply  convex polyhedron in $\SI^3$ is  either a 
$3$-hemisphere; a {\em bihedron}, i.e.,  the closure of a component of 
$\SI^3$ with two distinct great $2$-spheres removed;   or a {\em tube 
over  a $k$-gon} for $k \geq 3$, i.e., the closure  of  
an appropriate component of $\SI^3$ with 
$k$ distinct great $2$-spheres all containing a common pair of antipodal 
points removed.  Any  other convex polyhedron  in $\SI^3$  is a  simply convex 
polyhedron.  Since  they form a bounded  subset of the open hemisphere 
identified with   an affine $3$-space,  the  usual theory of Euclidean 
convex polyhedra applies.  
 
The only convex polyhedron in $\SI^3$ with  three sides is a tube over 
a triangle for  which we reserve  the term {\em trihedron}\/; for a  
simply convex polyhedron with four side, we reserve the term  
{\em tetrahedron}.  
 
\typeout{<<rfcvs.eps>>}  
 
\begin{defn}\label{defn:cone} 
Given a subsets $K$ of  $\SII$ and the origin $O$  in $\mathcal{H}$, we say 
that the union of all segments starting  from $O$ and ending at points 
of $K$ is the {\em cone over}\/ $K$.  
\end{defn} 
 
\begin{rem}\label{rem:cone} 
It is easy to see that a cone over a compact convex subset $K$ of $\SII$ 
is a compact convex set. We also have the correspondence that $K$ is a 
simply convex $i$-ball if  and only if the  cone over $K$ is a  simply 
convex $(i+1)$-ball. Moreover $K$ is a polyhedron  of dimension $i$ if 
and only if the cone over $K$ is a polyhedron of dimension $i+1$.  
\end{rem} 
 
A {\em convex segment}\/ in $\che M_h$ is  a subset $A$ such that $\dev|A$ 
is an imbedding onto a convex segment in  $\SI^n$.  Given a subset $A$ 
of $\che M_h$, we say that $A$ is convex if two points of $A$ can be 
connected by a convex segment in $A$.  A {\em tame}\/ subset of $\che M_h$ 
is a  convex  subset of $M_h$   or a convex  subset  of a compact 
convex subset of $\che M_h$. The closure in  $\che M_h$ of a convex subset 
of $M_h$ is always  tame  since for any  convex subset   $A$ of 
$M_h$, $\dev|A$ is  a $\bdd$-isometry onto $\dev(A)$, and, hence, 
$\dev$ restricted to the closure $\clo(A)$ of $A$ is a $\bdd$-isometry 
onto  $\clo(\dev(A))$. Since $\dev|A$ is a  $\bdd$-isometry for a tame 
set $A$, we also  see that $\dev|A$ for  a tame set  $A$ is  always an 
imbedding onto $\dev(A)$.  

Note that $M$ is {\em convex}\/ if $M_h$ is a convex subset of $\che M_h$.
Unfortunately, this definition is not equivalent to one in
the introduction; however, it is almost the same (see Theorem A.2 
in \cite{psconv}) as convex subsets of $\SI^n$ are almost always 
in an open hemisphere. (Recall that an open hemisphere corresponds 
to an affine patch.) For affine manifolds, the definitions are equivalent.
In this paper, there is no confusion since we study only
real projective manifolds coming from affine manifolds.

\begin{defn}\label{defn:i-ball} 
We define an {\em $i$-ball}\/ $A$ in $\che M_h$ to be 
a compact  subset of  $\che M_h$ such that $\dev|A$  
is a homeomorphism to an $i$-ball (not necessarily convex) in  
a great $i$-sphere and  its manifold interior $A^o$ is a subset of  
$M_h$. A {\em great $i$-sphere}\/ $A$ in $\che M_h$ is a compact  
subset of $\che M_h$ such that $\dev|A$ is a homomorphism onto 
a great $i$-sphere in $\SI^n$. 
\end{defn}  
As in \cite{psconv}, a topological $i$-ball is defined as above 
but $A^o$ is not required to be in $M_h$. 
 
We define {\em  convex segments}, {\em  convex sets}, {\em tame sets}, an 
{\em $i$-balls}\/ in $\che M$ in the exactly the same manner 
as in $\che M_h$. 
 
A convex {\em polyhedron}\/ in $\che M_h$ is a topological 
$i$-ball $A$ such that $\dev(A)$ is a polyhedron. 
A convex {\em polyhedron}\/ in $\che M$ is defined similarly. 
 
A compact convex subset $A$ of $\che M_h$ is tame, and hence  
is either a topological $i$-ball or a great $i$-sphere for 
some $i$.  
 
When  the dimension  of  $M$ equals  $3$,  we  will  use terms in   
the following manner.   A {\em convex polyhedron}\/ $A$ in   
$\che M_h$  is a convex $3$-ball in $\che M_h$ whose image is  
a convex polyhedron in $\SI^3$. A {\em convex polygon}\/ in $\che M_h$  
is a convex topological $2$-ball of $M_h$ whose image is a convex polygon  
in a great sphere in  $\SI^3$. These include lunes and triangles. 
A {\em side}\/ of the polyhedron or polygon is defined  
obviously as above. (Note that hemispheres, bihedra, trihedra,  
tetrahedra in $\che M_h$ are polyhedra defined in  
this sense.) These definitions also apply to subsets of 
$\che M$ in the obvious manner. 
 
\begin{rem}  
Note that all definitions  here are well-defined  under the change  of 
development  pairs and   our  definition   of convexity   is  directly 
influenced  by the fact that  we require $\dev$ restricted to relevant 
convex sets to be imbeddings to $\SI^n$.  
\end{rem}  
 
The  following theorem and lemma are  true both in  $\che M$ and $\che 
M_h$ as there are no difference in proofs.  Let $F_1$ and $F_2$ be two 
convex $n$-balls  in  $\che  M_h$.  We say   that  $F_1$ and  $F_2$ {\em 
overlap}\/ if  $F_1^o \cap F_2 \ne  \emp$, which is equivalent to $F_1 
\cap F_2^o \ne \emp$ or $F_1^o \cap F_2^o \ne \emp$.  
 
\begin{prop}\label{prop:overlap} 
If $F_1$ and $F_2$  overlap, then $\dev|F_1  \cup F_2$ is an imbedding 
onto  $\dev(F_1)  \cup \dev(F_2)$  and    $\dev| F_1 \cap  F_2$   onto 
$\dev(F_1) \cap \dev(F_2)$.  Moreover,  $F_1 \cup F_2$ is an $n$-ball, 
and $F_1 \cap F_2$ is a convex $n$-ball.  
\end{prop} 
\begin{proof} See Proposition 3.9 in \cite{psconv}. 
\end{proof} 
 
\begin{rem}\label{rem:overlap} 
If $F_1$ and $F_2$ be convex open balls in $M_h$ and $M$ is affine, the 
above   theorem  was proved  by  Carri\`ere   in  Proposition 1.3.1 of 
\cite{Card}. The above theorem implies this result since  
any convex open $n$-ball in $M_h$ has the  closure which  
is a compact $n$-ball.  
\end{rem} 
 
The above proposition follows from: 
\begin{prop}\label{prop:extmap} 
Let  $A$ be a $k$-ball in  $\che M_h$\cO and   $B$ an $l$-ball.  Suppose 
that $A^o \cap B^o  \ne \emp$\cO $\dev(A)  \cap \dev(B)$ is a  compact 
manifold    in    $\SI^n$ with interior    equal   to  $\dev(A)^o \cap 
\dev(B)^o$\cO  and   $\dev(A)^o \cap \dev(B)^o$  is arcwise-connected. 
Then $\dev| A \cup B$ is a homeomorphism onto $\dev(A) \cup \dev(B)$.  
\end{prop} 
\begin{proof} See Proposition 3.10 in \cite{psconv}. 
\end{proof} 
 
\chapter{$(n-1)$-convexity: previous results } 
\label{ch:psconv} 
 
In  this  chapter, we summarize the facts  that  we  need from  
the monograph \cite{psconv}:   We  will define $m$-convexity  of 
real  projective    manifolds,  and  explain   that  the    failure of 
$(n-1)$-convexity of a  real projective manifold implies the existence 
of   $n$-crescents in the completion $\che   M_h$ of the holonomy cover 
$M_h$  of $M$.  Finally,   we explain the   decomposition of  $M$ into 
$(n-1)$-convex  manifolds  and  concave  affine  $n$-manifolds   along 
$(n-1)$-dimensional     manifold  convex  from    the   $(n-1)$-convex 
sides. An intermediate step of splitting  the real projective manifold 
along   two-faced $(n-1)$-manifolds, if they exist, is explained.  
(Two-faced  submanifolds are closed  submanifolds which    
arise from  collections of $n$-crescents meeting at their  
boundary exclusively.)  
 
Let $M$ be   a real  projective $n$-manifold  and  $M_h$ its  holonomy 
cover.  An {\em  $m$-simplex}\/ $T$ in $\che  M_h$ is a tame subset of 
$\che   M_h$  such  that $\dev|T$   is  an  imbedding  onto  an affine 
$m$-simplex in an affine patch in $\SI^n$.  
 
\begin{defn}  
We say that $M$  is {\em $m$-convex}, $0 <  m <  n$, if the  following 
holds:  If $T \subset \che M_h$ be an $(m+1)$-simplex with sides $F_1, 
F_2, \dots, F_{m+2}$ such that $T^o \cup F_2  \cup \cdots \cup F_{m+2} 
\subset M_h$, then $T \subset M_h$.  
\end{defn} 
 
\begin{lem}\label{lem:dconv} 
$M$ is    $m$-convex   if  and  only  if   $\che   M_h$    includes no 
$(m+1)$-simplex $T$ with a  side $F_1$ such that  $T \cap \hideal{M} = 
F_1^o \cap \hideal{M} \ne \emp$.  
\end{lem} 
\begin{proof} 
See Proposition 4.3 in \cite{psconv}. 
\end{proof} 
 
\begin{lem}\label{lem:dconv2} 
$M$ is $m$-convex if and only if for an $(m+1)$-simplex $T$ with sides 
$F_1, \dots, F_{m+2}$, every real projective immersion $f: T^o \cup F_2 
\cup \cdots \cup F_{m+2} \ra M$ extends to one from $T$.  
\end{lem}  
\begin{proof} 
See Proposition 4.2 in \cite{psconv}. 
\end{proof} 
 
It is easy to see that $i$-convexity implies $j$-convexity whenever  
$i \leq j < n$ (see Remark 2 in \cite{uaf}).  
A compact real projective manifold with totally geodesic  
boundary is convex if and only if it is $1$-convex 
(see Theorem A.2 of \cite{psconv}). 
 
For example, Hopf manifolds and half-Hopf manifolds are not 
$(n-1)$-convex. $\mathcal{E}_2$ is $(n-1)$-convex but not  
$(n-2)$-convex for $n \geq 3$ (see Example \ref{exmp:hopf}). 
 
A bihedron in $\che M_h$ is said to be an {\em $n$-crescent}\/ if one of 
its sides is  a subset of  $\hideal{M}$ and  the other side  is not. An 
$n$-hemisphere $A$ in $\che  M_h$ is said  to be an {\em $n$-crescent}\/ 
if an $(n-1)$-hemisphere  in $\delta A$ is  a subset of $\hideal{M}$ 
and the boundary itself is not a subset of $\hideal{M}$.  
{\em Bihedral} $n$-crescents are the former ones, and  
{\em hemispheric} $n$-crescents are the latter ones. 
 
Given a bihedral $n$-crescent $R$, the interior of the side in  
$\hideal{M}$ is denote by $\alpha_R$ and the other side  
by $\nu_R$. Given a hemispheric $n$-crescent $S$,  
the union of all open $(n-1)$-hemispheres in  
$\hideal{M} \cap \delta S$ is denoted by $\alpha_S$ and  
its complement in $\delta S$ by $\nu_S$. Note that $\alpha_S$  
is homeomorphic to a connected open manifold and $\nu_S$ is  
a tame set homeomorphic to an $(n-1)$-ball. 
(Note that we define $n$-crescents and their objects in $\che M$  
in the same manner.) 
 
Assume  that  $M_h$ is    not  projectively diffeomorphic  to an  open 
$n$-bihedron or an  open $n$-hemisphere.  We  proved in \cite{psconv} 
the following theorem:  
 
\begin{thm}\label{n-1conv} Let $M$ be a compact radiant affine manifold 
with empty or totally geodesic boundary. Suppose that $M$  is  
not $(n-1)$-convex.  Then  $\che M_h$ includes an $n$-crescent.  
\end{thm}  
 
A {\em nice}\/ $(n-1)$-submanifold of $M$ is a totally geodesic closed 
real projective $(n-1)$-submanifold  in $M^o$ each component of  which 
is projectively diffeomorphic to a quotient of a connected open subset 
of   an affine patch  $\bR^{n-1}$ in   the  real projective space  
$\bR P^{n-1}$ by a  properly discontinuous  and  free action of   
a group of projective transformations.  
 
If two hemispheric $n$-crescents $R$ and $S$ overlap, then they  
are equal. If $R$ and $S$ meet but do not overlap, then $R\cap M_h$  
and $S\cap M_h$ meet in common components of $\nu_R \cap M_h$  
and $\nu_S \cap M_h$. Such components are called a {\em copied} 
components, and the union of all such components for 
all hemispheric $n$-crescents is a properly imbedded submanifold 
in $M_h$ and hence covers a closed submanifold in $M^o$, 
which is said to be a {\em two-faced submanifold arising from  
hemispheric $n$-crescents}. 
 
We say that $R \sim  S$ for two  bihedral $n$-crescents $R$ and $S$ if 
$R$ and  $S$ overlap. This generates  an equivalence relation which we 
denote by $\sim$ again. We defined $\Lambda(R)$ to be $\bigcup_{S \sim 
R}     S$    and  discussed   the   properties    of  $\Lambda(R)$  in 
\cite{psconv}. (See Chapter \ref{ch:tdim} discussing this  
in a more special setting.)  
The so-called two-faced $(n-1)$-submanifold  
arising from bihedral $n$-crescents is  
a properly imbedded submanifold in $M$ covered by the submanifold 
of $M_h$ that is the union of common components  
of the topological boundaries $\Bd \Lambda(R) \cap M_h$ and  
$\Bd \Lambda(S) \cap M_h$  
for all pair of non-equivalent $n$-crescents $R$ and $S$. 
Of course, such a set may be empty and a two-faced  
submanifold doesn't exist for our real projective manifold $M$. 
 
A two-faced submanifold is  a   nice  submanifold.  We  will     
repeat  the construction in  \cite{psconv}   in  Chapter \ref{ch:tdim}    
in  the radiant affine case, as needed by this paper.  
 
\begin{defn}\label{defn:cafm} 
A {\em concave affine $n$-manifold $M$ of type I}\/ is a real projective 
manifold  such  that   $M_h$ is   a  subset   of  a  hemispheric 
$n$-crescent  in $\che{M_h}$.  A  {\em concave affine $n$-manifold $M$ 
of type  II}\/ is a  real  projective  manifold with   possibly concave 
boundary such that  $\che  M_h$ includes no  hemispheric $n$-crescents 
and $M_h$ is a subset of $\Lambda(R)$ for a bihedral $n$-crescent 
$R$ in $\che M_h$.  
\end{defn}  
 
\begin{thm}\label{thm:mainpsII}  
Suppose that  $M$ is compact  and $\che M_h$ includes  
a hemispheric $n$-crescent. 
Then $M$ includes a   compact concave  affine $n$-submanifold $N$ 
of type I  or $M^o$ includes the two-faced $(n-1)$-submanifold 
arising from hemispheric $n$-crescents.  
\end{thm}  
 
\begin{thm}\label{thm:mainps} 
Suppose that $M$ is compact and $\che M_h$ includes  
bihedral $n$-crescents but no hemispheric $n$-crescents.  
Then $M$ includes a   compact concave  affine $n$-submanifold $N$ 
of type II  or $M^o$ includes the two-faced $(n-1)$-submanifold 
arising from bihedral $n$-crescents. 
\end{thm} 
 
A real projective manifold has {\em convex  boundary} if each of its 
boundary points has a neighborhood where the chart restricted to it is 
an imbedding onto a convex set. 
A real projective manifold has {\em concave  boundary} if each of its 
boundary points has a neighborhood where the chart restricted to it is 
an imbedding onto the complement of an open convex set inside a simply 
convex open  ball.   
 
The  manifold interior of   a concave  affine   
submanifold  always admits a compatible affine structure  
(see \cite{psconv}).  
 
The following corollary was stated and proved in \cite{psconv} 
in a slightly different form. The proof in Chapter 10 
of \cite{psconv} easily shows:  
\begin{cor}\label{cor:spl1}  
Suppose that $M$ is compact but not $(n-1)$-convex.  Then  
\begin{itemize} 
\item  after  decomposing $M$ along  the two-faced $(n-1)$-submanifold 
$A_1$  arising  from hemispheric   $n$-crescents{\rm ,}  the resulting 
manifold $N$  decomposes into concave affine manifolds  of type  I and 
real  projective $n$-manifolds   with  totally geodesic   boundary the 
Kuiper completion of whose holonomy  cover includes no hemispheric 
$n$-crescents.  
\item We let $N^1$ be the disjoint union of the resulting manifolds of 
the decomposition other than concave  affine  ones of type I.   
After splitting $N^1$ along the two-faced $(n-1)$-manifold $A_2$  
arising from bihedral $n$-crescents{\rm ,} the  resulting manifold  
$N^{1 \prime}$ decomposes into  concave affine manifolds  of type II  
and real projective $n$-manifolds with convex boundary which is  
$(n-1)$-convex. The Kuiper completions of the holonomy covers 
of $(n-1)$-convex pieces include no $n$-crescents.  
\end{itemize} 
\end{cor}  
 
Note that  $A_1$ and $A_2$ could be  empty.  If $A_1  = \emp$, then we 
define  $N = M$  and if $A_2 = \emp$,  then  we define $N^{1 \prime} = 
N^1$.  The  proof depends  on the following  proposition  which we will be 
using later.  
 
\begin{rem}\label{rem:boundary} 
The manifolds $A_1$ and $A_2$ are nice. So are the boundary components   
of a concave affine manifold of type I since the manifold-boundary  
points of the universal cover of a concave affine manifold of type I 
lie in $\nu_R$ for a bihedral $n$-crescent $R$. 
\end{rem} 
 
\begin{prop}\label{prop:lfin} 
Let $A$ be a submanifold and a closed subset  of a regular cover $\hat 
N$ of a compact $n$-manifold $N$ so  that for each deck transformation 
$\vth$, $\vth(A)$ either equals $A$ or is disjoint  from $A$.  Suppose 
that  the  collection composed  of  elements $\vth(A)$   for $\vth \in 
\Aut_N(\hat  N)$  for the deck-transformation group  $\Aut_N(\hat N)$ is 
locally finite, i.e., each  point of $\hat N$  has a neighborhood that 
meets only finitely many members of the collection.  
Then $p| A: A \ra p(A)$ for the covering map  $p: \hat N  \ra N$  
is a covering  map onto a compact manifold $A$.  Furthermore,  
if  $A$ is of  dimension $< n$, and is  an open manifold, then  
$p(A)$ is a closed submanifold of same dimension. 
If $A$ is of dimension $n$, and has boundary $\delta A$, then $p(A)$ 
is a compact submanifold  of  dimension $n$ with  boundary  
$p(\delta A)$.  
\end{prop} 
\begin{proof} 
The readers  can easily supply  the proof.  
\end{proof}

\chapter[Radiant vector fields]
{Radiant vector fields, generalized affine suspensions,  
and the radial completeness}  
\label{ch:racom} 
 
We  will discuss general facts  about  radiant affine $n$-manifolds in 
this chapter: Given a radiant affine manifold $M$, the radiant 
vector field in $\brno$ induces a radiant vector field on $M_h$  
and $M$. As examples of radiant affine $n$-manifolds,  
we define generalized affine  suspensions over  real  projective surfaces and 
orbifolds (see Carri\`ere \cite{Carlet}  and Barbot \cite{barbot1}).  
We show that a radiant affine manifold is a generalized affine suspension  
if and only if it admits a total cross-section to the radial flow. 
The results of Barbot \cite{barbot1} that we need are stated here. 
 
From an estimation of the radial flow in terms of $\bdd$, 
it follows that the radial flow extends to the Kuiper completion as well,   
and the point on $\che{M_h}$ corresponding to the origin  under $\dev$  
is  unique. The ideal sets are also radial-flow-invariant,  
and the completion of the holonomy cover of a radiant affine  
$n$-manifold, if not $(n-1)$-convex, must  contain a radiant  
$n$-crescent,  i.e.,  a radiant $n$-dimensional bihedron.  
Finally, we say about an intersection property of two radiant  
$n$-bihedra. 
 
Let $e_0, e_1, \dots, e_n$ be the orthonormal basis of $\bR^{n+1}$ and 
$x_0,  x_1, \dots, x_n$   the associated coordinates, and  $\SI^n$ the 
unit  sphere in  $\bR^{n+1}$ with origin  $O$.  Let  $\phi$ denote the 
radial  projection     from $\bR^{n+1} -     \{O\}$ to   $\SI^n$.  For 
convenience, we identify $\bR^n$ with the interior $\mathcal{H}^o$ of the 
upper hemisphere $\mathcal{H}$ in $\SI^n$  by sending $x  \mapsto x + e_0$ 
and  then  sending it  to $\SI^n$ by   $\phi$.  This identification is 
compatible with  standard  affine and  real projective  structures  on 
$\bR^n$ and $\mathcal{H}^o$ as  given in the introduction, i.e., geodesic 
structures are  preserved. (Note that  the origin  $O$ of $\bR^n$   
corresponds to $e_0$.)    
 
We define the radiant vector  field $\bV$  on  $\bR^n$ to be given  by 
$\sum_{i=1}^n  x_i \partial/\partial x_i$,  where $\bV(O)$ is the zero 
vector.    
The  vector field $\bV$  on  $\mathcal{H}$ generates a flow. It is easy  to  
compute that $\bV$  is a vector field on $\mathcal{H}$ of bounded  
$\bdd$-length and $\bV$ extends to a continuous vector field on $\SI^n$  
by a zero vector field on $\SI^n - \mathcal{H}$.

Since $\bV$  is  a  $\bdd$-bounded vector  field,  the flow $\Phi'_t$ on   
$\SI^n$ generated by $\bV$ satisfies the  following inequality:  
\begin{equation}\label{eqn:flowquasi}  
C(t)^{-1} \bdd(x, y)  \leq \bdd(\Phi'_t(x), \Phi'_t(y)) \leq  
C(t)\bdd(x, y)  
\end{equation} 
for $x, y \in \mathcal{H}^o - O$ and a positive constant $C(t)$  
smoothly depending only on $t$. (See Lemma 3 of Section 2.1 of  
Abraham-Marsden\cite{AM}.) 
 
Let $M$ be a compact radiant affine $n$-manifold with a developing map 
$\dev:M_h \ra \bR^n$ and the holonomy homomorphism $h: \pi_1(M)/\pi_1(M_h) 
\ra   \Aff(\bR^n)$ and the holonomy  group $\Gamma$ equal to $h(\pi_1(M))$,  
fixing the origin $O$. As $\bV$  on  $\mathcal{H}$  is $\Gamma$-invariant,   
the  induced vector field  $\bV_h$  in  $M_h$   by  $\dev$ is invariant  
under the  deck transformations; hence,  
there exists an induced  vector field  
$\bV_M$ on  $M$ called the {\em radiant vector field}\/ of $M$,  
which defines a flow $\Phi: \bR \times  M \ra M$.  
Since  $M$  is a closed manifold, $\Phi$ is a complete flow.
The vector field $\bV_h$ induces the flow  
\[\Phi_h: \bR \times M_h \ra M_h\]   
so that the following diagram commutes for every $t\in \bR$  
\begin{eqnarray}\label{eqn:flow}  
M_h & \stackrel{\Phi_{h, t}}{\longrightarrow}  & M_h \nonumber \\ \dev 
\downarrow       &     &         \downarrow      \dev    \\      \brno 
&\stackrel{\Phi'_t}{\longrightarrow} & \brno. \nonumber  
\end{eqnarray}   
A {\em radial line}\/  in $\brno$ is  a  component of a complete  line 
with $O$ removed.  
 
\begin{lem}\label{lem:radcom}  
The images $\dev(M_h)$ and $\dev(\tilde M)$  miss  $O$, and  
for   each point $x$ of $M_h$, there exists a unique imbedded   
geodesic $l$, maximal in $M_h$, passing through $x$  such  that  
$\dev| l$   is a diffeomorphism  to a radial line in $\bR^n$.  
\end{lem} 
\begin{proof}  
The first statement is proved by Theorem 3.3  of \cite{FGH}.  The last 
statement follows from the completeness of the flow $\Phi_h$.  
\end{proof} 
In the  above  theorem, the line $l$  is  said to  be the  radial line 
through $x$.  
 
We  now define generalized affine  suspensions over   
real  projective manifolds and 
orbifolds (see Carri\`ere \cite{Carlet}  and Barbot \cite{barbot1}).  A 
{\em real projective orbifold}\/  is  simply an orbifold  with geometric 
structure modeled  on $(\rp^n,  \PGL(n+1, \bR))$. (For  definitions of 
orbifold and geometric structures on orbifolds, see Ratcliffe \cite{RT}.)

Let $\Sigma$ be a compact real projective $(n-1)$-manifold  
with empty or totally geodesic boundary with a projective  
automorphism $\phi$. Since $\SII$ is  
a real projective $(n-1)$-sphere, there exists a 
developing  map $\dev:\tilde \Sigma \ra  \SII$ and  a holonomy homomorphism 
$h:   \pi_1(\Sigma)  \ra   \Aut(\SII)$.  Choosing   an 
arbitrary Euclidean metric in  $\bR^{n}$, we define an immersion $\dev': 
\Sigma \times \bR  \ra \bR^{n}$ by simply   mapping $(x,   t)$  to 
$e^t u(x)$ where $u(x)$ is the unit vector at the origin in the  
direction of $\dev(x)$ in  $\bR^{n}$. Since $\Aut(\SII)$ can be identified 
with $\SL_\pm(n, \bR)$, there is a natural quotient map  
$\GL(n, \bR) \ra \Aut(\SII)$; we choose any lift 
$h': \pi_1(\Sigma) \ra \GL(n, \bR)$ of $h$, and define 
a corresponding action of $\pi_1(\Sigma)$ on $\tilde \Sigma 
\times \bR$ by $\vth(x, t) = (\vth(x), t + \log||h'(\vth)(u(x))||)$.  
 
Note that there is a unique lift $h'':\pi_1(\Sigma)  
\ra \SL_\pm(n, \bR)$ and $h'(\vpi) = k(\vpi)h''(\phi)$ for  
a homomorphism $k: \pi_1(\Sigma) \ra \bR^+$. 
Clearly $k$ induces a homomorphism $k_\#: H_1(\Sigma) \ra \bR$  
by taking logarithms. We choose $h'$ to satisfy 
\begin{equation}\label{eqn:kaut} 
k_\# \circ \phi_\# = \phi_\# 
\end{equation}  
for induced homomorphism $\phi_\#: H_1(\Sigma) \ra H_1(\Sigma)$ 
of $\phi$.  
 
Letting $\tilde\Sigma\times\bR$ have the  affine structure induced from 
the   immersion  $\dev'$, we see  that $\pi_1(\Sigma)$ 
defines a properly discontinuous  and free affine action  of $\Sigma 
\times \bR$ preserving   each fiber  homeomorphic  to  $\bR$, and  the 
quotient space is homeomorphic to $\Sigma \times \bR$, i.e., a trivial 
$\bR$-fiber bundle over $\Sigma$.  We identify the quotient space with 
$\Sigma\times  \bR$, and  choose a  section   $s: \Sigma \ra  \Sigma 
\times \bR$ so that $s(\Sigma)$ becomes a compact imbedded surface.  
 
The projective automorphism 
$\phi$ lifts  to a  projective   automorphism $\tilde \phi$  
of $\tilde \Sigma$. Since $\tilde \phi$ is a projective  automorphism,  
there exists an element $\rho$ 
in $\Aut(\SII)$ satisfying $\dev \circ \tilde \phi =  \rho \circ \dev$.   
We may choose any element $\rho'$ of  $\GL(n, \bR)$ which induces $\rho$, 
and  $\rho'$  defines  an affine  automorphism  $\rho''$ of   
$\tilde \Sigma \times   \bR$    given    by   
$\rho''(x, t) = (\tilde \phi(x),  t  + \log||\rho'(u(x))||)$.   
 
Given a  deck transformation  $\vth$ of 
$\tilde \Sigma$,   there  exists a  deck transformation  $\vpi$ 
satisfying $\phi \circ \vth  = \vpi \circ  \phi$.  
$\vpi$ is to be denote by $\phi^*(\vth)$ and  
$\phi^*:\pi_1(\Sigma) \ra \pi_1(\Sigma)$ defines a group 
automorphism. By equation \ref{eqn:kaut}, we see that 
$\rho'' \circ \vth = \vpi \circ \rho''$ on $\tilde \Sigma\times \bR$. 
Therefore, it follows  that $\rho''$ induces an  affine automorphism  
$a_{\rho'}$ of $\Sigma \times \bR$.  
 
We let $e^r  \Idd$ for $r \in \bR$ denote the  
{\em  dilatation}\/ multiplying each vector in $\brto$   
by  a factor  $e^r$, which induces  an affine automorphism $D_r$  
on $\tilde \Sigma \times \bR$ also called a {\em dilatation}\/  
given by $D_r(x,  t) = (x,  t + r)$.  Since $D_r$ commutes 
with any deck transformation of $\tilde \Sigma \times \bR$, it follows that 
$D_r$  defines  an  affine  fiber-preserving   automorphism $D'_r$  of 
$\tilde \Sigma \times \bR$.  
 
Any other choice  $\rho'_1$ in $\GL(n, \bR)$ of $\rho'$  equals  
$e^r \Idd \circ  \rho'$  for some $r$, and given the  affine  
automorphism $a_{\rho'_1}$ of $\Sigma \times \bR$ corresponding to   
$\rho'_1$, we see that $a_{\rho'_1}$ equals $D_r \circ a_{\rho'}$ for   
some $r$.  Hence by  choosing $r$ sufficiently large  $>1$, and   
positive, we can make 
$a_{\rho'_1}(s(\Sigma))$  and   $s(\Sigma)$          disjoint      and 
$a_{\rho'_1}(s(\Sigma))$ to lie   in the radially  outer-component  of 
$\Sigma   \times  \bR  -   s(\Sigma)$.     Since $a_{\rho'_1}$  is   a 
fiber-preserving  diffeomorphism, $a_{\rho'_1}(s(\Sigma))$ is  
another total cross-section.  
We let $N$ denote the compact $(n+1)$-manifold in  
$\Sigma \times \bR$ bounded  by $s(\Sigma)$ and  
$a_{\rho'_1}(s(\Sigma))$,  and identify $s(\Sigma)$  and  
$a_{\rho'_1}(s(\Sigma))$ by $a_{\rho'_1}$ to 
obtain a compact   radiant  affine $n$-manifold homeomorphic to  
the mapping torus $\Sigma  \times_\phi  \SI^1$,  i.e.,  
$\Sigma \times   I/\sim$ where $\sim$ is defined by  
$(x, 0)  \sim (\phi(x), 1)$.  We call 
the  resulting affine $n$-manifold  the {\em generalized affine suspension  
over}\/ $\Sigma$ {\em using}\/ the projective automorphism $\phi$.  
(Note the term ``affine suspension'' is reserved for the case  
$\Sigma$ and $\phi$ are both affine.) 
 
If we let $\phi$  be the identity  automorphism  of $\Sigma$,  then  
the generalized affine suspension is a so-called  Benz\'ecri suspension.  
But even  when 
$\phi$  is  of finite  order,  we  will also   call  it {\em  Benz\'ecri 
suspension}\/ over the  projective  orbifold  $\Sigma/\langle\phi\rangle$   
in this paper (sometimes, over the manifold $\Sigma$ also).
 
\begin{prop}\label{prop:Benzrad}  
A generalized affine suspension of a real projective $(n-1)$-manifold  
with totally geodesic or empty boundary is a radiant affine  
$n$-manifold with totally geodesic or empty boundary.   
\end{prop} 
\begin{proof} Straightforward. 
\end{proof} 
 
When $n = 3$, the boundary components are homeomorphic to tori  
or Klein bottles since they are tangent to the radial flow, and hence  
have zero Euler characteristic. 
 
\begin{exmp}\label{exmp:affsusp}  
Let $T$ be an $(n-1)$-dimensional Hopf manifold given  
as the quotient of $\bR^{n-1} - \{O\}$ by the cyclic group  
generated by $g$ sending $x \ra 2x$ for each vector $x$.  
Clearly, $\mathcal{E}_2$ is a Benz\'ecri suspension of $T$. 
 
Given a closed projective surface $\Sigma$, we can obtain  
Benz\'ecri suspensions easily. When $\Sigma$ is a quotient  
of a standard ball in $\bR P^2$; i.e., $\Sigma$ has a projective 
structure induced from hyperbolic structure, then  
the Benz\'ecri suspension can be obtained from  
the upper part of the interior of the null cone  
in the Lorentzian space $\bR^{1, n-1}$ by an action of 
a group generated by linear Lorentz transformations and  
a homothety, i.e.,  
a linear map of form $s\Idd$ for $s > 0$. In this case,  
our Benz\'ecri suspension carries Lorentzian flat conformal  
structure as well.  
 
Given an arbitrary closed real projective surface $\Sigma$ 
of negative Euler characteristic,
$\Sigma$ decomposes along disjoint simple closed geodesics  
into convex real projective surfaces and annuli with geodesic 
boundary (see \cite{cdcr1}, \cite{cdcr2}, and \cite{cdcr3}).  
Thus the Benz\'ecri suspension also decomposes  
along corresponding totally geodesic tori into convex  
radiant affine manifolds and Benz\'ecri suspensions of the annuli. 
\end{exmp} 
 
\begin{thm}\label{thm:Sorb} Let $\Sigma$ be a compact real projective  
surface with totally geodesic boundary. If the Euler characteristic of  
$\Sigma$ is negative, then  a generalized affine suspension   
of $\Sigma$ using an automorphism $\phi$ is homeomorphic to  
a Seifert space with zero Euler number with base orbifold  
homeomorphic to $\Sigma/\langle \phi\rangle$  
where $\phi$ is a finite order automorphism. Moreover,  
the generalized affine suspension has a finite cover which is  
a Benz\'ecri suspension over $\Sigma$.  
\end{thm} 
\begin{proof} 
Since by following Theorem  \ref{thm:finite},  $\phi$  is of  finite 
order, it follows from the fact that $M$ is a mapping torus of $\rho$  
over $\Sigma$ that $M$ is homeomorphic to a Seifert manifold over  
$\Sigma/\langle \phi\rangle$. 
 
Since a finite power of $\phi$ is the  identity  automorphism of 
$\Sigma$, it follows that $M$ is finitely   covered  by a Benz\'ecri 
suspension, which is homeomorphic    to a trivial circle  bundle  over 
$\Sigma$. Since the Euler number of the trivial circle bundle is zero, 
the conclusion follows.  
\end{proof}

A projective automorphism of $\SI^2$ is said to be {\em hyperbolic} 
if it is represented by a diagonal matrix with distinct positive 
eigenvalues in $\GL(3, \bR)$. It is said to be {\em quasi-hyperbolic}  
if it is represented by a nondiagonal matrix with two  
distinct positive eigenvalues (see \cite{cdcr2} for details). 
Let $S$ be a compact real projective surface with empty or  
geodesic boundary and negative Euler characteristic. 
Let $(\dev, h)$ be its development pair. 
A {\em strong tight curve} in a real projective surface $S$  
is a simple closed geodesic $\alpha$ such that its lift  
$\tilde \alpha$ in $\tilde S$ is a simply convex line and $h(\vth)$ for  
the deck transformation $\vth$ corresponding to $\tilde \alpha$ and  
$\alpha$ is hyperbolic or quasi-hyperbolic, and  
$\dev(\tilde \alpha)$ is a geodesic connecting the fixed point  
of the largest eigenvalue to that of the smallest eigenvalue of  
$h(\vth)$. 
 
\begin{rem} 
If the  Euler characteristic of $S$ is  zero, then the generalized affine 
suspension of $\Sigma$ is  not necessarily Seifert. Actually there are 
many nontrivial automorphisms which are not finite order  
up to isotopy.  
\end{rem} 
 
\begin{thm}\label{thm:finite} 
Let $\rho: S  \ra S$ be  a projective automorphism of a compact real 
projective surface $S$ of  negative  Euler characteristic.   Then 
$\rho$ is of finite order.  
\end{thm} 
\begin{proof} 
A {\em purely convex} real projective surface $L$ is a compact  
convex surface with totally geodesic or empty boundary  
and negative Euler characteristic such that it includes  
no compact annulus with geodesic boundary freely homotopic 
to a boundary component of $L$. 
The  real projective   surface  $S$  canonically decomposes  
into maximal purely  convex real projective  subsurfaces of  negative  
Euler characteristics and maximal  annuli  or M\"obius bands   with  
geodesic boundary  \cite{cdcr2}, \cite{cdcr3}.  By uniqueness of the  
decomposition in \cite{cdcr3}, $\rho$  acts on the union 
of the purely convex surfaces which have negative Euler characteristics 
and are   mutually disjoint.  Letting  $S'$  denote  the disjoint 
union, we see  that each boundary component of  $S'$ is a  strong 
tight curve by Proposition 4.5 in \cite{cdcr2}.  
 
Assume that $\rho$ is of infinite order.  By taking  a power of $\rho$ 
if necessary, we assume that $\rho$  acts on  every  component and 
boundary component of $S'$.  We assume without loss of generality 
that $S'$ is connected. Since $\rho$ acts on a boundary component 
$k$ of  $S'$, $\rho$  lifts to  an automorphism  $\tilde \rho$ of 
$\tilde S'$  acting   on the  image $\tilde  k$   of  a lift   of 
$k$. Since $\dev\circ \tilde \rho$ is another projective immersion 
$\tilde S' \ra \SI^2$, it follows that $\dev\circ \tilde \rho = 
h(\tilde \rho) \circ \dev$ for an element $h(\tilde \rho)$ of  
$\Aut(\SI^2)$. 
 
Letting $\vth$ be the deck transformation of $\tilde S'$ 
corresponding to $\tilde k$ and  $k$, we see  that $\vth$ and   
$\tilde \rho$ commute as $\rho$ is a  topological automorphism  
of $S'$. Hence given a development pair $(\dev, h)$ of $\Sigma'$,  
$h(\vth)$ and $h(\tilde \rho)$ act on the convex  domain $\dev(S')$  
and an open geodesic $\dev(\tilde  k)$.  The boundary  of  $\dev(S')$  
is the union of the  closure  of $\dev(\tilde  k)$ and  an open  arc  
$\alpha$ connecting the endpoint of $\dev(\tilde k)$.  
Since $S'$ is convex and has a negative Euler characteristic, 
$h(\vth)$ is either hyperbolic or quasi-hyperbolic by 
Proposition 4.5 of \cite{cdcr2}. The obvious $h(\vth)$-invariant  
subsets are simply convex triangles with vertices fixed points of 
$h(\vth)$ when $h(\vth)$ is hyperbolic and are lunes when  
$h(\vth)$ is quasi-hyperbolic (see Figures 2 and 3 of \cite{cdcr3}). 
 
The closure of $\dev(S')$ is an $\langle h(\vth) \rangle$-invariant set. 
Since $k$ is strong, we can show that $\alpha$ is a subset of an open  
$\langle h(\vth)\rangle$-invariant triangle or lune as above. 
Since both $h(\vth)$ and $h(\tilde \rho)$ act on the common arc  
$\alpha$ in the open invariant triangle or lune, we can  
verify by explicit calculations that  
both $h(\vth)$ and $h(\tilde \rho)$ lie in a common  
$1$-parameter  subgroup $H$ of $\Aut(\SI^3)=\SLt_\pm(3, \bR)$. 
 
For any regular neighborhood of $\dev(\tilde k)$ in   
$\dev(\Sigma')$, there exists an $H$-invariant arc  
$\tilde  \alpha_k$  connecting the endpoints  of $\dev(\tilde k)$. 
Now, $\tilde \alpha_k$  corresponds to a   simple   closed arc   in 
$\alpha_k$ freely  homotopic  to $k$  in   $\Sigma^{\prime o}$. If  we 
choose $\alpha_k$ for each component  $k$, then there exists a compact 
surface $\Sigma''$ in $\Sigma^{\prime  o}$ bounded by $\alpha_k$s; and 
$\rho$  acts on $\Sigma''$.  As  $\Sigma''$ is a compact subsurface in 
$\Sigma^{\prime  o}$, it is  a compact metric  space under the Hilbert 
metric  of $\Sigma^{\prime o}$. Since   the  group of isometries of  a 
compact metric space is compact, the closure  of  
$\langle\rho\rangle$ is a compact Lie group acting on $\Sigma''$.  
If the closure of $\langle\rho\rangle$ is discrete, then 
$\rho$ is of finite order. Since the identity component  of 
this Lie group thus does not consist of a point,  
the Lie group includes the compact  Lie  subgroup $\SI^1$   
homeomorphic  to  a  circle  acting on $\Sigma''$ nontrivially.   
Since $\Sigma''$  is homotopy equivalent  to 
$\Sigma'$, this is absurd as there is no nontrivial $\SI^1$-action on 
a surface  of negative Euler  characteristic.  Therefore $\rho$  is  
of finite order.  
\end{proof} 
 
\begin{prop}\label{prop:crssect}  
A compact radiant  affine  $n$-manifold $M$ with  totally geodesic 
or empty boundary admits a total cross-section $\Sigma$ to the radial  
flow if and only if it is affinely  diffeomorphic to a generalized  
affine suspension over a compact real projective $(n-1)$-manifold  
$\Sigma'$  with totally geodesic or empty boundary. Moreover,  
in the above case, $\Sigma$ with the induced real projective  
structure from the affine space by the radial flow is real  
projectively diffeomorphic to $\Sigma'$.  
\end{prop} 
\begin{proof}  
Suppose that  $M$ has a  total cross-section $\Sigma$ to  the  radial 
flow.  Then  $\Sigma$ obviously has  a real projective structure since 
the radial flow gives a  local projection of the  affine space to  the 
real projective space. The affine    transition maps preserving    the 
radial flow now become real projective transition maps.  
 
The existence  of the total cross-section shows that there   
exists a cover $M'$ of $M$ homeomorphic to $\Sigma \times \bR$.  
Let $\dev: \tilde M' \ra \brno$ for the universal cover  
$\tilde M'$ of $M'$ denote the developing  map  
and $h:\pi_1(M') \ra \GL(n, \bR)$ the holonomy homomorphism,  
considering $\GL(n, \bR)$ as the subgroup of $\Aut(\mathcal{H})$  
acting on $\mathcal{H}$ fixing $O$ where $\mathcal{H}$ 
is the $n$-hemisphere  
whose interior is identified with $\bR^n$. 
 
The inverse image $\tilde \Sigma$ of $\Sigma$ in $\tilde M'$ is  
obviously the universal cover of $\Sigma$. 
$\tilde M'$  is   obviously   homeomorphic to   
$\tilde \Sigma \times  \bR$, and $\pi_1(\Sigma)$ can be identified with 
with $\pi_1(M')$ since the deck transformations of $M'$  
act on $\tilde \Sigma$ as well. 
The  radial projection map $\Pi: \brno \ra \SII$  
composed with $\dev|\tilde \Sigma$ is the developing map for 
$\Sigma$  considered  now as    a real  projective manifold,  and    
$k: \pi_1(\Sigma)  \ra \Aut(\SII)$, which is obviously obtained from   
$h: \pi_1(\Sigma) \ra \GL(n, \bR)$ by  restriction to $\SII$   
is  the holonomy homomorphism for $\Sigma$.  
 
From  the real projective $(n-1)$-manifold $\Sigma$, we will now  
construct an affine  structure for $\Sigma \times \bR$.   
We  may   construct  a developing map  
$\dev':\tilde \Sigma \times \bR \ra \brno$  
as in  the construction of  the generalized affine suspension  using  
$\Pi\circ \dev|\tilde \Sigma$ by   radial    extension.   
Letting    $h'$   denote   the  composition 
\[\pi_1(\Sigma \times  \bR)    \stackrel{i}{\ra}     \pi_1(\Sigma) 
\stackrel{k}{\ra} \Aut(\SII),\] where  $i$  is  the    obvious 
identification, we choose the lift  $\rho$ of $h'$ into $\GL(n,\bR)$ 
so that $\rho$ equals above $h$ on $\pi_1(\Sigma)$. Then we 
define  the action   of  the deck transformation $\vth$ in  
$\pi_1(\Sigma \times \bR)$  on  $\tilde \Sigma\times \bR$ by 
$\vth(x,  t) =(i(\vth)(x), t+ \log||\rho(\vth)(u(x))||)$    
where $u(x)$ is the  unit vector in $\bR^n$ in the direction  
of  $\dev(x)$.  Denote by $G$ the group of deck transformations  
on $\tilde \Sigma\times \bR$ thus defined.  
 
Since $\dev'$ and $\dev$ are  radially   complete  
maps (see Lemma \ref{lem:radcom}),   
it  is easy to  see that $\dev|\tilde M'$  lifts to   
an affine diffeomorphism  $k:\tilde M'  \ra \tilde \Sigma\times \bR$,  
so  that $\dev'\circ k =  \dev$. Then by our choice 
of $\rho$ and an analytic continuation from a small open set,  
we see that $k \circ \vth = \vth' \circ k$ where $\vth$ and 
$\vth'$ are corresponding deck transformations  of $M'_h$ and $\Sigma' 
\times \bR$.  This  means that $M'$  and $(\Sigma' \times  \bR)/G$ are 
affinely diffeomorphic.  Since  $M$  is obtained  by an  action of  an 
infinite cyclic automorphism group  from $M'$, transferring the affine 
action of the cyclic group to $(\tilde \Sigma\times \bR)/G$ shows that  
$M$ is affinely diffeomorphic to a generalized affine suspension  
over $\Sigma'$.  
 
The converse part is obvious since a generalized affine suspension  
over a surface $\Sigma$ using an automorphism $\gamma$ is   
a  mapping torus  and is homeomorphic   to a bundle over $\SI^1$   
with fibers diffeomorphic to $\Sigma$ with  fiber map induced  
from  that of  the projection $\Sigma \times \bR \ra \bR$;  
a fiber is a total cross-section obviously.  
\end{proof}  
 
Let us state Barbot's results \cite{barbot1}  
which we will use in this paper, telling us  
when a radiant affine $3$-manifold is a generalized affine suspension. 
A closed orbit of the radial flow of a radiant affine $3$-manifold  
is of {\em saddle type}\/ if the differential of the return map  
associated with the orbit and the local section has two 
distinct eigenvalues one of which is greater than $1$ and  
the other less than $1$. 
 
\begin{thm}[Barbot]\label{thm:barbot}  
Let $M$ be a closed radiant affine $3$-manifold.  
If one of the following holds, then $M$ admits a total cross-section. 
\begin{itemize}  
\item $M$ includes a totally geodesic surface tangent to 
the radial flow. 
\item $M$ is a convex affine manifold. 
\item $M$ has a closed orbit of the radial flow that is  
not of saddle type. 
\end{itemize}  
\end{thm}  
 
\begin{thm}[Barbot-Choi (Appendix \ref{app:radd})]
\label{thm:barbot2} 
If a compact radiant affine $3$-manifold $M$ has nonempty totally 
geodesic boundary, and each boundary component is affinely 
homeomorphic to a quotient of a convex cone in $\bR^2 -\{O\}$ or  
$\bR^2 -\{O\}$ itself by a group of affine motions,  
then $M$ admits a total cross-section. 
\end{thm} 
A totally geodesic submanifold of an affine $3$-manifold carries  
a natural induced affine structure as $2$-manifolds. 
 
Actually, we can dispense with the assumption on boundary  
components using Barbot's results or from our proof of the Carri\`ere 
conjecture. But we do not do so since it requires far more  
work. 
 
We will now begin to discuss about the ideal sets of  
the Kuiper completions of the holonomy covers of radiant affine  
manifolds. Since $M$ is  a radiant affine  manifold, we  divide  
$\hideal{M}$ into two  parts.  We denote  by $\hIideal{M}$   
the set $\dev^{-1}(\SII)$ in $\che M_h$, and is, by the following  
lemma, a subset  of $\hideal{M}$. We let  
$\hfideal{M} = \hideal{M} \cap \dev^{-1}(\brno)$. 
The two are said to be {\em infinitely}\/ and {\em finitely}\/ ideal  
subsets of $\che{M_h}$ respectively.  
 
\begin{lem} Let $A$ be a subset of $\che M_h$ such that  
$\dev(A)$ is a subset of $\SII$. Then $A$ is a subset of $\hideal{M}$.  
\end{lem} 
\begin{proof}  
Since $\dev(M_h)$ is the subset of $\mathcal{H}^o$, $A$ is necessarily a 
subset of $\hideal{M}$.  
\end{proof} 
 
\begin{prop}\label{prop:Hopf}  
Assume that $n \geq 3$. 
If $\hfideal{M} = \emp$, $M$ is finitely covered by  
a Hopf manifold or a half-Hopf manifold. 
In fact, $M$ is a generalized affine suspension of a real projective  
$(n-1)$-manifold with a universal cover homeomorphic to  
$\SI^{n-1}$ or an $(n-1)$-hemisphere. 
\end{prop} 
\begin{proof}  
Assume that $n \geq 2$. 
We start with the case when $M$ is a closed manifold. 
By Lemma \ref{lem:radcom}, $\dev(M_h)$ is a subset of  
$\bR^n - \{O\}$. Let $\gamma:I = [0, 1] \ra \bR^n - \{O\}$ be  
a $\bdd$-rectifiable path with endpoints $p = \gamma(0)$ and  
$q =\gamma(1)$, where  
$p$ belongs to $\dev(M_h)$. Since $\hfideal{M}$ is empty,  
we see that we can lift $\gamma$ starting from $0$ and  
can continue without any problem. Thus, it follows that  
$\dev:M_h \ra \bR^n-\{O\}$ is a covering map.  
 
Now assume that $n\geq 3$. In this case $\dev:M_h \ra \brno$ is  
a diffeomorphism and $M_h = \tilde M$.  
Since the unit sphere $\SI^{n-1}$ is  
a compact subset of $\brno$, the set $L = \dev^{-1}(\SI^{n-1})$  
is a compact subset of $M_h$. Since the action of 
deck transformation is properly discontinuous, there exists  
a deck transformation $\vth$ such that $\vth(L)$ is  
disjoint from $L$. Now, $\vth(L)$ and $L$ bound a subset  
of $M_h$ diffeomorphic to $\SI^{n-1} \times I$ for  
an interval $I$. Clearly $N = M_h/\langle\vth\rangle$ is 
a Hopf manifold and covers $M$ finitely.   
 
Our manifold $M$ is covered by the Hopf manifold $N$  
finitely. Obviously there exists a map $f: N \ra \SI^1$  
as $N$ is a mapping torus of $\SI^{n-1}$. Regarding  
$\SI^1$ as a quotient space of the real line $\bR$, 
the radiant vector field takes positive values under 
the closed and non-exact $1$-form $df$ on $N$. 
 
Since the radiant vector field on $N$ is pushed to 
a radiant vector field on $M$, a $1$-form $\kappa$ on $M$  
obtained by averaging $df$ over the inverse images  
of evenly covered neighborhood takes positive  
values for the radiant vector-field. 
Since $\kappa$ is an average of closed forms locally,  
$\kappa$ is closed. 
Hence by Tischler \cite{Tis}, we see that $M$ admits  
a total cross-section $\Sigma$ to the radial flow.  
(See the proof of Lemma \ref{lem:form} of \ref{app:radd} for more details.) 
Hence $M$ is a generalized affine suspension over $\Sigma$  
by Proposition \ref{prop:crssect}.  
The universal cover of $N$ is homeomorphic to $\brno$, which is  
homeomorphic to the product of $(n-1)$-sphere with the real line  
$\bR$. Since the universal cover of $N$ is also homeomorphic  
to $\tilde \Sigma \times \bR$ for a cover of $\tilde \Sigma$  
of $\Sigma$, $\tilde \Sigma$ is homotopy equivalent to 
$\SI^{n-1}$. Since $\tilde \Sigma$ is a total cross-section,  
it is a closed real projective manifold. It is now obvious that 
$\tilde \Sigma$ is projectively diffeomorphic to $\SI^{n-1}$
as a developing map must be such a diffeomorphism. 
 
Suppose now that $M$ has a nonempty boundary. Let  
$S$ be a boundary component, and let $S_h$ be a component  
of $p^{-1}(S)$ in $M_h$. Then as there are no ideal  
points of $M_h$, the argument in the first paragraph  
shows that $\dev|S_h$ maps to $P -\{O\}$ as a covering  
map where $P$ is a two-dimensional subspace.  
As $-\Idd$ commutes with all holonomy elements, and  
$-\Idd$ acts on any cover of $P-\{O\}$, we see that we can obtain  
a double $\hat M$ of $M$ which is closed and has  
a radiant affine structure. Therefore our proposition  
follows from the first part. 
\end{proof} 
 
We will now show that $\che{M_h} \cap \dev^{-1}(O)$  is a unique  
point in $\hideal{M}$.  
 
\begin{lem}\label{lem:extphi} 
\begin{enumerate} 
\item We have the inequality  
\begin{equation}\label{eqn:flowquasiM}  
C(t)^{-1}\bdd(x, y) \leq  \bdd(\Phi_{h,  t}(x),  \Phi_{h, t}(y))  \leq 
C(t) \bdd(x, y)  
\end{equation}  
for each $x, y \in M_h$ and a positive number $C(t)$ depending only on 
$t$ smoothly.  
\item  Each $\Phi_{h, t}$ extends to  a homeomorphism $\Psi_{h, t}$ of 
$\che{M_h}$ for $t  \in  \bR$,  and  the above  inequality  holds  for 
$\Psi_{h, t}$ also.  
\item Obviously, $\Psi_{h, t}$ acts on $M_h$ and hence on $\hideal{M}$ 
and on $\hfideal{M}$ and fixes each point of $\hIideal{M}$.  
\end{enumerate} 
\end{lem} 
\begin{proof} 
(1) For every pair of points $x$ and $y$,  $\bdd(x, y)$ is the infimum 
of $\bdd$-lengths  of paths $\alpha$ connecting $x$  and $y$ on $M_h$. 
For  each  path $\alpha$  connecting  $x$  and  $y$, $\Phi_{h, t}\circ 
\alpha$ connects $\Phi_{h, t}(x)$  and $\Phi_{h, t}(y)$.   
For the  positive number   $C(t)$ from equation \ref{eqn:flowquasi}, 
the $\bdd$-length of $\Phi_{h, t}\circ\alpha$ is bounded below 
by  $C(t)^{-1}$ times  that  of $\alpha$  and bounded  above by $C(t)$ 
times that of $\alpha$. Hence, we get the inequality.  
 
(2)  Since by (1),  $\Phi_{h, t}$ is  a  quasi-isometry $M_h \ra M_h$, 
$\Phi_{h, t}$ extends to a  continuous map $\Psi_{h,  t}: \che M_h \ra 
\che  M_h$, where obviously the  inverse of $\Psi_{h, t}$ is  
$\Psi_{h,-t}$.  
 
(3) Since $\Psi_{h,  t}$  acts on $M_h$,  $\Psi_{h,  t}$ acts   on its 
complement $\hideal{M}$. Since  $\Psi_{h, t}$ is a continuous extension 
of $\Phi_{h, t}$, the following diagram also commutes  
\begin{eqnarray}\label{eqn:flow2}  
\che{M_h}  &     \stackrel{\Psi_{h, t}}{\longrightarrow}  &  \che{M_h} 
\nonumber    \\  \dev  \downarrow  &    &  \downarrow \dev   \\ \mathcal{H} 
&\stackrel{\Phi'_t}{\longrightarrow} & \mathcal{H}. \nonumber  
\end{eqnarray}   
Since $\Phi'_t$ acts on $\brno$, it follows that $\Psi_{h, t}$ acts on 
$\hfideal{M}$.  
 
For each point $x$ of $\hIideal{M}$, there  exists a sequence $x_i \in 
M_h$ converging to $x$, where a radial line $l_i$ passing through 
$x_i$ for each $i$. $l_i$ has two endpoints $p_i$ and $q_i$  
so that $\dev(p_i) \in   \SII$ and $\dev(q_i) = O$ since  
$\dev|l_i$   is   a radial line  connecting $O$ and  a  point of $\SII$.  
Since $l_i$ is a  flow line of $\bV_h$,   
$\Psi_{h, t}$ acts on $l_i$ and hence on its  closure $\clo(l_i)$; thus,  
$\Psi_{h, t}$ fixes each $p_i$.  Since $p_i$ converges to $x$  as well,  
and $\che{M_h}$ is a metric space, $x$ is fixed by $\Psi_{h, t}$ for  
each $t \in \bR$.  
\end{proof}  
 
\begin{rem}\label{rem:notPhi} 
By an abuse of notation, we will denote $\Psi_{h, t}$ by $\Phi_{h, t}$ 
from now on.  
\end{rem} 
  
\begin{defn}\label{defn:radseg} 
A {\em radial}\/ segment in $\che M_h$ is a convex segment in $\che M_h$ 
such that the image of its endpoint consists of the origin and a point 
of $\SII$.  
\end{defn} 
One sees that the closure of any radial line is a radial segment.  
 
\begin{prop}\label{prop:origin} 
$\che M_h$ includes  a point $p$  of a radial segment in  $\che{M_h}$ 
satisfying $\dev(p)  = O$.  Furthermore, given   a point $q$  of $\che 
M_h$ such that $\dev(q) = O$, then $p = q$.  
\end{prop} 
\begin{proof} 
An  endpoint $p$  of  a radial maximal   line $l_1$ clearly satisfies 
$\dev(p) = O$ by Lemma \ref{lem:radcom}.  
 
Since $q$ belongs  to $\che{M_h}$, there exists  a  path $\alpha$ with 
endpoint $q$ and a  point $x$ of $M_h$.  
Choose a path  $\beta$ from a  point $y$  of  $l_1$ ending  at $x$ and 
produce a path $\gamma$ from  $y$ to $q$  obtained by joining $\alpha$ 
and $\beta$ at $y$. Since $\dev\circ\gamma$ is an arc in $\bR^n$ bounded  
away from $\delta \mathcal{H}$, $\dev \circ \Phi_{h, t}\circ\gamma$  
is a path connecting $\dev\circ\Phi_{h, t}(y)$ and $\dev(q)$, and  
has $\bdd$-length less than or equal to a constant $c(t)$ times  
that of $\gamma$. As $t \ra -\infty$, it is easy to 
see that $c(t) \ra 0$. Hence as $t \ra -\infty$, the $\bdd$-length of  
$\Phi_{h,t}\circ  \gamma$ goes   to zero, meaning that   
$\Phi_{h,  t}(y)$ converges to  $q$.  Clearly, $\Phi_{h, t}(y)$ converges  
to  $p$ since $\Phi_{h, t}(y)$ belongs to $l_1$; thus, we obtain $p = q$.  
\end{proof}  
 
We denote the unique point $p$ such that $\dev(p) = O$ by $O$, and say 
it  is the  {\em  origin}\/ of  $\che{M_h}$   or $M_h$; 
$O$ is fixed by every deck transformation of $M_h$.  
 
\begin{lem}\label{lem:comprad}  
Through each   point of $\che{M_h}$  a radial   segment passes.  Every 
radial segment $l$ is either a  subset of $\hideal{M}$, or $l$ 
with its endpoints removed belongs to $M_h$. 
The endpoints of $l$ consists of $O$ and a point of $\hIideal{M}$.  
\end{lem} 
\begin{proof} 
Let $x$  be a point of $\che{M_h}$.  If $x =  O$, then  we are done by 
Proposition  \ref{prop:origin}. Assume  $x  \ne O$,  which  means that 
$\dev(x) \ne  O$, and there exists a  sequence $\{x_i\}$  of points of 
$M_h$   converging to $x$ with  respect   to $\bdd$ such that $\bdd(O, 
\dev(x_i)) \geq C$ for  some  $C > 0$.   Since  $x_i \in M_h$,   there 
exists  a radial  segment $l_i$   containing  $x_i$ such that   $l^o_i 
\subset M_h$. Then since $x_i \ra x$ with  respect to $\bdd$, there 
exists a path $\alpha_{ij}$ in $M_h$ with endpoints $x_i$ and $x_j$ so 
that for each  $\eps > 0$, there exists  $N$ such that  if $i,j >  N$, 
then the $\bdd$-length$(\alpha_{ij})$ are less than $\eps$.  
 
Let  $\SI^2_1$ denote the $\bdd$-sphere of radius  $\pi/4$ with center 
$O$. For each $i$,  let $y_i$ denote the $\bdd$-midpoint  of $l_i$, and 
$\Pi$ denote the radial projection  of $\mathcal{H} - \{O\}$ to 
$\SI^2_1$, the radial completeness of $M_h$ shows that $\Pi\circ 
\dev   \circ  \alpha_{ij}$  lifts  to   a  path $\beta_{ij}$  in $M_h$ 
connecting  $y_i$ and $y_j$.  Since  $\dev(x_i)$ is uniformly  bounded 
away  from $O$, it follows that  the  $\bdd$-length of $\beta_{ij}$ is 
bounded  above  by $c(C)$ times  that  of  $\alpha_{ij}$, where $c(C)$   
is a positive constant depending only on $C$ easily  calculable from   
spherical  geometry. Hence $y_i$ converges to $y$ in $\che M_h$ and  
$l_i$  converges to a radial segment $l$ which contains $x$ and $y$. 
 
Let $l$ be an arbitrary radial segment. Then  $\Phi_{h, t}$ for $t \in 
\bR$ acts transitively on $l^o$. Since $\Phi_{h, t}$ acts on $M_h$ and 
on $M_{h\infty}$, it follows that either $l^o \subset M_h$ or $l$ is a 
subset of $M_{h\infty}$.  
\end{proof} 
 
\begin{defn}\label{defn:radset} 
A {\em radiant}\/ set  in $M_h$ is a  $\Phi_{h, t}$-invariant  subset of 
$M_h$,  and   a {\em  radiant}\/  set  in  $\che{M_h}$  is    
a $\Phi_{h,t}$-invariant closed subset of $\che{M_h}$ that is  
not a subset of $\hIideal{M} \cup \{O\}$. 
\end{defn}   
Obviously, a radiant set in $\che M_h$ is a union of maximal radial 
segments. A {\em radiant}\/ set in $\mathcal{H}$ is a $\Phi'_t$-invariant  
subset of $M_h$ that is not a subset of $\delta \mathcal{H} \cup \{O\}$.  
A closed radiant set in $\mathcal{H}$ is a union of radial segments.  
 
\begin{lem}\label{lem:radiint}  
Let $S$ be a closed radiant subset of $\che{M_h}$. Then $S \cap M_h$ and  
$S \cap \hfideal{M}$ are closed radiant sets also.  
\end{lem}  
 
\begin{exmp}\label{exmp:suspcompl} 
Let $N$ be a Benz\'ecri suspension of a closed real projective  
surface $\Sigma$. Then $N_h$ may be identified with  
$\Sigma_h \times (0,1)$, where $\Sigma_h$ is a holonomy cover  
of $\Sigma$. Then the Kuiper completion of $N_h$ can be  
easily identified with $(\che \Sigma_h \times [0,1])/\sim$ where  
$\sim$ is the equivalence relation identifying all points of 
$\che \Sigma_h \times 0$ with one point. The previous results  
in this section are easily verified for this example. 
\end{exmp} 
 
A {\em  radiant}\/ $n$-bihedron $L$ is  an $n$-bihedron in $\che M$ such 
that the  boundary of $\dev(L)$ contains $O$  and $\dev(F)$ for a side 
$F$ of $L$ lies in $\SII$, which is obviously a radiant set.  We define  
$\alpha_L$ to be the interior of $F$.   
 
\begin{lem}\label{lem:cosub}  
Suppose that  $M$ is not  convex. Then the  boundary $\delta B$ of a 
convex compact $n$-ball $B$ in $\che M_h$ {\rm (}resp. $\che M$\/{\rm )}  
meets $M_h$ {\rm (}resp. $\tilde M$\/{\rm ).}  Furthermore, a radiant  
$n$-bihedron is a bihedral $n$-crescent, and $\che M_h$  
{\rm (} resp. $\che M$ {\em )} cannot contain any $n$-hemisphere.  
\end{lem} 
\begin{proof}  
If $\delta B$ does not meet $M_h$, then $M_h$ is a subset of $B$, 
$\delta   B    \subset   \hideal{M}$, and  $M_h$ is boundaryless.  
Therefore, we obtain $M_h  = B^o$, which is a contradiction since  
$B^o$ is convex.  
 
Let $L$ be  a  radiant $n$-bihedron. Since the  side  $F$ of $L$ with 
$\dev(F) \subset \SIT$ is  a subset of $\hIideal{M}$, and by  
the above statement, the intersection of  the other side of $L$ with 
$M_h$ is not empty, it follows that $L$ is a bihedral $n$-crescent.  
 
If $\che{M_h}$ includes  an $n$-hemisphere $L$,  
then  $\dev|L : L \ra \mathcal{H}$ is an imbedding onto $\mathcal{H}$.   
Since $\dev(M_h)$ is a subset of $\mathcal{H}^o$, $\delta L$ does not   
meet $M_h$. The  first part of this  lemma  says that this cannot  
happen.  
\end{proof} 
 
\begin{prop}\label{radbih}  
Let $M$ be a compact radiant affine manifold {\rm (}$n \geq 2${\rm )}. 
If $M$ is  not  $(n-1)$-convex, then  $\che M_h$  includes a radiant  
$n$-bihedron.  
\end{prop} 
\begin{proof} 
By Theorem 4.6 of \cite{psconv}, $\che M_h$ includes  
an $n$-crescent $L$. $L$ is a bihedral $n$-crescent  
by Lemma \ref{lem:cosub}. 
 
By  Lemma \ref{lem:radcom}, $\dev(L)^o$ is  disjoint  from $O$.  Let 
$L'$ be the union of $\Phi_{h, t}(L^o)$ for all  $t \in \bR$, which is 
an open subset of $\tilde M$. Since by equation \ref{eqn:flow},  
$\dev(L')$   is  the  union  of $\Phi'_{h,t}(\dev(L^o))$ for all $t$,  
$\dev(L')$ is a radiant open half-space in $\mathcal{H}^o$.  
For any  $t$ and $t'$ such that $\Phi_{h,t}(L^o)$ and $\Phi_{h, t'}(L^o)$  
are  overlapping,  $\dev| \Phi_{h,t}(L^o)\cup\Phi_{h,  t'}(L^o)$ is a   
diffeomorphism onto $\dev(\Phi_{h,t}(L^o))\cup \dev(\Phi_{h,t'}(L^o))$  
by  Remark \ref{rem:overlap}. By choosing a  sequence $\{t_i\}$ such that  
$\Phi_{h,t_i}(L^o)$  and    $\Phi_{h, t_{i+1}}(L^o)$   always  overlap,   
and by  an induction, we see that $\dev| L'$  is a diffeomorphism  to 
an open half-space. For the closure $L''$ of $L'$ in $\che M$,  
$\dev| L''$ is  an imbedding onto the closure  of the half-space   
in $\mathcal{H}$ (see Chapter \ref{ch:prel}).  
Hence, $L''$ is a radiant $n$-bihedron.  
\end{proof} 
 
\begin{prop}\label{prop:raecr} 
Suppose   that  $M$  is  not  convex   and  $M_h$   includes an   open 
$n$-bihedron. Then $\che{M_h}$  includes a unique radiant $n$-bihedron 
including the open bihedron.  In fact, given any $n$-crescent, there 
exists a unique radiant $n$-bihedron including it.  
\end{prop}  
\begin{proof} 
The existence part is analogous  to the  proof  of the above  theorem; 
that is, we radially  extend it. The  uniqueness follows from Proposition 
\ref{prop:overlap}.  
\end{proof} 
 
Finally, we say about the intersection properties of radiant  
$n$-bihedra (see Appendix \ref{app:dipping} for definition of  
transversal intersection): 
\begin{cor}\label{cor:radtrs} 
Suppose that $M$ is not convex. Let $R_1$ and $R_2$ be two overlapping  
radiant $n$-bihedra in $\che M_h$. Then $R_1 = R_2$ or $R_1$ and $R_2$  
as $n$-crescents intersect transversally. 
\end{cor} 
\begin{proof}  
By Theorem \ref{thm:transversal}  
we need only to show that $R_1$ cannot be  
a proper subset of $R_2$ and vice versa. This is so since $\dev(R_1)$  
being radiant cannot be a proper subset of $\dev(R_2)$ from elementary 
geometry. 
\end{proof}

\chapter[Three-dimensional radiant affine manifolds]
{Three-dimensional radiant affine manifolds and  
concave affine manifolds}  
\label{ch:tdim} 
 
We will now  assume that  $M$  is a three-dimensional compact  radiant 
affine manifold with totally geodesic or empty boundary and is  
not two-convex from now  on.   
Since $\che  M_h$ does not  contain any hemispheric $n$-crescents by  
Lemma  \ref{lem:cosub},  Corollary 1.2 of \cite{psconv} shows that after   
splitting   along  the two-faced submanifolds  arising from bihedral    
$n$-crescents,  the  resulting manifold  $N$  decomposes into concave    
affine manifolds with concave boundary and two-convex  affine manifolds  
with  convex boundary.   The concave boundary may not be totally geodesic.   
However we will show in the next chapter that 
the  concave  boundary  of  the  concave  affine    
manifold is totally geodesic in the radiant affine $3$-manifold;  
our decomposition takes place along totally geodesic closed surfaces  
in $M$. (For full treatment of the decomposition theory, see \cite{psconv}.)

We   suppose   that $\che   M$   includes a bihedral 
$3$-crescent  since   otherwise    $M$   will   be   $2$-convex  by 
\cite{psconv} and we can proceed to the next chapter.   
By the previous chapter, any     bihedral  $3$-crescent is  
a subset of a unique  radiant 
$3$-bihedron. Let $B$ denote the set of all bihedral $3$-crescents and 
$\mathcal{B}$  that of all radiant $3$-bihedra.   We recall the equivalence 
relation    on  $B$    defined  in   \cite{psconv}.    Given  radiant 
$3$-crescents $R$ and $S$, we say $R$ and $S$ are simply equivalent if 
they overlap, which generates an  equivalence relation $\sim$ on  the 
set of all bihedral $3$-crescents. So we define as in \cite{psconv} 
for a bihedral $3$-crescent $R$: 
\[ \Lambda(R) = \bigcup\limits_{S \sim R} S, \quad 
\delta_\infty \Lambda(R) = \bigcup\limits_{S  \sim  R}   \alpha_S,\quad 
\Lambda_1(R) = \bigcup\limits_{S \sim  R} S -  \nu_S.  \]  
Given radiant $3$-bihedra $R$ and $S$, we say that $R$ and $S$  
are simply equivalent if they overlap, which again generates  
an equivalence relation $\sim'$ on $\mathcal{B}$. We define 
for a radiant bihedron $R$: 
\[ \Lambda'(R) = \bigcup\limits_{S \sim' R} S, \quad 
\delta_\infty \Lambda'(R) = \bigcup\limits_{S  \sim'  R}   \alpha_S,\quad 
\Lambda'_1(R) = \bigcup\limits_{S \sim'  R} S -  \nu_S.  \]  
We  easily see 
that given a $3$-crescent $T$, $\Lambda'(T') \subset \Lambda(T)$ for 
the unique radiant $3$-bihedron $T'$ including $T$ since any radiant 
bihedron equivalent  to $T'$  is  equivalent  to  $T$; conversely,  
by  Proposition \ref{prop:raecr}, we  have $\Lambda'(T') = \Lambda(T)$.   
Since radiant $3$-bihedra are radiant, $\Lambda'(R)$  is a radiant set;  
$\Lambda'(R) \cap  M_h$  and $\Lambda(R) \cap    M_h$  are radiant sets.  
 
Since any  $3$-crescent $R$ is included  in a unique  radiant bihedron 
$R'$, a  radiant bihedron is  a $3$-crescent by Lemma \ref{lem:cosub}, 
and $\alpha_R$  being in $\hideal{M}$  cannot meet  $R^{\prime o}$, it 
follows   that   $\alpha_R  =   \alpha_{R'}$   with $\alpha_R  \subset 
\hIideal{M}$.   Hence, $\delta_\infty   \Lambda(R) = \delta_\infty 
\Lambda'(R')$ and $\Lambda_1(R) =  \Lambda'_1(R')$.  (Because of these 
properties,  we will drop  primes from $\Lambda'(R)$ and $\sim'$,  
$\delta_\infty \Lambda'(R),$ and  $\Lambda'_1(R)$  respectively  
for a   given radiant bihedron $R$.)  
 
Let   us recall   properties   of  $\Lambda(R),   \Lambda_1(R),$   and 
$\delta_\infty \Lambda(R)$  from  Chapter 7 of \cite{psconv},  
the proofs of which are given there. $\Lambda(R)$  and $\Lambda_1(R)$   
are path connected. Given a natural projective structure  
induced from $3$-crescents,  
$\Lambda_1(R)$ is a real projective $3$-manifold with boundary  
$\delta_\infty  \Lambda(R)$ a totally geodesic open surface. 
Since   Corollary \ref{cor:radtrs}   shows   that for  two  overlapping 
$3$-crescents $R_1$  and   $R_2$, $\alpha_{R_1}$   and  $\alpha_{R_2}$ 
extend  each  other into   a  larger surface,  there exists a 
unique  great      sphere $\SI^2$   including  $\dev(\delta_\infty 
\Lambda(R))$ and a unique component $A_R$  of $\SI^3 - \SI^2$ such 
that  $\dev(\Lambda(R))  \subset  \clo(A_R)$  and   $\dev(\Lambda(R) - 
\clo(\delta_\infty    \Lambda(R)))     \subset A_R$.     
In our case $\SI^2$ must equal $\SIT$ and $A_R$ equal $\mathcal{H}^o$ 
by our fixed choice. For  a deck transformation   
$\vth$     acting    on   $\Lambda(R)$,    $A_R$   is 
$h(\vth)$-invariant.     
$\Lambda(R)  \cap M_h$ is a closed subset 
of $M_h$. The topological boundary $\Bd  \Lambda(R)  \cap   M_h$   
is a properly imbedded topological   
surface in $M_h^o$, and  $\Lambda(R)  \cap M_h$ is   
a topological submanifold of $M_h$ with concave boundary  
$\Bd \Lambda(R) \cap M_h$.  
 
Given a  deck transformation $\vth$,   we  have  $\vth(\Lambda(R))  = 
\Lambda(\vth(R))$, $\vth(\Lambda_1(R)) =   \Lambda_1(\vth(R))$,    and 
$\vth(\delta_\infty \Lambda(R)) = \delta_\infty \Lambda(\vth(R))$.  
 
\begin{prop}\label{prop:equipr} 
Given two radiant bihedra $R$ and $S$, precisely one of the following 
possibilities holds\/{\rm :}  
\begin{itemize} 
\item $\Lambda(R)= \Lambda(S)$ and $R \sim S$.  
\item $\Lambda(R) \cap M_h$ and $\Lambda(S) \cap M_h$ are disjoint and 
$R \not\sim S$.  
\item $\Lambda(R) \cap \Lambda(R) \cap M_h$ equals the union of common 
components of $\Bd \Lambda(R) \cap M_h$  and $\Bd \Lambda(S) \cap M_h$ 
where  $R  \not\sim S$, which furthermore equals the  union of common 
components of  $\nu_T \cap M_h$ and $\nu_{U}  \cap M_h$ for two  
radiant bihedra $T \sim R$ and $U \sim S$.  
\end{itemize} 
\end{prop} 
 
We recall  from  \cite{psconv}  that the  {\em   copied}\/ components 
arising  from  bihedral $3$-crescents  are  common  components of $\Bd 
\Lambda(R) \cap  M_h$ and $\Bd \Lambda(T)\cap M_h$  for some  pair of 
inequivalent radiant bihedra $R$ and $T$.  
 
\begin{lem}\label{lem:presplit} 
The collection whose elements are copied components are locally finite, 
the union $A$  of all  copied  components is a properly  imbedded 
submanifold in $M_h$, and $p|A: A \ra  p(A)$ covers a compact submanifold 
$p(A)$ of codimension one.  
\end{lem} 
 
$A$ is said to be the  {\em pre-two-faced submanifold}, and $p(A)$ the 
{\em two-faced submanifold}. If  $C$  is  a  copied component of  
$\nu_{R'}  \cap M_h$ in  $A$, then $p(C)$ is a component of  the closed   
codimension-one  submanifold $p(A)$. Since for a radiant  $3$-bihedron $R'$,  
$\nu_{R'} \cap  M_h$ is  foliated by radial  lines,   
the following lemma shows that $p(A)$ consists  of  
components   that   are totally geodesic two-sided tori or one- or  
two-sided Klein bottles.  
 
\begin{lem}\label{lem:convtori} 
Suppose that a  closed projective surface $S$   is covered by  an open 
subset $C$ of $\bR^2 - \{O\}$ foliated by radial lines. Then $S$ is an 
affine torus{\rm ,} or an affine Klein bottle covered by a convex open 
cone in $\bR^2 - \{O\}$ or covered by $\bR^2 - \{O\}$ itself.  
\end{lem}  
\begin{proof}  
Suppose that $C$ is a proper subset of $\bR^2 - \{O\}$.  Given on  
$\bR^2$ polar  coordinates $(r, \theta)$ where $r,  \theta \in  \bR$,  
$C$ is given by $\theta_1 <  \theta < \theta_2$ where  
$\theta_2 - \theta_1 <  2\pi$. Suppose that $\theta_2 -  \theta_1  > \pi$.    
Let $l_1$ and  $l_2$ be the radial  lines  corresponding to $\theta_1$ and 
$\theta_2$ respectively.  Then a radial  line  $m_1$ in  the  opposite 
direction to $l_1$ and $m_2$ to $l_2$ both lie in $C$.  
 
Since   $p(C)$     is  a  closed    surface,    a   subgroup  $G$   of 
$\pi_1(M)/\pi_1(M_h)$ acts on $C$  so that $p$ induces a homeomorphism 
$C/G \ra  p(C)$.  Since  $C$ is  invariant under  the affine action of 
$G$, it follows that $m_1 \cup m_2$ is $G$-invariant, and  
$(m_1  \cup m_2)/G$ is a simple closed  curve or the union of two simple   
closed curves. As in Lemma 2 of \cite{cdcr1}, this means that a double  
cover of the surface $p(C)$ is homotopy equivalent to a simple closed   
curve in the cover, and such  a closed surface   does  not exist.  Hence,   
$\theta_2 -\theta_1 \leq\pi$, and $C$ is convex. Since $G$ preserves  
the foliation, $S$ is foliated; $S$ is homeomorphic to a torus or  
a Klein bottle.  
 
If $C$  equals $\bR^2 - \{O\}$,  then $S$ is  again foliated by radial 
lines, and is homeomorphic to a torus or a Klein bottle.  
\end{proof}   
 
We split $M$ and $M_h$ by $p(A)$ and  $A$ respectively, and we obtain 
an $n$-manifold $N$  and its cover  $M'$ including a copy of 
$M_h - A$ as a dense subset such that the covering map  
$p| M_h - A: M_h - A \ra M -p(A)$ naturally extends  to a covering  
$p: M'  \ra N$ and $\dev| M_h   - A$ extends to a developing map  
$\dev: M' \ra \SI^3$. As explained in Chapter 8 of \cite{psconv}, 
given a component $N_i$ of $N$, we choose a component $P_i$ of  
$M_h - A$ covering a component of $M - A$ dense in $N_i$.  
Then the disjoint union of $P_i$ is a holonomy cover of $M - A$,  
and the disjoint union of $L_i$ where $L_i$ is the component of $M'$  
where $P_i$ is dense in is a holonomy cover of $N$. 
(Note that the cover is not necessarily connected here.) 
We denote the disjoint union by $N_h$.  
 
Since $M'$ admits an immersion to $\SI^3$, we may pull-back 
the standard metric $\mu$ to $M'$ and obtain a completion $\che M'$  
of $M'$. We can show that the completion $\che M'$ of $M'$  
includes radiant bihedra with same interiors as those of $\che M_h$;   
there is  a  one-to-one correspondence of radiant bihedra in this way.  
Using this we can show easily that  
in $\che M'$, there are no more copied components of  
$\Bd \Lambda(R) \cap N_h$ for every crescent $R$ in $\che M'$.  
(See Chapter 8 of \cite{psconv} for details.)  
Hence, in $\che N_h$ there are no copied components. 
 
Let $\mathcal{A} = \bigcup_{R \in B'} \Lambda(R)  \cap N_h$ where $B'$ is 
the  set of    representatives   of the  equivalence    classes  of 
radiant bihedra in $\che N_h$ under $\sim$.  
Since we did the splitting, given  two radiant bihedra $R$ and 
$S$,   either $\Lambda(R) =   \Lambda(S)$  or $\Lambda(R)\cap M_h$ and 
$\Lambda(S) \cap M_h$ are  disjoint.  Since deck transformations send 
radiant bihedra to radiant bihedra, the deck transformation group acts 
on $\mathcal{A}$.   Since the collection consisting  of elements of  form 
$\Lambda(R')$ is locally finite (see \cite{psconv}),  
$\mathcal{A}\cap N_h$ is a submanifold of $N_h$. It follows that $\mathcal{A}$  
covers a compact codimension $0$ submanifold $L$  of $N$  
with boundary $\delta L$ in $N^o$ by  Proposition  
\ref{prop:lfin}; $L$ is turns out to be a concave affine manifold of  
type II.  (See Chapters 7 and 10 of \cite{psconv}  for more details.) 
  
Proposition \ref{prop:totgeo} shows that $\Bd L$ is totally  
geodesic since each component of $\Bd L$ is covered by  
a component of $\Bd \Lambda(R) \cap N_h$ for  
some $R$.  
 
\begin{prop}\label{prop:totgeo} 
$\Bd  \Lambda(R) \cap  N_h$  is  a   totally geodesic submanifold   in 
$N_h^o$.  
\end{prop} 
\begin{proof} See Chapter \ref{ch:totgeo} for the proof. 
\end{proof} 
 
Since    $\Lambda(R)  \cap N_h$ is  radiant,  $\Bd 
\Lambda(R) \cap N_h$  is  also foliated  by  radial lines.   
 
\begin{prop}\label{prop:bdL} 
Each component of  $\Bd \Lambda(R) \cap  N_h$ is a component of $\nu_T 
\cap N_h$ for a radiant bihedron $T$, $T \sim R$.  
\end{prop} 
\begin{proof} 
Let   $x$ be a point    of a component   $F$  of $\Bd \Lambda(R)  \cap 
N_h$ so that $x \in T$ for $T \sim R$. Since $x$ is a boundary point, $x 
\in \nu_T$. The totally geodesic surface $F$ and  
the component $F'$  of $\nu_T \cap N_h$ containing 
$x$ are tangent   at $x$ since otherwise,  $F$  intersects $T^o$ which  is 
absurd.  Since $F$  and $F'$ are   both totally geodesic and  properly 
imbedded, we obtain $F = F'$.  
\end{proof} 
 
\begin{thm}\label{thm:dec2conv} A compact radiant affine $3$-manifold 
$M$ decomposes into two-convex radiant affine manifolds and radiant  
concave affine manifolds of type II along totally geodesic affine tori  
or Klein bottles covered by convex open cones in $\bR^2$  
or covered by $\bR^2 -\{O\}$. The Kuiper completion of the holonomy  
cover of each of the two-convex radiant affine manifolds includes 
no $3$-crescents or radiant bihedra.  
\end{thm} 
 
\chapter{The decomposition along totally geodesic surfaces} 
\label{ch:totgeo} 
 
The aim of this chapter is to prove Proposition  
\ref{prop:totgeo} needed in Chapter \ref{ch:tdim}. 
We will use the notations of the chapter: 
We prove  that $\Bd \Lambda(R) \cap N_h$  is totally geodesic  
by  showing that at  each tiny ball  
neighborhood $B(x)$ of $x \in \Bd \Lambda(R) \cap N_h$, 
there exists a  unique supporting  plane to  the  convex ball $B(x)  - 
\Lambda(R)$.  If  this is not true, we  will show that there exist two 
transversally intersecting crescents nearby $x$, which will be shown to 
be a contradiction to certain  maximal  properties of $\Lambda(R)$  by 
Lemma \ref{lem:twocresc}.  
 
We begin the proof.   Let $x \in \Bd \Lambda(R)  \cap N_h$  and 
the tiny ball neighborhood $B(x)$ of $x$.  
Then  since  $x  \in  N_h^o$,  $B(x)^o$  is an open 
neighborhood of $x$ and $B(x)^o - \Lambda(R)$ is a convex open set $K$ 
by definition of concave  boundary (see Chapter 3 of  \cite{psconv}). 
Note that the boundary of $K$  in $B(x)^o$ equals $\Bd \Lambda(R) \cap 
B(x)^o$.   Since $\Lambda(R) \cap N_h$ is  a radiant set, $K$ and $\Bd 
K$ are foliated by maximal radial lines in $B(x)^o$.  
 
We claim   that every point  of $\Bd  \Lambda(R)  \cap B(x)^o$  has an 
identical    supporting  hyperplane in   $B(x)^o$   to  $K$. Then $\Bd 
\Lambda(R)  \cap B(x)^o$   equals   this  hyperplane  and is   totally 
geodesic, proving the proposition. 
 
Suppose that the claim is not true --(*). Then we will show that there 
exist two radiant $3$-bihedra $T$ and $T'$  
($T, T' \sim R$) meeting transversally, satisfying  
$T^o \cap T^{\prime o} \cap B(x)^o \ne \emp$,  
$\nu_T \cap \nu_{T'}$ meets $B(x)^o$, and $\nu_T$ contains 
a point of $B(x)^o \cap  \Bd \Lambda(R)$; our claim  will be proved by 
showing that this cannot happen by Lemma \ref{lem:twocresc}.  
 
Firstly, suppose that  at a point  $y$ of $\Bd \Lambda(R) \cap B(x)^o$ 
there are at least two distinct supporting  hyperplanes $P_1$ and $P_2$ 
in $B(x)^o$ to $K$.  Since $K$ is foliated  by maximal radial lines in 
$B(x)^o$,  $\dev(P_1)$ and  $\dev(P_2)$  respectively are included  in 
hyperplanes in $\bR^3$ passing through $O$.  Since $P_1$ and $P_2$ lie 
in the closure  of the complement of  $K$,  it follows that  $P_1, P_2 
\subset  \Lambda(R)  \cap B(x)^o$.   Since $\Lambda(R)$  is a union of 
radiant  $3$-bihedra, it follows that $\Lambda(R) \cap   N_h$  is  
foliated  by radial lines, $P_1$ and $P_2$ must be foliated by subarcs  
of radial lines in $B(x)^o$, and they meet at  a subarc  of a radial   
line in $B(x)^o$ passing through $y$.  
 
Let $\SI^{2}_r$ be the sphere with center $O$ and radius  
$r=\bdd(O,\dev(x))$. Then $\dev^{-1}(\SI^{2}_r)\cap B(x)^o$ is  
a smooth  hypersurface $\Sigma_x$ in $B(x)^o$ 
$\mu$-orthogonal to radial lines.  
 
\begin{figure}[t] 
\centerline{\epsfxsize=3in \epsfbox{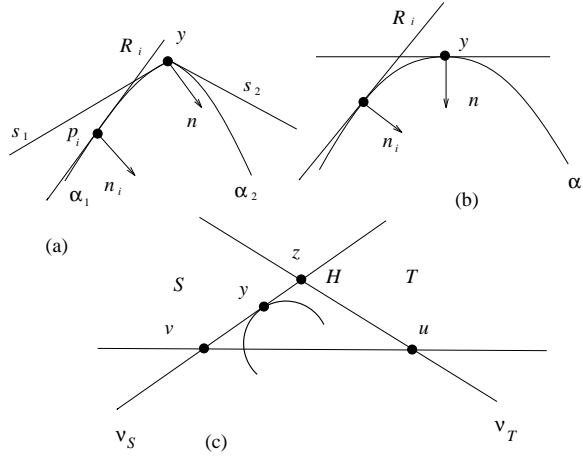}}  
\caption{\label{fig:crossect} The cross-section $\Sigma_x$ of $B(x)$ 
meeting with crescents.}  
\end{figure}

It is obvious that $\Sigma_x$ has an induced projective structure from 
the imbedding $\SI^{2}_r \hookrightarrow  \SI^3$, i.e., $\Sigma_x$ has 
a   flat  projective  structure of  an   open   subset of  the  sphere 
$\SI^{2}_r$ with the projective structure induced from its Riemannian  
metric. Let $s_1$   and $s_2$ be arcs  in  $P_1 \cap  \Sigma_x$  and  
$P_2 \cap \Sigma_x$  respectively ending at $y$, oriented away from $y$,   
and  at points of $\Bd B(x)^o$ as shown in Figure \ref{fig:crossect}(a). 
$\Sigma_x  \cap K$  is a convex open subset with a boundary point $y$.  
The boundary $\alpha$ of $\Sigma_x \cap K$ in $\Sigma_x$ forms an open arc,  
and $\alpha - \{y\}$ has two  components $\alpha_1$ and $\alpha_2$,  
which we orient away  from $y$. We may  assume  that $s_1$  are  tangent   
to   $\alpha_1$  and  $s_2$  to $\alpha_2$. (See Figure  
\ref{fig:crossect} (a).)  
 
\typeout{<<rf15.eps>>}  
 
Choose   a sequence $\{p_i\}$  of  points  of $\alpha_1$ converging to 
$y$. Since $\alpha_1$ is a  subset of $\Lambda(R)$, $p_i \in R_i$ 
for a radiant $3$-bihedron  $R_i$, $R_i \sim R$; as $p_i \in \Bd 
\Lambda(R)$,  it follows that $p_i \in  \nu_{R_i}$ for  each $i$.  Let 
$n_i$ be the outer normal vector to $\nu_{R_i}$  at $p_i$ with respect 
to the spherical metric $\bdd$, which is a vector tangent to the spherical  
surface $\Sigma_x$ at $p_i$.   Recall that $\nu_{R_i} \cap B(x)$ is a 
compact $(n-1)$-ball with boundary in $\Bd B(x)$ and $R_i^o \cap B(x)$ 
is  a  component of  the   complement of this   disk (see  Lemma 3.13  of 
\cite{psconv}); since the endpoint $y$ of $s_1$ is not a  subset of  
$R_i^o  \cap  B(x)^o$, it follows  easily that the   angle $\theta_i$  
(i.e.   the spherical angle)   between the tangent vector to the segment   
on $\Sigma_x$ connecting $y$ and  $p_i$ oriented towards $p_i$ at $p_i$  
and $n_i$ is greater  than or equal to $\pi/2$ by a geometric  
consideration (see Figure \ref{fig:crossect}). Choosing   a subsequence  
if necessary,  we  may assume without loss of generality that $n_i$  
converges to a unit  vector $n$ at $y$. By Corollary 3.16 of  
\cite{psconv} (or Section   6.2 of \cite{cdcr1}), we see that there  
exists a radiant $3$-bihedron $T$ in $\che M_h$ with $y\in \nu_T$,  
$n$ the outer normal vector to $T$, and $R \sim T$.  
 
Since $n$ is a limit of  $n_i$, the angle  between $n$ and the tangent 
vector to $s_1$ at $y$ is greater than or equal  to $\pi/2$.  Since $T 
\cap B(x)$ is a closure of a component of $T  - P$ for an $(n-1)$-ball 
$P = \nu_T \cap B(x)$, and $\alpha_1$ is in $\Bd \Lambda(R) \cap 
B(x)^o$, $s_1$ cannot  point  towards $T^o$;  otherwise, a  small open 
arc in $\alpha_1$ near $y$ is included in $T^o$.  Therefore, $s_1$ and 
$n$ are perpendicular and $s_1 \subset \nu_T$. 
 
Similarly, we have $s_2  \subset T'$  for a radiant  $3$-bihedron with   
$y \in \nu_{T'}$ and a normal vector $n'$ perpendicular to the tangent  
vector to $s_2$ at $y$.  
 
Since $s_1$ and $s_2$ make  an angle less  than $\pi$ and greater than 
$0$, it follows that  $n$ and $n'$ make  an angle less than  $\pi$ and 
greater than $0$.  Since $\nu_{T}$  and  $\nu_{T'}$ both  pass $y$, it 
follows that $\nu_T \cap \nu_{T'}$ meets $B(x)^o$, and  
$T^o \cap T^{\prime o} \cap   B(x)^o \ne \emp$,  and $T \sim T'$.   
Since the normal vectors  are not parallel, it follows that $T$ and $T'$  
intersect transversally by Theorem \ref{thm:transversal} and obviously both  
$T$ and $T'$ contain $y$ in $B(x)^o \cap \Bd \Lambda(R)$.  
 
Secondly, suppose now that at every  point $y$ of $\Bd \Lambda(R) \cap 
B(x)^o$,  there exists   a unique supporting   hyperplane  to  $K$  in 
$B(x)^o$ at $y$, and $\Sigma_x \cap \Bd \Lambda(R) \cap B(x)^o$ is a 
differentiable convex  arc,   say  $\alpha$.  Since by our   assumption   
(*) $\alpha$  is not projectively geodesic in   the sphere $\Sigma_x$ with 
the  induced projective  structure, there   is a  point  $y$ in 
$\alpha$ with a compact neighborhood in $\alpha$ where there exists an 
infinite  collection of distinct  supporting  hyperplanes to points of 
the neighborhood.  Hence, we can choose a sequence $\{p_i\}$ of points 
of $B(x)^o \cap \Bd \Lambda(R)$ converging  to $y$ of $B(x)^o \cap \Bd 
\Lambda(R)$    so that supporting  hyperplanes  at  $p_i$ are mutually 
distinct.  
 
We have $p_i \in R_i$ for $R_i \sim R$ for each $i$.  Since $p_i$ does 
not belong to $R_i^o$, we have $p_i  \in \nu_{R_i}$ for $R_i \sim R$. 
Since $\nu_{R_i} \subset  \Lambda(R)$, it follows that $\nu_{R_i} \cap 
B(x)^o$ is the unique supporting hyperplane  at $p_i$; thus, $R_i$ are 
mutually distinct.  Similarly to the above argument, $y$ is a point of 
a radiant $3$-bihedron $S$. Let  $n_i$  denote  the normal vector to  
$R_i$ at  $x_i$  and $n$  that to $S$   at $y$. By the uniqueness  
of tangent lines, it  is easy  to see that $n_i$ converges to $n$.  
This, the convexity of $\alpha$, and an elementary spherical geometry  
show that for $i$ sufficiently large $\nu_S \cap \nu_{R_i}$ meets  
$B(x)^o$, and we have $S^o \cap R_i^o \cap B(x)^o  \ne \emp$, and  
$S \not\subset R_i$ and $R_i \not\subset S$ (see Figure \ref{fig:crossect}  
(b)).  Since  $S$ and $R_i$ overlap, $S$ and $R_i$ intersect transversally,   
by Corollary \ref{cor:radtrs} and $R_i$ contains $p_i$ in  
$B(x)^o \cap \Bd \Lambda(R)$. Therefore, our proof is completed by 
the following lemma.\qed 
 
\begin{lem}\label{lem:twocresc} 
We have two   radiant $3$-bihedra $S$  and  $T$ equivalent to  $R$ and 
overlapping with a tiny ball $B(x)$ such that $\nu_S \cap \nu_T$ meets 
$B(x)^o$, and $S$  and $T$ meet  transversally, and $S^o \cap T^o 
\cap B(x)^o  \ne \emp$.  Then  $S$ cannot contain  a  point of $B(x)^o 
\cap \Bd \Lambda(R)$,  
\end{lem} 
\begin{proof} 
Suppose  not.  Let $y$  be  a point  of  $\nu_S$ in   $B(x)^o \cap \Bd 
\Lambda(R)$. Since $S$ and $T$ meet  transversally, $\nu_S \cap \nu_T$ 
is a segment $H$ with interior in $\nu_S^o$ and $\nu_T^o$ and boundary 
in $\delta  \nu_S$ and  $\delta \nu_T$.   
 
If   $y$ does not  belong  to  $H$, choose a  point   $z$  of $H  \cap 
\Sigma_x^o$, a  point $u$ of $\nu_T \cap   \Sigma_x$ nearby $z$, and 
$v$  of  $\nu_S \cap  \Sigma_x$   near $y$  so   that the projective 
geodesic  $\ovl{zv}$ on $\Sigma_x$ contains  $y$  in the interior.  
We obtain the triangle  $\triangle(vzu)$ in  $\Sigma_x$ with  
edges $\ovl{vz}$  and $\ovl{zu}$ lying in $\nu_S$  and $\nu_T$   
respectively.   If $y$ belongs to $H$,  then choose a  point $z$  
to be  $y$ and a point $u$ of $\nu_T \cap \Sigma_x$ nearby $z$, and  
a point $v$ of $\nu_S \cap \Sigma_x$ near $z$  so  that   
$\triangle(vzu)$ lies in $\Sigma_x$  and its edges $\ovl{vz}$  and 
$\ovl{zu}$ in $\nu_S$ and $\nu_T$ respectively. By our construction,  
$\ovl{vu}$ is   not  a subset   of   $(S   \cup  T)   \cap  B(x)^o$  
(see   Figure \ref{fig:crossect} (c)).  
 
We denote by $P$ the subspace of $M_h$ that is the union of all radial 
lines   passing   through    $\triangle(vzu)$.    Again   by  equation 
\ref{eqn:flow}, we see that  $\dev| P$ is an  imbedding onto  a convex 
subset bounded  by totally  geodesic  surfaces  in  $\mathcal{H}^o$.   In 
particular, $\clo(P)$  is a radiant tetrahedron in  $\che M_h$ with  a 
vertex $O$ bounded by three  sides $F_1$, $F_2$, and $F_3$ corresponding 
to  $\ovl{vu}$, $\ovl{vz}$, and  $\ovl{zu}$  respectively and the side 
$F_4$ in $\hIideal{M}$.  
 
Since $S$ and $T$   are overlapping, $\dev| S  \cup T$ is  an imbedding  
onto $\dev(S)  \cup \dev(T)$ where  $\dev(S)$ and $\dev(T)$ are two  
radiant $3$-bihedra  with sides in  $\SIT$ and $\dev(\nu_S)$ and  
$\dev(\nu_T)$ meeting transversally.   
Let $H'$  be the $(n-1)$-hemisphere in $\mathcal{H}$ including $F_1$.  
Then $H'$ is a side  of a unique radiant $3$-bihedron  $U$ in $\mathcal{H}$  
including $\dev(\clo(P))$ and whose other side is a hemisphere in  
$\dev(\clo(\alpha_S)) \cup \dev(\clo(\alpha_T))$.  
We can show easily  that $\dev| S \cup T  \cup \clo(P)$ is an 
imbedding;  hence, $S    \cup  T  \cup \clo(P)$   includes   a radiant 
$3$-bihedron  $U'$ mapping onto $U$.   Since  $\clo(P) \subset U'$ and 
$\ovl{vz}$ is  transversal to $F$, it  follows  that $\ovl{vz} - \{v\} 
\subset U^{\prime  o}$,  which means $y \in   U^{\prime o}$; since $U' 
\sim R$ obviously, this contradicts $y \in \Bd \Lambda(R)$.  
\end{proof}

\chapter{$2$-convex radiant affine manifolds}\label{ch:2conv} 
 
We will now show that $2$-convex radiant affine compact 3-manifolds  
with empty or totally geodesic boundary further decompose along totally  
geodesic  affine  tori or Klein  bottles  into convex radiant  affine   
manifolds and concave-cone affine manifolds which are either Seifert  
spaces or finitely covered by a bundle over a circle with fiber  
homeomorphic  to  tori. (We assume that the Kuiper completion of 
the holonomy cover of the radiant affine $3$-manifold does not 
include any $3$-crescents or radiant bihedra.) As   we  said  in   
the introduction, this major step will be accomplished in Chapters  
\ref{ch:2conv}-\ref{ch:obtain}.  
 
For this purpose, we will obtain a crescent-cone (Theorem  
\ref{thm:crscone}) in this chapter,  which will play  
the role of a crescent for real projective  surfaces to give us  our 
desired decomposition. 
First, we will  define various   objects  in $\che  M_h$.   
Next,   we  find  a  radiant tetrahedron $F$ which detects     
the nonconvexity of $M_h$  and name various  parts to this tetrahedron.   
 
Denoting  by $F_1,F_2, F_3$,  and $F_4$ the sides of the tetrahedron   
$F$, we  find a sequence of points in the face $F_3$ of $F$ which  
leaves every compact 
subset of $M_h$ and equidistant from $U$ and  $D$, the upper and lower 
subsets of $F_1  \cup F_2$. We do the  standard pull-back argument  in 
\cite{cdcr1};  we  pull the     points   of  the  sequence  by    
deck transformations  to a fixed fundamental domain   of $M_h$ and pull $F$ 
and  the named  subsets  of $F$  along  with the  points. Then by 
choosing subsequences, we assume that the each sequence of pulled-back 
images of the  objects converges.  We will show  that the ``limit'' of 
the  sequence of images  of $F$   is a nondegenerate    $3$-ball, 
i.e., a  radiant 
trihedron or  a radiant tetrahedron, and identify all  the ``limits'' of the 
sequences of images of every object. ({\em Note that we  will use the term 
``nondegenerate''  merely  to indicate it  is not  a lower dimensional 
object. The term is redundant actually but is used for emphasis}.)  
 
Using the {\em claim}, i.e., Proposition \ref{prop:MainEqu},  
in the introduction to be proved in Chapters \ref{ch:claim},  
\ref{ch:tric}, and \ref{ch:rtrh}, we verify at the end of the chapter: 
\begin{thm}\label{thm:crscone} 
Let $M$ be a compact $2$-convex but nonconvex radiant affine    
$3$-manifold with totally geodesic or empty boundary.  
Assume  that the Kuiper completion $\che M_h$ of  
its holonomy cover $M_h$ does not  include radiant bihedra.  
Then  $\che M_h$ includes a crescent-cone.  
\end{thm} 
 
We will use the following hypothesis of above theorem to be referred to  
in Chapters \ref{ch:2conv}-\ref{ch:obtain}.  
\begin{hyp}\label{hyp2} 
We  assume  that  $M$ is a  $2$-convex  but  non-convex radiant affine 
$3$-manifold with totally geodesic boundary  or empty boundary, and 
$\che M_h$ does not include $3$-crescents or radiant bihedra.  
\end{hyp}  
We  fix an  identification of  $\bR^3$ with the  interior  of a closed 
$3$-hemisphere $\mathcal{H}$, and denote by $\SIT$  the sphere at infinity,  
fix the  development pair  $\dev:  \che   M_h  \ra    \SI^3$ whose  
image lies in $\mathcal{H}$ as $M$ is affine and  $h:\pi_1(M)/\pi_1(M_h)\ra  
\Aut(\mathcal{H})$ whose image fixes the origin $O$. 
 
We recall  some terminology:  A  {\em  lune}\/ in  $\SI^3$ is a closed 
convex subset of a great $2$-sphere $A$ in $\SI^3$ that is the closure 
of a component of $A$ with two distinct  great circles in $A$ removed; 
a {\em triangle}\/  is the closure of  a component of $A$ with 
three great circles meeting in  general position removed.  
A {\em trihedron}\/ in $\SI^3$ is a closed convex subset of 
$\SI^3$ that is the closure of 
a component of $\SI^3$ with three great $2$-spheres meeting in 
general position removed. 
The  boundary of a trihedron  is the union  of three 
lunes with  disjoint interiors.  A {\em tetrahedron}\/ in $\SI^3$ is a 
closed convex subset of $\SI^3$ that is the  closure of a component of 
$\SI^3$ with four great $2$-spheres meeting  in   general 
position removed.  
The boundary of  a tetrahedron  in  $\SI^3$ is the union of 
four triangles with disjoint interiors. Recall how the corresponding  
objects in $\che M_h$ and $\che M$ are defined. 
 
\begin{rem}\label{rem:cone3} 
A convex set in $\clo(\mathcal{H})$ is  radiant if it is a cone 
over a  convex   set in $\SIT$.   It   is easy to  see   that  given a 
hemisphere, a lune,  or a simply  convex triangle in  $\SIT$, the cone 
over such sets are a  radiant  bihedron, a  radiant  trihedron, or a 
radiant tetrahedron respectively (see Remark \ref{rem:cone}).  
\end{rem} 
 
\begin{figure}[t] 
\centerline{\epsfxsize=2.5in \epsfbox{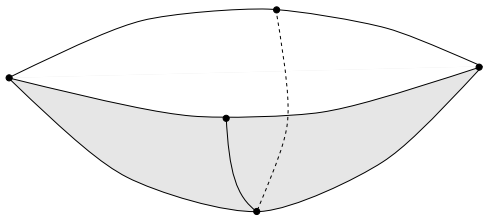}}  
\caption{\label{fig:crescent-cone} An example of a crescent-cone $T$, 
as a trihedron, where the grey area indicates the side in  
the closure of $T \cap \hideal{M} - \hIideal{M}$.  The top side is  
in $\hIideal{M}$. The picture is only schematic  
as the geodesics are bent.}  
\end{figure}

\begin{defn}\label{defn:crscone} 
A {\em radiant finite lune}\/ in $\che M_h$ is a  lune $A$ such that the 
boundary of $\dev(A)$ is  the union of two edges, one of which passes 
through the origin and the other in $\SIT$.  An {\em infinite lune}\/ is 
a  lune in $\che M$ included  in $\hIideal{M}$.  A {\em crescent-cone}\/ 
is a trihedron $A$ in $\che M_h$ such the boundary of $A$ is the union 
of  three lunes, one  of which  is infinite,  the  second is a radiant 
finite lune and lies in $\hideal{M}$ and the third is a radiant finite 
lune meeting $M_h$. (See Figure \ref{fig:crescent-cone} and  
Example \ref{exmp:crscone}.)  
\end{defn} 
 
\begin{exmp}\label{exmp:crscone} 
Let $x, y, z$ be the standard coordinates of $\bR^3$.  
The holonomy cover of $\mathcal{E}_2$ for the three-dimensional case 
equals $U - l$ where $U$ is the upper half-space given by $z \geq 0$  
and $l$ the $z$-axis. The closure of the subset of $U$ given by  
$ax + by > 0$ for $a, b \in \bR$ not both zero, is  
a crescent-cone. In fact the closure of the component 
missing $l$ of $U^o$ with a subspace of $\bR^3$ meeting $U^o$
removed is a crescent-cone in the Kuiper completion 
$\hat{\mathcal{E}}_2$ of the holonomy cover $U-l$ of $\mathcal{E}_2$. 
When $N$ is a Benz\'ecri suspension of a closed real projective surface  
$\Sigma$, the Kuiper completion of the holonomy cover $N_h$ of $N$  
is homeomorphic to $(\che \Sigma_h \times I)/\sim$ 
for an interval $I = [0, 1]$ and the Kuiper completion $\che \Sigma_h$ 
of a holonomy cover $\Sigma_h$ of $\Sigma$
where the equivalent relation $\sim$ identifies
$\Sigma \times \{0\}$ to a point. Here crescent-cones are  
precisely the closures of crescents times $(0,1)$. Recall that  
crescents correspond to annuli in $\Sigma$. Hence, crescent-cones  
correspond to the Benz\'ecri suspensions of the annuli. 
\end{exmp} 
 
\typeout{<<rfcr.eps>>}  
 
We will now begin the pull-back process as in \cite{cdcr1}.  
(To see examples, see Example \ref{exmp:pullback}). 
As a first step, we have the following: 
\begin{prop} If $M$ is not convex, then there exists a radiant  
tetrahedron $F$ with two sides  $F_1, F_2$ not in $\hideal{M}$   
and  $F_4$ in $\hIideal{M}$ and  the  remaining one $F_3$ such   
that  $F_3^o$ meets both $M_h$ and $\hideal{M}$ and  
$F_1 - (\{O\} \cup F_4)$ and $F_2 - (\{O\} \cup F_4)$ 
are subsets of $M_h$.  
\end{prop} 
\begin{proof} 
Let $\SI^2_1$ denote the sphere of $\bdd$-radius $\pi/4$ in $\mathcal{H}$ 
with center $O$. Consider $C = \dev^{-1}(\SI^2_1) \cap M_h$.  
Then $C$ is connected since $C$ meets every flow line in $M_h$. 
Since $\SI^2_1$ has a flat real projective structure induced  
from the radial projection $\brto \ra \SIT$, the surface $C$ 
has a real projective structure with geodesic boundary. 
 
Suppose that $C$ is not convex.  
$\dev|C: C \ra \SI^2_1$ can be considered a developing map  
of the real projective surface $C$. By inducing the metric 
from the Riemannian metric $\mu$ on $\SI^3$ restricted on  
$\SI^2_1$ to $C$ and completing the induced path-metric 
on $C$ (see \cite{psconv}). Let $\che C$ denote  
the Kuiper completion of $C$ and $\ideal{C}$ the set of 
ideal points. Since $C$ is not convex, there exists  
a real projective triangle $H$ in $\che C$ with an edge  
$e$ of $H$ satisfying 
$H \cap \ideal{C} = e^o \cap \ideal{C} \ne \emp$ 
by Theorem A.2 in \cite{psconv}. 
 
Since the inclusion map $C \ra M_h$ is a quasi-isometry, 
it follows that the $\che C$ may be regarded as a subset of 
$\che M_h$ so that $\ideal{C}$ is included in $\hfideal{M}$. 
 
Let $J$ be the set $H \cap C$. The interior of the union $J'$ of  
all radial lines in $M_h$ through $J$ is obviously  
a radiant open tetrahedron in $M_h$, and the closure $F$  
in $\che M_h$ is a radiant tetrahedron with a side  
$F_4 \subset \hIideal{M}$. 
(To see this, apply $\Phi_{h, t}$ to an ``$\eps$-thickened''   
$J$ for $t \in (-\infty, \infty)$ and take a union and apply  
Proposition \ref{prop:extmap}.) 
Since two edges $e_1$ and $e_2$ of $H$ other than $e$  
are in $M_h$, it follows that two radiant sides $F_1$ and  
$F_2$ of $F$ with $O$ and $\hIideal{M}$ removed are  
subsets of $M_h$ respectively. Since $e$ is not a subset of $M_h$,  
it follows that the face $F_3$ corresponding to $e$ meets  
$M_h$ but is not a subset of $M_h$. Therefore, we found 
the desired radiant tetrahedron $F$.  
 
If $C$ is convex, 
then $C$ is projectively diffeomorphic to a quotient of  
a convex domain in $\SI^2_1$ by Proposition A.2 of \cite{psconv}. 
As $C$ admits a developing map $\dev|C$ to $\SI^2_1$, it follows  
that $C$ is projectively diffeomorphic to a convex domain  
in $\SI^2_1$ by $\dev|C$. By a classification of convex sets  
in $\SI^2_1$ (see \cite{psconv}), $\dev(C)$ is either 
$\SI^2_1$ itself, a hemisphere, a convex subset  
of a $2$-hemisphere of $\SI^2_1$. The final possibility  
implies that $M_h$ is convex. 
Hence, $C$ is projectively diffeomorphic to $\SI^2_1$ 
or a hemisphere, and $M_h$ to $\brto$ or a half-space  
with the origin removed by $\dev$. Thus $\che M_h$ can be identified  
with $\mathcal{H}$ or a closed half-space by $\dev$ where  
$\hideal{M}$ is identifiable to $\{O\}$ union with 
$\SIT$ or a hemisphere in $\SIT$. 
We can easily find a tetrahedron in $\che M_h$ 
showing that $M_h$ is not $2$-convex (see Lemma \ref{lem:dconv}  
and Chapter 4 of \cite{psconv}).  
\end{proof} 
 
\begin{figure}[t] 
\centerline{\epsfxsize=4in \epsfbox{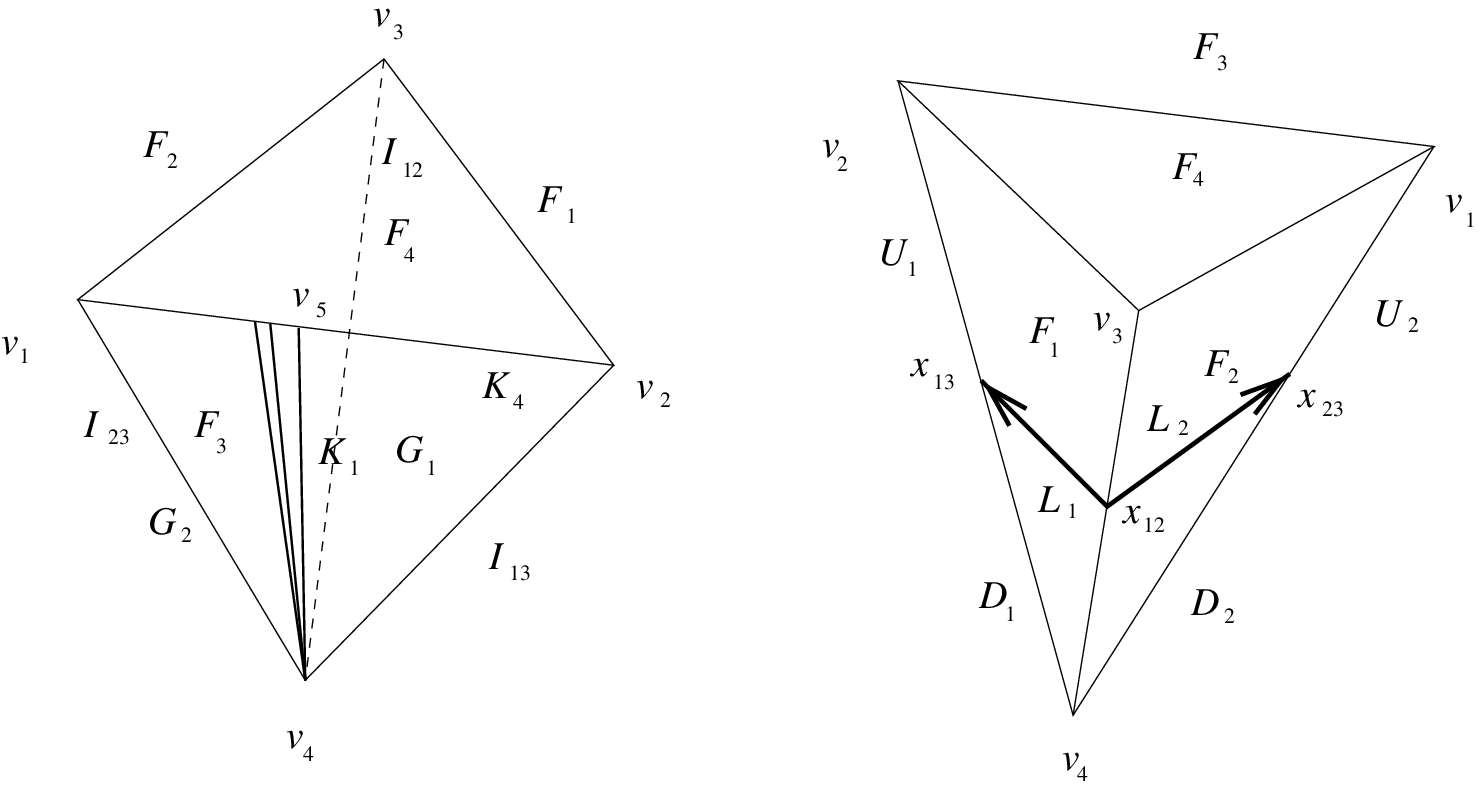}}  
\caption{\label{fig:stetra} $F$, its sides, vertices, $L_1$, and $L_2$ 
viewed from two different points. 
(For convenience,    we orient $L_1$  towards 
$x_{13}$ and  $L_2$  towards  $x_{23}$  and  draw arrows   to indicate 
directions.)}  
\end{figure} 
\typeout{<<rfs1.eps>>}  
 
{\em We now  produce notation which  will be used  extensively in this 
paper}\/;  {\em  they have fixed  meanings in  Chapters 
\ref{ch:2conv}-\ref{ch:rtrh}}.  Let $v_i$ 
denote  the vertex of $F$  opposite $F_i$ for each  $i =  1, 2, 3, 4$, 
with $v_4  =  O$; let  $I_{ij}$  for  $i \ne  j$  denote the   edge of 
intersection of sides $F_i$ and $F_j$. From above $F_3 \cap \hideal{M}$ is 
the  union of  rays from  the origin  and the  edge in $\hIideal{M}$. 
Hence $F_3  - \hideal{M}$ has two  components $G_1$  and $G_2$ meeting 
with $F_1$ and $F_2$ at edges  $I_{13}$ and $I_{23}$  respectively.   
$\clo(G_1) \cap  \hideal{M}$ is  the union of 
two segments  $K_1$ and $K_4$ where $K_4$  equals $\clo(G_1) \cap F_4$ 
and  $K_1^o$  is a   radial line  in  $F_3^o$.  Let  $v_5$ denote the 
intersection point of $K_1$ and $K_4$.  The boundary of $\clo(G_1)$ is 
the union of $I_{12}, K_1$ and $K_4$.  Let  us draw a segment $L_1$ on 
$F_1$   with  endpoints  $x_{13}$  in   $I_{13}^o$  and  $x_{12}$   in 
$I_{12}^o$, and a  segment $L_2$ in  $F_2$ with  endpoints $x_{23}$ in 
$I_{23}^o$ and $x_{12}$. We may assume  that $x_{ij}$ are midpoints of 
$I_{ij}$ with  respect to $\bdd$ for  each pair $i, j$, $i=1,\dots, 4$ 
and $j=1,\dots,4$. The closures of 
components of $F_1 - L_1$ are a quadrilateral and a triangle, which we 
denote by  $U_1$  and $D_1$  respectively.    Similarly we denote  the 
closures of components of $F_2 - L_2$ by $U_2$ and $D_2$ respectively, 
and finally denote by $U$ the set $U_1 \cup  U_2$ and $D$ the set $D_1 
\cup D_2$.  (See Figure \ref{fig:stetra}.)  
 
Recall the metric $d_M$ on $M_h$ induced from the path-metric $d_M$  
on $M$ obtained from the Riemannian metric on $M$. 
For each  point $x$ of $G_1$,  $d_M(x,  D\cap M_h)$ is  realized by a 
point  $y \in D  \cap M_h$ as $D \cap  M_h$ is a  properly imbedded in 
$M_h$; similarly, $d_M(x, U \cap M_h)$ is realized by a  point $y \in U 
\cap M_h$.  It follows from these facts that  $d_M(\cdot, D \cap M_h)$ 
and $d_M(\cdot, U \cap M_h)$ are continuous functions on $G_1$.  
Let $\mathcal{E}$ be the  set of points $x  \in G_1$  satisfying $d_M(x, D 
\cap M_h) = d_M(x, U \cap M_h)$.

For a point $x$ in the interior of the line $D \cap G_1 \cap 
M_h -L$, we obtain $d_M(x,  D \cap  M_h) < d_M(x,   U \cap M_h)$  (see 
Figure  \ref{fig:stetra}).  Similarly, for a  point $x$  in that of $U 
\cap G_1 \cap M_h  - L$, we  have $d_M(x,  U\cap  M_h) <  d_M(x, D\cap 
M_h)$.  Since $\mathcal{E}$ separates $G_1$ into two nonempty subsets where 
$d_M(x, D \cap M_h) > d_M(x, U\cap M_h)$ holds and where $d_M(x, D 
\cap M_h) < d_M(x, U \cap M_h)$ holds  respectively, $\mathcal{E}$ is not a 
compact subset of $G_1$.  Therefore, there exists a sequence of points 
$\ini^i$ converging  to  a  point  $\ini^\infty$ of   $\clo(G_1)  \cap 
\hideal{M} = K_1 \cup K_4$ (under the metric $\bdd$). 
Since $U \cap M_h$ and $D \cap M_h$ are properly imbedded submanifolds,  
we have a sequence $\{\upp^i\}$   in $U$ and $\{\dow^i\}$ in $D$ so that  
\[d_M(\ini^i, \upp^i) = d_M(\ini^i, \dow^i).\] 
 
Since  $\upp^i$  and $\dow^i$   belong to  the compact  subset  $F$ of 
$\che{M_h}$, we may assume without loss of generality that $\upp^i \ra 
\upp^\infty$ and $\dow^i \ra \dow^\infty$ for points $\upp^\infty$ and 
$\dow^\infty$ in $F$. Choose the  closure $\mathcal{W}$ of a fundamental  
domain  of $M_h$  and a deck transformation $\deck^i$ such that  
$\deck^i(\ini^i) = \fin^i$ for $\fin^i \in   \mathcal{W}$. 
 
We define  
\begin{eqnarray}  
F^i = \deck^i(F), & F_j^i = \deck^i(F_j), & G_1^i = \deck^i(\clo(G_1)), 
\nonumber 
\\ I_{jk}^i    =\deck^i(I_{jk}), & L_l^i    =  \deck^i(L_l), &   
K_1^i   = \deck^i(K_1), \\   
K_4^i  = \deck^i(K_4), & U_l^i  =  \deck^i(U_l), & \hbox{ and }  
D_l^i = \deck^i(D_l) \nonumber  
\end{eqnarray}  
for each $j, k = 1,2,3$ and $l=1,2$, and  define  $x_{jk}^i$ to be 
$\deck^i(x_{jk})$ and $v_l^i$ to be $\deck^i(v_l)$.  
 
We will  need the  material on the  Hausdorff convergence   of compact 
subsets of $\SI^n$ (see Appendix \ref{app:seqballs}).  
\begin{hyp}\label{hyp:fundseq} 
By choosing  subsequences,   we may  assume  the  following conditions. 
(Recall that $h(\deck_i)\circ \dev = \dev \circ \deck_i$.)  
\begin{itemize} 
\item $\fin^i$ converges to a point  $\fin^u$ of $\clo(\mathcal{W})$, i.e., 
$\fin^u \in M_h$.  
\item  
\[ 
\begin{array}{c} 
\begin{array}{cc}  
\dev(F^i) = h(\deck^i)(\dev(F)), & \dev(F_j^i) = h(\deck^i)(\dev(F_j)), 
\\  \dev(I_{jk}^i)   =  h(\deck^i)(\dev(I_{jk})),&   \dev(L_l^i)     = 
h(\deck^i)(\dev(L_l)),  
\end{array} \\ 
\dev(G_1^i) = h(\deck^i)(\dev(\clo(G_1))), \\  
\begin{array}{cc}  
\dev(K_1^i)        =    h(\deck^i)(\dev(K_1)),   &      \dev(K_4^i)  = 
h(\deck^i)(\dev(K_4)),  \\   \dev(U_l^i) =   h(\deck^i)(\dev(U_l)),  & 
\hbox{and } \dev(D_l^i) = h(\deck^i)(\dev(D_l))  
\end{array} 
\end{array} 
\] 
for each  $j, k  =  1,2,3$ and  $l =   1,2$ converge geometrically  to 
compact convex sets  
\begin{equation}\label{eqn:objconv} 
F^\infty, F_j^\infty,    I^\infty_{jk},       L_l^\infty,  G^\infty_1, 
K_1^\infty, K_4^\infty, U_l^\infty, \hbox{ and } D_l^\infty  
\end{equation} 
in $\SI^3$ respectively.  
\item The sequence of outer-normal vectors $n^i$ at $\fin^i$ converges to 
a unit vector at $\fin^u$. 
\item $h(\deck^i)(\dev(x_{jk}))$ and  $h(\deck^i)(\dev(v_l))$ converge 
to points $x_{jk}^\infty$ and $v_l^\infty$ respectively for each $j, k 
= 1,2,3$ and $l = 1, 2, 3, 4, 5$.  
\end{itemize} 
\end{hyp} 
 
It follows easily that the dimension  of $F^\infty$ is $\leq 3$, those 
of $F_j^\infty, G_1^\infty,   U_l^\infty,$ and $D_l^\infty$  $\leq 2$, 
and those of $I_{jk}^\infty$, $K_1^\infty$, $K_4^\infty$, $L_l^\infty$ 
and  $U_l^\infty$ $\leq  1$. These are all balls of such dimensions 
by Proposition 2.8 of \cite{psconv}. 
(See   the  Appendix \ref{app:seqballs}).  
 
For example, $I^\infty_{jk}$ is either a segment of $\bdd$-length   
$\leq \pi$ or a point since the limit $I^\infty_{jk}$ is convex.  
$I^\infty_{jk}$ for $j, k$ not equal to $4$ is a radial segment always. 
Since $F_l^\infty$ is a cone over $I^\infty_{l4}$ for $l =  1, 2, 3$  
by Lemma \ref{lem:conelim}, it follows that $F_l^\infty$ is either a radial  
segment, a radiant triangle, or a radiant lune for $l =1,2,3$.  

\begin{figure}[t] 
\centerline{\epsfxsize=4in \epsfbox{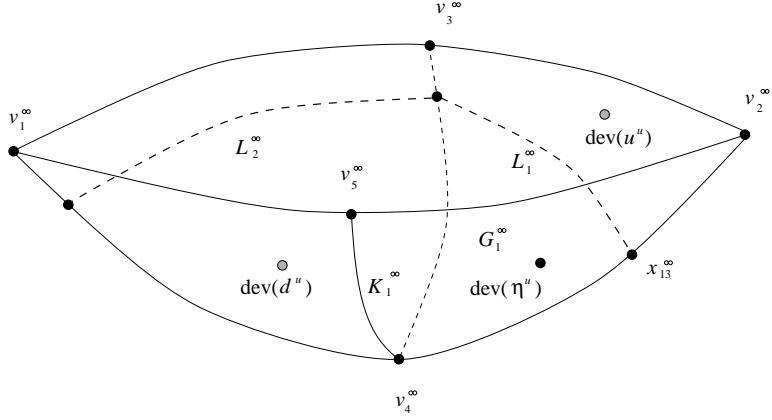}}  
\caption{\label{fig:typ} An example of a limit.}  
\end{figure} 
\typeout{<<rftyp.eps>>}

\begin{exmp}\label{exmp:pullback} 
We may do this process in $\mathcal{E}_2$ 
(see Example \ref{exmp:hopf}) in the three-dimensional case.  
Its holonomy cover equals $U - l$ where $U$ is the upper half-space  
and $l$ the $z$-axis. Then we can easily construct each of  
the objects in the process above. A total cross-section  
to the radial flow comes from in an affine plane parallel to  
the $xy$-plane. Hence components of the inverse image of  
the total cross-sections are affine planes parallel to  
the $xy$-plane. We may choose $L_1$ and $L_2$ to lie  
in a component $J$  of such an inverse image. 
In this case, the set to be denoted by $\mathcal{E}$ 
may be considered also to be the intersection of $G_1$ with $J$. 
(Keep this example in mind for the rest of the process.) 
 
Similar examples can be obtained from Benz\'ecri suspensions  
of nonconvex real projective surfaces. 
\end{exmp} 
 
\begin{figure}[t] 
\centerline{\epsfxsize=2.0in \epsfbox{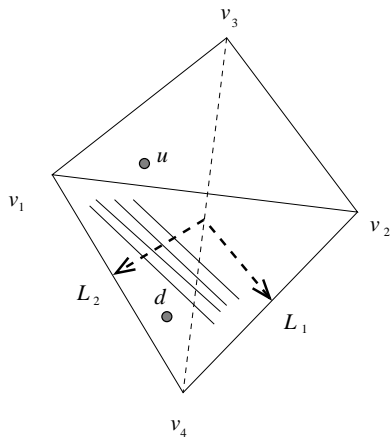}}  
\caption{\label{fig:nond} The segments indicate $\vth^{i, -1}(s_i)$.}  
\end{figure}

\begin{prop}\label{prop:Fnond} 
$F^i$ includes a common convex open ball $\mathcal{P}$ for sufficiently  
large $i$, and $F^\infty$ is a   $3$-ball\/{\rm  ;} that is,  
$h(\deck^i)(\dev(F))$ does  not degenerate into a convex  
set of dimension less than or equal to $2$.  
\end{prop} 
\begin{proof} 
Let  $B(\fin^\infty)$ be a tiny  ball of $\fin^\infty$, and $2\gamma$ be 
the   positive   number    equal   to   
$d_M(\fin^\infty,    \Bd B(\fin^\infty))$.  
For $i$ sufficiently large,  $d_M(\fin^\infty, \fin^i) < \gamma$. 
We assume this holds for $i$ in this proof. 
 
Suppose we have either $d_M(\ini^i, \upp^i) < \gamma$ 
or $d_M(\ini^i,  \dow^i) < \gamma$ for infinitely  many $i$,  either of 
of which means that $d_M(\fin^i, \deck^i(\upp^i)) < \gamma$ {\em and}\/ 
$d_M(\fin^i,     \deck^i(\dow^i))   <   \gamma$    by   equation 
\ref{eqn:eqdst}.     Since   for   each  $i$,    $\deck^i(\upp^i)$ and 
$\deck^i(\dow^i)$ belong to  $B(\fin^\infty)$, there exists a  segment 
$s_i$, $s_i \subset B(\fin^\infty)$ with  these points as endpoints by 
convexity of  $B(\fin^\infty)$; the segment  $\deck^{i, -1}(s_i)$ is a 
segment in the compact set $F$ with endpoints  
$\upp_i$  and $\dow_i$.   By choosing a subsequence if necessary,   
we  \awlg that the sequence  of  $\deck^{i,-1}(s_i)$ converges  
to a  segment  $s$ with  endpoints in $U$  and $D$ 
respectively. Such a segment $s$ must  have a point $t$ belonging 
to $F_1 \cap M_h$ if endpoints are in $U_1$ and $D_1$ respectively; to 
$F_2 \cap M_h$ if endpoints  in $U_2$ and  $D_2$ respectively; and  to 
$F^o \cup (F_3^o  - \hideal{M})$ otherwise by considering their 
images under $\dev$ and  an elementary geometry  argument (see Figure  
\ref{fig:nond}). 
 
Since  $B(\fin^\infty)$ is  a compact subset   of $M_h$, the images of 
$\deck^{i, -1}(s_i)$ may  not intersect  a  compact subset of  $F \cap 
M_h$   infinitely many  times by Lemma    \ref{lem:finiteness}.  But a 
compact neighborhood  of $t$ in   $M_h$ intersects infinitely  many of 
these, which  is a contradiction.    Therefore, for sufficiently large 
$i$, we have $d_M(\ini^i, \upp^i) \geq  \gamma$ and $d_M(\ini^i, \dow^i) 
\geq \gamma$; that is, we have $d_M(\fin^i,  \deck^i(F_1 \cup F_2)) \geq 
\gamma$. If $i$ is sufficiently large so that  
$d_M(\fin^u, \fin^i) < \gamma/2$, 
then $\deck^i(F_1 \cup  F_2 \cup  F_4)$ does not meet the interior  
of the $\gamma/2$-ball $A$ of $\fin^u$ in $B(\fin^u)$.  
 
Since  $\fin^i$    converges     to    $\fin^u$,   and 
$\deck^i(F_1 \cup    F_2  \cup F_4)$  does  not   intersect the 
interior of $\dev(A)$, Lemma  \ref{lem:common} implies the conclusion. 
\end{proof} 
 
\typeout{<<rftyp.eps>>}  
 
\begin{lem}\label{lem:finiteness} 
Let $\phi^i$,  $i=1,2, \dots$,  be   a   sequence  of     distinct 
deck transformations  of $M_h$, and   $K$ a compact  subset  of $M_h$. 
Then $\phi^i(K)$ intersects a compact ball  neighborhood of a point of 
$M_h$ for only finitely many $i$.  
\end{lem} 
\begin{proof} This follows since the action of the deck transformation  
group is properly discontinuous.  
\end{proof} 
 
Recall that a sequence of compact convex triangles in $\SI^2$ converges 
to a point, a segment, a   hemisphere,  a lune, or a   simply convex 
triangle (see  Chapter 2 of \cite{psconv} and the Appendix    
of  \cite{cdcr1}.)  Since $F^\infty$ is a nondegenerate $3$-ball  
by Proposition \ref{prop:Fnond}, and $F^\infty$ is  a cone over  
$F_4^\infty$  by Lemma \ref{lem:conelim} 
(see Definition  \ref{defn:cone} and  Remark \ref{rem:cone}),  
$F_4^\infty$ is  either a hemisphere, a lune,  
or a simply convex triangle.  By   Remark \ref{rem:cone3},  
$F^\infty$  is either a radiant bihedron, a radiant trihedron,  
or a radiant tetrahedron.  
 
Similarly    to   the proof of Theorem 4.6 in \cite{psconv},  
Theorem \ref{thm:seqconv} implies that there exists a convex $3$-ball   
$F^u$ in $\che M_h$ such that $\dev(F^u)$ equals  $F^\infty$   
since we have a sequence of $3$-balls $F^i$ in $\che M_h$   
including the common  ball $\mathcal{P}$ by above Proposition  
\ref{prop:Fnond}.  
 
Since  $\dev| F^u$ is  an imbedding  onto  $F^\infty$,  $F^u$ includes 
subsets   \[ F^u_i, I^u_{jk},    L_l^u,  U_l^u, D_l^u \]  respectively 
corresponding to  
\[ F_i^\infty, I_{jk}^\infty, L_l^\infty, U_l^\infty, D_l^\infty \]  
for each $i,  j, k =  1, 2, 3, 4$, $j  \ne k$, and $l  = 1, 2$.  Also, 
there are points $x^u_{jk}$  corresponding to $x^\infty_{jk}$ for each 
$j, k   = 1, 2,   3$, $j \ne k$  and  points  $v^u_j$  corresponding to 
$v_j^\infty$ for $j = 1, 2, 3, 4, 5$. (See Figure \ref{fig:typ}.) 
  
From the knowledge we have of objects with superscripts $\infty$, we  
gain informations about ones with superscripts $u$. For example,  
$I_{jk}^u$ is a segment of $\bdd$-length $\leq \pi$ or a point,  
and $F_{i}^u$ for $i =1,2,3$ is a radial segment, a radiant triangle, or 
a radiant lune. 
 
For the purposes of the later chapters, we state here: 
\begin{rem}\label{rem:endL} 
Notice that $I_{12}^u, I_{13}^u$, and $I_{23}^u$ are radial segments. 
We have $I_{12}^u = I_{13}^u$ or $I_{12}^u \cap I_{13}^u = \{O\}$  
respectively.  Similarly, we have $I_{12}^u = I_{23}^u$ or  
$I_{12}^u \cap I_{23}^u =\{O\}$. 
The following facts are useful: $L_1^u$  either equals $O$, is a 
segment with   one endpoint in   $I_{12}^u$  and  the another  one  in 
$I_{13}^u$, or  a  point in  $I_{12}^u$ where  $I_{12}^u  = I_{13}^u$. 
$L_2^u$ equals $O$, is  a segment with one  endpoint in $I_{12}^u$ and 
the another   one  in $I_{23}^u$,   or a  point in   $I_{12}^u$  where 
$I_{12}^u =  I_{23}^u$.  The  endpoints  of  $L_1^u$  and $L_2^u$  in 
$I_{12}^u$ are identical.  These follow from looking at  $L_1^\infty$, 
$L_2^\infty$, and so on. (Here, the endpoint of a point is itself.) 
\end{rem}

If  $F^\infty$ is  a radiant   bihedron,  then  $F^u$  has a  side  in 
$\hideal{M}$, and since $F^u$  is not  equal to  $\che M_h$,  it follows 
that   a side   of   $F^u$ intersects  $M_h$   and  $F^u$   must be  a 
$3$-crescent.  Since  $3$-crescents do not exist in $\che{M_h}$  
by  Hypothesis~\ref{hyp2},  we obtain: 

\begin{prop}\label{prop:Fu} 
$F^u$ is either a radiant trihedron or a radiant tetrahedron.  
$F^u$ includes $\mathcal{P}$. \qed 
\end{prop} 
 
\begin{proof}[The proof of Theorem \ref{thm:crscone}] 
Since $(F_1 \cup F_2) \cap M_h$ equals $(U \cup D) \cap M_h$, 
by Proposition \ref{prop:MainEqu}, i.e. the claim in the introduction,  
we can choose a sequence  of  points 
$\ini_i$  in $\mathcal{E}$ so  that $d_M(\ini_i, (F_1  \cup F_2) \cap M_h)$ 
converges  to $+\infty$. Using    the  pull-back argument as   in 
Chapter \ref{ch:2conv}, we assume  Hypothesis  \ref{hyp:fundseq}.  By 
Theorem \ref{thm:seqconv}, we obtain  a  radiant tetrahedron $F^u$  in 
$\che{M_h}$ with sides $F_1^u,   F_2^u,  F_3^u,$ and $F_4^u$ or   a 
radiant trihedron $F^u$ with boundary $F_1^u \cup F_2^u \cup F_3^u  
\cup F_4^u$. But as $F_1^i$ and $F_2^i$ are ideal  sequences by Proposition 
\ref{prop:MainEqu}, it  follows  from Theorem  \ref{thm:seqconv}   that 
$F_1^u \cup F_2^u$ is a subset of $\hideal{M}$.

\begin{defn}\label{defn:pcrc}  
A  radiant tetrahedron  $T$ in $\che  M_h$   which has  three sides in 
$\hideal{M}$ is said to be a {\em pseudo-crescent-cone}.  
\end{defn} 
 
By Lemma \ref{lem:cosub}, there is a side  of  $F^u$ meeting $M_h$.   
Since $F_4^i$ is a subset of $\hIideal{M}$, $F_1^u \cup F_2^u \cup F_4^u$  
is a subset of $\hideal{M}$; $F_3^u$ is the only totally geodesic disk  
in $\delta F$ intersecting with $M_h$.  Thus  if $F^u$  is a  radiant  
trihedron, then $F^u$ is a crescent-cone (see  Definition   
\ref{defn:crscone}). If $F^u$ is  a radiant   tetrahedron,  then   
$F^u$ is a pseudo-crescent-cone  by Definition \ref{defn:pcrc}.  

In Chapter  \ref{ch:pcrc},  we will show that $\che  M_h$ cannot  
have a  pseudo-crescent-cone  since $M$ is not convex.  
\end{proof}

\chapter{The claim and the rooms}\label{ch:claim} 
 
We will be proving the {\em claim}\/ in the introduction,  
which will be proved in this chapter and Chapters  
\ref{ch:tric} and \ref{ch:rtrh} completely: 
We assume otherwise and find  contradictions.  
The purpose of this chapter is to list properties  
which will be used in Chapters \ref{ch:tric} and \ref{ch:rtrh}:  
Propositions \ref{prop:U_lD_l},    
\ref{prop:room},\ref{prop:notsamel}, and \ref{prop:roomG1}, and 
Remark \ref{rem:UDdisk}.  
 
\begin{prop}\label{prop:MainEqu} 
The supremum of $d_M(x, D\cap M_h)$ on $\mathcal{E}$ is infinite.  
\end{prop} 
 
We will prove this by contradiction: 
we assume the following hypothesis from now on  
(in this chapter and Chapters \ref{ch:tric} and \ref{ch:rtrh}):  
\begin{hyp}\label{hyp:supr} 
Let $C_{ud}$ denote the finite supremum of  $d_M(x, D\cap M_h)$ for $x 
\in \mathcal{E}$.  
\end{hyp} 
 
To show contradiction, let $\ini^i$ be a sequence of points 
converging  to  a  point  $\ini^\infty$ of   $\clo(G_1)  \cap 
\hideal{M} = K_1 \cup K_4$ (under the metric $\bdd$)  
as in the above chapter satisfying the Hypothesis \ref{hyp:fundseq}. 
$d_M(\ini^i, D\cap M_h)$ is  bounded  above; we have a sequence  
$\{\upp^i\}$   in $U$ and $\{\dow^i\}$ in $D$ so that  
\begin{equation}\label{eqn:eqdst} 
d_M(\ini^i, \upp^i), d_M(\ini^i, \dow^i) < C_{ud} + 1 \hbox{ for all } 
i, \hbox{ and } d_M(\ini^i, \upp^i) = d_M(\ini^i, \dow^i).  
\end{equation}  
 
We will use the notations of previous chapters but we will have  
some additional requirements: 
Let $\mathcal{W}'$  the  $C_{ud}+1$ $d_M$-neighborhood   of  $\mathcal{W}$,    
to  which  $\deck^i(\upp^i)$ and $\deck^i(\dow^i)$ belong.  
In addition to Hypothesis \ref{hyp:fundseq}, we will also require: 
\begin{itemize} 
\item   $\deck^i(\upp^i)$   and $\deck^i(\dow^i)$  converge  to points 
$\upp^u$ and $\dow^u$ in $\mathcal{W}'$ respectively.  
\end{itemize} 
From the previous chapter, $\dev(F^i)$ converges to $F^\infty$,   
$F^i$ shares a common open ball $\mathcal{P}$, and $\che M_h$  
includes a radiant tetrahedron or radiant trihedron $F^u$ 
including $\mathcal{P}$ so that $\dev(F^u) = F^\infty$. All associated 
objects are the same as in the chapter. 
 
Since $M$ has empty or totally geodesic boundary,  
each  point $x$ of $M_h$ has an open coordinate chart $U$ such that  
a lift  $\phi: U \ra \SI^n$ of a chart $U \ra \rp^n$  
so that $\phi(U)$ is a convex set in $\mathcal{H}^o$.   
Hence,   $U$  includes   a compact    convex   $n$-ball 
neighborhood $B$  of $x$   such that $\phi|B$   is  an imbedding  onto 
$\phi(B)$.  Since  $\dev$ is a continuation  of charts, $\dev| B$ is a 
restriction of a chart, and  hence, $\dev| B$  is an imbedding onto  a 
convex compact  $n$-ball.   We   say  that  $B$  a  {\em  tiny    ball 
neighborhood}\/ of $x$. We can easily show that $B$ can be chosen  
so that $\dev(B)$ is a $\bdd$-ball of certain radius, perhaps  
intersected with a closed affine half-space in $\mathcal{H}^o$ if $B$ meets  
the totally geodesic boundary $\delta M_h$. 
 
Next, we give some information on $\fin^u, \dow^u,$ and $\upp^u$. 
Recall  that every  point $x$  of $M_h$  has  a tiny-ball neighborhood 
$B(x)$.  Since  $\mathcal{W}'$ is  compact and $M_h$ has  
totally geodesic boundary, using Lebesgue number,  
we  may assume without loss of  generality that each  point  
$x$  of  $\mathcal{W}'$ has  a tiny-ball neighborhood $B(x)$  such that the   
image $\dev(B(x))$  is a $\bdd$-ball of constant radius $\eps$, $\eps > 0$,  
perhaps intersected with a closed half-space of $\mathcal{H}^o$.  
 
\begin{lem}\label{roomexists} 
Points $\fin^u$, $\upp^u$, and  $\dow^u$ have tiny-ball  neighborhoods 
$B(\fin^u)$,     $B(\upp^u)$   and  $B(\dow^u)$.    $\dev(B(\fin^u))$, 
$\dev(B(\upp^u))$,   and $\dev(B(\dow^u))$ are  uniformly bounded away 
from $O$ and $\SIT$.  
\end{lem}  
\begin{proof} 
$\dev(M_h)$  is disjoint from $O$ by  Lemma \ref{lem:radcom}.  Since 
$\dev(M_h)$ is  a subset of  $\mathcal{H} -(\SIT \cup  \{O\})$,  
$\dev(\mathcal{W'})$ is  a compact subset of   
$\mathcal{H} - (\SIT  \cup \{O\})$, and the last statement follows.  
\end{proof} 
 
By an upper component of $F_l^u -L_l^u$, we mean one further away from $O$ 
than $L_l^u$; and by a lower component, we mean the one closer to $O$ than  
$L_l^u$.  
\begin{prop}\label{prop:U_lD_l} 
\begin{itemize} 
\item Suppose  that $F_l^u$ is  a radiant triangle  or a radiant lune. 
$U_l^u$ equals  
\begin{itemize} 
\item either the  union of the upper component  of $F_l^u - L_l^u$ and 
$L_l^u$, or  equals $L_l^u$ itself when $F_l^u  -  L_l^u$ has only the 
lower component\/{\rm ;} and  
\item $D_l^u$ equals either the union of the lower component of $F_l^u 
- L_l^u$  and $L_l^u$, or  equals $L_l^u$ itself  when $F_l^u - L_l^u$ 
has only the upper component.  
\end{itemize} 
\item If $F_l^u$ is a radial segment, then  
\begin{itemize}  
\item  either $U_l$  is  the closure of $F_l^u  -  \{p\}$ for a  lower 
endpoint $p$ of $L_l^u$  if  $p$ does not equal the upper  endpoint $p$   
of $F_l^u$, or $U_l$ equals the upper endpoint of $F_l^u$ 
otherwise \/{\rm ;} and  
\item either  $D_l$ is  the closure  of $F_l^u  - \{p\}$ for   
an upper endpoint $p$ of $L_l^u$ if $p$ does not equal $O$,  
or $D_l$ equals $\{O\}$ otherwise.  
\end{itemize}  
\end{itemize}  
\end{prop} 
\begin{proof} 
Since  $U_1^i$  is the  closure of   the upper component   of $F_l^i - 
L_l^i$, the geometric limit $U_l^\infty$ is as stated above.  
\end{proof} 
 
\begin{lem}\label{lem:finduuu} 
We have  $\fin^u \in F^u,$  and $\dow^u, \upp^u  \in  F^u$.  Moreover, 
$\fin^u \in G_1^u$, $\upp^u \in  F^u_l$, and $\dow^u \in F^u_{l'}$ for 
some $l$, $l' = 1,2$. 
\end{lem} 
\begin{proof} 
We know that $\upp^i \in F^i$ for each $i$. Since $\upp^i \ra \upp^u$, 
it follows  that  $\upp^i \subset  \inte  B(\upp^u)$ whenever $i$ is 
sufficiently large for  a tiny  ball $B(\upp^u)$   of $\upp^u$.   
Hence,  $F^i$  and $B(\upp^u)$ overlap and  
$\dev| F^i    \cup B(\upp^u)$ is an   
imbedding   onto $\dev(F^i) \cup \dev(B(\upp^u))$ by Proposition 
\ref{prop:extmap}. Since $F^i$ and $F^u$ both include  
$\mathcal{P}$ for $i$ sufficiently large, $\dev| F^i \cup F^u$ is  
an imbedding onto $\dev(F^i) \cup \dev(F^u)$. 
By choosing $B(\upp^u)$  carefully,    we may  assume   that 
$(\dev(F^i) \cup \dev(B(\upp^u))) \cap  (\dev(F^i) \cup \dev(F^u))$ is 
a  connected   submanifold with  interior  equal  to  $(\dev(F^i) \cup 
\dev(B(\upp^u)))^o \cap (\dev(F^i)   \cup \dev(F^u))^o$, which can  be 
accomplished  by   making these sets star-shaped  from a point of 
$\dev(F^i)$.  
Since Proposition \ref{prop:extmap} shows that 
$\dev| F^i \cup F^u \cup B(\upp^u)$ is  an imbedding, and  
$\dev(\upp^i)  \ra \dev(\upp^u)$, and 
$\dev(\upp^u) \in  \dev(F^u)$,  we conclude that 
$\upp^u \in  F^u$.  Similarly, we  obtain $\fin^u \in  F^u$ and 
$\dow^u \in F^u$.  
 
Since  $\dev(\upp^i) \in \dev(F_m^i)$ for  each  $i$ and some $m$, $m=1,2$, 
depending on $i$, we have  $\dev(\upp^u) \in \dev(F_l^u)$ for 
some  $l$, $l=1,2$.  Since  $\upp^u  \in F^u$, we have  
$\upp^u   \in F^u_l$. Similarly, we obtain $\fin^u \in G_1^u$,  
and $\dow^u \in F^u_{l'}$ for some $l'$, $l' = 1,2$.  
\end{proof} 
 
\begin{lem}\label{lem:seqconvFi} 
Suppose that a point  $x$ of $M_h$  belongs to $F^u$.  Then $\dev| F^u 
\cup B(x)$ is an imbedding  for a choice of the tiny ball $B(x)$ of $x$,  
and $\dev| F^u \cup B(x) \cup F^i$ is an imbedding onto  
$F^\infty \cup B(x) \cup \dev(F^i)$ for $i$ sufficiently large.  
\end{lem} 
\begin{proof} 
We choose  $B(x)$ so that  $\dev(B(x)) \cup F^\infty$ is a star-shaped 
set  from a point of  $\dev(\mathcal{P})$ which is the  common open ball in $F^i$ 
from the proof of Proposition \ref{prop:Fnond}.  The lemma now follows from 
Proposition \ref{prop:extmap}.  
\end{proof}  
 
\begin{prop}\label{prop:room} 
The  set $L_1^u  \cup L_2^u$ is   a subset of  $\hideal{M}$.  The set 
$F_l^u \cap B(\upp^u)$   is a  subset of  $U_l^u$  with nonempty  
relative interior in $F_l^u$ if $\upp^u \in F_l^u$  
and $F_l^u \cap  B(\dow^u)$ is that of  
$D_l^u$ with nonempty relative interior in $F_l^u$ 
if $\dow^u \in F_l^u$. In  particular,   
$U_l^u$ includes a relatively open subset of $F_l^u$ if  
$\upp^u  \in   F_l^u$, and $D_l^u$ includes a relatively open subset  
of $F_l^u$ if $\dow^u \in F_l^u$.  
\end{prop} 
\begin{proof} 
The set $L_1^i  \cup L_2^i$ equals  $\deck^i(L_1 \cup  L_2)$, and $L_1 
\cup L_2$ is a compact subset of $M_h$. If $L_1^u \cup L_2^u$ contains 
a point $x$ of $M_h$, then $x$ belongs to $F^u$.  Choose the tiny ball 
$B(x)$ as  in Lemma \ref{lem:seqconvFi} so  that $\dev|  F^u \cup B(x) 
\cup F^i$ is a  homeomorphism onto their  images.  Since a sequence of 
points $p_i$, $p_i \in \dev(L_1^i \cup L_2^i)$, converges to $\dev(x)$ 
in   $\SI^3$, Lemma \ref{lem:seqconvFi} shows that $\dev(L_1^i)   \cup 
\dev(L_2^i)$ intersects  $\dev(B(x))$ for sufficiently large  $i$, and 
hence so does  $L_1^i \cup L_2^i$ with  $B(x)$. This contradicts Lemma 
\ref{lem:finiteness}, and the first statement follows.  
 
Since $B(\upp^u)$   does not meet   $L_1^u \cup   L_2^u$, so does  not 
$B(\upp^u) \cap F_l^u$.  Since $\upp^u$ belongs to $U_l^u$, $B(\upp^u) 
\cap F_l^u$ is  included in $U_l^u$   by the connectedness of  $U_l^u$ 
implied by  Proposition \ref{prop:U_lD_l}. The same argument works for 
$B(\dow^u)$. The last statement follows easily from this.  
\end{proof} 
 
\begin{prop}\label{prop:notsamel}  
If $\upp^u  \in F_1^u$, then  $\dow^u \in  F_2^u$, and if  $\upp^u \in 
F_2^u$, then $\dow^u \in F_1^u$.  Moreover,  no $L_l^u$ passes through 
$F_l^{u,   o}$  for  each $l$, $l =1,2$.   Finally, $F_l^u$ is a {\rm 
(}nondegenerate{\rm )} radiant triangle or lune for $l$, $l=1,2$.  
\end{prop} 
\begin{proof} 
Recall that $F_l^u$ is either a radial segment, a radiant triangle,  
or a radiant lune for $l =1,2$.  
Suppose that $F_l^u$ is a radiant triangle or lune  
and $\upp^u$ and $\dow^u$ belong to $F_l^u$  for a  given $l$.   
Then  both $U_l^u$  and $D_l^u$  are convex 
$2$-balls with nonempty relative interiors and are  the closures of the 
components of $F_l^u  - L_l^u$.  In  order for this to happen, $L_l^u$ 
passes through   $F_l^{u, o}$ with  endpoints  in the  radial sides of 
$F_1^u$ (see Remark \ref{rem:endL}).  Since this implies that every segment 
from $O$ to a point of $I^u_{l4}$ intersects $L_l^u$, $F_l^{u, o}$ is a 
radiant   set, and $F_l^{u, o} \cap \hideal{M}$ is a radiant set,  
it follows that $F_l^u\subset \hideal{M}$. This contradicts the fact   
that $U_l^u$ contains a  point $\upp^u$, a point of $M_h$.  
(See Figure \ref{fig:typ}.) 
 
We can  show similarly that the situation  where  $F_l^u$ is a radial 
segment and  $\upp^u,  \dow^u \in F_l^u$  does  not  occur also.  This 
implies the first part of the lemma.  
 
If $L_l^u$ passes through $F_l^{u, o}$ for some $l$, then $F_l^u$ is a 
subset of $\hideal{M}$ as shown in the above paragraph. Hence, $F_k^u$ 
for $k  \ne l$, $k  = 1, 2$, must  contain both $\upp^u$ and $\dow^u$, 
points  of $M_h$. This  is a contradiction  by  the first statement of 
this  lemma, and the  second part of the lemma is proved.  
 
Suppose  that   $F_1^u$ is a  radial    segment.  Since $F^u$   is  a 
``nondegenerate''  ball,  $F_2^u$  has to be a radiant lune 
and $F_1^u$ is an edge of $F_2^u$. By above  paragraph,  
$F_2^u$ cannot contain both $\upp^u$  and $\dow^u$.    
But if  $\upp^u$ belongs to $F_1^u$  and $\dow^u$ to $F_2^u$, 
then $\upp^u$ also belongs to $F_2^u$, which is a contradiction. 
Similarly a contradiction  follows if $\dow^u$ belongs  
to $F_1^u$ and $\upp^u$ to $F_2^u$.  
 
Analogously, $F_2^u$ cannot be a radial segment also.  
\end{proof}  
 
\begin{rem}\label{rem:UDdisk} 
By above,  $F_l^u$  are disks  always for $l=1,2$.    
If  $\upp^u  \in F_l^u$, then $U_l^u$ includes  an open subset     
of $F_l^u$   by Proposition \ref{prop:room}, and $U_l^u$ is   
a nondegenerate disk. If $\dow^u  \in F_l^u$, then $D_l^u$ is  
a nondegenerate disk. In particular, $U^u$ and $D^u$ must have   
nonempty interior in the manifold-boundary  $\delta F^u$ of  
the $3$-ball $F^u$.  
\end{rem} 
 
\begin{prop}\label{prop:roomG1} 
There   exists a   compact  nondegenerate $2$-disk $G_1^u$ in $F^u$ 
corresponding to  $G_1^\infty$, $F_3^u$ is a nondegenerate $2$-disk, 
$\delta  G_1^u$ includes  segments $K^u_1$ and $K^u_4$ corresponding   
to $K_1^\infty$ and  $K_2^\infty$. Moreover $K_1^u$ and $K_2^u$ are  
subsets of $\hideal{M}$.  
\end{prop} 
\begin{proof} 
A sequence $\fin^i$ converges to $\fin^u$, which  has a tiny ball 
$B(\fin^u)$. We  assume that $\fin^i \in \inte B(\fin^u)$ for each $i$,  
and $\dev| F^i \cup  B(\fin^u)$ is an imbedding onto  its image by  
Proposition \ref{prop:overlap}.  Let  $J^i$   be  the  maximal totally    
geodesic subsurface in $M_h$ including the totally geodesic subsurface  
$G^i_1 \cap M_h$ so that $J^i \cap   B(\fin^u)$ is a  convex ball   
$D^i$ with boundary  in $\delta B(\fin^u)$.  
 
Since  $K_1^i$ and $K_4^i$  are  subsets of  $\hideal{M}$, they  don't 
intersect $B(\fin^u)$, and only the  interior of the edge $I^i_{13}$ 
of  $G^i_1$ may  meet  with  $B(\fin^u)$.   Let  $l^i$ be the  segment 
$I^i_{13} \cap  B(\fin^u)$ which    has  endpoints  in $D^i  \cap    \Bd 
B(\fin^u)$   if $I^i_{13} \cap B(\fin^u)$  is  not empty; otherwise, let 
$l^i    = \emp$.  Since   no  other   segment  of  $G^i_1$  meets with 
$B(\fin^u)$, and $\fin^i$  belongs to $G^i_1$,  it follows that $G_1^i 
\cap B(\fin^u)$  is the closure  of a component  $D^{\prime i}$ of the 
convex disk  $D^i$ with $l^i$ removed, and the closure of  
$D^{\prime i}$ is a compact convex disk bounded by $l^i$ and an arc  
$\alpha^i$ in $\Bd B(\fin^u)$.  
 
Since    $\bdd(\fin^u, \Bd B(\fin^u))$  is   bounded  from below by   
a positive constant $c$,  we   obtain $\bdd(\fin^u, \alpha^i) \geq   c$. 
Since $\fin^i \ra  \fin^u$, we obtain  $\bdd(\fin^i, \alpha^i)  > c/2$ 
for $i$ sufficiently large.    Since $D^{\prime i}$  contains $\fin^i$ 
for   each  $i$, an   elementary geometry  shows  that  $D^{\prime i}$ 
includes a  $\bdd$-disk of radius $c/2$  for each $i$  (using the same 
idea as in the proof of Lemma  \ref{lem:common})  
implying that $\dev(G_1^i)$ is not  degenerating  
into a segment, the first statement. The rest follows  from  
Theorem \ref{thm:seqconv}.  
\end{proof}

\chapter{The radiant tetrahedron case} 
\label{ch:tric} 
 
The  aim  of  Chapters  \ref{ch:tric} and    
\ref{ch:rtrh} is to show  that Hypothesis 
\ref{hyp:supr} leads to contradictions under the following assumptions:  
\begin{itemize} 
\item $F^u$ is a radiant tetrahedron.  
\item $F^u$  is  a   radiant  trihedron.    
\end{itemize} 
Since these two  are the only possible shapes  of $F^u$ by Proposition 
\ref{prop:Fu}, we     will  have  completed   the   proof   of Proposition 
\ref{prop:MainEqu}.  
 
\begin{figure}[b] 
\centerline{\epsfxsize=3.5in \epsfbox{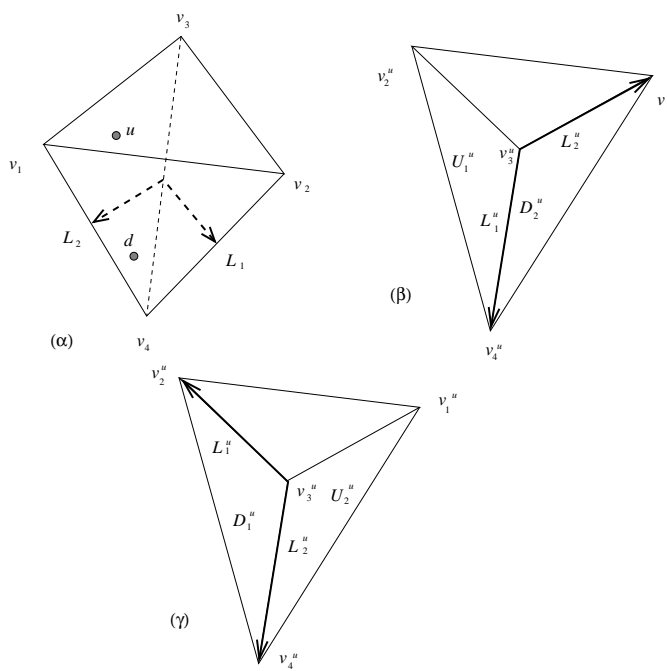}}  
\caption{\label{fig:radseg} The objects of $F$, and two examples of 
$F^u$}  
\end{figure} 
\typeout{<<rfr1.eps>>}  
 
We will show that the first case does  not happen in this chapter by 
singling  out two cases by  considering $L_l^u$  and $U_l^u$ and $D_l^u$ 
using mainly Propositions  \ref{prop:room}  and  \ref{prop:notsamel}. 
Then we show that  the two cases cannot  happen by showing  that their 
existence  implies the existence of radiant bihedra, which was ruled 
out by our hypothesis \ref{hyp2}.  We obtain the radiant bihedra 
by  putting  our tetrahedra  in standard positions  and estimating the 
eigenvalues  of  holonomy   action  asymptotically from  the positions  
of $L_1^u$ and $L_2^u$. (The proof in this case will be a prototype  
of all later cases which will occur.)  
 
\begin{figure}[t] 
\centerline{\epsfxsize=3.5in \epsfbox{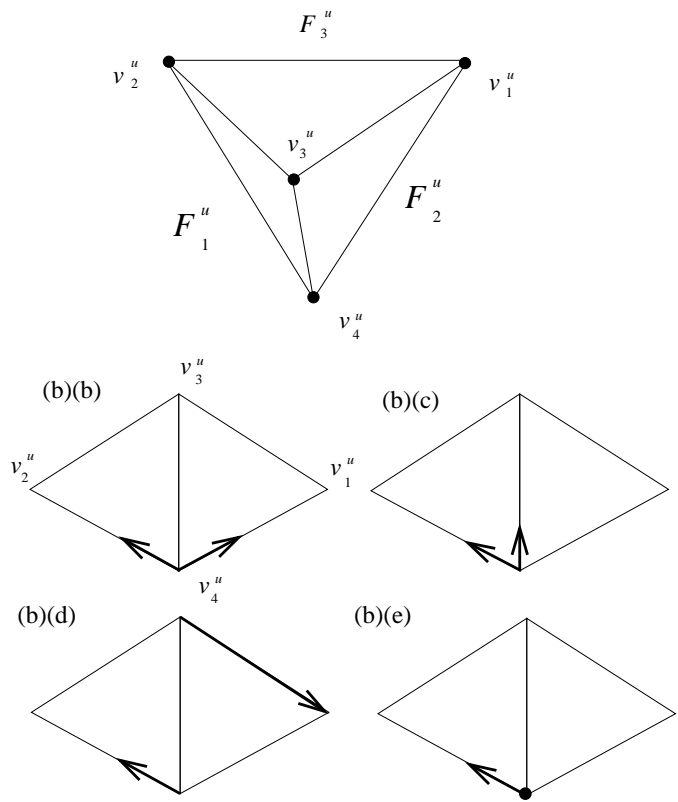}}  
\caption{\label{fig:cases} The arrows indicate the orientations.}  
\end{figure} 
\typeout{<<rf13.eps>>}  
 
\begin{figure}[t] 
\centerline{\epsfxsize=4.5in \epsfbox{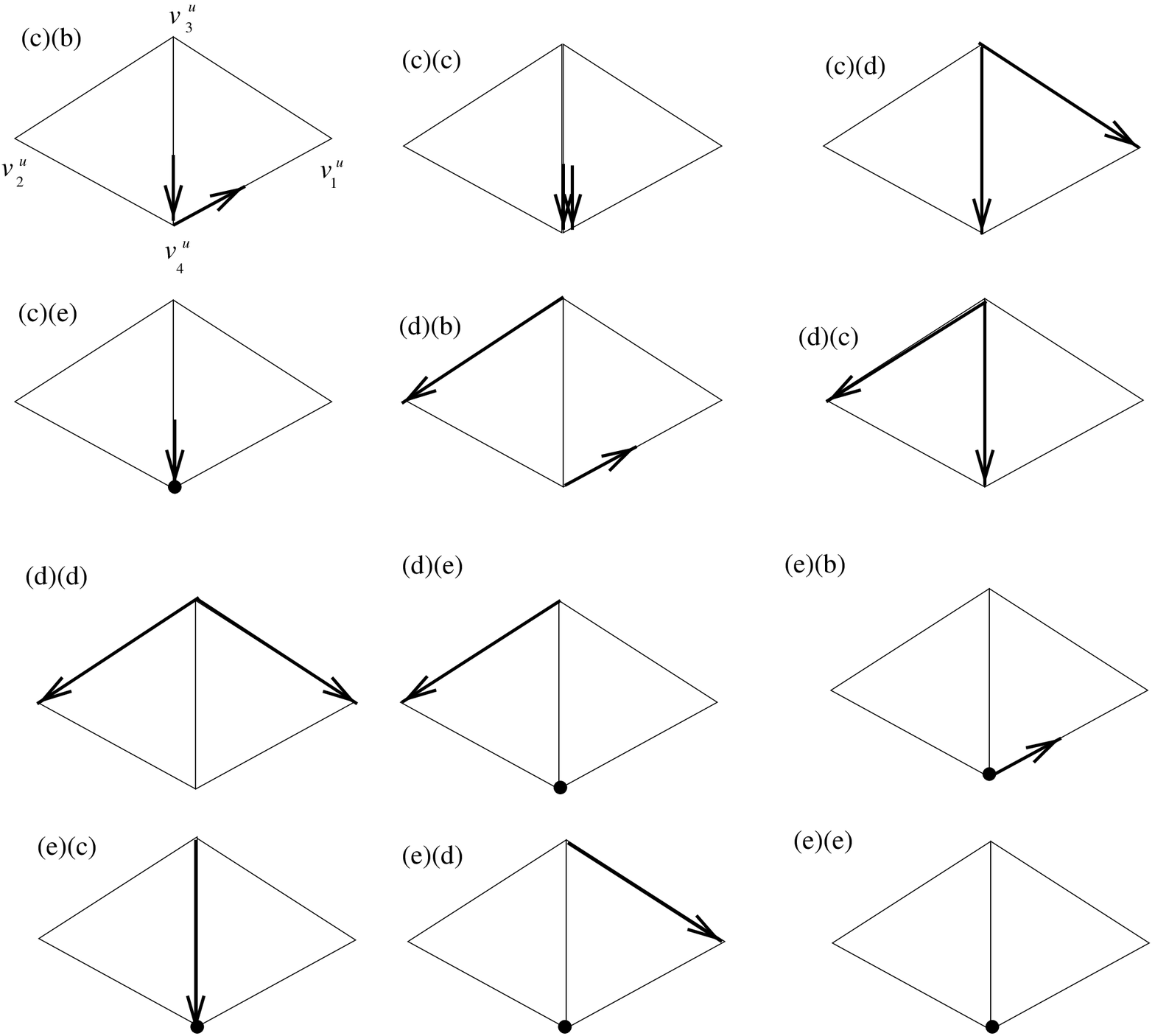}}  
\caption{\label{fig:cases2}}  
\end{figure} 
\typeout{<<rf14.eps>>}  
 
\begin{rem}\label{rem:ovlterm} 
Whenever points $x$ and $y$  both belong to $F^i$  (resp. $F^u$) for a 
fixed   $i$ and $\dev(x)$ and  $\dev(y)$   not antipodal, we denote by 
$\ovl{xy}$ the  unique segment in  $F^i$ (resp. $F^u$)  connecting $x$ 
and $y$.  
\end{rem} 
 
To begin, $F^\infty$ is a triangle cone with vertices $v^u_i$ for $i = 
1, 2, 3, 4$, each $F_l^u$ is a nondegenerate  triangle which forms sides 
of  $F^\infty$, and each $I_{jk}^u$ is  a segment forming the sides of 
$F_l^u$ for some $l$.  Since $L_1^u$ is a segment  or a point with one 
endpoint in $I_{12}^u$ and the another one  in $I_{13}^u$ (see Remark 
\ref{rem:endL}),  $L_1^u$  satisfies   one  of the following  mutually 
exclusive statements:  
\begin{itemize} 
\item[(a)] $L_1^u$ passes a point of $F_1^{u, o}$.  
\item[(b)] $L_1^u$ equals  a segment  in  $I^u_{13}$ with  an endpoint 
$O$.  
\item[(c)]  $L_1^u$  equals a  segment in $I^u_{12}$  with an endpoint 
$O$.  
\item[(d)] $L_1^u$ equals $I^u_{14}$.  
\item[(e)] $L_1^u$ equals the point $O$.  
\end{itemize} 
Similarly, $L_2^u$ satisfies  one of the following mutually  exclusive 
statement:  
\begin{itemize} 
\item[(a)] $L_2^u$ passes a point of $F_2^{u, o}$.  
\item[(b)]  $L_2^u$ equals a   segment in $I^u_{23}$  with an endpoint 
$O$.  
\item[(c)] $L_2^u$  equals a  segment  in $I^u_{12}$ with  an endpoint 
$O$.  
\item[(d)] $L_2^u$ equals $I^u_{24}$.  
\item[(e)] $L_2^u$ equals the point $O$.  
\end{itemize} 
 
We list  out all  possibilities   of $L_1^u$  and $L_2^u$ in   Figures 
\ref{fig:cases}  and \ref{fig:cases2}. By Proposition \ref{prop:notsamel}, 
$L_l^u$ can be as in (b)-(e) and $L_1^u$ and $L_2^u$ share an endpoint 
$x_{12}^u$.  By  these two conditions  and Remark \ref{rem:UDdisk}, we 
can have only one of the following possibilities:  
\begin{itemize} 
\item[(i)] $F_1^u  = U^u$, $F_2^u =  D^u$, $L_1^u = \ovl{v_3^uO}$, and 
$L_2^u = \ovl{v_1^uv_3^u}$ ((c)(d) in Figure \ref{fig:cases2}), or  
\item[(ii)] $F_1^u = D^u$, $F_2^u  = U^u$, $L_1^u = \ovl{v_2^uv_3^u}$, 
and $L_2^u =\ovl{v_3^uO}$ ((d)(c) in Figure \ref{fig:cases2}).  
\end{itemize} 
 
\begin{defn}\label{defn:compare} Given two sequences of  
positive numbers $a_i$ and $b_i$ with 
\[\lim_{i \ra \infty} a_i/b_i \in [0, \infty],\] 
we say  $a_i \gg b_i$  if $a_i/b_i  \ra \infty$. We  say 
$a_i \sim b_i$ if $|\log(a_i/b_i)|$ is bounded.  
\end{defn} 
Given  any pair of  sequences of positive  numbers, up  to a choice of 
subsequences, one of the three must hold.  
 
\begin{lem}\label{lem:eigen} 
Let  $s$ be a segment in  $\SI^n$ of $\bdd$-length  $< \pi$, and $p$ a 
point of  $s^o$. Let  $\{\vth_i\}$  be a  sequence of  real projective 
transformations acting  on $s$ fixing each  endpoints  $x$ and  $y$ of 
$s$. Then $\vth_i(p)$  converges to $x$  if and only if $\lambda_i \gg 
\nu_i$ for the eigenvalue $\lambda_i$ and $\nu_i$ corresponding to $x$ 
and $y$ of $\vth_i$. Moreover,  if $\vth_i(p)$ is convergent, then the 
limit is in $s^o$ if and only if $\lambda_i \sim \nu_i$.  
\end{lem} 
\begin{proof}  
The proof is obvious. 
\end{proof}  
 
There exist projective automorphisms $\help$ and $\help^i$ such that 
$\help(\dev(v_j))$ and $\help^i(\dev(v_j^i))$ for $j = 1, 2, 3, 4$ are 
in standard positions, i.e., at  
\[ [1, 0, 0, 0], [0, 1, 0, 0], [0, 0, 1, 0],  
\hbox{ and } [0, 0, 0, 1]\] the vertices  of the standard tetrahedron 
$T_s$ in  $\SI^3$, respectively,  such  that $\help^i$ form  a  bounded 
sequence in  $\Aut(\SI^3)$, which can be obtained  since $\dev(v_j^i)$ 
converges   to   the   vertices     of  $F^\infty$  a    nondegenerate 
tetrahedron. By choosing subsequences, we may assume that $\help^i \ra 
\help^\infty$ for some  real projective automorphism $\help^\infty$ in 
$\Aut(\SI^3)$.    Since  $\dev(F^i)  \ra   F^\infty$, it follows  that 
$\help^\infty(\dev(F^u))$ is in a standard position, i.e.,  
$\help^\infty(\dev(F^u)) = T_s$. 
 
We also note that  
\begin{equation}\label{eqn:deckstd}  
\help^i_\& =\help^i \circ h(\deck^i) \circ \help^{-1}  
\end{equation}   
fixes  the vertices  of   the  standard tetrahedron   $T_s$.    Such a 
projective automorphism has  a diagonal matrix expression with entries 
nonnegative; let $\lambda_1^i, \dots, \lambda_4^i$ be the eigenvalues 
of $\help^i_\&$ corresponding to  
\[[1, 0, 0, 0], [0, 1, 0, 0], [0, 0, 1, 0], \hbox{ and } [0, 0, 0, 1]\] 
respectively. 
For convenience, we assume  that $\help(x_{12})$, $\help(x_{23})$, and 
$\help(x_{13})$  are   at  $[0,0,1,1]$,  $[1,0,0,1]$,  and $[0,1,0,1]$ 
respectively.  
In  table  \ref{tab:guide}, we summarized    the relation between  all 
objects, which we will continue to use  
for situations in Chapters \ref{ch:tric}, and \ref{ch:rtrh}.  
 
\begin{table} 
\caption{\label{tab:guide} The lower table shows how to obtain various 
objects from another objects in the upper table.}  
\centerline{  
\begin{tabular}{|c|c|c|c|c|} \hline 
(1) & (2) & (3) & (4) & (5) \\  \hline\hline $\dev(F)$ & $\dev(F^i)$ & 
$\help(\dev(F))$  & $\help^i(\dev(F^i))$ &$\help^\infty(\dev(F^u))$ \\ 
\hline    $\dev(F_j)$   &    $\dev(F_j^i)$  & $\help(\dev(F_j))$     & 
$\help^i(\dev(F_j^i))$   &$\help^\infty(\dev(F_j^u))$    \\     \hline 
$\dev(L_l)$        &     $\dev(L_l^i)$   &     $\help(\dev(L_l))$    & 
$\help^i(\dev(L_l^i))$    &$\help^\infty(\dev(L_l^u))$    \\    \hline 
$\dev(x_{jk})$    &   $\dev(x_{jk}^i)$  &   $\help(\dev(x_{jk}))$    & 
$\help^i(\dev(x_{jk}^i))$  &$\help^\infty(\dev(x_{jk}^u))$ \\   \hline 
$\dev(\clo(G_1))$      &    $\dev(G_1^i)$    &      
$\help(\dev(\clo(G_1)))$      & 
$\help^i(\dev(G_1^i))$    &$\help^\infty(\dev(G_1^u))$     \\   \hline 
$\dev(K_1)$     &    $\dev(K_1^i)$        &    $\help(\dev(K_1))$    & 
$\help^i(\dev(K_1^i))$ &$\help^\infty(\dev(K_1^u))$ \\ \hline  
\end{tabular}} 
\medskip 
\centerline{  
\begin{tabular}{|c|c|c|c|c|} \hline  
$h(\deck^i)$ &  $\help$ & $\help^i_\&$  & $\help^i$ & convergence 
\\ \hline\hline (1) $\lora$ (2) & (1) $\lora$ (3) &  (3) $\lora$ (4) & 
(2) $\lora$ (4) & (4) $\lora$ (5) \\ \hline  
\end{tabular}}  
\end{table} 
 
\begin{figure}[t] 
\centerline{\epsfxsize=4.1in \epsfbox{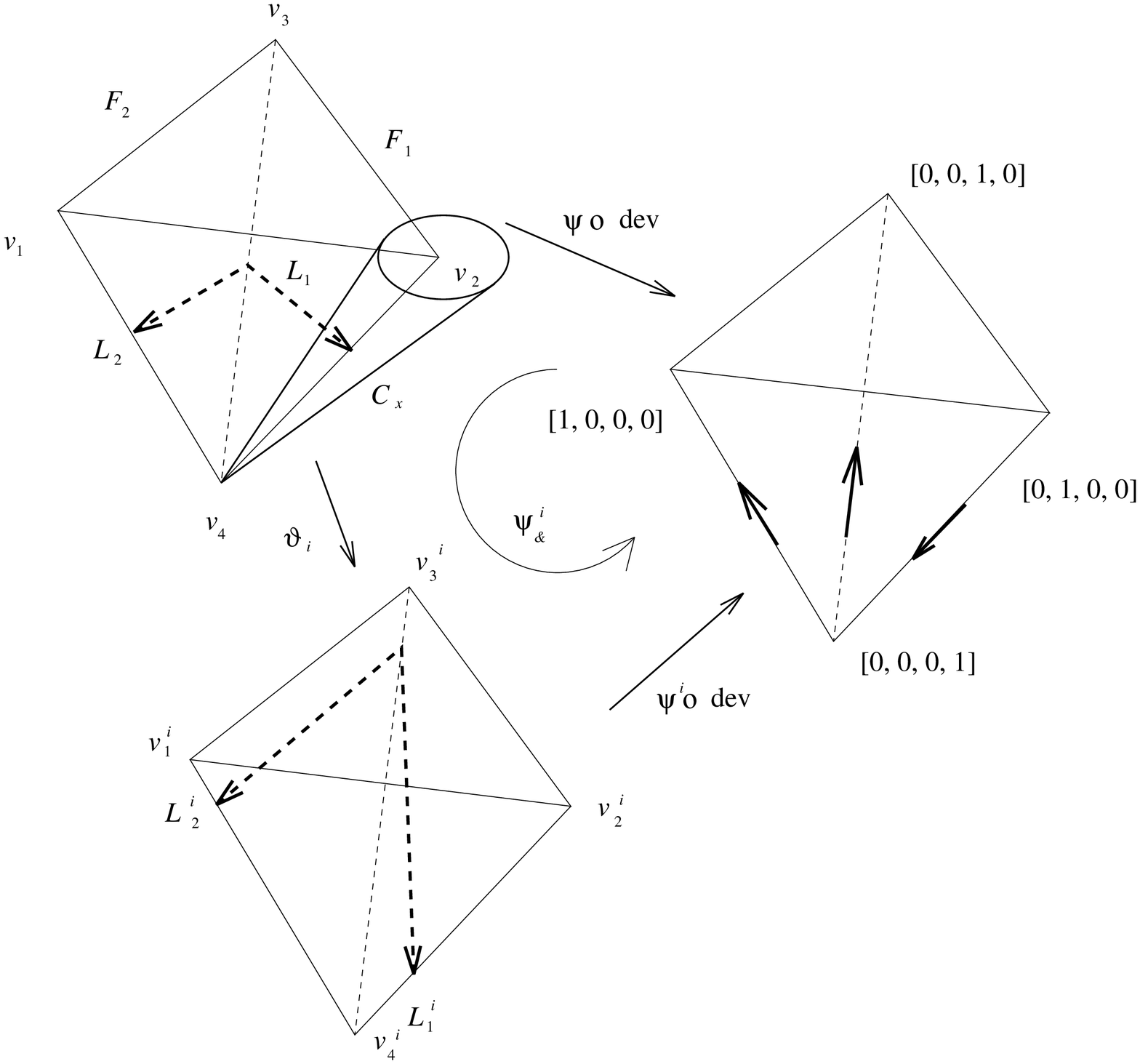}}  
\caption{\label{fig:caseI(i)} Case (i).}  
\end{figure} 
\typeout{<<rf23.eps>>}  
 
We will show  (i) cannot happen. Suppose  (i)  happened.  Then $L_1^u$ 
equals  $I_{12}^u$ and $L_2^u$  equals $I_{24}^u$.  Hence, $x^u_{13} = 
O$, $x^u_{12}   =   v^u_3$,  and  $x^u_{23}  =   v^u_1$.  (See  Figure 
\ref{fig:caseI(i)}.)  
Lemma \ref{lem:eigen} and Table \ref{tab:guide} imply that  
\begin{equation}\label{eqn:eigeninq} 
\lambda_4^i \gg \lambda_2^i,  \lambda_1^i \gg \lambda_4^i, \lambda_3^i 
\gg \lambda_4^i  
\end{equation}  
since \[x_{13}^\infty  = O, x^\infty_{23}  = v_1^\infty,  \hbox{ and } 
x_{12}^\infty   =    v_3^\infty\]   respectively.    Hence,  
$\lambda_2^i$ form the {\rm least eigenvalue sequence}\/ of  
$\help^i_\&$.  
 
Choose a small neighborhood $B(x)$ of a point  $x \in I_{13}^o$.  Then 
since  $B(x) \in M_h$, every ray  through $B(x)$ belongs to $M_h$; let 
$C_x$ be the cone formed from the union of these rays.  
Since $\help\circ  \dev(\clo(C_x))$ includes  an open 
cone including $\ovl{[0,0,0,1][0,1,0,0]}$, the condition on 
the eigenvalues implies that  
$\help^i_\& (\help \circ \dev(\clo(C_x)))$  
converges to a radiant  bihedron including  
$\help^\infty(F^\infty) = T_s$.  
As $\help^i_\& (\help \circ \dev(\clo(C_x)))$  
equals $\help^i \circ \dev(\deck^i(\clo(C_x)))$,  
and $\help^i$  is a convergent   sequence,       
$\dev(\deck^i(\clo(C_x))))$ converges   to    
a radiant bihedron $C^\infty$, which includes $F^\infty$. 
 
\begin{figure}[b] 
\centerline{\epsfxsize=4.1in \epsfbox{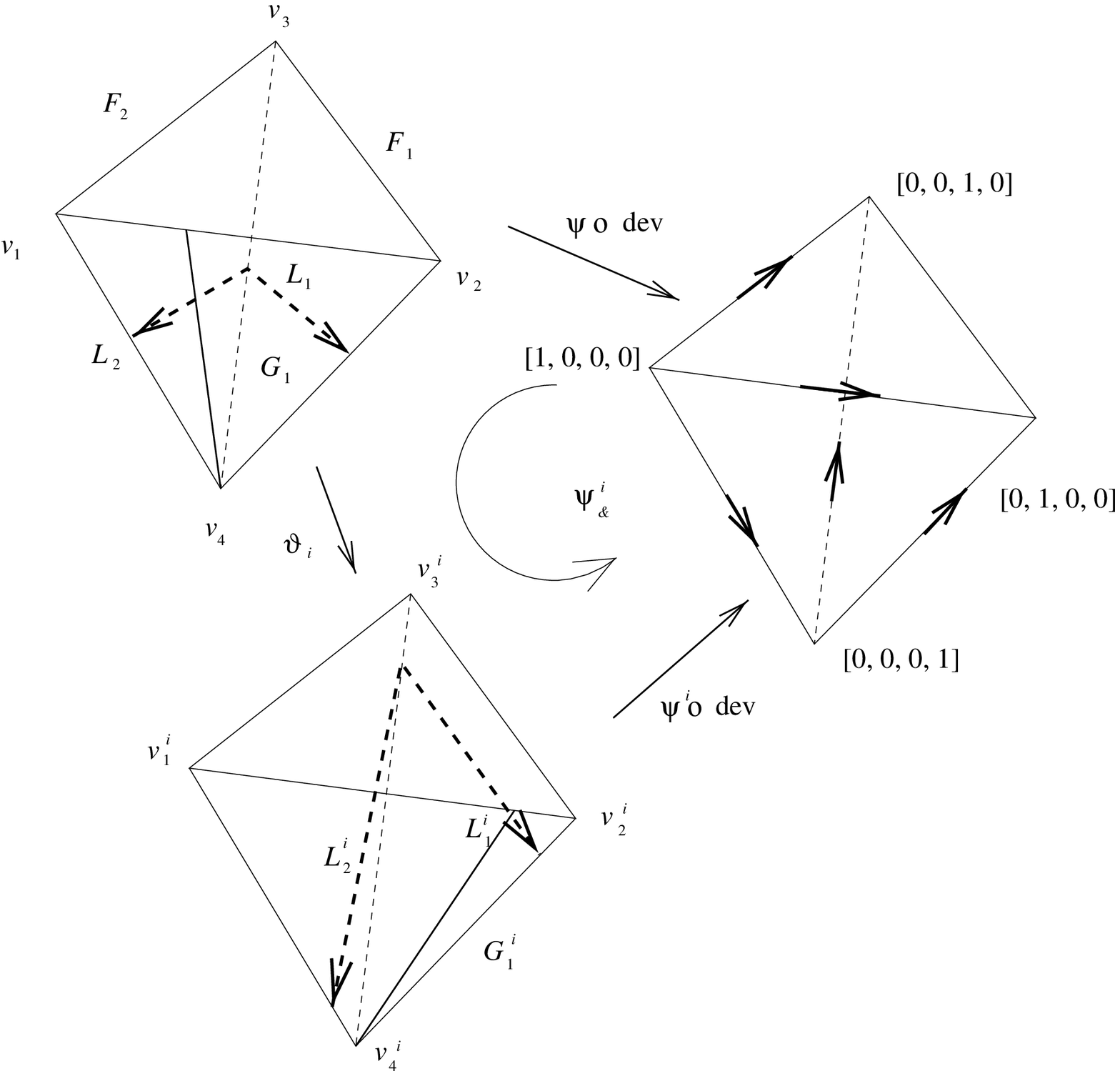}}  
\caption{\label{fig:caseI(ii)} Case (ii).}  
\end{figure} 
\typeout{<<rf24.eps>>}  

Since $\deck^i(\clo(C_x))$ overlaps with $\deck^i(F)$ for each $i$,  
by the dominating part of Theorem \ref{thm:seqconv}, there  
exists a radiant bihedron including $F^u$. This contradicts  
our hypothesis \ref{hyp2}, and the case (i) cannot happen.  
  
We  will  now show  that  (ii)  does   not  happen. Suppose that  (ii) 
happened.  Then Lemma  \ref{lem:eigen} and Table \ref{tab:guide} imply 
that   $\lambda_2^i \gg  \lambda_4^i$  since  $x_{13}^\infty = v_2^u$; 
$\lambda_4^i  \gg  \lambda_1^i$  since $x^\infty_{23}   = v_4^u$;  and 
$\lambda_3^i   \gg     \lambda_4^i$       since   $x_{12}^\infty    = 
v_3^\infty$. Hence, $\lambda_1^i$  is the least eigenvalue sequence of 
$\help^i_\&$, we obtain $v_5^i \ra  v_2^u$ as  
$\lambda_2^i \gg \lambda_1^i$, and, in  particular,   
$h(\deck_i)\dev(\clo(G_1))$ converges  to a segment,  
which is a  contradiction by   Proposition \ref{prop:roomG1}. (See Figure  
\ref{fig:caseI(ii)}.)

\chapter{The radiant trihedron case} 
\label{ch:rtrh} 
 
When $F^u$  is a radiant   trihedron, we divide our cases to (A), (B),  
and (C) depending on the  locations of the vertices $v_1^u$, 
$v_2^u$, and $v_3^u$ (see below). These  cases  will be ruled out  
in a  manner similar to the above section  for triangle cones; i.e.,  
we  use information  on  
$L_l^u$   and $U^u$  and   $D^u$ to restrict our configurations, and   
using an  asymptotic eigenvalue  estimation for holonomy actions, 
we can either obtain a radiant bihedron, which  is a contradiction as   
before, or obtain other contradictions straight away.  
 
Assume that  $F^u$ is a trihedron  with a face lune  in $\hIideal{M}$. 
$\dev(F^u)$  is a cone over  $\dev(F_4^u)$,  which must be a lune. 
$\dev(v_1^u)$,  $\dev(v_2^u)$,  and  $\dev(v_3^u)$ are   points in the 
boundary   of  $\dev(F_4^u)$.    Since   $\dev(F_4^i)$   converges  to 
$\dev(F_4^u)$ geometrically,  an elementary   geometric  consideration 
shows that two distinct points among $\dev(v_1^u)$, $\dev(v_2^u)$, and 
$\dev(v_3^u)$ form the vertices of the lune $\dev(F_4^u)$.  
 
We claim   that  $v_1^u  \ne  v_2^u$.   Otherwise,  $\dev(v_1^i)$  and 
$\dev(v_2^i)$ converge   to    a  common  point   respectively,    and 
$\dev(G_1^i)$  converges    to  a segment    since  $\dev(G_1^i)$ is a 
cone     over     $\ovl{\dev(v_1^i)\dev(v_5^i)}$           and 
$\ovl{\dev(v_1^i)\dev(v_5^i)}$     is        a     subsegment       of 
$\ovl{\dev(v_1^i)\dev(v_2^i)}$.  Since  $G_1^u$   is a   nondegenerate 
$2$-disk by Proposition \ref{prop:roomG1}, this is a contradiction.  
By Proposition \ref{prop:notsamel}, the  sides $F_1^u$ and $F_2^u$   
are not segments. Therefore, $v_3^u$ is distinct from $v_1^u$ and  
$v_2^u$.  We obtain:  
\begin{prop}\label{prop:distv} 
The vertices $v_1^u, v_2^u,$ and $v_3^u$ are mutually distinct.  
\end{prop}

\begin{figure}[b] 
\centerline{\epsfxsize=2.7in \epsfbox{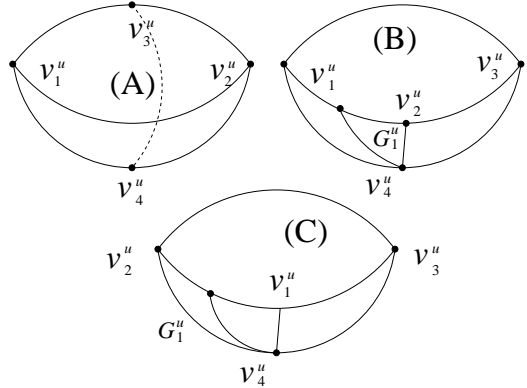}}  
\caption{\label{fig:caseABC} Cases A, B, C.}  
\end{figure} 
\typeout{<<rf17.eps>>}  
 
So, it follows that $v_1^u,  v_2^u, v_3^u$ must  be situated as in one 
of the following configurations (see Figure \ref{fig:caseABC}):  
\begin{itemize} 
\item[(A)] $v_1^u$  and $v_2^u$  form  antipodal vertices of the  lune 
$F^u_4$ with $v_3^u$ in  the interior of  a segment connecting $v_1^u$ 
and $v_2^u$.  
\item[(B)] $v_1^u$ and $v_3^u$ form antipodal vertices of $F^u_4$ with 
$v_2^u$ in the interior of the segment connecting $v_1^u$ and $v_3^u$.  
\item[(C)] $v_2^u$ and $v_3^u$ form antipodal vertices of the lune 
$F^u_4$ with $v_1^u$ in the interior of the segment connecting $v_2^u$ 
and $v_3^u$.  
\end{itemize}

We will show that  none of the above case  is possible, starting with 
the case  (A), using the strategy similar to what was in Chapter  
\ref{ch:tric}.  
 
From  our description  (A), we see  that  $v_5^u$  is  a point of  the 
segment $I^u_{34}$ of $\bdd$-length $=\pi$, $v_3^u$ in the interior of 
the segment   $I^u_{14}\cup I^u_{24}$   of $\bdd$-length  $=\pi$,  and 
$v_4^u$ at the origin $O$.  
 
As in the above section, Proposition \ref{prop:notsamel},  Remark 
\ref{rem:endL}, Remark \ref{rem:UDdisk}, and the endpoint matching 
condition show that we have only two possibilities:  
\begin{itemize} 
\item[(i)] $F_1^u  = U^u$, $F_2^u = D^u$,  $L_1^u = \ovl{v_3^uO}$, and 
$L_2^u = \ovl{v_1^uv_3^u}$, or  
\item[(ii)] $F_1^u = D^u$, $F_2^u  = U^u$, $L_1^u = \ovl{v_2^uv_3^u}$, 
and $L_2^u =\ovl{v_3^uO}$.  
\end{itemize} 
(See Figures \ref{fig:cases}  and \ref{fig:cases2} also, which are 
still true schematically but not geometrically.)  
 
From now on, $\ovl{xyz}$ for three points $x, y, z$ of $\SI^3$ so that  
$z$ is antipodal to $x$ denotes the unique convex segment with  
endpoints $x$ and $z$ passing through $y$. Such a segment has  
$\bdd$-length equal to $\pi$. 
Let $T_s$ be the {\em standard trihedron}\/ with vertices $[1, 0, 0, 0]$ 
and $[-1, 0, 0, 0]$ containing $[0,  1, 0, 0]$ and  $[0, 0, 1, 0]$ and 
$[0,  0, 0, 1]$ in its  edges respectively.  Now  we choose projective 
automorphisms $\help$ and $\help^i$ so that  
\begin{itemize} 
\item $\help(v_j)$ is in standard positions  
\[[1, 0, 0, 0], [0, 1, 0, 0], [0, 0, 1, 0], [0, 0, 0, 1]\] respectively 
for $j=1,2,3,4$,  
\item $\help^i(\dev(v_j^i))$  are in  standard positions for $j=1,3,4$ 
respectively, and  
\item $\help^i(\dev(I_{34}^i))$ is  on  the segment   
$\ovl{[-1,  0, 0,0][0, 1,  0, 0][1, 0,  0, 0]}$ with endpoints 
$[1,0,0,0]$  and $\help^i(\dev(v_2^i))$. 
\end{itemize} 
We choose $\help^i$ to form a bounded sequence in $\Aut(\SI^3)$. 
(We assume that $\help^i$  converges to $\help^\infty$ in 
$\Aut(\SI^3)$.)  Since  $\help^i$  is bounded,  $\help^i(\dev(v_2^i))$ 
converges to $[-1, 0, 0, 0]$.  Thus,  
\begin{equation}\label{eqn:deck} 
\help^i_\& =\help^i \circ h(\deck^i) \circ \help^{-1}  
\end{equation}  
acts on  $T_s$ and fixes the vertices  of the standard trihedron $T_s$ 
and $[0, 0, 1, 0]$ and $[0, 0,  0, 1]$; we  let $\lambda_k$ for $k=1, 
3, 4$  denote the eigenvalue  of $\help^i_\&$ associated with \[[1, 0, 
0, 0], [0, 0, 1, 0], \hbox{ and } [0, 0, 0, 1]\] respectively.  
$\phi^i_\&$ has the matrix form 
\[ \begin{bmatrix} \lambda_1^i & a^i & 0 & 0 \\  
0 & \lambda_2^i & 0 & 0 \\ 0 & 0 & \lambda_3^i & 0 \\  
0 & 0 & 0 &\lambda_4^i \end{bmatrix} \hbox{ for } a^i \in \bR.\] 
 
For  convenience, we let $\help(\dev(x_{12}))$, $\help(\dev(x_{13}))$, 
and  $\help(\dev(x_{23}))$ equal  $[0, 0,   1, 1]$,  $[0,  1, 0,  1]$, 
and $[1,0,0,1]$ respectively. We also note that  
$\help^\infty(\dev(F^u)) = \help^\infty(F^\infty) = T_s$. 
 
\begin{figure}[t] 
\centerline{\epsfxsize=4.5in \epsfbox{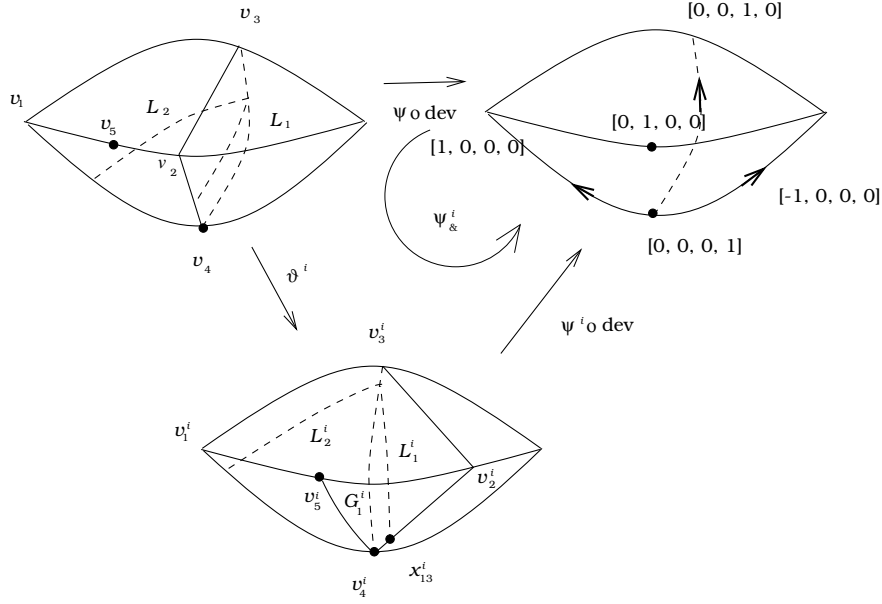}}  
\caption{\label{fig:case(A)(i)} The case (A)(i).}  
\end{figure} 
\typeout{<<rf11.eps>>}  
 
(i) Notice that $\help^i_\&$ is acting on  the lune $B$ with  
vertices $[1,0, 0, 0]$ and $[-1, 0,  0, 0]$ and containing   
$[0, 1, 0, 0]$ and $[0,0, 0, 1]$  in its edges where  
$\help^i_\&([0,  1, 0, 0])$ converges to $[-1, 0, 0, 0]$ and  
$\help^i(\dev(x_{13}^i))$  converges to $[0, 0, 0,1]$.  
 
From the  description of   the endpoints  of  $L_1$ and  $L_2$,  Table 
\ref{tab:guide} on $x_{23}$ and $x_{12}$ and Lemma \ref{lem:eigen}, we 
obtain that $\lambda_1^i, \lambda_3^i \gg \lambda_4^i$.  
 
In   the  affine  $3$-space $\mathcal{H}^o$ with origin $[0,0,0,1]$, 
$\help^i_\&$ is a linear map represented by a $3\times 3$-matrix.  
 
We  introduce the coordinates on  the affine plane containing $B^o$ so 
that $[0,0,0,1]$ has coordinates $(0, 0)$ and $\ovl{[1,0,0,0][0,0,0,1]} 
=  \help(I_{23})$  and   $\ovl{[0,1,0,0][0,0,0,1]}    = \help(I_{13})$ 
correspond to    the   $x$-axis  and   the    $y$-axis   respectively. 
$\help(\dev(x_{13}))$ and  $\help(\dev(x_{23}))$   have    coordinates 
$(1,0)$ and $(0,1)$ respectively.  
 
If we restrict $\help^i_\&$ to the vector subspace  
corresponding to the lune $B$, then $\help^i_\&$ has a matrix expression  
\[M^i_1 =  \begin{bmatrix} a_i & b_i \\ 0 & d_i \end{bmatrix}, 
a_i, d_i >  0. \] Since $\help^i(\dev(x_{23}^i))  \ra [1, 0, 0, 0]$, by 
equation     \ref{eqn:deck}   $a_i\ra  \infty$,    and     since 
$\help^i(\dev(x_{13}^i))  \ra [0, 0, 0, 1]$,  we obtain $b_i, d_i \ra 0$.   
Since $[0, 0, 1,   0]$ is a fixed  point  of $\help^i_\&$, it   follows  
that $\help^i_\&$ has a $3\times 3$-matrix expression:  
\[M^i_2 = \begin{bmatrix} a_i & b_i & 0 \\ 0 & d_i & 0 \\ 0 & 0 & e_i  
\end{bmatrix} \quad a_i, b_i, e_i > 0.\]  
The eigenvalues  of $M^i_1$ are  $a_i$ and $d_i$ and the corresponding 
eigenvectors are $(1,  0)$ and $(-b_i/(a_i -  d_i), 1)$.   Since $b_i, 
d_i \ra  0$  and  $a_i  \ra \infty$,   it  follows that 
the sequence of the  second eigenvectors converges to $(0,  1)$.  
Thus, for sufficiently large $i$, 
the matrix $M^i_2$ has a fixed point $s_i$ near $[0,  1, 0, 0]$ on the 
segment   $\ovl{[1, 0, 0, 0][0,1,0,0][-1,0,0,0]}$, 
with the associated eigenvalue $d_i$.  Since $a_i, e_i \ra \infty$ 
as $\lambda^i_1, \lambda^i_3   \gg \lambda^i_4$,  
and $d_i \ra 0$, it follows that $d_i$ is the least eigenvalue  
sequence of $M^i_2$. 
 
Choose  a point $x$ on the edge  $I_{13}$  of  $F$ and  a tiny-ball 
neighborhood $B(x)$ of  $x$ and form  a radiant cone  $C_x$  
containing $B(x)$.  
By the condition on  the eigenvalues and eigenvectors of  
$M^i_2$, $\help^i_\&(\help(\dev(\clo(C_x))))$    
converges to  a radiant bihedron including $T_s$. 
The trihedron $T_s$ equals $\help^\infty(\dev(F^u))$. 
Since $\help^i\circ \dev(\deck^i(\clo(C_x)))$ converges  
to a radiant bihedron including $T_s$,  
$\dev(\deck^i(\clo(C_x)))$ converges to  
a radiant bihedron including $F^\infty$. 
As in the case (i) of Chapter \ref{ch:tric},  
since $\deck^i(\clo(C_x)))$ always overlaps $\deck^i(F)$, 
the dominating part of Theorem \ref{thm:seqconv} shows  
that $\che   M_h$  includes a  radiant bihedron,  
contradicting Hypothesis \ref{hyp2}. 
 
\begin{figure}[b] 
\centerline{\epsfxsize=4.5in \epsfbox{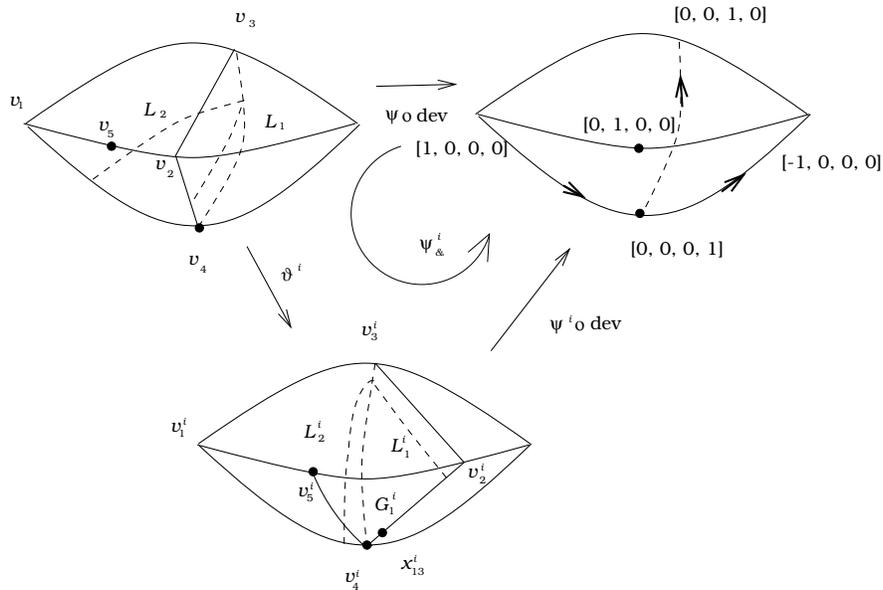}}  
\caption{\label{fig:case(A)(ii)} The case (A)(ii).}  
\end{figure} 
\typeout{<<rf12.eps>>}  
 
(ii) Here $L_1^u =  \ovl{v_2^uv_3^u}$  and $L_2^u = \ovl{v_3^uv_4^u}$.   
Since  $x_{23}^u = O$  and $x_{12}^u  = v_3^u$, we 
have that     $\lambda_3   \gg \lambda_4    \gg  \lambda_1$  by  Table 
\ref{tab:guide} and Lemma \ref{lem:eigen}.  
 
Our transformation  $\help^i_\&$  acts on the  lune  $B$ with vertices 
$[1, 0, 0, 0]$ and  $[-1, 0, 0, 0]$ and  containing $[0, 1, 0, 0]$ and 
$[0, 0, 0, 1]$ in its edges. We introduce coordinates as above and let 
$M^i_1$  denote the $2\times  2$-matrix  corresponding to $\help^i_\&$ 
acting on the vector subspace corresponding to $B$:  
\[M^i_1 = \begin{bmatrix} a_i & b_i \\ 0 & d_i \end{bmatrix} \quad 
a_i, d_i  > 0.\] Since $\help^i(\dev(x^i_{23}))$  converges to $[0, 0, 
0,   1]$, we obtain   $a_i     \ra    0$.      
Since the      direction      of 
$\ovl{O\help^i(\dev(x_{13}))}$       converges       to      that   of 
$\ovl{O[-1,0,0,0]}$, it follows that $b_i/d_i \ra -\infty$. Here $a_i, 
d_i > 0$ and $b_i < 0$ for $i$ sufficiently large.   
Moreover, $b_i \ra  -\infty$ by the condition 
that $\help^i(\dev(x_{13}^i))$ converges to $[-1, 0, 0, 0]$.  
 
For any pair $(c, d)$ of numbers, $d > 0$, $M^i_1$  maps $(c, d)$ to 
$(ca_i + d b_i, d d_i)$.  Since we have  
\[ \left| \frac{c}{d} \frac{a_i}{d_i}\right| /\left| \frac{b_i}{d_i} \right|  
=  \left|\frac{c}{d}\right|   \left|\frac{a_i}{b_i}\right| \ra 0,\]  
the ratio of the   coordinates  
$(ca_i +  d b_i)/(d d_i)  = (c/d)a_i/d_i  +  b_i/d_i$  
converges  to  $-\infty$.  Thus, for every point $x$ of $B^o$,  
the  sequence of rays $\help^i_\&(\ovl{Ox})$ 
converges    to   $\ovl{O[-1,  0,      0,   0]}$, meaning  that 
$\help^i_\&\circ\help(\dev(\clo(G_1)))$ converges to  
$\ovl{O[-1,   0,   0,  0]}$. This shows that 
$h(\deck^i)(\dev(\clo(G_1)))$ converges to a segment, 
contradicting Proposition \ref{prop:roomG1}.  
 
We now rule out the case (B). In this case, $F_2^u$ is a lune opposite 
to $v_2^u$ and $F_1^u$ is a  triangle, and so is  $F_3^u$.  As in case 
(A), Proposition \ref{prop:notsamel}, Remarks \ref{rem:UDdisk} and 
\ref{rem:endL}, and the endpoint matching condition show that  
\begin{itemize}  
\item[(i)] $F_1^u = D^u$, $F_2^u  = U^u$, $L_1^u$ is the segment 
$I^u_{14}$ and  $L_2^u$  is a subsegment of $I_{12}^u   \cup I_{23}^u$ with 
endpoints $v_3^u$ and $x_{23}^u \in I_{23}^u$.  
\item[(ii)] $F_1^u = U^u$,  $F_2^u = D^u$, $L_1^u$  equals $I_{12}^u$, 
and $L_2^u$ equals $I_{24}^u$.  
\end{itemize}

\begin{figure}[b] 
\centerline{\epsfxsize=3.1in \epsfbox{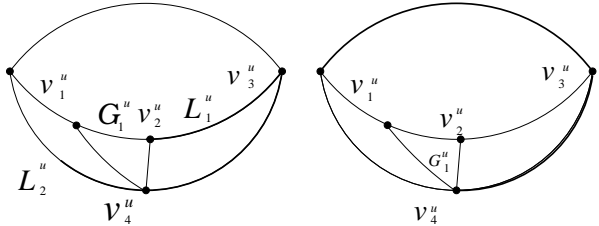}}  
\caption{\label{fig:case(B)(i)(ii)} Cases (B)(i) and (B)(ii).}  
\end{figure} 
\typeout{<<rf19.eps>>}

By the convergence condition, i.e., equation (\ref{eqn:objconv}),  
we may find a unique element $\help$ of  $\Aut(\SI^n)$ and  
a sequence of uniformly bounded transformations $\help^i$ of  
$\Aut(\SI^n)$ such that  
\begin{itemize}  
\item  $\help(\dev(v_1)), \help(\dev(v_2)),   \help(\dev(v_3)),$   and 
$\help(\dev(v_4))$ are at \[[1, 0,  0, 0], [0, 1, 0,  0],  
[0, 0, 1,0], \hbox{ and } [0, 0, 0, 1]\] respectively,  
\item $\help^i(\dev(v_1^i)) = [1, 0,  0, 0]$, $\help^i(\dev(v_2^i)) 
=[0, 1, 0, 0]$, and $\help^i(\dev(v^i_4)) = [0,0,0,1]$.  
\item $\help^i(\dev(I^i_{24}))$ is a segment in   
$\ovl{[1, 0, 0, 0][0,0, 1, 0][-1, 0, 0, 0]}$ with endpoints  
$[1, 0,0,  0]$ and  $\help^i\circ\dev(v_3^i)$,  the latter of   
which  forms a sequence converging to $[-1,0,0,0]$.  
\end{itemize}  
We assume as before  that   $\help^i$ converges to $\help^\infty$   in 
$\Aut(\SI^3)$.

As in  the previous chapters,   we define $\help^i_\&$ to be  $\help^i 
\circ h(\deck^i) \circ \help^{-1}$.  Then  $\help^i_\&$ acts on  $T_s$ 
and fixes each of $[1, 0, 0, 0], [0,  1, 0, 0],  [-1, 0, 0, 0], [0, 0, 
0, 1]$ but not $[0,0,1,0]$. Thus, $\help^i_\&$ has a matrix form:  
\begin{equation}\label{eqn:helpis}  
\begin{bmatrix}  
\lambda_1^i & 0  & a^i &  0 \\  0 & \lambda_2^i  &  0 & 0 \\  0  & 0 & 
\lambda_3^i & 0 \\ 0 & 0 & 0 & \lambda_4^i \end{bmatrix}  
\end{equation}   
Obviously, we may assume $\lambda_1^i,  \lambda_2^i,  \lambda_4^i  >   0$.    
As we  require $\lambda_1^i\lambda_2^i\lambda_3^i\lambda_4^i= 1$,  
it follows that $\lambda_3^i > 0$.  
 
For    convenience, we        assume    that    $\help(\dev(x_{12}))$, 
$\help(\dev(x_{23}))$,  and $\help(\dev(x_{13}))$ equal   $[0,0,1,1]$, 
$[1,0,0,1]$, and $[0,1,0,1]$ respectively.  
Again, $\help^\infty(\dev(F^u)) = T_s$.

\begin{figure}[t] 
\centerline{\epsfxsize=3.6in\epsfbox{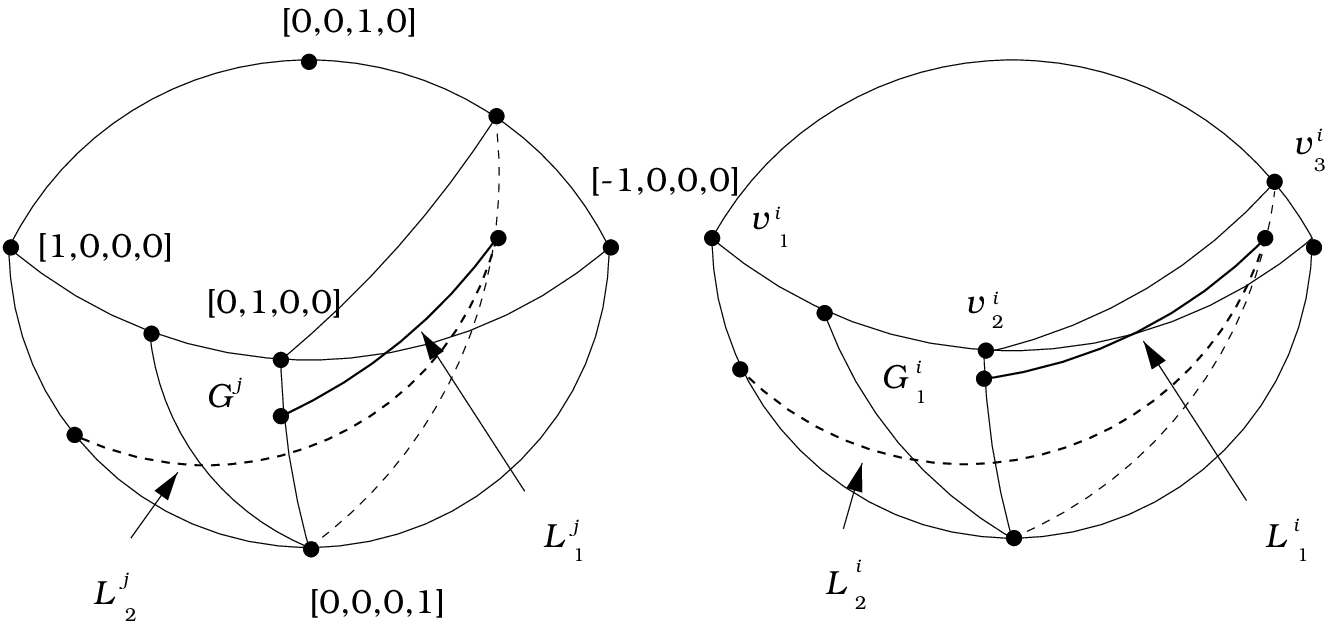}}  
\caption{\label{fig:caseBi} The case (B)(i).}  
\end{figure}

We begin with (B)(i).  
By condition (i), we see that $\help^i \circ \dev(x^i_{13})$ converges 
to $[0, 1, 0, 0]$,  $\help^i\circ \dev(x^i_{12})$ to  $[-1, 0, 0, 0]$, 
and $\help^i \circ \dev(L_1^i)$ converges  to the segment  
$\ovl{[0, 1,0, 0][-1, 0, 0, 0]}$.  $\help^i\circ \dev(x^i_{23})$  
converges to 
a point of $\ovl{[1, 0, 0, 0][-1, 0,  0,  0]}$,   and 
$\help^i\circ\dev(L_2^i)$ converges to a subsegment of 
$\ovl{[1,0,0,0][0,0,0,1][-1,0,0,0]}$    with  an 
endpoint $[-1,0,0,0]$    and   containing $[0,0,0,1]$. (See  Figure 
\ref{fig:caseBi}.)  
 
\typeout{<<rf18.eps>>}

From Table \ref{tab:guide}, we will get  estimations of eigenvalues of 
$\help^i_\&$ using  Lemma  \ref{lem:eigen}. By considering $x_{13}^u$, 
we  obtain $\lambda_2^i \gg \lambda_4^i$, and from  the following lemma  
we get $\lambda_1^i \gg \lambda_4^i$: Hence,   
$\help^i_\& \circ \help \circ \dev(x_{23})$ converges to $[1, 0, 0, 0]$ 
and so does $\help^i \circ \dev(x_{23}^i)$; 
$\help^i \circ  \dev(L_2^i)$ converges to 
$\ovl{[1, 0, 0, 0][0, 0, 0,1][-1, 0, 0, 0]}$.  
 
\begin{lem}\label{lem:eignsim} 
$\lambda_1^i \gg \lambda_2^i$ or $\lambda_1^i \sim \lambda_2^i$.  
\end{lem}  
\begin{proof}  
Since $\help^i\circ \dev(F_3^i)$ equals the triangle  
\[T= \triangle([1, 0, 0, 0][0, 1, 0, 0][0, 0, 0, 1]),\]  
$\help^i_\&$  acts on $T$. If    $\lambda_2^i  \gg    \lambda_1^i$,     
then    $\help^i_\&  \circ \help(\dev(\clo(G_1)))=\help^i(\dev(G_1^i))$    
converges to  the segment $\ovl{[0, 1,   0,    0][0, 0,  0,  1]}$. 
This contradicts Proposition \ref{prop:roomG1}.    
\end{proof}

Let us denote by $B$ the lune with  vertices $[1, 0,  0, 0]$ and  
$[-1,0, 0, 0]$ and two sides passing through $[0, 0, 0,  1]$ and  
$[0, 0, 1,0]$. Let $\SI^2_B$ denote the great sphere including $B$.  
Then $\help^i_\&$ restricts to a projective automorphism on  
$\SI^2_B$ with the matrix form  
\begin{equation} 
A_i = \left[ \begin{array}{ccc}  
\lambda_1^i & 0 & a^i  \\  
0 &\lambda_4^i & 0 \\ 0 & 0 &\lambda_3^i  
\end{array} \right] 
\end{equation}  
where $[1,0,0,0]$, $[0,0,0,1]$, and $[0,0,1,0]$  
form the order of the coordinate system. 
 
A {\em dilatation} is defined to be a projective 
automorphism of $\SI^3$ with a matrix of form  
\begin{equation} A_\lambda =  
\begin{bmatrix}  
1 & 0 & 0 & 0 \\ 
0 & 1 & 0 & 0 \\ 
0 & 0 & 1 & 0 \\ 
0 & 0 & 0 & \lambda \end{bmatrix}, \lambda > 0 
\end{equation} 
with respect to the standard basis $[1,0,0,0], [0,1,0,0], 
[0,0,1,0],$ and $[0,0,0,1]$  
 
\begin{lem}\label{lem:indconv} 
Let $\phi^i$ be a sequence of projective maps acting  
on the standard bihedron $T_s$. Let $S$ be a segment  
with vertices respectively at the middle of two radial  
edges $e$ and $f$ of a standard tetrahedron. Suppose that  
$\phi^i(S)$ converges to a compact set $S^*$. 
Assume that $S^*$ is dilatation-invariant.  
If $S'$ is another segment with vertices in $e^o$ 
and $f^o$. Then $\phi^i(S')$ also converges to $S^*$.  
\end{lem} 
\begin{proof} 
Suppose that $S'$ equals $A(S)$ for a dilatation $A$. Then  
as $A$ commutes with $\phi^i$, $\phi^i(A(S)) = A(\phi^i(S))$  
converges to $S^*$ for each $\lambda$. 
 
Suppose now that $S'$ is arbitrary. 
By choosing $\lambda \gg 1$ and $\delta \ll 1$,  
we can show that our $S'$ lies in the quadrilateral in a side of  
a standard tetrahedron bounded by $A_\lambda(S)$ and $A_\delta(S)$ 
and segments in $e$ and $f$ respectively.  
Since both sequences of images of $A_\lambda(S)$ and  
$A_\delta(S)$ under $\phi^i$ converges  
to $S^*$, it follows that $\phi^i(S')$ converges to $S^*$. 
\end{proof} 
 
Recall that $\help^i_\& (\help(\dev(L_2)))$ converges to  
\[\ovl{[1,0,0,0][0,0,0,1][-1,0,0,0]},\] a dilatation-invariant set. 
This is true no matter how we chose $x_{ij}$ in $I_{ij}$: 
Let $x'_{ij}$ denote another choice in $I_{ij}$.  
Let $L'_1$ and $L'_2$ be the resulting segments for our  
new choice. By above lemma \ref{lem:indconv}, 
the limit of the sequence of images $\help(\dev(L'_2))$ under $\help^i_\&$  
is same as that of $\help(\dev(L_2))$. (For  $L_1$ and $L'_1$,  
we can say similar things.) 
 
On $\SI^2_\infty$, $\help^i_\&$ restricts to a projective map  
with the matrix form  
\begin{equation} 
B_i =  
\left[ \begin{array}{ccc} \lambda_1^i & 0 & a^i  \\  
0 &\lambda_2^i & 0 \\ 0 & 0 &\lambda_3^i  
\end{array} \right] 
\end{equation}  
where $[1,0,0,0]$, $[0,1,0,0],$ and $[0,0,1,0]$ form the order of  
the basis. 
 
Since $\lambda_2^i \gg \lambda_4^i$, we see that   
$B_i = K_iA_i$ (strictly as $3\times 3$-matrices) for $K_i$  
equal to  
\begin{equation} 
\left[ \begin{array}{ccc} 1 & 0 & 0  \\ 0 & k_i& 0 
\\ 0 & 0 &1  
\end{array} \right] 
\end{equation}  
where $k_i$ converges to $\infty$.   
Let $L_2''$ be the segment in $F_4$ with vertices  
in the interiors of $I_{14}$ and $I_{34}$ respectively.  
Then by above equation,  
we easily see that the limit of the sequence 
$\help^i_\&(L_2'')$ is a subsegment of  
\[\ovl{[1,0,0,0][0,1,0,0][-1,0,0,0]}\] or the singleton $\{[0,1,0,0]\}$ 
up to choosing a subsequence. 
This follows since under the sequence of projective automorphisms  
of $\SIT$ with {\em matrices} $A_i$, this happens already for $L_2''$ 
by the conclusion of the above paragraphs and $K_i$ only  
enforces this for $L_2''$ to make the limit to be as stated 
since $K_i$ pushes everything near $[0,1,0,0]$. 
(We are applying facts obtained for $\SI^B$ to $\SIT$ in  
a somewhat convoluted way.) 
 
A {\em cone-neighborhood} of a subset $A$ of $\mathcal{H}$ is  
a radial-flow invariant neighborhood of $A$ in $\mathcal{H}$. 
 
The triangle $\help\circ\dev(F_2)$ has vertices  
$[1,0,0,0], [0,0,1,0]$, and $[0,0,0,1]$. 
Consider a cone-neighborhood of  
$\help\circ\dev(F_2 \cap M_h) - \{O\}$  
and a quadrilateral disk $J$ in $\SI^2_\infty$ with vertices  
the vertices of $L_2''$ and $[1,0,0,0]$ and $[0,0,1,0]$.  
Reflect this quadrilateral disk by  
an order two $\bdd$-isometric reflection  
$R_{2}$ fixing points of the triangle $\help\circ\dev(F_2)$. 
Then the union of the quadrilateral disk and its  
reflected image is another quadrilateral disk, say $J'$.   
If we choose the vertices of $L_2''$ sufficiently close to  
the edge $\ovl{[1,0,0,0][0,0,1,0]}$,   
the interior of the cone over $J'$ is a subset of  
the cone-neighborhood of $\help\circ \dev(F_2 \cap M_h) -\{O\}$. 
 
Since $F_2 \cap M_h$ is a subset of $M_h$, there  
exists a cone-neighborhood $C$ of $F_2 \cap M_h$ mapping into  
the above cone-neighborhood of $\phi\circ\dev(F_2) -\{O\}$ under  
$\phi\circ \dev$ if we choose the cone-neighborhood  
sufficiently close to $\help\circ\dev(F) - \{O\}$. 
Therefore, $\clo(C)$ includes a compact $3$-ball $J''$ mapping to  
the cone over $J'$ under $\phi\circ \dev$. 
 
By the condition on the limit of the sequence of 
images of $L_2''$ under $\help^i_\&$,  
we can show that $\help^i_\&(J)$ converges to  
the lune in $\SI^2_\infty$ with vertices $[1, 0, 0,0]$ and  
$[-1,0,0,0]$ with edges passing through $[0,1,0,0]$ and  
$[0,0,1,0]$ respectively.  
As $R_{2}$ commutes with $\help^i_\&$, we see that  
the limit of $\help^i_\&(J')$ is a $2$-hemisphere 
with center $[0,0,1,0]$.  
Therefore, $\help^i_\&(\help \circ \dev(J''))$ converges to  
a radiant bihedron including $T_s$. 
The trihedron $T_s$ equals $\help^\infty(\dev(F^u)) 
= \help^\infty(F^\infty)$ 
the limit of $\help^i_\& (\help \circ \dev(F))$.  
Thus, $\dev(\deck^i(J''))$ converges to a radiant  
bihedron including $F^\infty$ up to a choice of subsequences.  
As $\deck^i(J'')$ always overlaps with $\deck^i(F)$, similarly  
to case (i) in Chapter \ref{ch:tric} the dominating part of  
Theorem \ref{thm:seqconv} shows that  
there exists a radiant bihedron in $\che M_h$, 
a contradiction. 
 
We will now study (B)(ii).  By condition (B)(ii), we see that  
$\help^i\circ\dev(x^i_{23})$   converges  to $[1,0,0,0]$, $\help^i 
\circ  \dev(x^i_{13})$ to   $[0, 0, 0,  1]$,  $\help^i\circ\dev(x^i_{12})$  
to $[-1, 0, 0, 0]$, and $\help^i \circ \dev(L_1^i)$ to the segment     
$\ovl{[0,  0, 0,  1][-1,   0,  0, 0]}$.  $\help^i\circ\dev(L_2^i)$   
converges   to   the segment $\ovl{[1,0,0,0][0,0,1,0][-1,0,0,0]}$.  
 
\begin{figure}[b] 

\centerline{\epsfxsize=3.5in\epsfbox{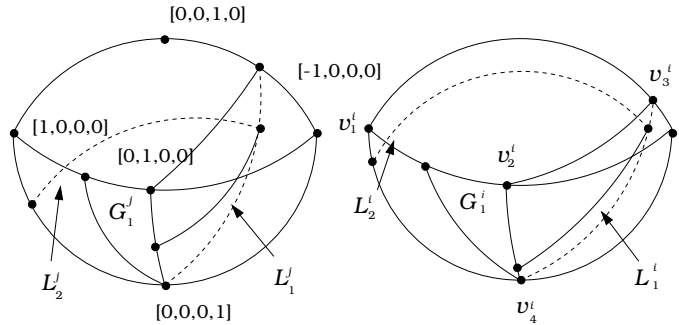}}  
\caption{\label{fig:caseBii} The case (B)(ii).}  
\end{figure} 
\typeout{<<rfc5.eps>>}  
 
From  the condition on   $x^i_{13}$,  we  have that $\lambda^i_4   \gg 
\lambda^i_2$. By the condition on $x^i_{23}$,  we get $\lambda^i_1 \gg 
\lambda^i_4$, and obtain  
\begin{equation}\label{eqn:Bii}  
\lambda^i_1 \gg \lambda^i_4 \gg \lambda^i_2.  
\end{equation}  
 
As above, we denote by $B$ the  lune with vertices $[1,  0, 0, 0]$ and 
$[-1, 0, 0, 0]$ and two sides passing through $[0, 0,  0, 1]$ and $[0, 
0, 1, 0]$.   
 
Let $\SI^2_B$ denote the great sphere including $B$.  
Then $\help^i_\&$ restricts to a projective automorphism  
on $\SI^2_B$ with the matrix form  
\begin{equation} 
A_i = \left[ \begin{array}{ccc}  
\lambda_1^i & 0 & a^i  \\  
0 &\lambda_4^i & 0\\  
0 & 0 &\lambda_3^i  
\end{array} \right] 
\end{equation}  
where $[1,0,0,0]$, $[0,0,0,1]$, and $[0,0,1,0]$  
form the order of the coordinate system. 
 
By assumption, $\help^i_\& (\help(\dev(L_2)))$ converges to  
\[\ovl{[1,0,0,0][0,0,1,0][-1,0,0,0]}.\] 
This is true no matter how we chose $x_{ij}$ in $I_{ij}$ 
as above by Lemma \ref{lem:indconv}. 
 
On $\SI^2_\infty$, $\help^i_\&$ restricts to a projective map  
with the matrix form  
\begin{equation} 
B_i =  
\left[ \begin{array}{ccc}  
\lambda_1^i & 0 & a^i  \\  
0 &\lambda_2^i & 0 
\\ 0 & 0 &\lambda_3^i  
\end{array} \right]. 
\end{equation}

Since $\lambda_4^i \gg \lambda_2^i$, we see that  
$B_i = K_iA_i$ for $K_i$ equal to  
\begin{equation} 
\left[ \begin{array}{ccc} 1 & 0 & 0  \\ 0 & k_i& 0 
\\ 0 & 0 &1  
\end{array} \right] 
\end{equation}  
where $k_i$ converges to $0$. Let $L_2''$ be the segment in $F_4$  
with vertices in the interiors of $I_{14}$ and $I_{34}$ respectively.  
A sequence of segments obtained by 
projective automorphisms of $\SIT$ with  
{\em matrices} $A_i$ applied to $(L_2'')$ converges to  
\[\ovl{[1,0,0,0][0,0,1,0][-1,0,0,0]},\] a fact following  
from seeing what happens in $\SI^2_B$ under $\help^i_\&$. 
Then by the above equation, we easily see that the limit of  
$\help^i_\&(L_2'')$ is the segment  
\[\ovl{[1,0,0,0][0,0,1,0][-1,0,0,0]}\] as $K_i$ can only 
help the segments to do so. 
 
Let $L''_1$ be the segment in $F_4$ with endpoints in endpoints in  
the interior of $I_{14}$ and $\ovl{[0,1,0,0][-1,0,0,0]}$  
and the endpoint in $I_{14}$ agrees with that of $L''_2$,  
which we denote by $x''_{14}$.  Then the limit of the sequence  
of images of  $x''_{14}$ under $\help^i_\&$ is clearly  
$[-1,0,0,0]$. Let $y$ be the other endpoint of $L''_1$.  
Since $\lambda^i_1 \gg \lambda^i_2$, we see that the limit of  
the sequence of images of $y$ under $\help^i_\&$ equals $[-1,0,0,0]$.  
Thus, that of $L''_1$ under  
$\help^i_\&$ converges to $\{[-1,0,0,0]\}$. 
 
Choose a cone-neighborhood  $C$ of  $I^o_{13}$ in $M_h$. We choose  
$L''_1$ and $L''_2$ sufficiently close to $[0,1,0,0]$ so that  
they lie in $\phi\circ\dev(\clo(C))$. 
Let $R$ be the $\bdd$-reflection so that the set of fixed points  
is precisely the great circle containing $[0,1,0,0]$ and $[0,0,0,1]$.  
Then $R$ commutes with $\help^i_\&$. $L''_1$, $L''_2$,  
$R(L''_1)$, and $R(L''_2)$ bound a quadrilateral disk $J'$ in $\SIT$.  
We of course assume that $R(L''_1)$ and  
$R(L''_2)$ also lies in $\phi\circ\dev(\clo(C))$. 
 
We see easily that the limit of the sequence of
the image of $J'$ under $\help^i_\&$ 
converges to the hemisphere with center $[0,1,0,0]$.  
There exists a cone in $J''$ in $\clo(C)$ which maps to  
the cone over $J'$ under $\help\circ \dev$.  
The sequence of images of $J''$ under  
$\help^i_\& \circ \help \circ \dev = \help^i\circ \dev \circ \deck^i$  
converges to a radiant bihedron including  
$\help^\infty(\dev(F^u)) = T_s$ by the above discussion. 
That of images of $J''$ under $\dev\circ\deck^i$ converges  
to a radiant bihedron including $F^\infty$. 
As above by  Theorem \ref{thm:seqconv}(4),  
there exists  a radiant  bihedron  in  $\che{M_h}$,  
a contradiction to Hypothesis \ref{hyp2}.  
We showed (B)(ii) does not happen.     
 
Now  we will show that  (C) does not  happen.  Since the arguments are 
similar to the  case (B), we  will only sketch  the argument:  
 
In this  case, $F_1^u$ is  a lune  opposite  $v_1^u$ and $F_2^u$  is  
a triangle, and so is  $F_3^u$. 
As in case (A) or (B), one of the following is true:  
\begin{itemize}  
\item[(i)]  We have $F_1^u = D^u$,  $F_2^u = U^u$,  and $L_1^u$ is the 
segment $I_{14}^u$ and $L_2^u$ is a segment $I_{12}^u$.  
\item[(ii)] We have $F_1^u  = U^u$, $F_2^u  =  D^u$, and $L_1^u$ is  a 
subsegment of $I_{12}^u  \cup I_{13}^u$ with  an  endpoint $v_3^u$ and 
$x_{13}^u$ and $L_2^u$ equals $I_{24}^u$.  
\end{itemize}  
 
\begin{figure}[t] 
\centerline{\epsfxsize=3.2in\epsfbox{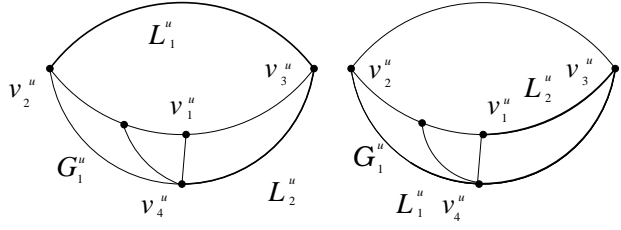}}  
\caption{\label{fig:case(C)(i)(ii)} Cases (C)(i) and (C)(ii).}  
\end{figure} 
\typeout{<<rfc7.eps>>}  
 
By the convergence condition equation \ref{eqn:objconv}, we may find a 
unique element $\help$ of  $\Aut(\SI^3)$  and a sequence of  uniformly 
bounded transformations $\help^i$ of $\Aut(\SI^3)$ such that  
\begin{itemize}  
\item  $\help(\dev(v_1)),     \help(\dev(v_2)), \help(\dev(v_3)),$ and 
$\help(\dev(v_4))$ are at  \[[0, 1, 0, 0], [1,  0, 0, 0],  
[0, 0, 1,0], \hbox{ and } [0, 0, 0, 1]\] respectively,  
\item $\help^i(\dev(v_1^i))  =  [0, 1, 0,   0]$, $\help^i(\dev(v_2^i)) 
=[1, 0, 0, 0]$, $\help^i(\dev(v_4^i))= [0, 0, 0, 1]$,  and
\item $\help^i(\dev(I^i_{14}))$ is a segment in   
$\ovl{[1, 0, 0, 0][0,0, 1, 0][-1, 0,  0,  0]}$, with  endpoints  
$[1, 0, 0,  0]$  and $\help^i\circ \dev(v_3^i)$, the latter of 
which forms a sequence  converging to $[-1, 0, 0, 0]$.  
\end{itemize}  
 
\begin{figure}[b] 

\centerline{\epsfxsize=3.4in\epsfbox{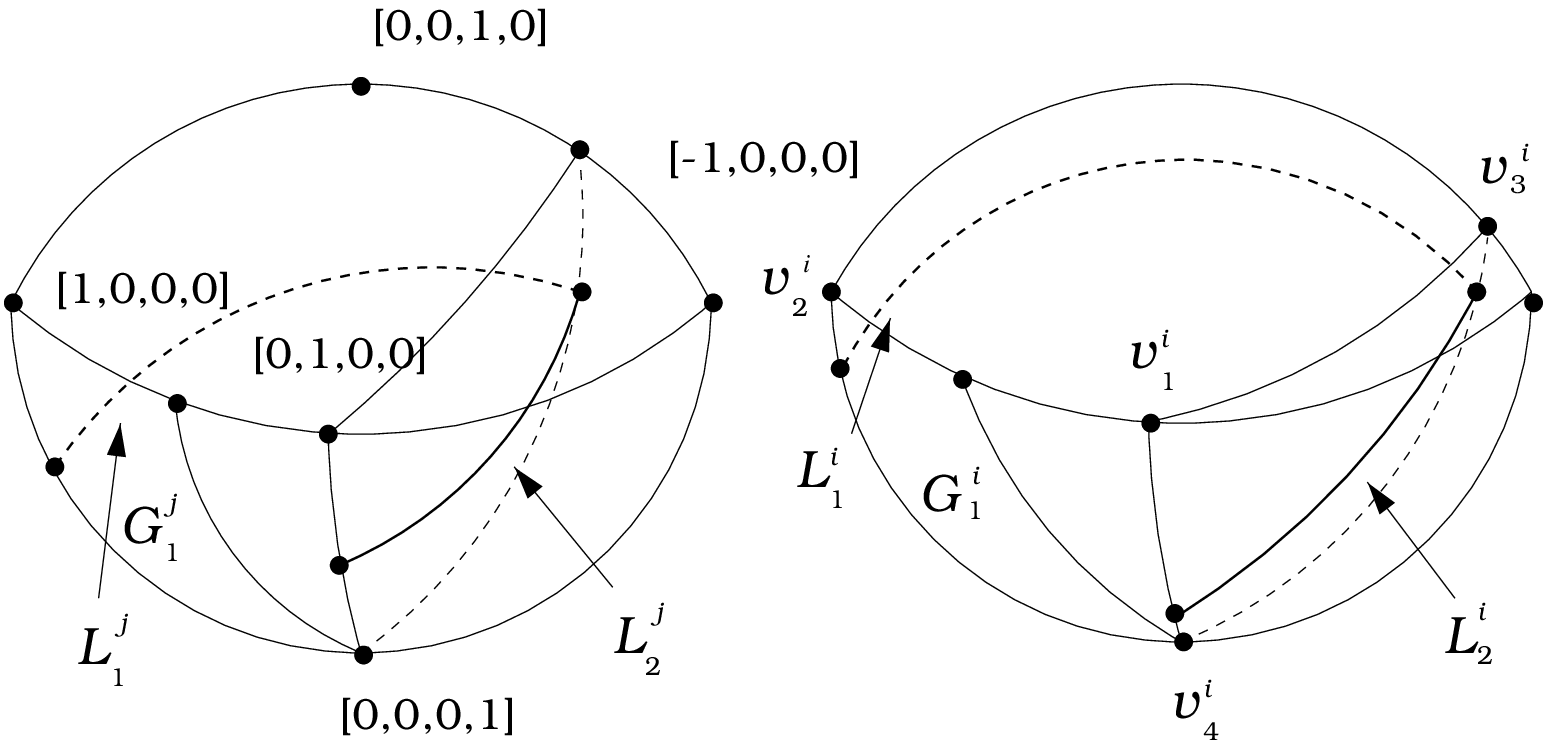}}  
\caption{\label{fig:case(C)(i)}  The  case  (C)(i). }  
\end{figure} 
\typeout{<<rfc6.eps>>}  
 
As before,   we define $\help^i_\&$ to be $\help^i  \circ h(\deck^i) 
\circ \help^{-1}$.  Then $\help^i_\&$ fixes  each of  $[1, 0, 0,  0],$ 
$[0,  1, 0, 0],$  $[-1, 0, 0, 0],$ and $[0,  0, 0, 1]$, and has a matrix 
form:  
\begin{equation}\label{eqn:helpis2}  
\left[   \begin{array}{cccc}   
\lambda_1^i &   0  &  a^i  &  0  \\   
0 & \lambda_2^i  & 0 &  0 \\   
0 & 0  &  \lambda_3^i &  0  \\  
0  & 0  & 0 & \lambda_4^i  
\end{array} \right] 
\end{equation}   
where $\lambda_j^i$ are positive as before.  
 
We begin with  (C)(i).  By condition (i), we  see that  $\help^i \circ 
\dev(x^i_{23})$   converges    to  $[0,   0, 0,    1]$,  $\help^i\circ 
\dev(x^i_{13})$  to  $[1,   0,  0,   0]$,  and  hence  $\help^i  \circ 
\dev(L_1^i)$  to the  segment
\[\ovl{[1, 0,   0, 0][0,  0, 0, 1][-1,  0, 0,  0]}.\]    
From this, we obtain the eigenvalue 
estimation  $\lambda_4^i  \gg  \lambda_2^i$    and $\lambda_1^i    \gg 
\lambda_4^i$; which implies $\lambda_1^i \gg \lambda_2^i$. Then 
$\help^i_\&(\help\circ\dev(G_1^u)) = \help^i \circ \dev(\deck^i(G_1^u))$  
converges to the segment $\ovl{[0,0,0,1][0,1,0,0]}$. As $\help^i$ 
is bounded, this is a contradiction to Proposition  
\ref{prop:roomG1}. 
 
We finish the series of the arguments with (C)(ii). By condition  (ii),   
$\help^i\circ   \dev(x^i_{23})$ converges to $[0,1,0,0]$, implying 
$\lambda_2^i \gg \lambda_4^i$.  
 
\begin{figure}[t] 
\centerline{\epsfxsize=3.4in\epsfbox{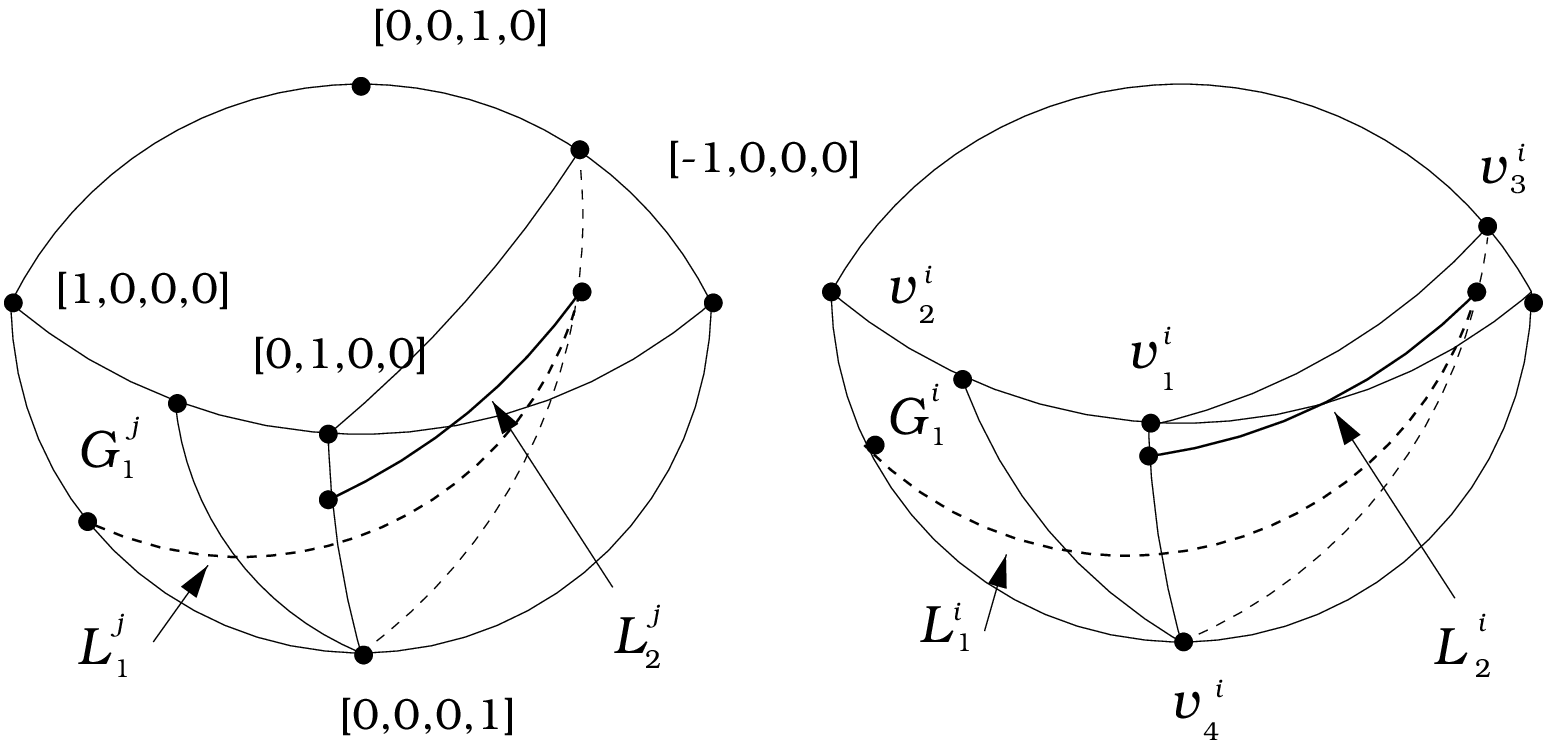}}  
\caption{\label{fig:case(C)(ii)} The case  (C)(ii).} 
\end{figure} 
\typeout{<<rf22.eps>>}  
 
As above, we denote by $B$ the  lune with vertices  $[1, 0, 0, 0]$ and 
$[-1, 0, 0, 0]$ and two sides passing through  $[0, 0, 0, 1]$ and $[0, 
0,  1, 0]$.   
 
Let $\SI^2_B$ denote the great sphere including $B$.  
Then $\help^i_\&$ restricts to a projective automorphism on  
$\SI^2_B$ with the matrix form  
\begin{equation} 
A_i = \left[ \begin{array}{ccc} \lambda_1^i & 0 & a^i   
\\ 0 &\lambda_4^i & 0 \\ 0 & 0 &\lambda_3^i  
\end{array} \right] 
\end{equation}  
where $[1,0,0,0]$, $[0,0,0,1]$, and $[0,0,1,0]$  
form the order of the coordinate system. 
 
By assumption, $\help^i_\& (\help(\dev(L_1)))$ converges to  
a subsegment of 
\[\ovl{[1,0,0,0][0,0,0,1][-1,0,0,0]}.\] 
This is true no matter how we chose $x_{ij}$ in $I_{ij}$ 
as before. 
 
On $\SI^2_\infty$, $\help^i_\&$ restricts to a projective map  
with the matrix form  
\begin{equation} 
B_i =  
\left[ \begin{array}{ccc} \lambda_1^i & 0 & a^i  \\ 0 &\lambda_2^i & 0 
\\ 0 & 0 &\lambda_3^i  
\end{array} \right]. 
\end{equation}

Since $\lambda_2^i \gg \lambda_4^i$, we see that   
$B_i = K_iA_i$ for $K_i$ equal to  
\begin{equation} 
\left[ \begin{array}{ccc} 1 & 0 & 0  \\ 0 & k_i& 0 
\\ 0 & 0 &1  
\end{array} \right] 
\end{equation}  
where $k_i$ converges to $\infty$.   
Let $L_1''$ be the segment in $F_4$ with vertices  
in the interiors of $I_{14}$ and $I_{34}$ respectively.  
Then as in case (B)(i), 
we easily see that the limit of the sequence $\help^i_\&(L_1'')$ 
is a subsegment of  
$\ovl{[1,0,0,0][0,1,0,0][-1,0,0,0]}$ or  the set of the point $[0,1,0,0]$. 
 
Exactly the same steps in (B)(i) will apply from here on and   
show that   there exists  a radiant    bihedron  in 
$\che{M_h}$, which is a  contradiction.    
Hence, the proof of  Proposition \ref{prop:MainEqu} is completed.

\chapter[Obtaining concave-cone affine manifolds]
{Obtaining concave-cone affine manifolds} 
\label{ch:obtain} 
 
We    will discuss the transversal intersection of two  crescent-cones,  
and define  the union $\Lambc(R)$   of a  collection   
of  overlapping  crescent-cones for  a 
crescent-cone $R$ in $\che M_h$.   The main part  of the  chapter   
is to derive  the properties of $\Lambc(R)$, particularly,   
the  equivariance and  local finiteness.   
We may find ``two-faced manifolds'' as in $n$-crescent   
cases  in   \cite{psconv}. After splitting off along these  
totally geodesic surfaces,   
we will be  able to show that the union of all subsets of form  
$\Lambc(R)$ for a crescent-cone $R$  covers a compact submanifold  
in of codimension $0$, i.e., the union of concave-cone 
affine submanifolds, 
using Proposition \ref{prop:lfin}. The complement  
of this compact manifold  is convex, and hence, we will achieve the   
decomposition. This will in turn prove Theorem A, and Corollaries
A and B. 
  
The main result of this chapter, implying Theorem A, is as follows: 
\begin{thm}\label{thm:dec2}  
Assume that $M$  is a $2$-convex but nonconvex radiant affine  
compact $3$-manifold with totally geodesic or empty boundary.   
The Kuiper completion of the holonomy cover of $M$  
does not include a radiant bihedron. 
Then  $M$ decomposes along disjoint totally  geodesic tori  
or  Klein bottles into  convex radiant affine  $3$-manifolds and  
concave-cone affine manifolds.  
\end{thm}  
 
We will be splitting and decomposing $M$ in the following  
steps. The completions of the holonomy cover of the resulting  
manifolds also do not include radiant bihedra and  
pseudo-crescent-cones. We need this to insure that the constructed 
manifold in each step satisfies the required assumptions  
for the next step.  
\begin{prop}\label{prop:subspl} 
Let $N$ be a radiant affine $3$-submanifold of $M$, closed in $M$  
or a radiant affine $3$-manifold obtained from $M$ by splitting  
along a properly imbedded totally geodesic surfaces. 
Then if the completion $\che N_h$ of the holonomy cover $N_h$ of  
$N$ includes a radiant bihedron, a crescent-cone, or a pseudo-crescent  
cone, then so does $\che M_h$. 
\end{prop} 
\begin{proof}  
This is similar to the proof of Proposition 8.9 of \cite{psconv}. 
Basically, we show that $N_h$ is a subset of $M_h$ when $N$  
is a submanifold of $M$, and if $R$ is a radiant bihedron in  
$\che N_h$, then the closure of $R^o$ in $M_h$ is a radiant  
bihedron in $\che M_h$. The same argument works if $R$ is  
a crescent-cone or a pseudo-crescent-cone. The details are omitted since  
they are more completely explained in \cite{psconv}. 
When $N$ is obtained from $M$ by splitting, a similar argument  
also applies. 
\end{proof} 
 
From now  on, we assume that   $\che M_h$ includes  crescent-cones  
but no pseudo-crescent-cones, which will be proved in
Chapter \ref{ch:pcrc}. We also assume that $\che M_h$ does not 
include any radiant bihedron because we decomposed using 
results of \cite{psconv} (see Theorem \ref{thm:dec2conv}). 
 
We will now construct an object in $\che{M_h}$ 
playing an analogous role to concave set in  \cite{cdcr1}; that is, we 
will  find an  equivariant and locally finite set from the suitable  
collections of crescent-cones.  
 
For a convex polyhedron or polygon $G$  in $\che{M_h}$, a side meeting 
$M_h$ is said  to be a {\em finite  side}\/ of  $F$.   A side  of $G$ in 
$\hideal{M}$  is said to  be   an {\em  ideal  side}.    If it is   in 
$\hIideal{M}$, we call it an {\em infinitely ideal  side}, and if not, 
a {\em finitely ideal side}.  
 
\begin{defn}\label{defn:crc}  
Let $T$ be a  crescent-cone in $\che{M_h}$, which is a trihedron 
with three lune sides. We denote by  $\nu_T$ the 
side of  $T$ meeting $M_h$, by  $\alfi_T$ the interior of the finitely 
ideal side of  $T$ and $\alin_T$ the  interior of the infinitely ideal 
side  of $T$ and $\alfiin_T$  the interior of  the segment that is the 
intersection of the finitely ideal  side and infinitely ideal side  of 
$T$. 
\end{defn}  
Incidentally, $\alfiin_T$ is an open line of $\bdd$-length $\pi$. 
 
Theorems    similar  to  the   transversal   intersection   theorem of 
\cite{cdcr1} hold for  crescent-cones (and  pseudo-crescent-cones).  Let 
$B$ be a convex $i$-ball in a convex $j$-ball $C$  for $i < j$.  Then 
$B$ is {\em   well-positioned}\/ if $B^o$  is  a subset of $C^o$   and 
$\delta B$ that  of $\delta C$.  (It is  fairly  easy to list  all 
types  of well-positioned convex  $i$-balls  in  a  lune, a triangle,  
a  tetrahedron, a trihedron or bihedron.)  
 
\begin{figure}[t] 
\centerline{\epsfxsize=3.3in\epsfbox{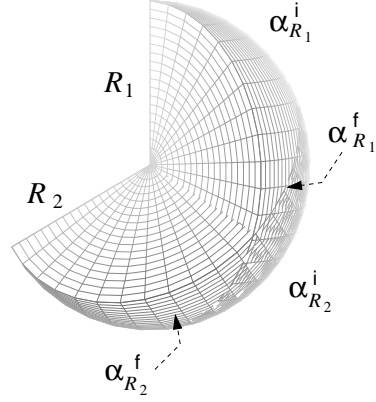}}  
\caption{\label{fig:trscc}     The     transversal   intersection   of 
crescent-cones.}  
\end{figure}

The definition of the transversal  intersection is quite difficult   
to write down but is fairly simple as Figure \ref{fig:trscc} shows.  
\begin{defn}\label{defn:trscc}  
Let $R_1$ and  $R_2$ be both crescent-cones.    We say that $R_1$  and 
$R_2$ intersect transversally if the  following statements hold ($i=1, 
j=2$, or $i=2, j=1$):  
\begin{enumerate}  
\item $R_1$ and $R_2$ overlap.  
\item $\nu_{R_1} \cap  \nu_{R_2}$ is a well-positioned radial segment 
$s$ in $\nu_{R_1}$ and similarly in $\nu_{R_2}$.  
\item $\nu_{R_i}  \cap  R_j$ equals  the   closure of  a  component of 
$\nu_{R_i}  - s$, and  it is a well-positioned  convex triangle $D$ in 
$R_j$. Two sides   of  $D$ other  that  $s$  are  finitely  ideal  and 
infinitely  ideal  respectively.  (The  three sides      of  $D$ are 
well-positioned in three sides of $R_j$ respectively.)  
\item $R_i \cap R_j$ is the closure of a component of $R_i - D$.  
\item Let $t$ be the  union of two sides of  $D$ other than $s$.  Then 
$t$ intersected   with the   infinitely   ideal side of $R_i$   is   a 
well-positioned segment, and $\alin_{R_i} \cap \alin_{R_j}$ equals one 
of  two components of $\alin_{R_i}   -  t$. Hence, $\alin_{R_i}   \cup 
\alin_{R_j}$  is an open  $2$-disk. The same statement holds for the 
finitely ideal sides $\alfi_{R_i}$ and $\alfi_{R_j}$.  
\item   $t \cap  \alfiin_{R_i}$   is  a  point.  $\alfiin_{R_i}   \cap 
\alfiin_{R_j}$  equals one of  two components  of $\alfiin_{R_i} - t$, 
and $\alfiin_{R_i} \cup \alfiin_{R_j}$ is an open line.  
\end{enumerate} 
\end{defn} 
 
\begin{rem}  
We now give another characterization of transversal intersection,  
which is easily shown to be equivalent.  
$R_1$ and $R_2$ overlap, $\dev(\alin_{R_1})$ and $\dev(\alin_{R_2})$  
are subsets of a common $2$-hemisphere in $\SIT$ bounded by  
a great circle including both $\dev(\alfiin_{R_1})$ and  
$\dev(\alfiin_{R_2})$, and $\dev(\nu_{R_1})$ and $\dev(\nu_{R_2})$  
meet transversally in a well-positioned segment. 
\end{rem} 
 
\begin{exmp} 
From Example \ref{exmp:crscone}, we can easily obtain  
two crescent-cones and in general they meet transversally. 
\end{exmp} 
 
\begin{thm}\label{thm:trscc}  
Let  $R_1$ and $R_2$ are overlapping  crescent-cones. Then  
either $R_1 \subset  R_2$ or $R_2 \subset R_1$  hold or  
$R_1$ and $R_2$ intersect transversally.  
\end{thm} 
\begin{proof}  
The proof is similar to that of Theorem \ref{thm:transversal}.  
Only technical differences exist although we have to use  
the following lemma \ref{lem:52} instead of Lemma 5.2 of  
\cite{psconv}. 
\end{proof} 
 
\begin{lem}\label{lem:52} 
Let $N$ be a closed real projective manifold with a developing map 
$\dev:\che N_h \ra \SI^3$ for the completion $\che N_h$ of 
the holonomy cover $N_h$ of $N$. Suppose that $\dev$ is  
an imbedding onto a union of two radiant bihedra $H_1$ and $H_2$  
in the $3$-hemisphere $\mathcal{H}$ meeting in a radiant bihedron  
or a radiant trihedron. Then $H_1 = H_2$, and $N$ is convex. 
\end{lem} 
\begin{proof} The proof is similar to that of Lemma 5.2 in \cite{psconv}. 
That is, we find a holonomy-group-invariant codimension-one  
submanifold in $H_1^o \cup H_2^o$ if $H_1$ is distinct from  
$H_2$. 
\end{proof}

Let $\mathcal{R}_1$ be a collection of all  crescent-cones in $M_h$.  We 
define a relation that $R \sim S$ for $R, S \in \mathcal{R}_1$ if 
they overlap, i.e., $R^o$ and $S^o$  meet. (Recall that since $R^o$  
and $S^o$ are open, this is equivalent to the statements that $R$  
and $S^o$ meet or  $R^o$ and  $S$ meet.)  Using  this, we define  
an equivalence relation  on $\mathcal{R}_1$, i.e., $R   \sim S$  
if  there exists a finite sequence $R^1, R^2, \dots, R^n$ in  
$\mathcal{R}_1$ such that $R^1 =  R$   and $R^n =  S$  where   
$R^i$  and $R^{i+1}$ overlap   for each $i=1,\dots, n-1$.  
 
We define the following sets as  in  \cite{cdcr1} and exhibit their 
equivariant properties:  
\begin{eqnarray}\label{eqn:Lambc}  
\Lambc(R)  = \bigcup\limits_{S \sim R}  S,  
&\dinLambc(R) = \bigcup\limits_{S \sim R} \alin_S,  
&\dfLambc(R) = \bigcup\limits_{S  \sim R} \alfi_S,\\  
\dfiLambc(R) = \bigcup\limits_{S \sim R} \alfiin_S,   
&\Lambc_1(R) = \bigcup\limits_{S \sim R} (S - \nu_S)  
&\hbox{ for } R \in \mathcal{R}_1 \nonumber  
\end{eqnarray} 
 
\begin{defn}\label{defn:ccaff} 
A {\em  concave-cone affine manifold} $N$ is  an  affine manifold  
such that $N_h$ is a subset of $\Lambc(R)$  for a crescent-cone  
$R$ in $\che N_h$.  
\end{defn}

We list the properties of $\Lambda(R)$.  
\begin{enumerate}  
\item   $\Lambc(R)$   and   $\Lambc_1(R)$  are 
path-connected.  
\item  $\Lambc(R) \cap M_h$  is a closed radiant  subset of $M_h$ with 
totally geodesic boundary, which is also a radiant set.  
\item $\Lambc_1(R)$ is a real projective $3$-manifold  
with boundary $\delta_\infty \Lambc(R)$ defined as 
$\dinLambc(R) \cup \dfLambc(R) \cup \dfiLambc(R)$ 
if given an obvious real projective structure induced from  
each crescent-cone. (Actually, the manifold-structure has corners 
\cite{KeSi}). 
\item $\dinLambc(R)$ is an open surface in $\hIideal{M}$.  
\item $\dfLambc(R)$ is an open surface in $\hfideal{M}$.  
\item  $\dfiLambc(R)$   is  an arc  in  $\hIideal{M}$.     The union 
$\delta_\infty \Lambc(R)$  is an open surface in 
$\hideal{M}$ which maps totally geodesic except at $\dfiLambc(R)$ 
under $\dev$. 
\item $\Lambc(R)$  is maximal, i.e., for  any triangle $T$  
in $\che M_h$ with sides $s_1, s_2$, and $s_3$,  
if $s_2, s_3  \subset \Lambc(R) \cap M_h$, then $T \subset \Lambc(R)$.  
\item $\Bd \Lambc(R) \cap M_h$  is a properly imbedded countable union 
of surfaces in $M_h^o$, each component  of  which is  a totally geodesic,    
radiant, properly imbedded open triangle or lune. 
\item For any deck transformation $\vth$, we  have $\Lambc(\vth(R))= 
\vth(\Lambc(R))$,      $\dinLambc(\vth(R))    =   \vth(\dinLambc(R))$, 
$\dfLambc(\vth(R))    = \vth(\dfLambc(R))$, and   $\dfiLambc(\vth(R))= 
\vth(\dfiLambc(R))$.  
\end{enumerate}  
The proofs of the above items can be obtained by adding radial  
dimension to the proofs for real projective surfaces  
(see \cite{cdcr1}). Since there is nothing surprising in this  
extension, we will not give proofs here. 
 
\begin{figure}[b] 
\centerline{\epsfxsize=3.3in \epsfbox{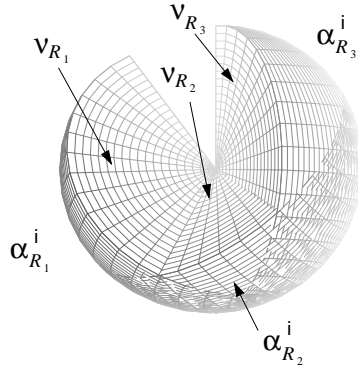}}  
\caption{\label{fig:Lam} A picture of $\Lambc(R)$ when there are only 
three crescent-cones. Imagine when there  are infinitely many.  $\dev$ 
restricted to it might be something like the universal covering  map  
of a solid torus, but the map is usually not a covering map.}  
\end{figure} 
 
\begin{exmp} 
In Example \ref{exmp:crscone}, we see that $\Lambc(R)$ for any  
crescent-cone $R$ equals the closure of $U-l$, i.e.,  
$\che \mathcal{E}_2$. $\dfLambc(R)$ equals the $xy$-plane 
with $O$ removed; $\dinLambc(R)$ equals the interior of  
the hemisphere that is the intersection of $\SIT$ with  
the closure of $U$ with the endpoints of $l$ removed. 
\end{exmp} 
 
\begin{prop}\label{prop:lfin2} 
The collection consisting of  elements of form $\Lambc(R)$ for  
some crescent-cone $R$ is locally finite, i.e.,  for any point $x$ 
of  $M_h$, there exists a neighborhood  which intersects only finitely 
many elements of form $\Lambc(R)$. {\rm (} Note if two sets of form  
$\Lambc(R)$ are equal, we don't count twice\/{\rm.)}  
\end{prop} 
\begin{proof}  
Suppose not. Let  $x$ be a  point of $M_h$, and  $B(x)$ a tiny ball of 
$x$ and $B'(x)$ a closed tiny-ball-neighborhood of $x$ in $B(x)^o$.   
Assume that our conclusion does not hold. Then there exists a sequence 
of crescent-cones $R_i$ such that $R_i$ meets $\inte B'(x)$ for  
each $i$ and $R_i$ is equivalent to a crescent-cone $S_i$ where 
$\Lambc(S_i)$   are   mutually   distinct.   
We see from the definition of $\Lambc(R)$ that $R_i$  
are mutually inequivalent.  
 
If $B'(x)$ is a subset of $R_i$ for some $i$, then $R_j$ and $R_i$  
overlap for each $j$, which is a contradiction. Hence,  
$\nu_{R_i} \cap B'(x)$ is not empty; let $p_i$ be a point  of  
$\nu_{R_i} \cap B'(x)$. We may assume without loss of generality that   
$p_i$ converges to a point  $p$ of $B'(x)$ and the $\bdd$-outer  normal  
vector $n_i$  to $\nu_{R_i}$  at $p_i$ converges  to a unit vector at  
$p$. Since $\alfi_{R_i}  \cup \alin_{R_i} \cup \alfiin_{R_i}$ does not 
meet  $B(x)$, Lemma \ref{lem:common}  shows that there  exists  
a common  open  ball in $B(x)$ in $R_i$ for $i$ sufficiently  large;  
that is, $R_i \sim  R_j$ for  $i, j$ sufficiently large, a contradiction  
to the above paragraph.  
\end{proof} 
 
As an aside, we have: 
\begin{prop}\label{prop:lfin3} 
Suppose  that  a   convex $3$-ball  $B$  in  $\che  M_h$ is  such that 
$\delta  B$ is a   union of a  convex  $2$-disk $\nu_B$  and a  disk 
$\alpha_B$ in $\delta B$ so that  $\nu_B$ meets $M_h$ and $\alpha_B$ 
is  a subset    of  $\hideal{M}$.  Suppose   that  $B$   satisfies  the 
equivariance property, i.e., for  each deck transformation either $B = 
\vth(B)$  or $B \cap  \vth(B) \cap M_h =  \emp$.   Then the collection 
consisting of  elements  of form   $\vth(B)$ for  deck transformations 
$\vth$ is locally finite.  
\end{prop} 
\begin{proof} 
If the collection is not locally finite,  then there exists a sequence 
$\vth_i(B)$ such  that $\vth_i(B) \cap \vth_j(B)  \cap M_h = \emp$ and 
$\vth_i(B)$  meets a tiny ball $B(x)$  of $x$ for   every $i$.  We get 
contradiction similarly to above by Lemma \ref{lem:common}. 
\end{proof} 
 
\begin{prop}\label{prop:splitting} 
Given   two  crescent-cones  $R$  and $S$,    we  have  one   of three 
possibilities\/{\rm :}  
\begin{itemize}  
\item $\Lambc(R) = \Lambc(S)$ and $R \sim S$,  
\item $\Lambc(R)$ meets $\Lambc(S)$  at the union of common components 
of $\Bd \Lambc(R) \cap M_h$ and $\Bd \Lambc(S) \cap M_h$ where $R$ and 
$S$ are not equivalent.  
\item $\Lambc(R) \cap M_h$ and $\Lambc(S) \cap M_h$ are disjoint.  
\end{itemize}  
\end{prop}  
\begin{proof}  
The proof is same as in the results in Chapter 7 of \cite{psconv}.   
The role played by $n$-crescent is played by crescent-cones.  
\end{proof} 
 
Suppose that $\Bd \Lambc(R)  \cap  M_h$ is  empty. Then  clearly  $M_h 
\subset  \Lambc(R)$ for  a single  crescent-cone $R$, and  $M$ is a 
concave-cone affine $3$-manifold and we are done.   We now assume that 
$\Bd\Lambda(R)\cap M_h$  is not empty  for  each crescent-cone $R$ in 
$\che M_h$.  
 
As in Chapter 7 of \cite{psconv}, a component   
$\Bd  \Lambc(R) \cap M_h$ which is also a component of  
$\Bd  \Lambc(S) \cap M_h$, for $R  \not\sim S$, is 
said to be a {\em copied}\/ component  of $\Bd \Lambc(R)  \cap M_h$.   
As in \cite{psconv}, we can show that a copied component  
is a common component of $\nu_T  \cap M_h$  for some crescent  
cone $T$, $T \sim R$ and $\nu_U \cap M_h$ for $U$, $U \sim S$.  
Two copied components are either disjoint or identical,  
and the collection of all copied components are 
locally  finite.  (This is proved similarly to the arguments in   
Propositions 6.4 and 7.6 in \cite{psconv}.) 
The union of  all copied components of  $\Bd \Lambc(R) \cap M_h$ where 
$R$ runs over all crescent-cones is denote by  $P_{M_h}$ and is said 
to  be  a  {\em  pre-two-faced surface}.   It  is a  properly imbedded 
$2$-manifold in $M_h$, and moreover by   Proposition   \ref{prop:lfin}, 
$p|P_{M_h}$,  for the covering map $p:M_h \ra  M$,  covers a properly 
imbedded compact  surface $A_{M}$. Since each component of  
$\Bd \Lambc(R)$ is a radiant  set, each component of $A_{M}$ is radially  
foliated. Hence $A_{M}$ is a finite union of disjoint totally geodesic  
tori or Klein bottle. $A_M$ is said to be the {\em two-faced submanifold}.

Note that the interior of every crescent-cone in $\che M_h$ is disjoint 
from $P_{M_h}$. Otherwise, a component of $P_{M_h}$ and $R^o$ for 
a crescent-cone $R$ meets. Since each component of $P_{M_h}$ is 
a components of $\nu_T \cap M_h$ for a crescent-cone $T$.  
We have $R \sim T$ and this means that $R^o$ is a subset of  
$\inte \Lambc(S)$ for $S \sim T$, a contradiction. 
 
To see some examples of these, the Benz\'ecri suspension 
of a real projective surface in Example 7.9 of \cite{psconv}  
will give us a two-faced surface and sets of form 
$\Lambda(R)$. Further examples may be obtained from  
Benz\'ecri suspensions of real projective surfaces.  
 
As in  \cite{psconv}, we  can  split $M$   along   $A_{M}$.   
Let  $N$ be the split manifold, and $M'$ the split manifold of  
$M_h$ by $P_{M_h}$.  
Then $M'$ contains a copy  of $M_h - A$ and $\dev| M_h$ extends to   
a projective  map   on $M'$, again denoted by $\dev$. By pulling-back  
the metric $\mu$ and getting a path-metric $\bdd$ on $M'$, we complete  
$M'$ to a space $\che M'$. 
 
\begin{prop}\label{prop:crccorr} 
The crescent-cones in $\che M'$ and $\che M_h$ are in  
one-to-one correspondence given by $R \leftrightarrow R'$  
if and only if $R^o = R^{\prime o}$ in $M_h - P_{M_h}$. 
\end{prop} 
\begin{proof} See Chapter 8 of \cite{psconv}. 
\end{proof} 
 
We can easily show by Proposition \ref{prop:crccorr} that $M'$   
contains  no copied components as in \cite{psconv}  
since we split $M$ to obtain $N$. Hence, given any two  
crescent-cones $R$ and  $S$ in $\che M'$, we have either  
$\Lambc(R) = \Lambc(S)$  or $\Lambc(R) \cap M'$ and  
$\Lambc(R)\cap M'$ are disjoint.  
 
The holonomy cover $N_h$ of $N$ is obtained by taking a disjoint  
union of single appropriate components of $M'$ for all components of $N$ 
(see Chapter 8 of \cite{psconv}). 
The developing map $\dev$ is obtained by restricting $\dev:M' \ra \SI^3$  
to $N_h$, and the holonomy homomorphism is given by  
$\vth \mapsto h(\vth)$ when $\vth$ is a deck transformation of  
$M'$ acting on a component of $N_h$. (Note here, we do not consider  
the action of deck transformations switching components.) 
$\dev|N_h$ is again a developing map, and we obtain  
the Kuiper completion $\che N_h$ of $N_h$ using this  
developing map. 
 
Obviously the split manifold $N$ is still $2$-convex.  
If not, $\che N_h$ includes a $3$-crescent  
by Theorem 4.6 of \cite{psconv}  
but $\che N_h$ cannot include radiant bihedra 
and hence $3$-crescents by Propositions  
\ref{prop:subspl} and \ref{prop:raecr}. 
 
Let  $A$ be $\bigcup_{R \in {\mathfrak  R}_1(N)}  \Lambc(R) \cap N_h$ where  
${\mathfrak R}_1(N)$ is  the set of representatives  of the equivalence  
classes in the set ${\mathfrak R}_1(N)$ of all  crescent-cones in $\che N_h$.  
By Proposition \ref{prop:lfin2},  the collection of such  subsets is 
locally finite, and $A\cap  N_h$ is a  submanifold of $N_h$ by above 
paragraph since each $\Lambc(R') \cap  N_h$ for a crescent-cone $R'$ 
in $\che N_h$ is a submanifold. Therefore  by Proposition \ref{prop:lfin},  
$P| A:   A  \ra p(A)$ is a covering map  onto a compact submanifold  
$p(A)$  in $N$ of codimension $0$ with nonempty boundary.  Since each  
component of $\Bd \Lambc(R) \cap N_h$ is totally geodesic  
(see the properties) and radiant, i.e., foliated by radial  lines,  
each component of $p(A)$ is a  compact submanifold with boundary  
whose components are totally geodesic tori or Klein  bottles.   
Thus, we  see  that $N$ decomposes into $p(A)$ and the closure of  
components of $N-p(A)$.   
 
The holonomy cover of a component $S$ of the closure of $N-p(A)$ 
can be considered a component $S_h$ of the closure of $N_h - A$ 
with the induced covering map $p: S_h \ra S$. 
 
\begin{prop}\label{prop:subcorr} 
Let $A$ be a closed submanifold in $\che N_h$. Then  
the collection of crescent-cones in the Kuiper completion 
of $\che A$ of $A$ are in one-to-one correspondence with 
those of $\che N_h$ in $\clo(A)$ by  
$R \leftrightarrow R'$ if and only if $R^o = R^{\prime o}$. 
\end{prop} 
\begin{proof} See Section 8 of \cite{psconv}. The proof  
can be obtained from the arguments there. 
\end{proof} 
 
Each component $S$ of the closure of $N-p(A)$ is obviously  
$2$-convex. Otherwise, $\che S_h$ must include a $3$-crescent  
by Proposition 4.6 of \cite{psconv} but this possibility 
was ruled out earlier by Proposition \ref{prop:subspl}. 
 
The components of the  closure  of $N-p(A)$ are convex since  
otherwise, the completion of the holonomy cover of at least one  
of them includes a crescent-cone by Theorem \ref{thm:crscone}. 
Hence $\che S_h$ includes a crescent-cone for some component  
$S$, and by Proposition \ref{prop:subcorr}, $\clo(S_h)$ in  
$\che N_h$ includes a crescent-cone in $\che N_h - A$, where  
$S_h$ is regarded as a subset of $N_h$. Since for each  
crescent-cone $R$ in $\che N_h$, $R \cap N_h$ is a subset of  
$A$ by the definition of $A$, $R \cap N_h$ is a subset of  
the intersection of $S_h$ and $A$. Since $R\cap N_h$ includes 
$R^o$, an open set and $S_h \cap A$ does not contain an open  
set, this is a contradiction. Thus each component $S$ is convex.

Let $K$  be a component  of $p(A)$. Then  $K$ is covered by  
$\Lambc(R) \cap M_h$ for some crescent-cone $R$. By a reason similar  
to that in Chapter 8 of \cite{psconv}, i.e., Lemma 8.1, we can show   
that $\Lambc(R) \cap M_h$ is a holonomy cover of $K$ with  
developing map $\dev| \Lambc(R) \cap M_h$. The surface  
$\Lambc(R) \cap  M_h$ is a subset of $\Lambc(R')$ in  
the completion $(\Lambc(R) \cap  M_h)\che{}$ of $\Lambc(R) \cap M_h$  
for the crescent-cone $R'$ in the completion with the 
identical interior as $R$ by Proposition \ref{prop:subcorr}.   
Thus, each component of $p(A)$ is a concave-cone affine manifold  
with  totally geodesic  boundary, and this proves  
Theorem \ref{thm:dec2}.  
 
Theorems \ref{thm:cafB}, \ref{thm:cafA}, and \ref{thm:dec2} complete  
the proof of Theorem A. (Of course, assuming the result of
Chapter \ref{ch:pcrc}.)
 
\begin{proof}[Proof of Corollary A] 
Suppose that $M$ is closed. If $M$ decomposes as  
in Theorem A, then since the decomposing surfaces 
are tangent to the radial flow, Theorems \ref{thm:barbot}  
shows that $M$ admits a total cross-section.  
 
If $M$ does not decompose, then Theorem  A    shows   
that $M$ is   either    a convex   radiant affine 
$3$-manifold or a generalized affine suspension as in (2) of  Theorem A.   
By Theorem \ref{thm:barbot}, $M$ is a  
generalized affine suspension.  
 
If $M$ has nonempty boundary, and $M$ decomposes  
nontrivially, then each boundary component of $M$  
is in one of the decomposed pieces. The decomposed  
manifolds are either convex, concave-affine or  
concave-cone affine. If a boundary component $A$ is  
in a convex piece, then $A$ is convex. If $A$ is in  
the other two types, then a component $A_h$ of the inverse  
image of $A$ in the holonomy cover of the piece $N_i$,  
where $A$ is, is a subset of a boundary of  
a crescent or a crescent-cone intersected with $M_h$.  
In these cases, $A_h$ is affinely homeomorphic to 
a convex cone in $\bR^2 -\{O\}$ and $\bR^2 -\{O\}$  
itself by Lemma \ref{lem:convtori}. Theorem \ref{thm:barbot2} 
shows that $M$ is a generalized affine suspension. 
 
Therefore, $M$ admits a total cross-section to the radial flow  
by Proposition \ref{prop:crssect}. 
The proof is completed by Theorem \ref{thm:Sorb}. 
\end{proof} 
 
\begin{proof}[Proof of Corollary B] 
Since a Benz\'ecri suspension over an orbifold is always finitely covered 
by a  Benz\'ecri suspension over  a surface  obtained using the identity 
automorphism,  and such a  Benz\'ecri suspension over a surface $\Sigma$ 
is a   trivial  circle bundle over $\Sigma$   (see  Chapter  
\ref{ch:racom}), it  follows   that  the  Euler   number is  zero    
(see  Scott \cite{Sc}).  Thus, if $M$   is  a Benz\'ecri  suspension,  
then   we  are done. If $M$ is a generalized affine suspension over  
a Euler  characteristic nonnegative  real projective surface,   
then $M$  is homeomorphic to  a 
bundle  over   a circle with fiber    homeomorphic to the   surface by 
considering only the topological aspects  of the generalized  
affine suspension. If the Euler  characteristic of a  fiber is  
positive,   then the fiber is either a  sphere,  $\rp^2$, or a disk,   
and  $M$ falls into  the first case.  If the Euler characteristic is zero,  
then $M$ is as in the second case of the conclusion of the corollary.  
\end{proof} 
 
\chapter[Concave-cone radiant $3$-manifolds]
{Concave-cone radiant affine $3$-manifolds and  
radiant concave affine $3$-manifolds}  
\label{ch:caff} 
 
In this chapter, we will first classify concave-cone affine manifolds, 
which is  an easy consequence of Tischler's  work \cite{Tis}.  Next, we 
classify radiant concave  affine manifolds. Let $M$ be a compact  
radiant concave affine manifold with empty or totally geodesic  
boundary. Let $M_h$ be its holonomy cover and $\che M_h$ the Kuiper 
completion. First,  we do the classification for the case when  
there exists a unique radiant bihedron in $\che M_h$.  We show 
that if  there are three   radiant bihedra in  $\che M_h$  in general 
position,  than $M$ is finitely  covered by a generalized  
affine suspension of a sphere.  We will do  this by  showing   that  
$\dev| \Lambc(R)^o \cap M_h$ is a diffeomorphism onto  
the complement of a closed set $K$ which is either 
a simply convex cone or the origin.  If $K$ is not the origin, then 
the deck-transformation group will act on the inverse  image $K'$ in  
$M_h$ of $-K$ as an infinite cyclic group, leading to a desired  
contradiction.  If $K$ is the origin only, then $M$ is shown to  be  
a generalized affine suspension  of a sphere or a real projective plane.   
Assume that there are   no three radiant bihedra  in  $\che  M_h$   
in general position,  which implies that the   developing  images  
of all  radiant bihedra contain a  common  complete real line   
passing  through $O$ in their respective boundaries.  
We show that $\Lambda^o(R)$  equals $M_h^o$.  
If $M$ has nonempty boundary, we easily obtain  
a contradiction by obtaining a deck-transformation-group  
invariant $3$-crescent.  If not, $M_h$ is a 
cover of the complement  of a complete real line in $\bR^3$; such a  
radiant affine $3$-manifold is shown not to exist by Barbot  
and Choi in Appendix \ref{app:radd}. 
 
\begin{thm}\label{thm:cafB}  
Let  $M$ be a compact concave-cone affine $3$-manifold with  empty or  
totally geodesic  boundary. Then $M$ is a generalized  affine suspension  
of affine tori, affine Klein bottles, affine annuli with  
geodesic boundary, or affine M\"obius bands with geodesic boundary.  
\end{thm}  
 
Note  that affine  tori, affine   Klein  bottles, affine  annuli  with 
geodesic boundary, or affine M\"obius bands with geodesic boundary are 
classified by  Nagano-Yagi essentially  (see Nagano-Yagi \cite{NY}).   
We use an   argument  communicated to us by T. Barbot.  
 
\begin{proof}  
We have  that $M_h \subset \Lambc(R)$ for  a crescent-cone $R$.  Since 
$\dfLambc(R)$ is totally geodesic  being  a union of totally  geodesic 
surfaces   $\alfi_T$,  $T  \sim      R$,  extending one another, 
$\dev(\dfLambc(R))$ is a subset of an affine plane $P$ passing through 
$O$, and $\dev(\Lambc(R))$ is a subset  of the closure of a half-space 
$H$ with boundary including $P$.  Given any deck transformation  
$\vth$ of $M_h$, since  $\vth(M_h)  =    M_h$,  we have  
$\vth(\Lambc(R)) =   \Lambc(R)$, and  
$\vth(\dfLambc(R))   =   \dfLambc(R)$; thus, $P$  is 
$h(\pi_1(M))$-invariant.  
 
Introduce affine  coordinates $x, y,$ and  $z$ on  $\bR^3$ so that the 
origin $O$ has coordinates $(0,0,0)$ and $P$  is the $xy$-plane.  Then 
the  $1$-form $dz/z$ is  invariant   under $h(\pi_1(M))$.  This   form 
induces a closed $1$-form on $M_h$ by $\dev$  which is invariant under 
the deck-transformation group,  and  the  form  descends to a  closed 
$1$-form transversal to the  radial flow. Hence by Tischler's argument 
(\cite{Tis}),   it follows that the radial   flow has a global section 
$S$ (see Lemma \ref{lem:form} of Appendix \ref{app:radd}). 
 
There exists a closed   hemisphere   $L$ which equals $\clo(H)    \cap 
\SIT$. Let $L'$ be the interior of $L$.  Let $\tilde S$ be a component 
of the  inverse image of $S$ in  $M_h$. Then $\tilde S$ is transversal 
to radial flows and a subgroup $G$ of  the deck-transformation group  
acts on $\tilde S$ so that $p$ induces a  homeomorphism  
$\tilde S/G \ra S$. For  the radial   projection   $\Pi$ from   
$\bR^3  - \{O\}$ to $\SIT$,  $\Pi  \circ  \dev  |\tilde  S$  maps   
into $L'$ as an immersion. Since $h(\pi_1(M))$ acts on $H$, it acts   
on $L'$.  By Berger \cite{B}, $L'$ can be identified  with an affine   
space and $h(\pi_1(M))$ acts as  affine transformations on $L'$.  
Moreover, $\Pi  \circ \dev| \tilde S$ is equivariant  with  respect  
to the  homomorphism   $i_{L'} \circ h$, where $i_{L'}$ is  
the restriction homomorphism to $L'$.  
 
This means that  $S$ has an  affine  structure, and since $\delta  S 
\subset \delta M$ and $\delta M$ is empty  or totally geodesic, it 
follows that $\delta S$ is empty or  geodesic. By the classification 
of affine surfaces, it follows that $S$  is an affine torus, an affine 
Klein bottle, an  affine annuli with  geodesic boundary,  or an affine 
M\"obius  band with geodesic boundary. Moreover,  the existence of the 
total cross-section tells us  that $M$ is a generalized affine suspension   
over these affine surfaces by Proposition \ref{prop:crssect}. 
\end{proof} 
 
\begin{thm}\label{thm:cafA} 
Let $M$ be a compact radiant concave affine manifold with empty or totally 
geodesic boundary.  Assume  that  $M$ is  not convex.  Then $M$ is   
a generalized affine suspension  of a real projective  hemisphere,  
a real projective plane, a real projective sphere, or  
a $\pi$-annulus {\rm (}or M\"obius band\/{\rm )} of type C.  
\end{thm} 
 
Since  $M$  is  radiant and  concave affine,    $M_h$ is  a  subset of 
$\Lambda(R)$ for a  radiant bihedron $R$  in $\che M_h$.  Recall  that 
$\Lambda(R)$  is the  union   of  all radiant  bihedra  equivalent  to 
$R$. Given any radiant  bihedron $S$ in $\che M_h$,  
$S^o$  meets a radiant  bihedron equivalent to $R$, and hence  
$S \sim R$.  
 
Suppose that there is  only one radiant   bihedron equivalent to  $R$, 
namely itself.  Then $M_h = R^o \cup (\nu_R^o \cap M_h)$.  If $M_h$ is 
closed, then $M_h  =  R^o$, which means  $\nu_R^o  \cap M_h =\emp$,  a 
contradiction to the definition of  $3$-crescent. Thus, $\delta M_h$ 
is not empty and equals $\nu_R^o \cap M_h$.  
 
Since $\delta M_h$ is   not empty,   
$M$ admits a total cross-section  $S$ to the radial flow by   
Theorem \ref{thm:barbot2}.  
Let $\tilde S$ be a  component of $p^{-1}(S)$. Defining 
$\Pi: \bR^3 - O \ra \SIT$ be the  radial flow, we  see that $\Pi \circ 
\dev|\tilde S$ is an immersion into $\dev(\clo(\alpha_R))$ equivariant 
with respect to  $h|  \pi_1(S)$.  Furthermore, $\Pi\circ  \dev| \tilde 
S^o$ is an imbedding onto the open hemisphere $\dev(\alpha_R)$.  Thus, 
$S$  is a  real projective surface    such that $\tilde  S^o$ is  real 
projectively     isomorphic to   an    open  hemisphere.     
By the following theorem \ref{thm:nohemrps},    
$M$ is a generalized affine suspension of an hemisphere, a $\pi$-annulus or  
a $\pi$-M\"obius band of type C.

\begin{thm}\label{thm:nohemrps} 
Let $S$ be a compact real projective surface with non-empty  
geodesic boundary. Then if $\tilde S^o$ is real projectively 
diffeomorphic to an open hemisphere, then $S$ is  
either a closed hemisphere, a $\pi$-annulus 
of type C or a $\pi$-M\"obius band of type C. 
\end{thm} 
\begin{proof} 
Since $\tilde S^o$ is dense in $\tilde S$,  
its closure in $\che S$ equals $\che S$. 
Since $\tilde S^o$ is tame, $\che S$ is tame, and  
$\dev:\che S \ra \SI^2$ is an embedding onto its images 
$\dev(\che S)$, a closed hemisphere $H$. Thus, we see that  
$\dev| \tilde S$ is an imbedding onto its image $\Omega$, 
a domain whose interior is an open hemisphere $H^o$.  
Moreover, $\ideal{S}$ clearly maps into $\delta H$. 
 
From   the classification  of  real   projective   surfaces   
(see   \cite{cdcr1}, \cite{cdcr2}, and \cite{cdcr3}),  
if the Euler characteristic $\chi(S)$ of $S$ is negative,  
either $S$ itself is convex, or $S^o$ contains a simple closed 
geodesic $c$ such that its lift $\tilde c$ to $\tilde S^o$ 
is of $\bdd$-length $< \pi$. The closure of  
the image of $\tilde c$ is a simply convex segment $l$ with  
endpoints in $\ideal{S}$, and $\dev|l$ is a diffeomorphism 
to its image. However, there exists no simply convex  
line in $H^o$ ending in $\delta H$ in $H$. 
 
If $S$ is convex, then Lemma 1.5 of \cite{cdcr1} shows that  
$\tilde S^o$ is simply convex. As $H^o$ is not 
simply convex, this is a contradiction. 
 
If $\chi(S) > 0$, then since $S$ has nonempty boundary, 
$S$ is diffeomorphic to a disk and hence $S$ is a closed  
hemisphere. 
 
If $\chi(S) = 0$, then assume that $S$ is orientable, 
by taking a double cover if necessary. Since $S$ is an annulus,  
$\tilde S \cap \delta H$ has two components on each of  
which the generator $\vth$ of $\pi_1(S)$ acts on. 
Since $\vth$ is orientation-preserving, we may assume  
that $\vth$ is of form (1)-(7) in Section 1.1 of  
\cite{cdcr1} and the eigenvalues of $\vth$ is positive.  
by taking a double cover of $S$ if necessary. 
Then since $\langle \vth \rangle$ acts without fixed points on $H^o$,  
$\vth$ is of type (2) or (4).  
(4) can be ruled out since one cannot get a compact quotient   
that is Hausdorff. Now, investigating when $\langle \vth\rangle$  
acts properly and freely, we can easily show that  
$S$ must be a $\pi$-annulus of type C. (We omit  
this somewhat tedious verification.) 
\end{proof}

We will now assume that there are more than one radiant bihedron.  There 
has to be two  overlapping  bihedra since $\Lambda(R)$  contains $\che 
M$.  
 
Let $R_1, R_2, R_3$ be radiant bihedra  equivalent to $R$. We say that 
$R_1, R_2,$ and $R_3$   are in general position  if $\dev(\nu_{R_1})$, 
$\dev(\nu_{R_2})$, and $\dev(\nu_{R_3})$ are in general position.

\begin{figure}[b] 
\centerline{\epsfxsize=3.6in\epsfbox{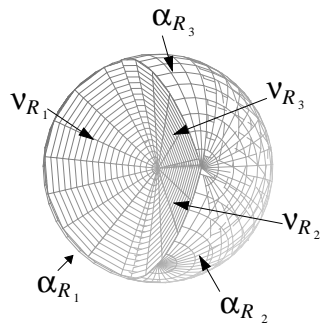}}  
\caption{\label{fig:big} Three radiant bihedra in general position.}  
\end{figure} 
 
\begin{prop}\label{prop:cafA1}  
Suppose that there exists mutually overlapping three radiant bihedra 
$R_1$, $R_2$, and   $R_3$     in  general position.   Then       
$\dev|\Lambda_1(R)\cap M_h$ is an imbedding onto $\bR^3 - O$ or  
the complement of a closed simply convex radiant set $K$ in $\bR^3$.  
\end{prop}  
\begin{proof}  
First,  $\dev| R_1  \cup R_2$  is   an imbedding  onto $\dev(R_1) \cup 
\dev(R_2)$ by Theorem \ref{thm:transversal}. Since 
$\dev(R_3) \cap (\dev(R_1)  \cup \dev(R_2))$ is a connected  $3$-ball, 
it follows that  $\dev| R_1 \cup R_2 \cup   R_3$ is an  imbedding onto 
$\dev(R_1)    \cup    \dev(R_2)   \cup     \dev(R_3)$  by  Proposition 
\ref{prop:extmap}.  
 
We assume that all crescents denoted by different symbols are  
distinct in $\che M$ in this proof. Let  $S$ be an arbitrary  
radiant bihedron equivalent to $R$ and overlapping some $R_i$ for  
$i = 1, 2,  3$. Assume without loss  of generality that $S$ 
overlaps  with $R_3$.  $R_1 \cap R_3$  and   $R_2 \cap R_3$  are  two 
distinct radiant trihedra.  From looking at images of these sets under 
the imbedding  $\dev| R_3$,  we see  that $R_3 \cap   S$ is  a radiant 
trihedron and overlaps  with at least  one of $R_1  \cap R_3$ and $R_2 
\cap R_3$.  Assuming  without loss  of generality  that  $R_3 \cap  S$ 
overlaps with  $R_2  \cap  R_3$, we  see  that  $S$, $R_2$, and  $R_3$ 
overlap mutually with one another.  Hence $\dev|  S \cup R_2 \cup R_3$ 
is an imbedding onto $\dev(S) \cup \dev(R_2) \cup \dev(R_3)$ as in the  
above paragraph.

Consider $\dev(R_1)^o \cap (\dev(R_2)^o \cup \dev(R_3)^o  \cup 
\dev(S)^o)$.  Since this also is a  connected open $3$-manifold by  
the following Lemma \ref{lem:ints}, Proposition \ref{prop:overlap} 
and Remark \ref{rem:overlap} show that  
$\dev| R_1^o \cup R_2^o \cup R_3^o \cup S^o$  is an imbedding.  
Since we may hence read the intersection pattern of  
$R_1^o, R_2^o, R_3^o,$ and $S^o$ by their image half-spaces  
in $\bR^3$, we obtain that $R_i, R_j,$ and $S$ are mutually  
intersecting and are in general position for some choice of  
pair $i, j$, $i, j = 1,2,3$ (i.e., check that the other possibilities  
cannot occur).

A collection of radiant bihedra $S_1, \dots, S_n$ is  
said to be a {\em chain} if $S_i$ overlaps with $S_{i+1}$ for 
$i=1, \dots, n-1$ with $n > 1$, or if $n=1$. 
We say that a collection $R_1, R_2, R_3, S_1, \dots, S_n$ is  
a {\em chained} collection if any two elements of  
the collection are connected by a chain of  
overlapping elements in the collection. 
We will now prove by induction that  
given a chained collection $R_1, R_2, R_3, S_1, \dots, S_n$  
\begin{enumerate} 
\item[(i)] $\dev| R_1^o \cup R_2^o \cup R_3^o  
\cup S_1^o \cup S_2^o \cdots \cup S_n^o$  
is an  imbedding, and 
\item[(ii)] given $i$, $i=1, \dots, n$,  
$S_i$ overlaps with $R_j$ for some $j$, $j=1,2,3$. 
\end{enumerate} 
 
For $n=1$, this was proved above.  
For $n=2$, let $R_1, R_2, R_3, S_1, S_2$ be a chained collection.  
Assume without loss of generality that $R_1, R_2, R_3, S_1$  
is a proper chained collection. 
By above, $S_1$ overlaps with $R_i$ and $R_j$ for 
some pair $i, j$, $i, j = 1,2,3$ and they are in general position. 
If $S_2$ overlaps with $S_1$, then an open 
trihedron $S_1^o \cap S_2^o$ must  
meet one of the open trihedra $S_1^o \cap R_i^o$  
and $S_1^o \cap R_j^o$ by geometry and Proposition  
\ref{prop:extmap}. Hence, we can assume that $S_2$ 
overlaps with at least one of $R_1, R_2,$ and $R_3$,  
and (ii) follows. 
 
Since the collection $R_1, R_2, R_3, S_1$ and  
the collection $R_1, R_2, R_3, S_2$ satisfy (i) and 
the intersection of the respective unions of 
the images of their interiors 
\begin{eqnarray}  
(\dev(R_1)^o \cup \dev(R_2)^o \cup \dev(R_3)^o) \cup  
(\dev(S_1)^o \cap \dev(S_2)^o)
\nonumber 
\end{eqnarray} 
is clearly a connected open subset, (i) follows  
for the collection $R_1, R_2, R_3, S_1, S_2$ 
by Remark \ref{rem:overlap}. 
 
Assume that (i) and (ii) are true if $n \leq k$ for  
an integer $k$. By above, we assume that $k \geq 2$. 
We will now prove for $n = k+1$.  
Let $R_1, R_2, R_3, S_1, \dots, S_{k+1}$ be  
a chained collection. We may assume without loss of generality 
that \[R_1, R_2, R_3, S_1, \dots, S_m\] is a proper chained  
subcollection with $m \leq k$ and $S_{k+1}$ overlaps with  
$S_l$ for $l \leq m$. 
Then by (ii) for $n \leq k$, $R_1, R_2, R_3, S_l, S_{k+1}$  
is a chained collection, and hence  
$\dev| R_1^o \cup R_2^o \cup R_3^o \cup S_l^o \cup S_{k+1}^o$  
is injective and $S_{k+1}$ overlaps with $R_j$ for some  
$j =1,2,3$. (ii) follows for $n= k+1$. 
 
We claim that $\dev$ restricted to  
$R_1^o \cup R_2^o \cup R_3^o \cup S_1^o \cup \cdots \cup S_{k+1}^o$  
is injective, which will prove (i) for $n= k+1$. If not, there  
exist points $p \in S_l$ and $q \in S_{m}$, some pair $l, m \leq k+1$,  
so that $\dev(p) = \dev(q)$. However since $R_1, R_2, R_3, S_l$  
and $R_1, R_2, R_3, S_m$ are both chained collections by above,  
it follows that $R_1, R_2, R_3, S_l, S_m$ is 
a chained collection, and $\dev$ restricted to  
$R_1^o \cup R_2^o \cup R_3^o \cup S_l^o\cup S_m^o$ is  
injective by (i) for $n=2$, a contradiction. 
 
From (i) for all $k$, it follows that $\dev|\Lambda_1(R) \cap M_h$ is 
injective.  
The complement of $\dev(\Lambda_1(R) \cap M_h)$ is  
disjoint from $\dev(R_1), \dev(R_2)$, and $\dev(R_2)$. Hence  
it must be simply convex. 
\end{proof} 
 
\begin{lem}\label{lem:ints}  
Let $L_1$  be an  open half space  in $\bR^3$,  and $L_2, L_3,  
L_4$ open  half spaces overlapping  with $L_1$ whose respective  
boundary planes $P_1, P_2, P_3, P_4$  pass through $O$.  If 
$P_1, P_2,$ and  $P_3$ are in general  position,  then  
$L_1 \cap  (L_2 \cup L_3 \cup L_4)$ is a connected open $3$-manifold.  
\end{lem}  
\begin{proof}   
We will show that $P_1 \cap (L_2 \cup L_3 \cup L_4)$ is  
a connected open set. Now our statement easily implies the lemma.  
 
Since  $P_1, P_2$, and $P_3$ are  in general position, it follows that 
$P_1 \cap P_2$ and  $P_1 \cap P_3$ are  distinct lines passing through 
the origin.  Thus,  $P_1 \cap  L_1$ and $P_1   \cap L_2$ are  distinct 
half-planes meeting in an open set.  Hence $P_1 \cap L_4$ must meet at 
least one of them in a nonempty open set, and 
$P_1 \cap (L_2 \cup L_3  \cup L_4)$ is  a connected  open set; this 
proves the above statement.  
\end{proof}  
 
Suppose that we are in the first case of the conclusion of  
Proposition \ref{prop:cafA1}. Then $\dev|\Lambda_1(R) \cap M_h$ is  
a homeomorphism onto $\brto$. Since  
$\Lambda_1(R) \cap M_h$ is dense in $M_h$,  
and $\dev(M_h)$ is a subset of $\brto$, it follows that  
$\Lambda_1(R) \cap M_h = M_h$, and 
$\dev: M_h \ra \mathcal{H}^o - \{O\}$ is a homeomorphism, 
implying $M$ is boundaryless.   
By the proof of Proposition \ref{prop:Hopf}, $M$ is  
a generalized affine suspension over a real projective sphere  
or a real projective plane. 
 
Suppose that we are in the second case; we will  show this leads to 
contradictions.  Since $K$ is simply convex, the set $-K$ consisting 
of $-v$ for $v \in  K$ is disjoint from  $K -  \{O\}$  and hence is  a 
proper subset of $\dev(\Lambda_1(R) \cap M_h)$.  
 
Since $\dev|\Lambda_1(R) \cap M_h$ is  
an  imbedding onto the image by Proposition 
\ref{prop:cafA1}, it  follows that there exists a set $K'$ in  
$\Lambda^o(R)$ such that $\dev|K'$ is a  diffeomorphism  onto  
$-K  -   \{O\}$.  For any deck transformation $\vth$, since    
$\vth(M_h^o)  = M_h^o$ and   $M_h^o$  is a   subset of $\Lambda(R)$,  
it  follows  that by Proposition \ref{prop:splitting},  
$\vth(\Lambda_1(R)) = \Lambda_1(R)$. Hence 
$h(\vth)(\dev(\Lambda_1(R))) = \dev(\Lambda_1(R))$,  $K$ is 
$h(\pi_1(M))$-invariant,   and the deck-transformation  group  
$G$ of $M_h$ acts on $K'$.  
 
Assume that $K'$ has non-empty interior. 
Since the deck transformation group acts on $K'$,  
the image $p(K')$ in $M$ is a  closed submanifold.  
Since $K'$ is homeomorphic to $D^2 \times \bR$ for a disk $D^2$ 
and  is a radiant set, the  boundary of $K'$ which is 
homeomorphic to  $\SI^1 \times \bR$  covers the boundary of  $K'/G$, and 
hence the  boundary of $K'/G$  is homeomorphic to a torus or a Klein  
bottle.  Since boundary is compressible, there exists an imbedded disk  
in $K'/G$ with boundary $\delta K'/G$ by the loop theorem. Since $K'$   
is  simply connected, $K'/G$ is homeomorphic  to a solid torus  or  
a solid Klein bottle. Hence $G$ is isomorphic to  
the group of integers $\bZ$. If $K'$ has empty interior,  
$K'$ is homeomorphic to $I \times \bR$ for an interval $I$ or  
to $\bR$. Since $G$ acts on $K'$ freely and properly,  
$K'/G$ is homeomorphic to an annulus, a M\"obius band, or a circle, and  
$G$ is isomorphic to $\bZ$. 
 
Recall that  the boundary of $M$ is  either empty or totally geodesic. 
Suppose the  boundary of $M$ is not  empty. Let $T$  be a component of 
$\delta M$, an affine torus  or an affine Klein bottle. Since the   
fundamental group  of  $\pi_1(M)$ equals $\bZ$, it follows that a simple   
closed  curve  $\alpha$ in $T$ is a boundary of an imbedded disk  $D$  
in $M$ by the loop theorem. Since the  holonomy of $\alpha$ is trivial,  
the  classification of affine structures on tori or Klein bottles  
(see Nagano-Yagi \cite{NY}) shows  
that  a component $T'$ of $p^{-1}(T) \subset \delta M_h$ is affinely  
isomorphic to a finite cover of $\bR^2 - \{O\}$.  
 
Since   $\delta M_h$ is disjoint from $\Lambda_1(R) \cap M_h$,   
$\Lambda_1(R) \cap M_h$ is dense  in $M_h$, and  
$\dev|  \Lambda_1(R) \cap M_h$ is an imbedding onto  
$\bR^3 -K$, it follows that $\dev(\delta M_h)$ is  
a subset of $K$. Since $K$ is a  simply convex subset,  
this contradicts the fact that an affine copy  
of a cover of $\bR^2  - O$ contains complete affine lines. 
(A {\em complete affine line} is a complete affine one-dimensional  
subset of $\bR^n$. It is real projectively homeomorphic  
to an open $1$-hemisphere, and has $\bdd$-length $\pi$.) 
 
Suppose next that $\delta M$  is empty.  
Let $t$ be a generator of $\pi_1(M)$.  Then since any $3\times 3$-matrix 
has a real eigenvalue, $h(t)$ must have a pair of antipodal fixed points  
in $\SIT$ (by taking  a double cover of $M$ if necessary). They  
correspond to   a pair of  antipodal radial  line in $\mathcal{H}^o$. At  
least one of them, say  $l$, is in $\dev(\Lambda_1(R) \cap M_h)$  
since $\dev(\Lambda_1(R) \cap M_h)$ is the    
complement of a simply convex  cone  in $\bR^3$. The radial  line $l'$  
in $M_h$  corresponding  to $l$ maps to a periodic orbit  $c$ in $M$.  
Suppose that $c$ is of saddle type in Barbot's  terminology  
(see Theorem \ref{thm:barbot}). 
Then  $h(t)$ has  three distinct eigenvalues, i.e.,  
diagonalizable. We may assume without loss of generality that $h(t)$  
has positive eigenvalues (by taking a double cover of $M$ if necessary).  
Then choosing  the fixed point  in $\SIT$ of attracting or repelling  
type, we obtain another radial line $m$ in $\Lambda_1(R) \cap M_h$  
fixed by $t$. Then $m$ corresponds in $M$  to a periodic orbit  
which is  not of saddle type, and by Theorem \ref{thm:barbot},  
$M$ is a generalized affine suspension of a projective surface $\Sigma$.  
If $c$ is not of saddle type to begin with, then $M$ is  
a generalized affine suspension of a projective surface $\Sigma$.   
 
Since  the fundamental group of $M$ equals $\bZ$, it follows that the  
projective surface $\Sigma$ is a disk or sphere. The surface cannot be  
a  disk since $M$ is closed. If the surface is a  sphere, then since the   
projective  structure on the  sphere is  unique up to isotopy, it  
follows that $M$ is obtained from a generalized affine suspension  
of the real projective sphere. Hence, our original $M$ is covered  
by the generalized affine suspension of the real projective sphere.  
However this means that $K$ is the origin only, a contradiction.  
 
Now we can  assume the following statement:  Whenever $R_1,  R_2$, and 
$R_3$  are  mutually overlapping   radiant  bihedra, $\dev(\nu_{R_1}), 
\dev(\nu_{R_2})$, and $\dev(\nu_{R_3})$ include a common complete real 
line.  We  will  show  that such a case is  impossible,  and this will 
complete the proof of Theorem \ref{thm:cafA}.  
 
By induction,  we will show  that $\dev(\nu_S)$  for   any $S \sim  R$ 
includes $l$   for a fixed complete   affine line $l$.   Let $S$  be a 
radiant bihedron so that $S \sim R$.  We  choose a sequence $R_1, R_2, 
R_3, \dots, R_n$ so that $R_i$ and $R_{i+1}$ overlap for $i =1, \dots, 
n-1$ and $R_n  = S$, and  $\dev(\nu_{R_i})$ contains a complete real 
line $l$ for $i=1,2$.  
 
We  may assume that $R_i$  are  distinct from  $R_{i+1}$ for each $i$ 
so that $\dev(\nu_{R_i})$ and $\dev(\nu_{R_{i+1}})$ are transversal.  
Assume that  $\dev(\nu_{R_i})$ contains  $l$  for $i=1, 2, \dots, k$   
for  $2 \leq k   < n$. We   will prove this  for  $k+1$. If 
$\dev(\nu_{R_{k+1}})$      does     not        contain    $l$,    then 
$\dev(\nu_{R_{k-1}})$,  $\dev(\nu_{R_k})$,  and  $\dev(\nu_{R_{k+1}})$ 
are  in general position. Since  $R_{k-1}$ and $R_k$ overlap and $R_k$ 
and $R_{k+1}$ overlap,  
$\dev| R_{k-1}   \cup  R_k  \cup   R_{k+1}$  is  an   imbedding   onto 
$\dev(R_{k-1})  \cup \dev(R_k)   \cup   \dev(R_{k+1})$ by  Proposition 
\ref{prop:extmap}.   
This means that $R_{k-1}, R_k$, and $R_{k+1}$ are mutually  
overlapping radiant bihedra in general position. 
This   case   was  ruled   out above.   Therefore, we conclude that 
$\dev(\nu_S)$ includes $l$ for any $S \sim R$.  
 
We will now prove that $M_h = \Lambda(R) \cap M_h$.  
\begin{lem}\label{lem:l}  
$\dev(M_h)$ misses $l$.  
\end{lem}  
\begin{proof}  
Suppose that $\dev(x)$ belongs to  $l$ for $x  \in M_h$.   
Since $M_h \subset \Lambda(R)$, there exists a radiant   
bihedron $S$ containing $x$ with $S \sim  R$.   
Since there is
more than one bihedral crescent in $\che{M_h}$, it follows that there 
is a bihedral crescent $S_1$ transversally meeting with $S$ such that 
$\dev|S \cup S_1$ is an  imbedding onto $\dev(S) \cup \dev(S_1)$,  and 
$\dev(S)$ is distinct from $\dev(S_1)$.  
 
Since $\dev(x) \in l$, we have $\dev(x) \in \dev(\nu_{S_1})$,  
and this means $x \in \nu_{S_1}$. Similarly, $x \in \nu_S$.  
 
Note that if two properly imbedded surfaces meet at a point  
tangentially, then they must agree completely. 
If $x \in \delta M_h$, then  since $x \in  \nu_S$ for $S \sim R$, it 
follows  that a component  of $\nu_S \cap  M_h$  equals a component of 
$\delta M_h$. Thus, if for another crescent $S'$, 
$x \in \nu_{S'}$, then $\nu_S \cap M_h$ and  $\nu_{S'} \cap M_h$ agree 
on this component.  By Corollary \ref{cor:radtrs}, it follows that that 
$S = S'$ actually; however, our choice of $S_1$ contradicts this.  
 
Suppose that $x \in M_h^o$.  Let $B$ be a small tiny ball neighborhood 
of $x$ in  $M_h^o$.  Then $\dev| \bigcup_{t  \in \bR} \Phi_{h, t}(B)$ is an 
imbedding   onto a simply  convex cone  $C'$ since  $\dev(x) \ne O$. 
Letting $C$ equal to this  union, we see  again that $\dev| C \cup S^o 
\cup S_1^o$    is an  imbedding  onto   $\dev(C)  \cup  \dev(S^o) \cup 
\dev(S_1^o)$.   By  the  geometry  of  the  situation,  there exists a 
radiant bihedron  $S'$ in the  image so that  its sides, $\dev(\nu_S)$, 
and $\dev(\nu_{S_1})$ meet  in   general position.  Let $S''$  be  the 
inverse image in  $C \cup S^o \cup S_1^o$.  Then the closure  $S_2$ of 
$S''$ in $\che{M_h}$  is a  radiant bihedron  such that  $S, S_1,$ 
and $S_2$ overlap mutually and $\dev(\nu_S), \dev(\nu_{S_1}),$     
and $\dev(\nu_{S_2})$ are in   general position, contradicting what   
we just assumed.  
\end{proof} 
 
\begin{lem}\label{lem:Mhsub} 
$M_h = M_h^o = \Lambda_1(R) \cap M_h$,  
and hence $M_h^o \subset \Lambda_1(R)$.  
\end{lem} 
\begin{proof} 
Since $S^o$ for any crescent $S$, $S  \sim R$ is  a subset of $M_h^o$, 
we obviously have $\Lambda_1(R)\cap M_h \subset M_h^o$.  
 
We  prove that $M_h^o \subset \Lambda_1(R)$.   Suppose not. Then there 
exists a  point $x$ of $M_h^o$ in  $\Bd \Lambda_1(R) \cap M_h$,  
which is not in $T^o$ for any $T \sim R$.   
Since $\Lambda(R) \cap M_h$ is closed, we  have $x \in  \nu_S$ for  
some $S \sim R$ and by the above lemma, $\dev(x)$ is not in $l$.  
  
Let $l'$ denote the inverse image of $l$  in $\nu_S$.  There has to be 
at least one  radiant bihedron $T$ overlapping  with $S$ and  distinct 
from $S$.  By  Corollary \ref{cor:radtrs},    $S$ and  $T$    intersect 
transversally. Since $\dev(\nu_S)  \cap \dev(\nu_T)$ contains  $l$, it 
follows that   $\nu_S \cap \nu_T$   contains $l'$, and there  exists a 
component $E$ of  $\nu_S^o -l'$ contained in  $T^o$ as we can see from 
the conditions for transversal intersections  
in Appendix \ref{app:dipping}.  While $E$ does not contain $x$ since  
$x$ does not belong  to $T^o$, we conclude  that a component,  
namely  $E$,  of $\nu_S^o - l'$  not  containing $x$   
is  an open lune completely contained in $M_h^o$.  
  
Assume without loss  of generality  that  $h(\pi_1(M))$ fixes $p$  and 
$-p$  (by taking a double  cover if necessary); that is, $h(\pi_1(M))$ 
acts on $l$ preserving the orientation of $\SI^3$.  We claim that $S$  
cannot overlap with its image under a deck transformation distinct  
from $S$. If our claim is false, then $S$  and $\vth(S)$ for  
a deck transformation $\vth$ intersect  transversally by  
Corollary \ref{cor:radtrs}. By an orientation   condition, $x$  
belongs to $\vth(S^o)$ or $\vth(x)$ belongs to $S^o$. In the latter  
case, $x$ belongs to $\vth^{-1}(S^o)$. Both of these are  
impossible. Therefore, for any  deck transformation $\vth$  
of $M_h$, we have either  $S = \vth(S)$, $S \cap \vth(S) \cap M_h =  
\emp$,  or $S \cap \vth(S)\cap  M_h$ is a union of common components 
of $\nu_S \cap M_h$ and $\nu_{\vth(S)} \cap M_h$ where $S$ and $\vth(S)$  
do not overlap as in \cite{psconv}.  
 
For a component  $A$ of $\vth(\nu_S)  \cap M_h$ for  
a deck transformation $\vth$ of $M_h$, we say  that  
$A$ is a {\em copied component for $\vth(S)$}\/ if $A$ is  
a component of $\vpi(\nu_S) \cap M_h$ for a deck transformation $\vpi$  
where $\vth(S)$ and $\vpi(S$) do not overlap.  
The union $\hat A$ of all copied components of $\vth(S)$ as $\vth$  
ranges over deck transformations $\vth$ of $M_h$ is  
a properly imbedded submanifold in $M_h$. The proof 
of this fact is similar to the proofs for pre-two-faced  
submanifolds for $3$-crescents and crescent-cones 
above (see also \cite{psconv}). Therefore, we see 
that $p(\hat A)$ must be a properly imbedded closed surface in $M$. 
 
We split off along $p(\hat A)$ to obtain the result $N$. As before  
the holonomy cover $N_h$ of $N$ is a disjoint union of components  
of $M_h$ split along $\hat A$. 
The Kuiper completion $\che N_h$ of $N_h$ is obtained by completing   
with respect to the metric induced by $\dev:N_h \ra \SI^3$ 
which is the immersion extended from $M_h - A'$. 
 
We can choose $N_h$ so that $N_h$ includes $S^o$. 
The closure $S'$ of $S^o$ in $\che N_h$ is a radiant bihedron in 
$\che N_h$. As in crescent cases \cite{psconv}, the splitting  
process insures that there is no copied  component for $\vth(S')$  
in $\che N_h$ where $\vth$ is a deck transformation of $\che N_h$. 
Lemma \ref{lem:invcr} shows that $E$ must be  a triangle,  
which is a contradiction.  
 
If there are  no  copied  components  in  $\nu_S \cap M_h$,  then  $S$ 
satisfies the premise  of Lemma \ref{lem:invcr}  showing that $E$ is a 
triangle, a contradiction. This shows that $M_h^o$ is a subset of 
$\Lambda_1(R)$. 
 
If $\delta M_h$ is not empty, then a point $x$ of $\delta M_h$ is  
in $\Bd \Lambda_1(R)$. $x$ must be a point of a radiant bihedron  
$S$, $S \sim R$, which has the same property as the radiant  
bihedron denoted by $S$ above. The same argument as above shows that  
this cannot happen. Hence it follows that $M_h = M_h^o$. 
\end{proof} 
 
\begin{lem}\label{lem:invcr} 
Suppose that given a radiant bihedron $S$,  either we have $S=\vth(S)$ 
or $S  \cap  \vth(S)  \cap  M_h = \emp$  for   any deck transformation 
$\vth$. Suppose that a complete affine line $l'$ in $\hideal{M}$ lies in  
$\nu_S$. Then   $\nu_S \cap M_h$   is a  union  of  two  radiant convex 
triangles respectively in the two components of $\nu_S -l'$.  
\end{lem} 
\begin{proof} 
Since the collection whose elements are of form  $\vth(S)$ is locally  
finite by Proposition \ref{prop:lfin3}, Proposition \ref{prop:lfin}  
implies that $S^o \cup (\nu_S\cap M_h)$ covers a compact three-dimensional  
submanifold $K$ in  $N$  with boundary  $p(\nu_S  \cap M_h)$.   Since  
$K$ is a radiant affine $3$-manifold  
it follows  that $K$ admits  
a total cross-section  $S$ by Theorem \ref{thm:barbot2}.   
Then by Theorem \ref{thm:nohemrps}, $K$ is a generalized  
affine suspension over a hemisphere, or a $\pi$-annulus  
(or M\"obius band) of type C.  
The first case  is not  possible  as  $\nu_S  \cap  M_h$  has  
at  least two components respectively in each side of  $\nu_S - l'$.  
The second case implies the conclusion of the lemma.  
\end{proof} 

Given a lune $L$ in  $\mathcal{H}$ with  boundary equal  to $l$ union  
a segment in $\SIT$ with endpoints $p$  and $-p$, we  say that  
$L -l$ is a {\em  proper lune}.  A subset $J$ of a lune in  
$\che{M_h}$  such that $\dev(J)$ is a proper lune in $\mathcal{H}$  
is said to be a {\em proper lune}.

\begin{lem}\label{lem:lift}  
$\dev|\Lambda_1(R):   \Lambda_1(R)  \ra \mathcal{H}-l$ lifts  to  an 
imbedding $\dev'$  to some cover $W$  of $\mathcal{H}  -l$, and the image 
$\dev'(M_h)$ is a union of open radiant bihedra in $W$.  
\end{lem} 
\begin{proof}  
Let $\SI^1$ be the great circle in $\SIT$  which is perpendicular to a 
segment  on $\SIT$ connecting $p$ and  $-p$.  Given a radiant bihedron 
$S$, $S \sim R$, let $l'_S$ denote the inverse image of $\dev^{-1}(l)$ 
in $S$.  Since for each point $x$ of $\Lambda_1(R)$,  
 $x \in T^o \cup \alpha_T$ for some $T \sim R$,  
it follows that $x$ belongs to a unique proper lune $J$.  
 
Let  $\mathcal{P}$ denote  the set of all  proper  lunes in $\Lambda_1(R)$. 
Given the topology of geometric convergence to $\mathcal{P}$ with respect to 
the metric  $\bdd$, $\mathcal{P}$  becomes  a differentiable $1$-manifold  
(without boundary).  
 
There exists an immersion $f: \mathcal{P} \ra \SI^1$ given by mapping $J$ 
to the intersection of $\SI^1$  with the closure  of $J$.  If $\mathcal{P}$ 
is simply connected, i.e., $\mathcal{P}$ is an open arc, then we let $W$ be 
the  infinite cyclic cover    of $\mathcal{H}  -l$.    If  $\mathcal{P}$  is 
homeomorphic to a circle, then let $W$ be the $n$-fold cyclic-cover of 
$\mathcal{H}-l$ where $n$ is the degree of $f$.  Then $\dev|\Lambda_1(R)$ 
lifts to an injective map into $W$ since no two  proper lunes map into 
same one in  $W$. Since $\dev|\Lambda_1(R)$ is an  open map, the first 
statement follows.  Since  $M_h = \Lambda_1(R) \cap M_h$,  
$\dev'(M_h)$ is a union of open radiant bihedra in $W$.  
\end{proof}

We let $q: W \ra \mathcal{H} -l$ denote the  covering map.  We denote by 
$W_\infty$ the part of $W$ covering $\SIT - l$. Using the radiant 
vector field on $W$ induced  from that of $\mathcal{H} -l$, we define  
$\Pi:\mathcal{H} -l \ra  \SIT -l$  the radial  projection and   
$\Pi_W: W \ra W_\infty$ that induced from  it  so that  
$\Pi  \circ  \dev| M_h =  q\circ \Pi_W \circ \dev' | M_h$.  
 
We see that $W_\infty$ is foliated by lines  corresponding to lines in 
$\SIT$ with endpoints $p$ and $-p$.  Since $M_h = \Lambda_1(R) \cap M_h$,  
$\Pi_W\circ \dev'(M_h)$ is a  union of  the lines  in this foliation   
$\mathcal{F}$.  Letting $\SI^1_W$ denote  the inverse image of $\SI^1$ in  
$W_\infty$,  we see that $\Pi_W\circ \dev'(M_h) \cap \SI_W$ is  
a connected open set $U$.  
 
Note that $\SI^1_W$ has a path-metric $\bdd$ induced  from  
the Riemannian metric pulled from $\SI^1$ by the covering map. 
The $\bdd$-length of $\SI^1_W$ is greater than equal to $\pi$  
since $M_h$ includes at least one radiant bihedron. Since  
there are more than one radiant bihedra, we have $\bdd$-length  
$> \pi$. 
 
Suppose that $U$ has an endpoint $q$ in $\SI^1_W$.  There exists a point 
$q'$  which  is at  the $\bdd$-distance $\pi$   away from $q$  in $U$. 
Since    $\dev'(M_h)$   is  the      union   of   open   hemispheres 
$\dev'(\alpha_S)$ for some $S \sim R$, it follows that there exists an 
open hemisphere $H$ in $W_\infty$ such that  $H \cap \SI^1_W$ equals the 
arc connecting $q$ and $q'$, and $H$ equals $\Pi_W(T^o)$ for a radiant 
bihedron $T$ in $\dev'(\che{M_h})$, where $T \sim R$. (To see this use 
the fact that $\dev'$ is an imbedding$: M_h^o \ra W$.)  
 
Since $T$ is a lune such that $\Pi_W(T^o)$ contains an end of $U$, and 
we assumed  that  $h(\pi_1(M))$ acts  on   $p$  and $-p$ and   
we assume that $h(\pi_1(M))$ acts on  the foliation $\mathcal{F}$  
in an orientation preserving  manner (by taking a 
double cover of  $M$ if necessary),  then the deck transformation group 
of $M_h$ acts on $T$.   At least one component  of $\nu_T^o - l'_T$ 
is an open lune in $M_h$ as in the proof of  Lemma \ref{lem:Mhsub}  
since  $T$ has to overlap with at least one radiant  bihedron  
different from itself.  By Lemma \ref{lem:invcr},  
we get  a contradiction as  before.  Therefore, $\dev'(M_h) \cap \SI^1_W$   
is  infinitely long  in both directions;  or equals $\SI^1_W$ itself and  
$W$ is a cyclic cover of $\mathcal{H}-l$.  
 
This means  that  $\dev'| \delta_\infty  \Lambda(R)$  is a map  onto 
$W_\infty$, and  $\dev'|M_h$ is an  imbedding onto $W$.   
This is a contradiction by the following lemma. 
 
\begin{lem}[Barbot-Choi (Appendix \ref{app:radd})] \label{lem:BC} 
A radiant affine $3$-manifold $N$ does not have    
the developing map  $\dev: N_h \ra \mathcal{H}^o$ that is  
a finite or infinite cyclic covering map onto  
$\mathcal{H}^o -l$.  
\end{lem} 
Obviously such $N$ must be a closed manifold. 
 
\chapter{The nonexistence of pseudo-crescent-cones}  
\label{ch:pcrc} 
 
As promised in Chapters  \ref{ch:2conv} and \ref{ch:obtain}, 
we will show that pseudo-crescent-cones do not exist 
assuming that $M$ is not convex. 
A radiant affine $3$-torus, as constructed by J. Smillie 
\cite{smillie}, which is obtained as a quotient  
of a radiant tetrahedron in $\mathcal{H}$, the Kuiper completion,
the closure of the tetrahedron, of the holonomy cover contains 
a pseudo-crescent-cone. 
 
In this chapter,  we will  work on the Kuiper completion $\che M$ of  
$\tilde M$ as we can lift any pseudo-crescent 
on  $\che M_h$ to $\che M$ by Lemma \ref{lem:pstart}.   
We will denote by $\dev$ the  developing map $\tilde  M \ra \SI^3$ and 
by $\bdd$ the metrics on $\che M$ and 
$\che M_h$ induced from $\SI^3$ by developing maps. 
(This is a slight abuse of notation.)  
 
First, we define the notion of infinite pseudo-crescent cones 
generalizing that of pseudo-crescent cones for convenience of proof 
(see Definition \ref{defn:pcrc2}).  We will 
show that  a subset of the  boundary  of a pseudo-crescent-cone  or an 
infinite pseudo-crescent-cone  cannot cover  a compact submanifold  in 
$M$ (see Proposition \ref{prop:decfund}) and some lemmas.   
 
We now assume that there exists a pseudo-crescent-cone or  
an infinite pseudo-crescent cone in $\che M$ and try  
to obtain contradictions. 
We discuss the transversal intersections of two pseudo-crescent-cones. 
We show that infinite  pseudo-crescent-cones are fictitious  as $M$ is 
assumed to be not convex.  This follows from the equivariance property 
of the maximal infinite  pseudo-crescent-cones.   Next, we define  the 
union of collections of  pseudo-crescent-cones and show that the union  
satisfies nice properties  such as  equivariance, and $\tilde M$ is  
a subset of one of such unions.  Next,  we show that the union  has   
triangles at  the  ideal boundary   which can  be linearly ordered.   
Since   the deck-transformation  group  acts  on the  set of linearly  
ordered triangles, if the action  is trivial, then $\pi_1(M)$ acts on  
a  pseudo-crescent-cone $R$  and  hence a component of  $\nu_R 
\cap \tilde M$ covers a closed surface, a contradiction by Proposition   
\ref{prop:decfund}.   If   the action is  nontrivial,  it follows by  
$3$-manifold  topology that $M$ is  a fiber bundle over a circle with   
fiber homeomorphic to a closed  surface, and, again, we obtain  
a contradiction by Proposition \ref{prop:decfund}. 
 
The sides of a compact radiant convex  $3$-ball $P$ consist of radial 
segments,  radiant triangles, and the  unique side $I_K$  in $\SIT$ so 
that $P$  is the cone over  $I_K$.  A  compact radiant convex $3$-ball 
$P$ in $\mathcal{H}$ is said to be an {\em infinite}\/ polyhedron if the 
collection of its   sides not in $\SIT$  is  locally finite except  at 
points of two radial edges of one side  $F$ of $P$ and the collection 
of the   sides  are countably  infinite.   $F$   is called  the   {\em 
fundamental side},  and an edge where  the ideal sides are not locally 
finite a {\em special}\/ edge of the fundamental side.  
 
Again a convex compact ball and its parts in $\che M$ are named by how 
their images are named in $\SI^3$.

\begin{defn}\label{defn:pcrc2}  
A radiant  convex  polyhedron $T$ in $\che  M$  that has only  a  side 
meeting $\tilde M$ and the rest  lying in $\ideal{M}$ is  said to be a 
{\em pseudo-crescent-cone}.   We define $\nu_T$ be  the side that meet 
$\tilde M$ and $\alin_T$ the interior of  the infinitely ideal side of 
$T$,  $\alfi_T$ the  interior  of the disk  that  is the  union of the 
finitely ideal sides of $T$, and $\alfiin_T$ the 
interior  of the  arc that is the   union  of  all segments in   the 
intersection  of  the infinitely   ideal side  and   the union  of the 
finitely ideal sides.  
\end{defn}

\begin{defn}\label{defn:ipcrc} 
A  radiant   convex infinite  polyhedron  $T$  in  $\che  M$  with the 
fundamental  side meeting  $\tilde  M$ and  the rest  of  the sides in 
$\ideal{M}$ is said to  be an {\em infinite pseudo-crescent-cones}. We 
define $\nu_T$ be the fundamental  side and $\alin_T$ the interior  of 
the infinitely ideal side  of $T$  and $\alfi_T$  the interior  of the 
disk   that  is the  union  of the  finitely  ideal  sides of $T$, and 
$\alfiin_T$ the open arc  that is the interior  of the arc that is the 
intersection  of  the  infinitely ideal   side  and the  union  of the 
finitely ideal sides.  
\end{defn}

\begin{figure}[t] 
\centerline{\epsfxsize=3.3in \epsfbox{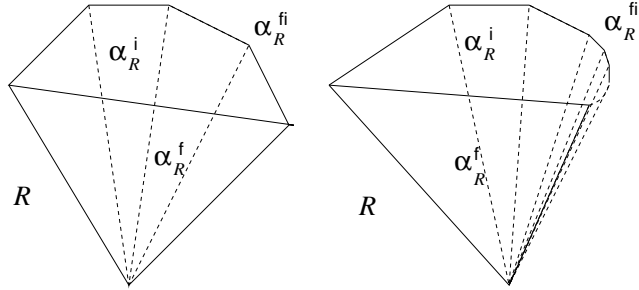}}  
\caption{\label{fig:pscc}  The  pictures  of a pseudo-crescent-cone  and 
an infinite pseudo-crescent-cone.}  
\end{figure} 
\typeout{<<rfpcc.eps>>}  
 
Pseudo-crescent-cones and infinite pseudo-crescent-cones on $\che M_h$  
are defined the same way as for $\che M$, and various parts of them  
are also defined as above. 
 
\begin{lem}\label{lem:pstart} 
If  $\che M_h$ includes  a pseudo-crescent-cone  {\rm (}resp. infinite 
pseudo-crescent-cone{\rm ),} then  $\che M$ includes a pseudo-crescent 
cone {\rm (}resp. infinite pseudo-crescent-cone{\rm ).}  
\end{lem} 
\begin{proof} 
Let $R$  be  a  pseudo-crescent-cone in $\che  M_h$.   Lift the simply 
convex open $3$-ball $R^o$ to $\che M$, and denote the image of the  
lift by $P$. Then since  $P$ is a  simply  convex open   $3$-ball also,  
$\dev|P$ is  an imbedding  onto a simply convex open $3$-ball   
$\dev(P) =\dev(R^o)$, and the  closure $\clo(P)$ of  $P$ in $\che M$  
is a simply convex $3$-ball. Since  the  covering map  
$c:\tilde   M \ra M_h$  is distance non-increasing with  
respect to $\bdd$, it extends to a map  
$\che c: \che M \ra \che M_h$. Since $c| P: P \ra  R^o$ is  
a $\bdd$-isometry, $\che c: \clo(P) \ra R$ is an embedding.  
 
We note that  under $\che c$, $\tilde   M$ maps into $M_h$ obviously.    
Hence  $\che  c^{-1}(\hideal{M})  \subset \ideal{M}$.  We also have  
$\che c^{-1}(M_h) \subset  \tilde M$:  Otherwise a point $p$ of $\ideal{M}$  
maps to a point $p'$ in $M_h$.  We may find a point $q$ in  $\tilde M$   
and a path  $\alpha$ in  $\tilde  M$ connecting $q$ to $p$. Then  
$c \circ  \alpha$ is a path  in $M_h$ with endpoint $p'$ and $q'$  
which lifts to $\alpha$. Since $c$ is a  covering map, this means 
$p \in \tilde M$.  
 
From above, it follows that the inverse image of $R \cap \hideal{M}$ 
under $\che c|\clo(P)$ is a subset  of $\ideal M$ and that of $R \cap 
M_h$ is a subset of $\tilde M$; that is, $\clo(P)$ is a pseudo-crescent.  
When $R$ is an infinite  pseudo-crescent cone, a similar argument applies.  
\end{proof} 
 
The converse of the above lemma is also true, but we won't prove it.

\begin{prop}\label{prop:decfund} 
Let  $F$  be a    component   of  $\nu_R    \cap  \tilde M$     for  a 
pseudo-crescent-cone or  infinite  pseudo-crescent-cone $R$.  Then $p| 
F$ is not a covering map of an imbedded closed surface in $M$.  
\end{prop}  
\begin{proof}  
The component $F$ is an open triangular component of   
$\nu_R \cap \tilde M$.  If $F$ covers a closed surface  $S$ in $M$,   
then there exists a subgroup $G$ of deck transformations acting  
on $F$ so that $p|F: F \ra S$ induces a homeomorphism $F/G \ra  S$.  
Up to choosing an  index  two subgroup, we may assume that   
$G$ preserves the  sides of the submanifold $F$, and $G$ acts on  
$R$ by Propositions \ref{prop:trspcc} and \ref{prop:trsipcc}  
and acts on $\nu_R\cap \tilde M$.

Since $F$ is a radiant  open set, $F$ is a  simply connected disk.  By 
Lemma \ref{lem:sinconv}, $F$ is not the  only component of $\nu_R \cap 
\tilde M$.   Since $G$ has to   act on $F$ and $\nu_R$,   and $F$ is a 
proper  cone in $\nu_R$,  it  follows that  $h(G)$  has at least three 
fixed points on the top side of $\dev(\nu_R)$.  Thus for every element 
$g$ of $G$, $h(g)$ restricts to a  dilatation, i.e., of form $s\Idd$, 
$s \in \bR^+$, in the hyperplane  containing $\dev(\nu_R)$, and so 
$\dev(F)/h(G)$ is not compact. Since the imbedding $\dev|F: F 
\ra \dev(F)$ induces a homeomorphism $F/G \ra \dev(F)/h(G)$, this is a 
contradiction.  
\end{proof}

\begin{lem}\label{lem:sinconv}  
Suppose  that  $N$  is a  real  projective $n$-manifold   and $\che N$ 
includes a convex $n$-ball $B$ such that $B^o  \subset  \tilde  N$ and 
$\Bd B \cap \tilde N$ is  an open {\rm (} connected\/ {\rm )}  
convex subset  $F$ of a side of $B$. Suppose that $F$ covers  
a compact $(n-1)$-dimensional submanifold in $N$.  Then $N$ is convex.  
\end{lem}  
\begin{proof}  
By  lifting to  a  finite cover, we  assume that  $p(F)$  and  $N$ are 
orientable.  
 
Since  $F$ is a  cover of an imbedded closed submanifold $F/G$ in $N$,  
for any   element $\vth$ of $\pi_1(N)$,  we  have either $\vth(F)  = F$  
or $\vth(F)$  and  $F$  are  disjoint, and, moreover, the collection 
$\{\vth(F)| \vth \in \pi_1(N)\}$ is locally finite in $\tilde N$.  
 
If $F$ is a  subset of $\delta \tilde  N$, then clearly $B^o$  
includes $\tilde N^o$. Since $\delta B \cap \tilde M = F$,  
it follows that $\tilde N = B^o \cup F$, which means that  
$\tilde N$ is convex.  
 
Assume now that  $F$ is a subset of  $\tilde N^o$.  Then $F$ separates 
$\tilde N$  into two parts,  one of  which  is $B^o$  and the other is 
$\tilde N - B$.  
 
If  $N - p(F)$ is connected,  then there exists  a simple closed curve 
$\alpha$ in $N$ intersecting with $p(F)$ only once.  This implies that 
the  homotopy class  of $\alpha$ is  not trivial,  and that  there are 
infinitely  many copies  of  $F$ in  $B^o$ obtained   by  applying the 
elements of $\pi_1(N)$ that  are   powers of the deck transformation 
corresponding to $\alpha$.  
 
Suppose that $N -p(F)$ has two components, and there are only finitely 
many copies of $F$ in $B$ under the  deck-transformation-group action. 
Then  we choose a copy  $F'$   so that a component  of   $B - F'$  not 
including $F$  includes no copies of  $F$. Let $B'$,  
$B' = B^{\prime o} \cup F$, denote the closure of this  component.  
Then $F'$  separates $\tilde N$  into a convex set $B^{\prime  o}$ and   
$\tilde N  - B'$.   Since   $B^{\prime o}$  is a component of  
$\tilde N -  p^{-1}(p(F))$, we  have that $B^{\prime  o}$ must cover a  
component of $N - p(F)$, and hence, $B^{\prime o} \cup F$ covers  the  
closure $N_1$   of a  component  of  $N -p(F)$. Each  
deck transformation acting on $B^{\prime o}$ acts  on $F$ nontrivially,  
and a deck transformation acting on   $F$ acts on $B^{\prime  o}$   
since  $N$ and $p(F)$ are orientable and  
$F$ is the only boundary component of $B'$, 
seen as a manifold. Since the  convex open ball $B^{\prime o}$ and 
$F$ are contractible, $N_1$ and $p(F)$ are  homotopy equivalent. If we 
double $N_1$  along $p(F)$,  the  the double  is  a  closed orientable 
$n$-manifold homotopy  equivalent to an  $(n-1)$-manifold $p(F)$. This 
is absurd, and therefore, there are infinitely many copies of $F$ in $B$.  
 
By the conclusion of the above  two paragraphs, given $j$, there exists 
$i$, $i > j$ with a deck transformation $\vth_i$ such that $\vth_i(F)$  
is  inside  the component $C_j$ of $B^o  - \vth_j(F)$ not containing $F$.  
Let $d_N$ denote the original Riemannian metric on $\tilde N$ induced  
from one  of $N$. Then  for each point $x$ of  $F$, 
$d_N(\vth_i(x), F) \ra \infty$ since such a path from $x$ to $F$ has  
to pass through an increasing  number of $\vth_j(F)$ as $i \ra \infty$,   
and  there exists a lower bound to the infimums of $d_N(y, z)$ if $y$  
and $z$ are in two different components of  $p^{-1}(F)$. (The  latter  
statement follows since any arc in $N$ with endpoints in $p(F)$ that is not  
homotopic to a path in $p(F)$ must have $d_N$-length longer than some  
uniform positive constant $\eps$.)  
 
Let $K$ be any compact subset of $\tilde N$ intersecting $F$.  Then it 
follows that $\vth_i(K)$ is eventually in  $B^o$.  From this we easily 
see   that any    two   points  $y$  and   $z$  can   be connected  by 
$\vth_i^{-1}(s)$   for some $i$ and  a segment $s$  in $B^o$ of  
$\bdd$-length $\leq \pi$ since $B^o$ is convex. Thus, $\tilde N$ is  
convex, and so is $N$.  
\end{proof}

\begin{figure}[b] 
\centerline{\epsfxsize=2.2in \epsfbox{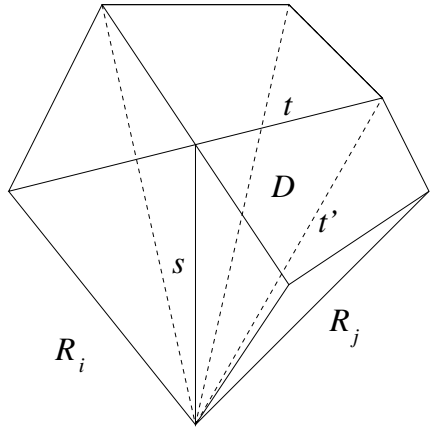}}  
\caption{\label{fig:trspcc}       The transversal  intersection     of 
pseudo-crescent-cones.}  
\end{figure} 
\typeout{<<rftpc.eps>>}  
 
\begin{defn}\label{defn:trspcc}  
Let  $R_1$  and $R_2$  be both pseudo-crescent-cones  or both infinite 
pseudo-crescent-cones.   We say    that   $R_1$ and  $R_2$   intersect 
transversally if the following statements  hold ($i=1, j=2$, or  $i=2, 
j=1$):  
\begin{enumerate}  
\item $R_1$ and $R_2$ overlap.  
\item $\nu_{R_1} \cap \nu_{R_2}$  is a well-positioned  radial segment 
$s$ in $\nu_{R_1}$ and similarly in $\nu_{R_2}$.  
\item $\nu_{R_i}   \cap  R_j$ equals the   closure  of a  component of 
$\nu_{R_i} - s$. Hence it is  a well-positioned convex triangle $D$ in 
$R_j$. Two  sides of  $D$  other  than  $s$   are finitely  ideal  and 
infinitely   ideal respectively.   The   three    sides of  $D$    are 
well-positioned in three sides of $R_j$   respectively; or a  finitely 
ideal side of $D$ equals the common side of two finitely ideal sides of 
$R_j$  and other two sides are well  positioned in $\nu_{R_j}$ and  
an infinitely ideal sides of $R_j$ respectively.  
\item $R_i \cap R_j$ is the closure of a component of $R_i - D$.  
\item Let $t$ be the infinitely ideal  side of $D$.  $\alin_{R_i} \cap 
\alin_{R_j}$ is  equal to one   of  two components  of $\alin_{R_i}  - 
t$. Hence, $\alin_{R_i} \cup \alin_{R_j}$ is  an open $2$-disk, i.e., 
$\alin_{R_i}$ and $\alin_{R_j}$ extend each other.  
\item Let $t'$ be the finitely ideal side of $D$.  $\alfi_{R_i} \cap 
\alfi_{R_j}$ is equal   to one of   two components of  $\alfi_{R_i}  - 
t'$. Hence, $\alfi_{R_i} \cup \alfi_{R_j}$ is an open $2$-disk.  
\item  $t  \cap  \alfiin_{R_i}$   is  a point.   $\alfiin_{R_i}   \cap 
\alfiin_{R_j}$  equals one of  two components of  $\alfiin_{R_i} - t$, 
and  $\alfiin_{R_i} \cup \alfiin_{R_j}$ is  an open  arc. 
\end{enumerate} 
\end{defn} 
 
\begin{prop}\label{prop:trspcc}  
Let $R_1$ and  $R_2$ are  both pseudo-crescent-cones overlapping  with 
each other. Then either $R_1 \subset R_2$ or $R_2 \subset R_1$ holds or 
$R_1$ and $R_2$ intersect transversally.  
\end{prop} 
\begin{proof}  
The proof is similar to  that of  Theorem \ref{thm:transversal}.    
Only technical differences exist but we have to use the following 
lemma \ref{lem:522} instead of Lemma 5.2 in \cite{psconv}. 
\end{proof}

\begin{lem}\label{lem:522}  
Let $N$ be a closed real projective manifold with a developing map 
$\dev:\che N_h \ra \SI^3$ for the completion $\che N_h$ of 
the holonomy cover $N_h$ of $N$. Suppose that $\dev$ is  
an imbedding onto a union of two {\rm (}infinite or otherwise\/{\rm )}  
radiant convex polyhedra $H_1$ and $H_2$ 
meeting in an {\rm (}infinite or otherwise\/{\rm )} 
radiant polyhedron. Then $H_1 \cup H_2$ is  
convex and so is $N$. 
\end{lem} 
\begin{proof} Suppose that $H_1$ and $H_2$ are finite. 
If $H_1 \cup H_2$ is not convex, then there exist 
sides of $H_1$ and $H_2$ meeting in an edge so that they  
can be extend into $H_1^o$ and $H_2^o$ as totally geodesic surfaces  
respectively. We find a holonomy-group-invariant codimension-one  
submanifold in $H_1^o \cup H_2^o$ from these. The rest of the proof  
is similar to that of Lemma 5.2 in \cite{psconv}. 
When at least one of $H_1$ and $H_2$ is infinite, a similar 
argument will work using support planes. 
\end{proof} 
 
\begin{prop}\label{prop:trsipcc} 
Let $R_1$ and $R_2$ be both infinite pseudo-crescent-cones overlapping 
with each other. Then either $R_1 \subset R_2$ or $R_2 \subset R_1$.  
The special edges of $R_1$ and $R_2$ must coincide in this case. 
\end{prop} 
\begin{proof} 
If the conclusion does not hold, then $R_1$ and $R_2$ must intersect 
transversally as the proof for Proposition \ref{prop:trspcc} also works  
here. Since $M$ is orientable, this means that a special edge 
$e$ of the fundamental side  of $R_1$ lies  in $\mathcal{L}_2$ the union  
of finitely ideal sides of  $R_2$ with  the special edge of  $R_2$ 
removed or conversely for $R_1$ and $R_2$. Assume the former case. 
The finitely ideal sides  of $R_1$ near the  special edge $e$ must lie 
in $\mathcal{L}_2$ from  the definition  of transversality.   However, the 
sides of $R_2$ are locally finite there while those  of  $R_1$ are 
not. This is a contradiction. The last statement follows obviously 
from the proof. 
\end{proof}

\begin{prop}\label{prop:asc} 
Let $R_i$  be   an  ascending  sequence  of  pseudo-crescent-cones  or 
infinite pseudo-crescent-cones\/{\rm ;} that is,  $R_1 \subset R_2  
\subset R_3 \subset \dots$. Then there exists a  pseudo-crescent-cone  
or an infinite pseudo-crescent-cone that includes  the  union of $R_i$  
as a dense subset.  
\end{prop}  
\begin{proof}  
By Theorem \ref{thm:seqconv}, there exists a $3$-ball $R$ including 
every $R_i$.  The fact that 
$\clo(\alfi_{R_i}   \cup \alin_{R_i}  \cup   \alfiin_{R_i})$  form  an 
increasing ideal sequence shows that a part of the boundary of $R$ 
is a subset of the ideal set. If $\che M$ is a subset of $R$, 
then $\tilde M$ equals $R^o$ and $M$ must be convex. Therefore, 
a side of the boundary of $R$ is a subset of $M$, and $R$ is
the desired object.  
\end{proof} 
 
We will now show  that infinite pseudo-crescent-cones  do not exist on 
$\che M$. Suppose  that $R$ is  an infinite  pseudo-crescent-cone.  We 
can introduce  an equivalence relation  on the collection $\mathcal{R}_2$ 
of    all    infinite      pseudo-crescent-cones:    We    say     two 
infinite pseudo-crescent-cones are equivalent if they overlap.   
The equivalence relation is generated  by this; i.e.,  $R \sim S$  
if there   
exists a finite sequence $R^1,  R^2, \dots, R^n$ in $\mathcal{R}_2$  
such that $R^1  = R$ and $R^n = S$  where $R^i$ and $R^{i+1}$  
overlap for each $i=1,\dots, n-1$.  
 
Propositions \ref{prop:trsipcc} and \ref{prop:asc} show that the 
equivalence class   of overlapping  infinite  pseudo-crescent-cones is 
totally ordered  and has  a unique  maximal  element.  Let $R$  be the 
maximal   element. For a deck transformation  $\vth$, either $R$ and 
$\vth(R)$ do not overlap or $R \subset \vth(R)$ or $\vth(R) \subset R$ 
by Proposition \ref{prop:trsipcc}.  If $R  \subset \vth(R)$, then since $R$ 
is maximal, we  have $R =  \vth(R)$; if $\vth(R)  \subset  R$, then we 
have $R = \vth(R)$.  If $R$ and $\vth(R)$ do not overlap but meet each 
other, then  $R \cap \vth(R)$   is the union  of  common components of 
$\nu_R \cap \tilde M$ and $\nu_{\vth(R)} \cap \tilde M$, which follows 
as  in Chapter 7 of \cite{psconv}.   As in Chapter 10 of \cite{psconv},  
we obtain the  so-called  two-faced  submanifolds.   
If a  two-faced  submanifold exists, then a  component of  
$\nu_R  \cap \tilde  M$ covers a closed surface in $M$,  
i.e., the two-faced submanifold. This was ruled  out by  
Proposition \ref{prop:decfund}, and we have either $R = \vth(R)$   
or  $R \cap \vth(R)  \cap \tilde M = \emp$.  (The detail of a similar  
argument was in Chapters \ref{ch:tdim} and \ref{ch:obtain} and we  
will see it again in this chapter.)  
 
By Proposition \ref{prop:lfin3}, the collection consisting of elements 
of form $\vth(R)$ is locally  finite in $\tilde M$.  The  equivariance 
and this show by Proposition \ref{prop:lfin} that $R$ covers a compact 
submanifold  in  $M$. Again  this  means  that a 
component of $\nu_R \cap \tilde M$ covers  a closed surface in $M$, 
a  contradiction.  Therefore, we    proved that there  exist  no 
infinite pseudo-crescent-cones:  
 
We restate our result as follows:  
\begin{prop}\label{prop:nipcc} 
Let $N$ be  a compact  radiant affine manifold  with  empty or totally 
geodesic  boundary.  If  $\che  N$  includes  infinite pseudo-crescent 
cones, then $N$ is convex.  
\end{prop} 
 
We will now show that pseudo-crescent-cones do not exist as well using  
a longer more involved argument. Let $\mathcal{R}_3$ be a collection of all   
pseudo-crescent-cones in $\che  M$. We define a relation that $R \sim S$  
for $R, S \in \mathcal{R}_3$ if they overlap, and define an  equivalence  
relation on $\mathcal{R}_3$ generated by this as usual.  
 
We define the  following  sets as  in  \cite{cdcr1}:  
\begin{eqnarray}\label{eqn:Lambp}  
\Lamp(R)  =     \bigcup\limits_{S   \sim   R}   S,     
&\dinLamp(R)   = \bigcup\limits_{S  \sim  R} \alin_S,  
&\dfLamp(R)  =  \bigcup\limits_{S \sim R} \alfi_S, 
\\  \dfiLamp(R) = \bigcup\limits_{S \sim R} \alfiin_S, 
&\Lamp_1(R) = \bigcup\limits_{S \sim R} S - \nu_S  
&\hbox{ for } R \in \mathcal{R}_3. \nonumber  
\end{eqnarray} 
 
\begin{figure}[b] 
\centerline{\epsfxsize=2.5in \epsfbox{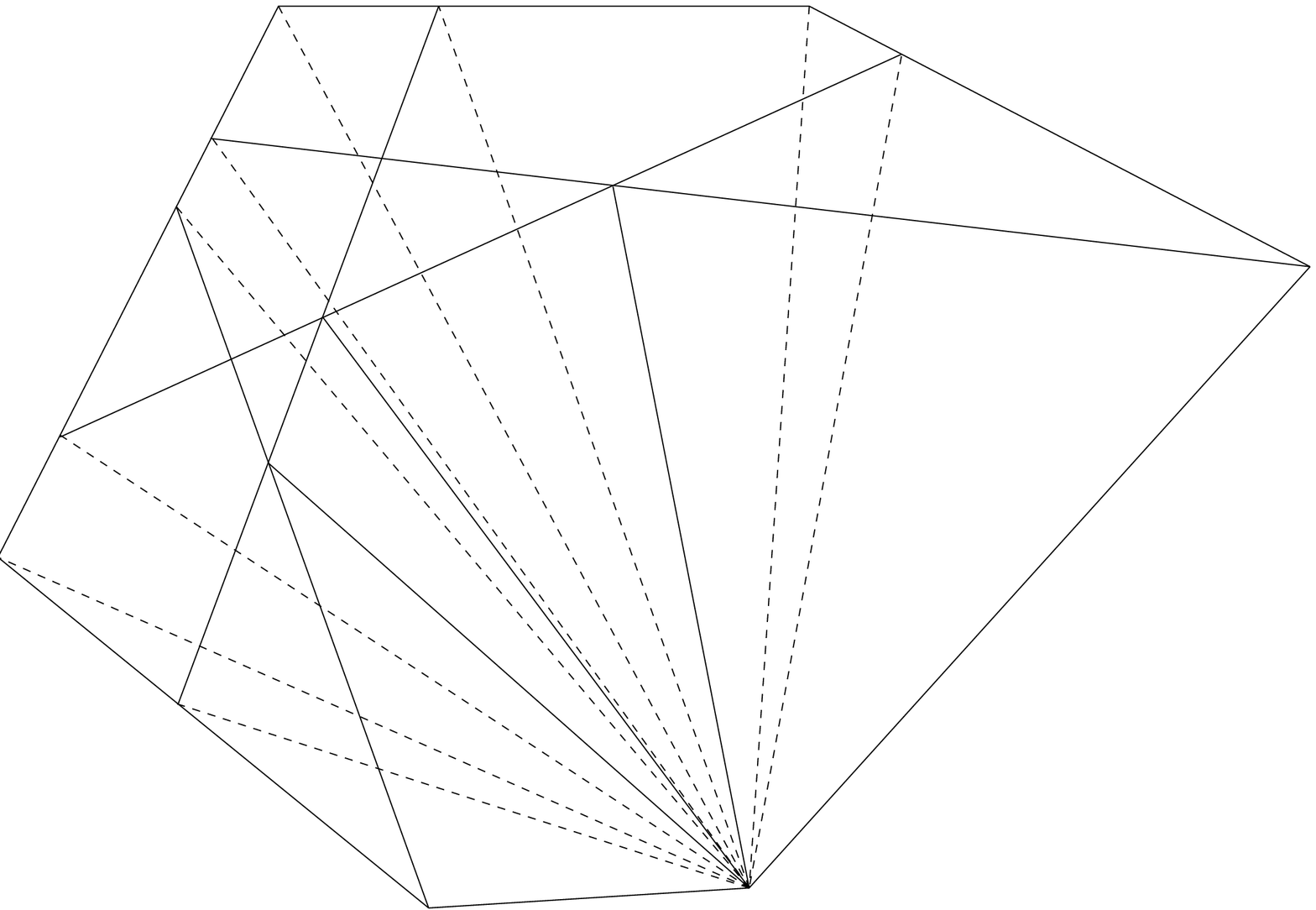}}  
\caption{\label{fig:Lamp} A picture of $\Lamp(R)$.}  
\end{figure} 
\typeout{<<rfLamp.eps>>}  
 
We begin by introducing the  equivariance properties of $\Lamp(R)$. We 
list properties of $\Lamp(R)$.  
\begin{enumerate}  
\item    $\Lamp(R)$ and  $\Lamp_1(R)$  are 
path-connected.  
\item $\Lamp(R) \cap  \tilde M$ is a closed  radiant subset of $\tilde 
M$ with totally geodesic boundary, which is also a radiant set.  
\item $\Lamp_1(R)$ is  a $3$-manifold with  boundary (and corners)  
$\delta_\infty \Lamp(R)$, defined as $\dinLamp(R) \cup 
\dfLamp(R) \cup \dfiLamp(R)$. The $3$-manifold has a real projective  
structure induced from the pseudo-crescent-cones. 
\item $\dinLamp(R)$ is an open surface imbedded in $\Iideal{M}$.  
\item $\dfLamp(R)$ is an open surface imbedded in $\ideal{M} - \Iideal{M}$.  
\item  $\dfiLamp(R)$ is an   arc   in  $\Iideal{M}$, and   the   union 
$\delta_\infty \Lamp(R)$  is  an  open surface  
imbedded in $\ideal{M}$.  
\item  $\Lamp(R)$ is {\em maximal}, i.e., for any triangle $T$ in  
$\tilde M$ with sides $s_1,s_2$, and $s_3$,   
if $s_2$ and $s_3$  are subsets of $\Lamp(R)$,   
then $s_1$ is a subset of $\Lamp(R)$.  
\item $\Bd  \Lamp(R) \cap \tilde  M$ is a properly  imbedded countable 
union of disks in $\tilde M^o$, each of which is a totally geodesic,  
radiant, properly imbedded open triangles.  
\item For any  deck transformation $\vth$, we  have  
\begin{eqnarray*}
&\Lamp(\vth(R)) = \vth(\Lamp(R)), &\dinLamp(\vth(R)) =\vth(\dinLamp(R)),\\ 
&\dfLamp(\vth(R)) =\vth(\dfLamp(R)), &\dfiLamp(\vth(R)) = 
\vth(\dfiLamp(R)).
\end{eqnarray*}  
\end{enumerate}  
The  proof    of these  facts   are similar    to  what is   given  in 
\cite{cdcr1} for dimension two and Chapter 7 of \cite{psconv}.   
Note  that for the closedness  in (2),  
we use the fact that a sequence  of pseudo-crescent-cones with   
a common open set in them give rise to  a pseudo-crescent-cone   or     
infinite pseudo-crescent-cone  by Theorem \ref{thm:seqconv}.  
Since infinite pseudo-crescent  cones do not exist, 
we obtain the closedness in (2).  
 
\begin{prop}\label{prop:splitting2} 
Given  two pseudo-crescent-cones $R$ and $S$, we have one of three 
possibilities\/{\rm :}  
\begin{itemize}  
\item $\Lamp(R) = \Lamp(S)$,  
\item $\Lamp(R)$ meets $\Lamp(S)$ at the union of common components of 
$\Bd \Lamp(R) \cap \tilde M$ and $\Bd \Lamp(S) \cap \tilde M$.  
\item  $\Lamp(R)  \cap  \tilde  M$ and   $\Lamp(S) \cap \tilde  M$ are 
disjoint.  
\end{itemize}  
\end{prop}  
\begin{proof}  
The proof is same as in Chapter 6 or 7 of \cite{psconv}.  
The role played by $n$-crescent is played by pseudo-crescent-cones.  
\end{proof}

\begin{prop}\label{prop:corner}  
$\dfLamp(R) \cup \dfiLamp(R) \cup\{O\}$ is  a  union of countable   
collection of triangles $T_i$. Each $T_i$ is a side of  
a pseudo-crescent-cone $S$, $S \sim R$ in the closure of $\alfi_S$. 
Precisely one of the following holds\/{\rm :}  
\begin{enumerate} 
\item The collection equals $\{T_0, T_1, \dots, T_{n-1}\}$ for $n \geq 
2$ and $T_i$ and $T_j$ meet in an edge  if $i-j = 1,  n-1$ mod $n$ and 
$T_i$ and $T_j$ are disjoint if otherwise and $i \ne j$. {\rm (}Indices are 
cyclic{\rm .)}  
\item The collection equals $\{T_0, T_1, \dots, T_{n-1}, \dots\}$ for  
$n \geq 2$ and $T_i$ and $T_j$ meet in an edge if $i-j =  \pm 1$  
and $T_i$ and $T_j$ are disjoint if otherwise and $i \ne j$.  
{\rm (} The collection could be finite or infinite {\rm ).} 
\item The collection equals $\{T_i| i \in \bZ \}$  and $T_i$ and $T_j$ 
meet in an edge if $i-j  = \pm 1$ and $T_i$  and $T_j$ are disjoint if 
otherwise and $i\ne j$.  
\end{enumerate} 
\end{prop}  
\begin{proof}  
This follows from Proposition \ref{prop:trspcc}.  
\end{proof} 
 
Let $R$  be a  pseudo-crescent-cone.   We orient each  radial segment 
toward the origin $O$.  The collection of  finitely ideal sides of $R$ 
in   $\fideal{M}$  can    be linearly  ordered   using   the boundary 
orientation induced from that of $R$, which in  turn was obtained from 
$\tilde M$: We give the orientation on $\tilde M$ and the ordering  
on $\{T_i\}$ so  that $T_{i+1}$ is to the right of $T_i$ for  
each $i=0, 1, \dots$.  
 
The {\em  right-most}\/ ideal side of  a pseudo-crescent-cone $R$ is the 
finitely ideal side of $R$ which  is right-most among the finitely 
ideal sides of $R$.  
 
\begin{lem}\label{lem:right}  
Each $T_i$ is the right-most  ideal side of the unique  
pseudo-crescent-cone maximal among pseudo-crescents having $T_i$ as  
the right-most ideal side. {\rm (}In the  second case above we  have  
to assume $i \geq 2.$\/{\rm )}  
\end{lem} 
\begin{proof} 
Clearly,  there exists  at least   one pseudo-crescent-cone  with  the 
right-most  side $T_i$ since we can choose such a pseudo-crescent-cone 
in a pseudo-crescent-cone $S$, $S \sim R$, including $T_i$ as  
a side in the closure of $\alfi_S$. 
Given  any two such  pseudo-crescent cones $R$ 
and $S$, since they overlap, we have either $R \subset S$ or  
$S \subset R$  or $R$ and $S$ intersect transversally by Proposition 
\ref{prop:trspcc}.  If $R$ and  $S$ intersect  transversally, we can see   
easily that  $R$ and $S$ cannot share the right-most ideal side. Hence,  
we have  $R \subset S$ or $S \subset R$.  
 
Also, given any ascending sequence $R_{i, j}$ of pseudo-crescent-cones 
such that  $T_i$ is the right-most  ideal side, we see that $\bigcup_j 
R_{i,  j}$  is    contained   in a pseudo-crescent-cone by Proposition 
\ref{prop:asc},  which has $T_i$ as   the right-most  ideal side.   The 
uniqueness of  the  maximal element follows  easily  from the previous 
paragraph.  
\end{proof}

We rule out the first possibility for the triangles.  
\begin{prop}\label{prop:infcorn}  
The collection $\{T_i\}$ is linearly ordered.  
\end{prop} 
\begin{proof} 
Suppose   not.     Then  we obtain by Lemma \ref{lem:right}  
a sequence  of pseudo-crescent-cones $P_i$ such that   
$P_i$ and $P_{i+1}$ overlap for $i= 0, \dots,  n-1$  mod  $n$  
where the   indices  are cyclic.  Orient 
$\nu_{P_i}$  by  the outer-normal to   $P_i$ for each $i$.  $\nu_{P_i} 
\cap P_{i+1}$ is  the closure of  the  left component of  $\nu_{P_i} - 
\nu_{P_{i+1}}$ and $\nu_{P_i} \cap P_{i-1}$ the closure of the right  
component of $\nu_{P_i} - \nu_{P_{i-1}}$. We see that  
$\nu_{P_i} \cap P_{i-1}^o$ is contained the right-most component of   
$\nu_{P_i}^o - \fideal{M}$ and $\nu_{P_i}\cap P_{i+1}^o$ is contained    
in    the  left-most component  of $\nu_{P_i}^o - \fideal{M}$.  
Hence, we  may choose a  closed simple arc $\alpha$ in $\tilde  M$  
so that $\alpha   \cap P_i$ is  a  connected arc with  two endpoints  
respectively in $\nu_{P_i} \cap P_{i-1}^o$ and $\nu_{P_i}\cap P_{i+1}^o$ 
for each $i$.   
 
Suppose that $\nu_{P_i}^o \cap \fideal{M}$  is  empty.  
Assume first that $\nu_{P_i}^o - P_{i-1}^o - P_{i+1}^o$ is not empty. 
Then there exists a segment $s$ in $\nu_{P_i}^o - P_{i-1}^o -  
P_{i+1}^o$.  The segment has a compact neighborhood  $N$ in   
$\tilde M$ which includes another segment $s'$ in $\tilde M$ outside $P_i$   
ending at $p' \in \nu_{P_{i-1}}$ and $q' \in \nu_{P_{i+1}}$  
transversally. We extend $s'$ in 
$P_{i-1}$ and $P_{i+1}$ to  a segment $s''$ with endpoints in  
$p'' \in \alfi_{P_{i-1}}$ and $q'' \in \alfi_{P_{i+1}}$ respectively 
where we assume that points occur on $s''$ in order
$p'', p', q', q''$. We are additionally required to choose $s''$ so that  
$q'' \in \delta \nu_{P_i} \cap T_i$. We may also assume that  
$s$ and $s'$ are the sides of a totally  geodesic convex disk  
$D$ with four edges where two other edges are in $\nu_{P_{i-1}}$  
and $\nu_{P_{i+1}}$ respectively.  
Then the union of all radial segments passing through $D$ is a convex  
$3$-ball $K$ bounded by five sides. The union of $K$, $P_{i-1}$,  
$P_i$, and $P_{i+1}$ contains a pseudo-crescent cone containing  
$P_i  \cup K$ properly and having $T_i$ as its right-most ideal side.  
To see this, we   need to   make  $K$ slightly larger and    
apply Proposition \ref{prop:overlap}  
as in the proof of Lemma \ref{lem:twocresc}.   
This is a contradiction  since $P_i$  is maximal.

If $\nu_{P_i}^o - P_{i-1}^o - P_{i+1}^o$ is empty, then  
$\nu_{P_i}^o$ is a subset of $P_{i-1}^o \cup P_{i+1}^o$.  
We clearly see that there exists a radiant triangle  
in $P_{i-1} \cup P_{i+1}$ near $\nu_{P_i}$ outside $P_i$ which  
bounds a pseudo-crescent cone $P'_i$ including $P_i$. We can  
choose $P'_i$ so that $T_i$ is the right-most ideal side of  
$P'_i$. This is a contradiction. Therefore,  
$\nu_{P_i}^o \cap \fideal{M}$ is not empty.

We may now choose a properly  imbedded radiant open triangle  
$\triangle_i$ in $P_i^o$ such that $\triangle_i$  meets  
$\alpha$ at a unique point (up to isotopying $\alpha$ in $P^o$),  
and $\clo(\triangle_i) - \triangle_i$ is  the  union of   
a radial segment in  $\nu_{P_i}  \cap \hideal{M}$,    
a radial segment in $\clo(\alfi_{P_i})$, and  a segment   
in  $\clo(\alin_{P_i})$ and  is a subset of $\ideal{M}$.    
Since   $\tilde M$ is  simply    connected, $\alpha$ must   
bound  a disk $D$.   Perturbing  $D$ and  $\alpha$ into 
general positions with respect to $\tri_i$ makes  $D \cap \tri_i$ to a 
properly  imbedded arcs  in $D$ with  endpoints  in $\delta D$ since 
$\tri_i$ is properly imbedded; thus, $\alpha$ intersects with $\tri_i$  
even  number of times by the  classification of $1$-manifolds,  
a contradiction.  
\end{proof}

\begin{claim} 
$\tilde M$ is a subset of $\Lamp(R)$, and $M$ is  
boundaryless. 
\end{claim}

\begin{proof} 
Suppose  that $\Lamp(R) \cap  \tilde M$ is a  proper subset of $\tilde 
M$; i.e., $\Bd \Lamp(R) \cap \tilde M$ is not  empty.  For any pair of 
pseudo-crescent-cones $R$ and $S$,  one of the following possibilities 
hold from Proposition \ref{prop:splitting2}:  
\begin{itemize} 
\item $\Lamp(R) = \Lamp(S)$,  
\item $\Lamp(R) \cap \Lamp(S) \cap  \tilde M  = \emp$ where   
$R\not\sim S$, or  
\item $\Lamp(R) \cap \tilde M$ and $\Lamp(S)\cap \tilde M$ meet at the 
union of common components  of $\Bd \Lamp(R) \cap  \tilde M$ and  $\Bd 
\Lamp(S) \cap   \tilde M$ where $R \not\sim  S$. As in \cite{psconv}, 
this set is also the union of common  components of $\nu_T \cap \tilde 
M$ for a pseudo-crescent-cone $T$,  $T \sim R$  and $\nu_U \cap \tilde 
M$ for $U \sim S$.  
\end{itemize}  
The  third   possibility  gives  rise  to   a   closed submanifold  of 
codimension $1$    as in  \cite{psconv}  or  in  above  Chapters 
\ref{ch:tdim} and \ref{ch:obtain},   analogously to  the two-faced 
manifolds.  This  means  that a  component of  $\nu_T  \cap \tilde  M$ 
covers a compact submanifold for a pseudo-crescent-cone $T$,  
which  is a contradiction by Proposition 
\ref{prop:decfund}, and the third possibility does not occur.  
 
Since for  every  deck transformation  $\vth$,  $\vth(\Lamp(R))   = 
\Lamp(\vth(R))$, and there is no  two-faced submanifold, we have that 
$\Lamp(R)   = \vth(\Lamp(R))$ or  $\Lamp(R)   \cap \vth(\Lamp(R)) \cap 
\tilde M = \emp$.  
 
The  collection consisting  of elements   of form $\vth(\Lamp(R))$ is 
locally  finite  since this can be proved by the  proof  of Proposition 
\ref{prop:lfin2} replacing crescent-cones by pseudo-crescent-cones.    
By Proposition  \ref{prop:lfin}, 
we have that $\Lamp(R) \cap \tilde M$ covers a compact submanifold $N$ 
in $M$.  Recall that $\Bd \Lamp(R)\cap \tilde  M$ is totally geodesic.   
Let  $F$ be  a component of $\Bd \Lamp(R)$. Then a point $x$  of $F$  
belongs to $\nu_T$  where $T$ is a radiant bihedron  equivalent to $R$.   
Since $\nu_T \cap \tilde M$ and $F$ must be tangent there,  and  
$F$  is  properly imbedded, it follows that  $F$ is a component of   
$\nu_T \cap \tilde M$.  Since $F$ covers a component of the boundary  
of $N$, $F$ covers a closed surface in $N$, which was ruled out by   
Proposition \ref{prop:decfund}, and our claim that $\tilde M$ is  
a subset of $\Lamp(R)$ is proved. 
 
Suppose that $\delta M$ is not empty. 
There exists a component $K$ of $\delta \tilde M$. Let  $x$ be  
a point of $K$. Then $x$ belongs to $\nu_T \cap \tilde M$ for some  
$T  \sim R$.  Let $K'$ be the component of $\nu_T  \cap \tilde M$  
containing $x$.  Then since $K$ and  $K'$ must be tangent at $x$,  
and both are properly imbedded, we have  $K = K'$. However,  
Proposition \ref{prop:decfund} shows that this is a contradiction. 
\end{proof} 
 
As before, $\vth(\Lamp(R)) = \Lamp(\vth(R))$ for 
a deck transformation $\vth$ includes $\tilde M$ and hence  
$\vth(\Lamp(R))$ meets $\Lamp(R)^o$.  
Hence \begin{eqnarray*} 
\vth(\Lamp(R)) = \Lamp(R), &\vth(\dfLamp(R)) = \dfLamp(R), \\ 
\vth(\dinLamp(R)) = \dinLamp(R),& \vth(\dfiLamp(R)) = \dfiLamp(R) 
\end{eqnarray*} 
hold for every deck transformation  $\vth$. We  have that     
the deck transformation   group $\pi_1(M)$ of $\tilde M$ acts  on   
the discrete ordered set  $\{T_i\}$ either  preserving the order  or   
reversing the  order since  adjacent triangles have to go adjacent ones.   
We have an exact sequence:  
\[ 1 \ra K \ra \pi_1(M) \ra Q \ra 1 \] 
where $Q$  is the group   of automorphisms  of linearly  ordered   set 
$\{T_i\}$ and $K$ is the kernel. Elements of $Q$ may preserve  
the order  or reverse it.  We assume  without loss of generality   
that  $Q$  preserves the order  (by choosing a double cover of $M$  
if necessary).  
 
Suppose that   $Q$ is a    trivial group;  that is, $\pi_1(M)$   acts 
trivially on $\{T_i\}$. Then $\pi_1(M)$ acts on each of $T_i$. For  
$i \geq 2$, choose a maximal pseudo-crescent-cone $R$ which contains  
$T_i$ as the right-most ideal side. Then for each deck transformation  
$\vth$, $\vth(R)$ is the maximal such pseudo-crescent-cone, and hence  
$R = \vth(R)$. 
Since $\pi_1(M)$ acts on $R$, the submanifold $R\cap \tilde M$  
covers  a compact submanifold of  $M$ with boundary union of closed 
surfaces by Proposition \ref{prop:lfin}. This shows  that a component   
$\nu_R \cap \tilde M$ covers a closed surface in $M$, a contradiction  
to Proposition \ref{prop:decfund}.  
 
Suppose now that the action is not trivial so that $Q$ is isomorphic to 
$\bZ$. By Theorem 11.6 in Hempel  \cite{Hemp}, we see that  
the Poincar\'e associate of $M$ is a fiber bundle over a circle,  
and $K$ is nothing  but the fundamental group  of  the fiber 
which must be a closed surface as $M$ is a closed $3$-manifold. 
Hence, the  subgroup $K$ of $\pi_1(M)$ isomorphic    
to the fundamental group of a surface acts trivially  on $\{T_i\}$.   
As above choose a  maximal  pseudo-crescent-cone $R_i$ which  contains  
$T_i$ as the right-most ideal side. Then $K$ acts on  
$\nu_{R} \cap \tilde M$, a union of disjoint open sets.  
 
Since   the  action of $\pi_1(M)$ on the collection $T_i$  
is  not trivial, the set of indices of $T_i$  
equals $\bZ$; so we choose a fixed $i$. $R_{i+1}$ exists, and  
since $R_i$ and $R_{i+1}$ are maximal and $R_i$ does not 
include $T_{i+1}$, either $R_i \subset R_{i+1}$  
or $R_i$  and  $R_{i+1}$ intersect transversally.  
In the   first case, since  $R_i$ is a proper  subset  of $R_{i+1}$,   
it follows from looking at their images under $\dev$ that  
$\nu_{R_i}^o  \subset R_{i+1}^o$. In  the   second   case, $R_i$ and      
$R_{i+1}$ intersect transversally,   and $\nu_{R_i}  \cap  R_{i+1}^o$    
is a radiant  open triangle adjacent to $T_i$. In both cases, let $E$  
denote the open triangle $\nu_{R_i} \cap R_{i+1}^o$.  
 
Since $K$ acts on  $T_i$, $K$ acts  on the component $L$ of $\nu_{R_i} 
\cap  \tilde M$ including $E$.  Since  $K$   is isomorphic  to  the 
fundamental group of a closed surface $F$ and $L$ is homeomorphic to an  
open disk, it follows that $L/K$ is a closed surface homeomorphic to  
$F$, contradicting Proposition \ref{prop:decfund} as above. Therefore,   
we completed to show that there exists no pseudo-crescent in $\che{M}$   
and hence in $\che{M_h}$ by Lemma \ref{lem:pstart}.

\begin{appendix} 
 
\chapter{Dipping intersections} 
\label{app:dipping} 
 
We will  gather   needed  facts  about  the  dipping   intersection of 
$n$-balls and   the transversal intersection  of  $n$-crescents  here, 
which were already stated and proved in \cite{cdcr1} and \cite{psconv}.  
 
Let  $D$ be a  convex  $n$-ball in $\che M_h$   such that $\delta D$ 
includes a totally geodesic convex $(n-1)$-ball $\alpha$.  We say that 
a  convex  $n$-ball $F$ is {\em   dipped into}\/ $(D,  \alpha)$ if the 
following statements hold:  
\begin{itemize} 
\item $D$ and $F$ overlap.  
\item $F \cap \alpha$  is a tame  $(n-1)$-ball $\beta$ with  $\delta 
\beta \subset \delta F$ and $\beta^o\subset F^o$.  
\item $F - \beta$ has two convex components $O_1$  and $O_2$ such that 
$\clo(O_1) = O_1 \cup \beta = F - O_2$ and $\clo(O_2) = O_2 \cup \beta 
= F - O_1$.  
\item $F \cap D$ is equal to $\clo(O_1)$ or $\clo(O_2)$.  
\end{itemize} 
We say that $F$ is  {\em dipped  into  $(D, \alpha)$ nicely}\/ if  the 
following statements hold:  
\begin{itemize} 
\item $F$ is dipped into $(D, \alpha)$.  
\item $F\cap D^o$ is identical with $O_1$ and $O_2$.  
\item $\delta(F\cap  D)  =  \beta \cup    \xi$ for  a  topological 
$(n-1)$-ball $\xi$ in the topological boundary $\Bd F$ of $F$ in $\che 
M$ where $\beta \cap \xi = \delta \beta$.  
\end{itemize} 
As a consequence, we have $\delta \beta \subset \Bd F$.

The direct generalization of Corollary 1.9 of \cite{cdcr1} gives:  
\begin{cor}\label{cor:dipping} Suppose that $F$ and $D$ overlap{\rm ,}  
and $F^o \cap  (\delta D - \alpha^o) =  \emp$.  Assume the following 
two equivalent conditions\/{\rm :}  
\begin{itemize} 
\item $F^o \cap \alpha \ne \emp$.  
\item $F \not\subset D$.  
\end{itemize} 
Then  $F$  is dipped  into  $(D, \alpha)$.   If $F \cap  (\delta D - 
\alpha^o) = \emp$  furthermore{\rm  ,} then $F$   is dipped into  $(D, 
\alpha)$ nicely.  
\end{cor} 
 
\begin{figure}[t] 
\centerline{\epsfxsize=3in \epsfbox{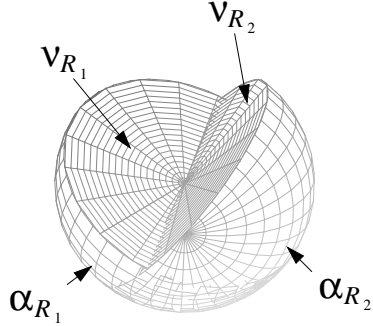}}  
\caption{\label{fig:tr}   A transversal intersection of $3$-crescents 
(stereographically projected to $\bR^3$).}  
\end{figure} 
 
Suppose that $R_1$ is  an $n$-crescent that  is an  $n$-bihedron.  Let 
$R_2$ be another  bihedral $n$-crescents with sets  $\alpha_{R_2}$ and 
$\nu_{R_2}$.      We    say that   $R_1$    and  $R_2$  intersect {\em 
transversally}\/ if the following conditions hold ($i=1, j= 2$; or $i= 
2, j= 1$):  
\begin{enumerate} 
\item $\nu_{R_1}  \cap \nu_{R_2}$ is an $(n-2)$-dimensional hemisphere 
if $n>2$ and consists of a single point $x$ if $n=2$.  
\item For the  intersection $\nu_{R_1}\cap \nu_{R_2}$  denoted by $H$, 
$H^o$ is a subset of the interior $\nu_{R_i}^o$, and $\nu^o_{R_i}$ and 
$\nu^o_{R_j}$ intersect transversally at points of $H^o$.  
\item  $\nu_{R_i} \cap R_j$  is  a tame $(n-1)$-bihedron with boundary 
the union of  $H$  and an $(n-2)$-hemisphere  $H'$  in the closure  of 
$\alpha_{R_j}$ with its interior $H^{\prime o}$ in $\alpha_{R_j}$.  
\item $\nu_{R_i}\cap R_j$ is the closure of  a component of $\nu_{R_i} 
- H$ in $\che M$.  
\item  $R_i  \cap  R_j$  is  the  closure  of  a component  of  $R_j - 
\nu_{R_i}$.  
\item  Both $\alpha_{R_i}   \cap \alpha_{R_j}$ and   $\alpha_{R_i}\cup 
\alpha_{R_j}$ are homeomorphic to open $(n-1)$-dimensional balls. 
(See Figure \ref{fig:tr}.)  
\end{enumerate}  
 
The above mirrors  the  property of  intersection of  $\dev(R_1)$  and 
$\dev(R_2)$  where $\dev(\alpha_{R_1})$ and $\dev(\alpha_{R_2})$   are 
included in a  common  great sphere $\SI^{n-1}$ of  dimension $(n-1)$, 
$\dev(R_1)$ and $\dev(R_2)$ are included  in a common   $n$-hemisphere 
bounded by $\SI^{n-1}$ and $\dev(\nu_{R_1})^o$ and $\dev(\nu_{R_2})^o$ 
meet transversally (see Proposition \ref{prop:overlap}).

\begin{thm}[Theorem 5.4 \cite{psconv}]\label{thm:transversal} 
Suppose  that $R_1$ and $R_2$  are overlapping.  Then either $R_1$ and 
$R_2$  intersect transversally, or $R_1 \subset  R_2$  or $R_2 \subset 
R_1$ hold. 
\end{thm} 
\begin{proof}  
The proof is in \cite{psconv}. Note we need the following lemma  
to prove this. 
\end{proof} 
 
\begin{lem}[Lemma 5.2 \cite{psconv}] 
Let $N$ be a closed real projective $n$-manifold. Suppose that  
$\dev: \che N_h \ra \SI^n$ is an imbedding onto the union of  
$n$-hemispheres $H_1$ and $H_2$ meeting each other in an $n$-bihedron  
or an $n$-hemisphere. Then $H_1 = H_2$, and $N_h$ is projectively  
diffeomorphic to an open $n$-hemisphere.  
\end{lem}

\chapter{Sequences of $n$-balls} 
\label{app:seqballs} 
First, we discuss for $\SI^n$ the  convergence of a sequence of convex 
$n$-balls, the limit of  a convergent sequence,  and the boundary  and 
the   interior  of  the  limit.     
Next, we consider sequences of convex  $n$-balls in  
a  Kuiper completion. Subsequences  
may not converge  since the Kuiper completion  is  not compact in  
general. However,  when   a sequence has   a  {\em core},  i.e.,  
a common convex open $n$-ball in  every sequence  element, a subsequence   
will have a ``limit." So, a certain criterion assuring the  existence  
of a core is presented first. Finally,  certain sequences with  cores  
are  shown to ``converge."  
The details of this chapter is given  
in Part I of \cite{psconv}. 
 
Since $\SI^n$ is a compact metric space, the collection of all compact 
subsets of   $\SI^n$   form a  compact metric   space   with  the {\em 
Hausdorff}\/ metric $\bdd^H$  defined as follows: for compact  subsets 
$A$ and   $B$,  given $\eps > 0$, we  define $\bdd^H(A,   B) <  \eps$    
if  $A$ is in $\eps$-$\bdd$-neighborhood of $B$ and $B$ is in that of  
$A$.  
 
A  sequence  of  compact convex   subsets of  $\SI^n$    always has  a 
convergent subsequence $\{D_i\}$ with respect  to $\bdd$. The limit is 
always a compact convex subset. The dimension of the limit is less  
than or equal to $\liminf_{i \ra \infty}\dim(D_i)$. We deduce from   
this that the limit of  a convergent  sequence of convex  $n$-balls   
is a convex $m$-ball for  $m  \leq n$.  As  in  \cite{cdcr1}, 
if $m=n$, we can deduce that $\delta  D$  is the  limit  of  
$\{\delta D_i\}$, and $\bigcup_{i =1}^\infty D_i^o \supset D^o$  
(see Part I of \cite{psconv}).

\begin{lem}\label{lem:conelim}  
Let $D_i$ be a sequence of $3$-balls in $\clo(\mathcal{H})$ which are cones 
over   $2$-disks $B_i$ in  $\SIT$.  Then  $D_i$  converges to a  limit 
compact  set $D_\infty$ if and  only   if $B_i$ converges  to a  limit 
compact set $B_\infty$ in $\SIT$.  Moreover, in this case, $D_\infty$ 
is a cone over $B_\infty$.  
\end{lem} 
\begin{proof} 
Straightforward.  
\end{proof} 
 
Recall that  $\mu$ denote the  Riemannian metric on $\che M_h$ induced 
from that of $\SI^n$ by $\dev$.  
\begin{lem}[Proposition 3.14 \cite{psconv}]\label{lem:common}  
Let $\{D_i\}$ be a sequence of convex $n$-balls in  $\che M_h$. Let $x 
\in M_h$, and  $B(x)$ a tiny ball  of $x$. Suppose that the  following 
properties hold\/{\rm :}  
\begin{enumerate}  
\item $\delta D_i$ includes an $(n-1)$-ball $\nu_i$.  
\item  $B(x)$ overlaps with  $D_i$ and  does  not meet $\delta D_i - 
\clo(\nu_i)$.  
\item A sequence $\{x_i\}$ converges to $x$ where  $x_i \in \nu_i$ for 
each $i$.  
\item The sequence $\{\mathbf{n}_i\}$  converges where ${\mathbf{n}}_i$ is the 
outer-normal vector to $\nu_i$ at $x_i$ with respect to $\mu$ for each 
$i$.  
\end{enumerate} 
Then there exist  a positive integer $N$  and a convex open disk  
$\mathcal{P}$ in $B(x)$ such that  
\[\mathcal{P} \subset D_i \hbox{ whenever } i > N.\] 
\end{lem} 

We say  that  a compact  subset   $D_\infty$ of $\SI^n$   is the  {\em 
resulting set}\/ of a sequence   $\{D_i\}$ of compact subsets  of $\che 
M_h$  if $\{\dev(D_i)\}$  converges to  $D_\infty$.  Let $\{D_i\}$ and 
$\{B_i\}$  be  sequences  of  convex  $n$-balls  with   resulting sets 
$D_\infty$ and $B_\infty$ respectively; let $\{K_i\}$ be a sequence of 
compact  subsets with the   resulting   set $K_\infty$.  We say   that 
$\{D_i\}$ {\em subjugates}\/ $\{K_i\}$ if $D_i \supset K_i$ for each $i$ 
and that  $\{B_i\}$  {\em  dominates}\/  $\{D_i\}$ if  $B_i$   and $D_i$ 
overlap for each $i$ and if  $B_\infty \supset D_\infty$. Moreover, we 
say that $\{K_i\}$ is {\em ideal}\/  if there is  a positive integer $N$ 
for  every compact subset $F$ of  $M_h$ such that $F  \cap K_i = \emp$ 
whenever $i > N$.  In particular, if $K_i$ is a subset of $\hideal{M}$ 
for each $i$, then $\{K_i\}$ is an ideal subjugated sequence.  
 
\begin{thm}[Proposition 3.15 \cite{psconv}]\label{thm:seqconv}  
Suppose that $\{D_i\}$ is a sequence of $n$-balls including a common  
open ball $\mathcal{P}$, $\{D_i\}$ subjugates a sequence of  
compact subsets $\{K_i\}$, and a sequence of $n$-balls $\{B_i\}$  
dominates $\{D_i\}$. Let $D_\infty, B_\infty$, and  
$K_\infty$ be the resulting sets of the above sequences  
$\{D_i\}, \{B_i\}$, and $\{K_i\}$. Then $\che M_h$ includes  
two convex $n$-balls $D^u$ and $B^u$ and a compact subset $K^u$ with  
the following properties\/{\rm :}  
\begin{enumerate} 
\item $D^u \supset \mathcal{P}$, and $\dev(D^u) = D_\infty$.  
\item $B^u \supset D^u$, and $\dev(B^u) = B_\infty$.  
\item $D^u \supset K^u$, and $\dev(K^u) = K_\infty$.  
\item If $\{K_i\}$ is ideal, then $K^u \subset \hideal{M}$.  
\end{enumerate} 
\end{thm}

\chapter[Radiant affine $3$-manifolds with boundary]
{Radiant affine $3$-manifolds with boundary, and certain   
radiant affine $3$-manifolds} \label{app:radd} 

\medskip
 
\centerline{By Thierry Barbot\footnote{ \large
Address: UMPA, \'Ecole Normale Sup\'erieure de Lyon, \hfill\break
46, all\'ee d'Italie, 
LYON, France \hfill\break
(e-mail: Thierry.Barbot{@}umpa.ens-lyon.fr)}
\footnote{supported by CNRS}
and Suhyoung Choi} 

\medskip

Let $G$ be a group acting on an analytic manifold $X$.  
An $(X,G)$-manifold is a manifold admitting an atlas  
with charts with value in $X$ and whose coordinate  
change mappings are restrictions of elements of $G$.  
It is well-known that equipping a manifold $M$  
with an $(X,G)$-structure is equivalent to giving  
a pair $(\dev, \rho)$, where $\dev$ is  
an immersion from the universal covering $\tilde{M}$  
of $M$ into $X$, and where $\rho$ is a homomorphism from the  
fundamental group $\Gamma$ of $M$ into $G$, such that  
\[ \forall \gamma \in \Gamma \;\;\;\;\;  
{\dev} \circ \gamma = \rho(\gamma) \circ \dev .\]  
Here, the action of $\Gamma$ on $\tilde{M}$ is of  
course the action by deck transformations.  
The map ${\dev}: \tilde{M} \rightarrow X$   
is the {\em developing map\/} of the  
structure, and $\rho: \Gamma \rightarrow G$ is  
the {\em holonomy homomorphism.\/}  
  
A radiant affine $n$-manifold is an $(X, G)$-manifold,  
where $X$ is the vector space $\bR^n$, and $G$ is the group  
$\GL(n,\bR)$ of linear transformations (see \cite{Carlet}).  
Such a manifold is naturally equipped with  
a transversely projective flow, the so-called {\em radial flow,\/}  
defined as follows:  
if $(x_{1},\ldots,x_{n})$ are local coordinates, the vector  
field generating the radial flow is:  
\[ X(x_{1}, \ldots,x_{n})= \Sigma_{i=1}^{n} x_{i}{\partial}_{x_{i}} \]  
Observe that this vector field does not depend on  
the coordinate systems as long as the origins   
are the same, and thus induces well-defined vector   
field on $\tilde M$ and on $M$, which is said   
to be a {\em radiant vector field}.  
The flow generated by the vector field is said to be  
a {\em radial flow}. The radial flow in $\bR^n$ has unique singularity   
at the origin $O$ but the radial flow on $M$ 
has no singularity since $\dev$ misses $O$ 
(see \cite{Carlet} and Lemma \ref{lem:radcom}). 
  
Let $N$ be a closed real projective $(n-1)$-manifold, i.e. 
an $(\bR P^{n-1}, \PGL(n, \bR)))$-manifold, where ${\bR}P^{n-1}$ is  
the real projective space of dimension $n-1$, and  
$\PGL(n,\bR)$ is the group of projective transformations.   
Let $\varphi$ be a projective automorphism of $N$.  
We can associate to the pair $(N,\varphi)$ a family  
of radiant affine closed $n$-manifold: i.e., generalized  
affine suspensions, homeomorphic to a topological suspension   
of $N$ by $\varphi$ (see Chapter \ref{ch:racom} and   
\cite{Ben,smillie,Carlet,resolu}).  
They can be characterized by the following property:   
{\em a closed radiant affine manifold is   
a generalized affine suspension if and only if its radial  
flow admits a total cross-section, i.e. there is a closed   
embedded submanifold transverse everywhere to the flow and   
which meets every orbit of the flow.\/}   
(See Proposition 3.2.)  
(Note that the term ``affine suspension is reserved for   
the case when $N$ and $\varphi$ are both affine.)  
A {\em Benz\'ecri suspension} is an affine suspension  
so that $\varphi$ is the identity or a finite-order  
automorphism of $N$. In this case all orbits are closed,   
and $N$ is Seifert-fibered.  
  
In this appendix, we study the following  
particular case which was left from Chapter \ref{ch:caff},   
i.e., Lemma \ref{lem:BC}, which is implied by  
the following theorem since developing maps of   
a universal cover is always obtainable from developing maps  
of the holonomy cover composed with the covering map   
to the holonomy cover from the universal cover.  

\begin{thm}\label{thm:noraff}  
There is no closed radiant affine $3$-manifold whose developing  
map from the universal cover is an infinite cyclic covering   
over ${\bR}^{3}$ minus a line.  
\end{thm}  
  
\vspace{.2cm}  
  
A subsurface $S$ in an affine $3$-manifold is {\em totally-geodesic}   
if every point of $S$ has a neighborhood $\mathcal{O}$ 
affinely diffeomorphic to   
an open subset of $\bR^3$ or of an affine half-space of   
$\bR^3$ so that $S \cap \mathcal{O}$ corresponds to a closed affine   
subspace of codimension-one intersected with the open set.  
A totally geodesic subsurface has a natural induced   
affine structure as a two-dimensional manifold.  
A boundary component of an affine $3$-manifold is   
{\em totally geodesic} if each boundary point has a neighborhood   
affinely diffeomorphic to an open subset of a half-space   
in $\bR^3$.   
  
\begin{thm}[Theorem B of Barbot \cite{barbot1}]  
Let $M$ be a closed radiant affine $3$-manifold.  
If $M$ includes a totally geodesic surface tangent   
to the radial flow, then $M$ admits a total cross-section\/{\rm ;}  
i.e., $M$ is a generalized affine suspension.  
\end{thm}  
  
The second result of this appendix (Theorem \ref{thm:barbot2}),  
we will prove is:  
\begin{thm}\label{thm:barbot2II}  
Let $M$ be a compact radiant affine $3$-manifold with a non\-emp\-ty   
totally geodesic boundary. Suppose that each component is   
affinely homeomorphic to the quotient of a convex cone  
or $\bR^2 -\{O\}$ by an affine action.  
Then $M$ admits a total cross-section  
to the radial flow, and hence is a generalized affine suspension.  
\end{thm}  
  
We cannot prove this theorem by a method of ``doubling'':  
Some radiant affine $3$-manifold $N$ may not   
be doubled; i.e., there may not be a radiant affine $3$-manifold  
homeomorphic to the topological double of $N$ in which $N$   
and a copy of $N$ are affinely imbedded, meeting at boundary.  
An example comes from a convex real projective surface $\Sigma$ with   
negative Euler characteristic and the boundary component   
with holonomy $\vth$ that has a nondiagonalizable   
$3\times 3$-matrix with two distinct positive eigenvalues   
(\cite{cdcr1} and \cite{cdcr3}.) These real projective surfaces   
exist, of course, as one can see that the construction   
of convex real projective surfaces in Goldman \cite{Gconv}   
can be modified to construct convex ones with this behavior.  
Such a holonomy $\vth$ does not commute   
with a projective automorphism in $\bR P^2$   
whose fixed points comprise a subspace that $\vth$ preserves.  
We can easily see that the Benz\'ecri suspension  
of $\Sigma$ cannot be doubled in the above sense.  
  
The proof of the theorem is essentially that of Theorem B   
in \cite{barbot1} where the totally geodesic surface now is   
in the boundary.  
  
We remark that the theorem is true without the assumption   
on boundary component which can be proved applying   
Barbot's method \cite{barbot1}. 
For simplicity of proof, we prove this weaker but   
sufficient version here.  
  
The proof of the first one goes as follows: we assume the existence   
of a radiant affine $3$-manifold whose developing map  
is an infinite cyclic covering over $\bR^{3} \setminus \{ x=y=0 \}$  
where $x, y$, and $z$ denote the standard coordinate   
functions of $\bR^3$.  
  
In the first section, we prove that the holonomy group  
is solvable: indeed, if not, it contains a hyperbolic element  
$\rho(\gamma)$ with one eigenvalue greater than $1$, another less  
than $1$ (and positive), and the last exactly equal to $1$.  
The contradiction nearly arises from the fact that   
such a linear transformation does not act properly discontinuously  
on $\bR^{3} \setminus \{ x=y=0 \}$, whereas $\gamma$ must  
act properly discontinuously on $\tilde{M}$.  
  
Since the holonomy group is solvable, the affine $3$-manifold  
is a generalized affine suspension (\cite{resolu}). Therefore,  
a short argument that no projective surface has a developing map 
which is an infinite cyclic covering over the sphere minus two points
completes the proof.
  
For the second theorem, we will only prove for the cases   
when the fact that $M$ has nonempty boundary makes any  
difference from the proof of Theorem B of \cite{barbot1}.  
  
If the holonomy of the boundary component   
contains a homothety, i.e., a linear transformation that is   
a positive multiple of the identity map, then all radial flow   
orbits are periodic and it follows easily that our manifold has   
a total cross-section. First, we show that if the boundary   
surface is not convex as an affine $2$-manifold, then   
our affine manifold is a half-Hopf manifold.  
Then we look at the holonomy group of the fundamental   
group of the boundary component, which we may assume is   
an affine torus, and classify them into six cases as   
in \cite{barbot1}. Only four cases are applicable since   
the boundary torus is convex. We will show that   
in each case, our radiant affine manifold $M$ is finitely   
covered by torus times an interval. Either $M$ is foliated 
by tori corresponding to orbits of certain group actions 
or decomposes into submanifolds which are affinely isomorphic to   
domains in $\bR^3$ quotient out by linear $\bZ + \bZ$-action.  
If $M$ is foliated, there is a quick away to show that $M$ is
a generalized affine suspension using Carri\`ere's volume 
argument \cite{Carlet}. Otherwise, we show that 
the third case is a generalized affine suspension  
and forth cases are impossible  
and in the remaining two cases, the pieces must be generalized  
affine suspensions. (This proof is mostly a generalization   
of that of Theorem B in Barbot \cite{barbot1})  
  
\section[Certain radiant affine $3$-manifolds]
{The nonexistence of certain radiant affine $3$-manifolds}  
  
Let $M$ be a closed radiant affine $3$-manifold.  
We denote by $\Phi^{t}$ its radial flow.  
We denote by $p: \tilde{M} \rightarrow M$ the  
universal covering (we don't worry about the choice of a base point).  
Let $\Gamma$ be the fundamental group of $M$; it acts  
naturally on $\tilde{M}$.  
  
Let ${\dev}: \tilde{M} \rightarrow \bR^{3}$ be the  
developing map, and $\rho: \Gamma \rightarrow \GL(3,\bR)$ be the  
holonomy homomorphism.  
They satisfy:  
\[ \forall \gamma \in \Gamma \;\;\;\;\;  
{\dev} \circ \gamma = \rho(\gamma) \circ \dev. \]  
  
As above, the radial vector field induces radial flows in $M$ and   
$\tilde M$ respectively. The orbits are said to be   
{\em rays} and $\dev$ restricted to each ray is a homeomorphism  
onto rays in $\bR^3$ by Lemma \ref{lem:inje}; i.e., an open   
half-line with an endpoint at $0$.  
  
\begin{lem}  
\label{lem:inje}  
Let $G$ be a Lie group acting on two manifolds $X$ and $Y$.  
Let $f: X \rightarrow Y$ be a function equivariant  
for the actions of $G$. Let $x$ be an element of $X$  
such that $f(x)$ is fixed by no element of $G$.  
Then, the restriction of $f$ to the $G$-orbit of $x$ is injective.  
\end{lem}  
  
\begin{proof}  
For every element $g$ of $G$ we have $f(gx)=gf(x)$.  
\end{proof}

We now assume that $\dev$ is an infinite cyclic covering map over   
$\bR^{3} \setminus \Delta$,  
where $\Delta$ is a line through the origin $O$.  
Our aim is to obtain a contradiction.   
  
Since we want to show that such a $M$ does not exist, we are free to  
replace $M$ by any finite covering of itself. For example,  
we can consider only the case where $M$ is oriented, i.e. the  
case where every element of the holonomy group is of  
positive determinant.  
  
Since $\dev$ is well-defined  
up to composition by a linear transformation, we can assume  
that $\Delta$ is the line $\{ x=y=0 \}$.  
Then, since $\Delta$ has to be $\rho(\Gamma)$-invariant, every   
element $\rho(\gamma)$ of the holonomy group is of the form:  
  
\[ \rho(\gamma)=  
\left( \begin{array}{cc}  
\bar{\rho}(\gamma) &   
\begin{array}{c}  
0 \\  
0 \end{array} \\  
\begin{array}{cc}  
\ast & \ast   
\end{array}  
& \lambda(\gamma)  
\end{array} \right)  
\]  
where $\lambda(\gamma)$ is a non-zero real number, and  
$\bar{\rho}(\gamma)$ an element of $\GL(2,\bR)$.  
Clearly, $\lambda$ and $\bar{\rho}$ are homomorphisms.

We discuss more on generalized affine suspensions   
(see also Chapter \ref{ch:racom}):  
Let $\varphi: N \rightarrow N$ a projective diffeomorphism of  
a real-projective $(n-1)$-manifold $N$.  
Recall that $\SI^{n-1}$ has a real projective structure induced   
from $\rp^{n-1}$ by the standard double cover, and   
the group $\Aut(\SI^{n-1})$ of projective automorphisms of   
$\SI^{n-1}$, which is isomorphic to the quotient group of   
$\GL(n, \bR)$ by homotheties. ($\SI^{n-1}$ with this   
structure is said to be a {\em real projective sphere}.)  
We can always lift the chart of $N$ to $\rp^{n-1}$ to $\SI^{n-1}$   
with respect to the standard double cover. Then   
the transition functions then lie in $\Aut(\SI^{n-1})$  
since a projective map defined on a small open   
subset of $\SI^{n-1}$ extends to one defined on $\SI^{n-1}$ always  
(see Chapter 2 of \cite{psconv}). We gather that   
$N$ has a natural $(\SI^{n-1}, \Aut(\SI^{n-1}))$-structure.  
  
Let $f_{i}: V_{i} \rightarrow U_{i} \subset \SI^{n-1}$ be  
a family of projective charts covering $N$. When  
$V_{i}$ meets $V_{j}$, we have an element $\bar{g}_{ij}$  
of $\Aut(\SI^{n-1})$ such that on $V_{i} \cap V_{j}$:  
\[ f_{i}=\bar{g}_{ij} \circ f_{j}. \]  
  
Let's choose representatives $g_{ij}$ of the $\bar{g}_{ij}$ in  
$\GL(n,{\bR})$. We impose the condition $g_{ij}g_{jk}g_{ki}=\Idd$  
if $V_{i} \cap V_{j} \cap V_{k}$ is not empty. Such a choice  
is always possible: take for example the unique representative  
of $\bar{g}_{ij}$ with determinant $\pm 1$. The set of the  
possible choices is parameterized by the cohomology group   
$H^{1}(N,{\bR})$ of $N$.  
For every $i$, let $W_{i}$ be the open cone in ${\bR}^{n}$  
with vertex at $0$, the union of the half lines belonging to  
$U_{i}$. Let denote by $W$ the quotient of the disjoint union of  
the $W_{i}$ by the relation identifying each element $x_{j}$ of  
$W_{j}$ with the element $g_{ij}(x_{j})$ of $W_{i}$ (when  
$g_{ij}(x_{j})$ belongs effectively to $W_{i}$, of course).  
This quotient is a noncompact radiant affine manifold, equipped   
with a complete radial flow $\hat{\Phi}^{t}$.   
The quotient of  
$W$ by the relation `being on the same orbit of $\hat{\Phi}^{t}$'  
is homeomorphic to $N$. The quotient map is a fibration by rays.  
Let $N_{0}$ be any section of this fibration.  
The manifold $W$ is diffeomorphic to $N \times {\bR}$.  
  
Remember the projective transformation $\varphi$ of $N$.  
It lifts to an affine diffeomorphism $\hat{\varphi}$ of $W$  
well-defined up to composition by $\hat{\Phi}^{t}$.  
If $T$ is big enough, $\hat{\Phi}^{T}\hat{\varphi}(N_{0})$ is   
a section of $\hat{\Phi}^{t}$ disjoint from $N_{0}$.   
Therefore, for every real positive $t$,  
$\langle \hat{\Phi}^{t}\hat{\varphi}\rangle$ acts freely and properly  
discontinuously on $W$. The quotient of this action is a closed  
radiant affine manifold homeomorphic to the topological suspension  
$N_{\varphi}$ of $\varphi: N \rightarrow N$.

Actually, such a lifting does not always exist for any  
choice of $W$, but for many of $W$ above the given  
$\Sigma$, we can perform such liftings.   
The condition is: let $\bar{\rho}: \pi_{1}(W) \rightarrow \GL(n,\bR)$   
be the holonomy homomorphism  
of $W$. Observe that $\pi_{1}(W)$ is isomorphic to $\pi_{1}(N)$.  
Let $\varphi_{\ast}$ be the automorphism of $\pi_{1}(N)$ induced  
by $\varphi$.  Then, $\varphi$ lifts if and only if  
$\det \circ \bar{\rho}$ is constant on the orbits of $\varphi_{\ast}$.  
For example, the choice of the $g_{ij}$'s of determinant  
$\pm 1$ works.   
  
Observe also that the construction is not uniquely defined:   
we made some choices. These choices are  
parameterized by an open subset of the first cohomology   
module $H^{1}(N_{\varphi},{R})$  
satisfying the above requirement.

By construction, the radial flow of a generalized affine suspension  
admits a closed total cross-section homeomorphic to $\Sigma$. 
Note that  
this section, equipped with the projective structure induced  
by the transverse projective structure of the radial flow  
is isomorphic to the initial projective surface $\Sigma$.  
(See Proposition 3.2.)

Returning back to our radiant affine $3$-manifold $M$:  
\begin{prop}  
\label{prop:solv}  
The holonomy group $\rho(\Gamma)$ is solvable.  
\end{prop}  
  
\begin{proof}  
Denote by $\Gamma'$ the first commutator subgroup of $\Gamma$.  
Since $\lambda$ and $\bar{\rho}$ are  
homomorphisms, for every element of $\Gamma'$ we have:  
\begin{itemize}  
\item $\bar{\rho}(\gamma)$ belongs to $\SL(2,\bR)$,  
\item $\lambda(\gamma)=1$.  
\end{itemize}  
Observe that by definition $\rho(\Gamma)$ is solvable if and only if   
$\bar{\rho}(\Gamma')$ is solvable.

Let $\mathcal{F}^{0}$ be the foliation of $\bR^{3} \setminus \Delta$  
whose leaves are the half-planes containing $\Delta$ in their  
boundaries. The leaf space of this foliation, i.e. the  
quotient of $\bR^{3} \setminus \Delta$ by the relation  
``being on the same leaf of $\mathcal{F}^{0}$'',  
is naturally identified with the double covering of  
the real projective line $\bR P^{1}$.   
Let $\mathcal{F}$ be the pull-back of $\mathcal{F}^{0}$ by  
$\mathcal{D}$. Since $\dev$ is an infinite cyclic covering,  
$\mathcal{F}$ is a foliation whose leaf space is naturally identified  
with the universal covering $\tilde{P}^{1}$ of $\bR P^{1}$.   
The action of $\Gamma'$ on the leaf space induced by the action  
of $\Gamma$ on $\tilde{M}$ is a lifting of   
the projective action of $\bar{\rho}(\gamma')   
\in \SL(2,\bR)$ over $\bR P^{1}$. According to Lemma \ref{lem:lift2}  
below, if $\bar{\rho}(\Gamma')$ is not solvable, there is   
an element $\gamma$ of $\Gamma'$ preserving a leaf  
$F$ of $\mathcal{F}$, and such that in an adequate coordinate system:  
\[ \rho(\gamma) =  
\left( \begin{array}{ccc}  
\lambda & 0 & 0 \\  
0 & \lambda^{-1} & 0 \\  
0 & 0 & 1  
\end{array}  
\right)  
\]  
for some real positive $\lambda$ different from $1$.  
  
We fix this coordinate system with coordinate functions   
denoted by $x, y,$ and $z$ still. Since the required coordinate   
change sends the standard coordinate vectors to  
the eigenvectors of $\rho(\gamma)$, and the $z$-axis is in  
the eigendirection, $\Delta$ is still given by $x = 0$ and $y =0$.  
Leaves of $\mathcal{F}^{0}$ are again given as zero sets   
of linear functions of the new coordinate functions $x$ and $y$ only.  
We may assume without loss of generality that $F$ maps   
in the plane given by $x = 0$ under $\dev$.  
  
Let $P$ be the inverse image by $\dev$ of the punctured plane   
$\{ z=0 \} \setminus \{ (0,0,0) \}$. Since $\dev$ is  
an infinite cyclic covering missing $\Delta$, it follows that $P$   
is connected, and the restriction of $\dev$ to $P$ is   
an infinite cyclic covering to the punctured plane.   
Moreover, $\gamma$  
preserves $P$. Since $\gamma$ preserves the leaf $F$ also,  
$\gamma$ preserves each connected component of   
${\dev}^{-1}(\{ z=x=0 \})$ and, by same reason,   
${\dev}^{-1}(\{ z=y=0 \})$.   
Let $C$ be a connected component of  
${\dev}^{-1}(\{ z=0 , x \geq 0 , y \geq 0 \})$. It  
is preserved by $\gamma$, and the restriction of $\dev$ to  
$C$ is a homeomorphism over $\{ z=0 , x \geq 0 , y \geq 0 \}   
\setminus \{ (0,0,0) \}$. The action of $\rho(\gamma)$ on   
$\{ z=0 , x \geq 0 , y \geq 0 \} \setminus \{ (0,0,0)$ is  
given by $(x,y,0) \mapsto (\lambda x , \lambda^{-1}y , 0)$.  
It is not properly discontinuous, since any path joining   
$\{ z=0 , x=0 , y \geq 0 \}$ to $\{ z=0 , x \geq 0 , y=0 \}$   
intersects all its iterates by $\rho(\gamma)$.   
This is a contradiction since the action of $\gamma$ on $C$ has   
to be properly discontinuous. It follows that $\bar{\rho}(\Gamma)$,   
and therefore $\rho(\Gamma)$, is solvable.  
\end{proof}  
  
For the proof of the following lemma \ref{lem:lift2},   
we must first recall some  
facts about the action of $\PSL(2,\bR)$ on $\bR P^1$,
and the universal covering group $\widetilde{\SL}(2,\bR)$ 
action on the universal covering space $\tilde P^1$ of $\bR P^1$.
Let $q:\widetilde{\SL}(2,\bR) \rightarrow \PSL(2,\bR)$  
denote the covering map. Every element $g$ of   
$\PSL(2,\bR)$ is either:  
\begin{itemize}  
\item {\em elliptic:\/} $g$ has no fixed point in $\bR P^{1}$. It  
is conjugate to a rotation,  
\item {\em parabolic:\/} $g$ has one and only one fixed point.  
This fixed point is of saddle-node type, i.e.  
attractive on one side, and repulsive on  
the other side,  
\item {\em hyperbolic:\/} $g$ has two fixed points: a repulsive one and  
an attractive one. It is conjugate to the element represented by:  
\[ \left( \begin{array}{cc}  
\lambda & 0 \\  
0 & \lambda^{-1}  
\end{array}\right).  
\]  
\end{itemize}  
  
An element of $\widetilde{\SL}(2,\bR)$ is said to be {\em elliptic,   
parabolic or hyperbolic}\/ according to the nature of its projection   
$q(g)$. If this projection is trivial, $g$ belongs to the center $H$ of  
$\widetilde{\SL}(2,\bR)$. The group $H$ is infinite cyclic. Let $h$   
be a generator of $H$. If $g$ is not trivial and admits fixed points  
on the universal covering $\tilde{P}^{1}$ of $\bR P^1$, it is  
parabolic or hyperbolic. In the first case, its fixed points are  
of saddle-node type; in the second case, they are attractive  
or repulsive.   
  
\begin{lem}  
\label{lem:lift2}  
Let $\Gamma$ be a subgroup of $\widetilde{\SL}(2, \bR)$.  
We assume that $q(\Gamma)$ is not solvable.  
Then, it contains  
a hyperbolic element that fixes a point of $\tilde{P}^{1}$.  
\end{lem}  
\begin{proof}  
Note that $\Gamma$ is not solvable since it is a cyclic extension  
of $q(\Gamma)$ which is not solvable. According to H\"older's theorem  
(see e.g. \cite{HH}, IV.3.1), a group acting freely on the real line  
is abelian. Therefore, the action of $\Gamma$   
on $\tilde{P}^{1}$ is not  
free: some element $\gamma_{0}$ of $\Gamma$ admits a fixed  
point $x_{0}$ in $\tilde{P}^{1}$. If $\gamma_{0}$ is hyperbolic,  
we are done. If not, $\gamma_{0}$ is parabolic.  
Then, the fixed points of $\gamma_{0}$ are the $H$-iterates of $x_{0}$.  
We denote by $x_{i}$ the image of $x_{0}$ by $h^{i}$.  
Observe that since $\tilde{P}^{1}$ is homeomorphic  
to the real line $\bR$, orienting $\tilde{P}^{1}$ is equivalent  
to equip it with an archimedian total order.  
We orient $\tilde{P}^{1}$ in such a way that  
$x_{1}$ is greater than $x_{0}$. Taking the inverse of $\gamma_{0}$   
if necessary, we can assume that all the $\gamma_{0}$-orbits   
in the open interval $]x_{0},x_{1}[$ go from $x_{0}$ to $x_{1}$.

The stabilizer in $\PSL(2,\bR)$ of a point in $\bR P^{1}$  
is isomorphic to the group of affine transformations of the line.  
It is therefore solvable. It follows that there is an  
element $\gamma$ of $\Gamma$ such that $\gamma(x_{0})$ is not  
one of the $x_{i}$'s. Let $\gamma_{1}$ be the conjugate   
$\gamma\gamma_{0}\gamma^{-1}$.  
It is parabolic and fixes $\gamma(x_{0})$. Therefore,  
it admits a fixed point $x'_{0}$ in $]x_{0},x_{1}[$.   
Then, $\gamma_{1}^{-1}\gamma_{0}(x_{0})=\gamma_{1}^{-1}(x_{0})$  
is less than $x_{0}$, and $\gamma_{1}^{-1}\gamma_{0}(x'_{0})$  
is greater than $x'_{0}$, since $x'_{0}$ is a fixed point  
of $\gamma_{1}^{-1}$ and $\gamma_{0}(x'_{0})$ is greater than  
$x'_{0}$. Therefore, the closed interval $[x_{0},x'_{0}]$ is  
contained in its image by $\gamma_{1}^{-1}\gamma_{0}$.  
It follows that $\gamma_{1}^{-1}\gamma_{0}$ is a hyperbolic  
element admitting a repulsive fixed point in  
$]x_{0},x'_{0}[$.  
\end{proof}

We know from Proposition \ref{prop:solv} that the  
holonomy group is solvable. It follows from  
Theorem A of \cite{resolu} that $M$ is  
affinely isomorphic to a generalized affine suspension.  
In particular, the radial flow admits a total cross-section.

Let $\Sigma$ be such a total cross-section, and $\tilde \Sigma$ a  
lifting of $\Sigma$ in $\tilde{M}$, i.e. a connected  
component of $p^{-1}(\Sigma)$. Let $\tilde{\Phi}^{t}$ be  
the lifting of $\Phi^{t}$ in $\tilde{M}$. Since  
$\Sigma$ is a total cross-section, it is a fiber of some fibration  
of $M$ over the circle. Hence, $\tilde M  \setminus \tilde \Sigma$   
is not connected. Every orbit of $\tilde{\Phi}^{t}$  
meets $\tilde{\Sigma}$. This orbit remains in the past in one  
connected component of $\tilde{M} \setminus \tilde{\Sigma}$,   
and in the future, it remains in  
the other connected component. In other words,  
every orbit of $\tilde{\Phi}^{t}$  
meets $\tilde{\Sigma}$ at one and only one point.   
The developing map sends injectively every orbit  
of $\tilde{\Phi}^{t}$ over a half-line in $\bR^{3}$   
(Lemma \ref{lem:inje} applied to the $\bR$-actions).  
Therefore, it induces a local homeomorphism $\dev'$ from  
$\tilde{\Sigma}$ on the sphere $\SI^{2}$ of half-lines.  
We denote by $\SI_{\ast}$ the sphere $\SI^{2}$ punctured at   
$(0,0,1)$ and $(0,0,-1)$.  
Since $\dev$ is an infinite cyclic covering over   
$\bR^{3} \setminus \Delta$, the map 
$\dev'$ is an infinite cyclic covering over   
$S_{\ast}$. Therefore, $\tilde{\Sigma}$ is the universal covering  
of $\Sigma$, and $\dev'$ is the developing map of a   
real projective structure on $\Sigma$.   
The holonomy homomorphism $\hat{\rho}$ of this structure is   
the composition of the restriction of $\rho$ to $\hat{\Gamma}$   
with the projection of $\GL(3,\bR)$ in $\Aut(\SI^2)$, where   
$\hat{\Gamma}$ is the group of elements of $\Gamma$ which   
preserve $\tilde{\Sigma}$.

In order to find a contradiction, i.e., in order  
to achieve the proof of Theorem \ref{thm:noraff},   
it suffices to show:  
  
\begin{prop}\label{prop:final}  
Given a real projective structure on the closed surface $\Sigma$,  
its developing map $\dev'$ can not be an infinite cyclic  
covering over $\SI_{\ast}$.  
\end{prop}  
\begin{proof}   
Suppose that $\Sigma$ is a closed surface with such a structure.  
We first complete $\tilde \Sigma$ by the path-metric induced from   
the Riemannian metric $\mu$ on $\SI^2$ by $\dev'$  
obtaining the Kuiper completion $\che \Sigma$ of $\Sigma$.   
Recall that $\ideal{\Sigma}$ denotes the set of ideal points   
$\che \Sigma - \tilde \Sigma$. The developing map $\dev'$   
also extends to an obvious distance decreasing map   
$:\che \Sigma \ra \SI^2$. In our situation, it is   
easy to see that there exist only two ideal points in   
$\che \Sigma$, mapping to $(0,0,1)$ and $(0,0,-1)$, and that   
$\che \Sigma$ is obtained from $\tilde \Sigma$ by adding these   
two points.  
  
Our surface $\Sigma$ is obviously not convex since   
$\tilde \Sigma$ is not convex. A {\em $2$-crescent} in $\che M$ is   
a convex hemisphere or lune $D$ in $\che M$ with interior in $\tilde M$   
and the interior of a convex segment in the boundary $\delta D$  
of $D$ includes the nonempty $\tilde M \cap \delta D$.   
Theorems 4.6 and 4.5 of \cite{psconv} show   
that $\che \Sigma$ includes a $2$-crescent (see also  
Section 5 of \cite{cdcr1}). By definition of $2$-crescents,  
there exists a nontrivial open arc in the boundary of   
the $2$-crescents that is in $\ideal{\Sigma}$, and hence   
the set of ideal points in a $2$-crescent is uncountable.  
However, $\che \Sigma$ contains only two ideal points.   
This is a contradiction.  
\end{proof}  
  
\section[Radiant affine $3$-manifolds with boundary]
{Radiant affine $3$-manifolds with boundary have 
total cross-sections}  
  
Now we begin the proof of Theorem \ref{thm:barbot2II}.  
Let $M$ be a compact radiant affine $3$-manifold with totally  
geodesic boundary. Since we may prove the result for   
a finite cover of $M$, we assume without loss of generality  
that the boundary components of $M$ are tori.

\begin{lem}\label{lem:form}   
A radiant affine $3$-manifold $M$ admits a total cross-section  
if and only if it has a closed $1$-form taking a positive   
value for each radiant vector.  
\end{lem}  
\begin{proof}  
The existence of a total cross-section and the flow   
is transversal to it shows that $M$ is diffeomorphic to   
a bundle over a circle so that the radiant vector field   
corresponds to the vector field transversal to each fiber.   
The differential of the fiber map gives us the closed form.

Given a closed form with above property, we can   
approximate it by a non-vanishing closed form with rational period.   
Such a closed form obviously gives a fibration   
$M \ra \SI^1$ (see \cite{Tis}).  
\end{proof}

\begin{lem}   
If $M$ is a radiant affine $3$-manifold, and a finite cover $N$  
of $M$ admits a total cross-section to the radial flow,   
then $M$ admits a total cross-section.  
\end{lem}  
\begin{proof}   
A total cross-section in $N$ corresponds to a closed   
$1$-form on $N$ which is positive for radial vectors.   
Clearly such a $1$-form descends to one on $M$ by averaging   
over the finite group action as the action preserves 
the flow direction.
\end{proof}

Let $\tilde M$ be the universal cover of $M$ with a development pair   
$(\dev, h)$. The radiant flow lines induce a {\em radiant foliation}\/   
on $\tilde M$. Let $Q$ be the space of   
leaves of the radiant foliation in $\tilde M$, which has   
a natural real projective structure. The group of  
deck transformations acts on $Q$ as a group of projective   
automorphisms (see Barbot \cite{barbot1} for details).   
As $\tilde M$ is simply connected, $Q$ is simply-connected also.

There is a quotient map $f: \tilde M \ra Q$ which is   
a fibration whose fibers are rays.   
We see that the developing map $\dev$ induces   
an immersion $\dev':Q \ra \SI^2$ where $\SI^2$ is the space of  
rays in $\bR^3$, and the deck-transformation   
group acts on $Q$ so that $h'(\vth) \circ \dev'   
= \dev' \circ \vth$ for a deck transformation   
$\vth$ and $h'(\vth)$ the induced projective map  
from $h(\vth)$.   
  
Choose a boundary component $K$ of $M$ and   
a component $\tilde K$ of $p^{-1}(K)$.  
As $K$ is tangent to the radial flow, $K$ has Euler characteristic   
zero. By taking a finite cover of $M$, $K$ is assumed to be   
a torus. A deck-transformation group $G$ which is a subgroup   
of $\pi_1(M)$ and isomorphic to $\bZ$ or   
$\bZ + \bZ$, acts on $\tilde K$. Let $c$ be the image of   
$\tilde K$ in $Q$, which is a simple geodesic in   
the boundary $\partial Q$ of $Q$.  
  
An affine automorphism $\varphi:\tilde M \ra \tilde M$ is   
a {\em homothety} if each ray is preserved and   
$\dev \circ \varphi = s \Idd \circ \dev$ for a positive scalar $s$.   
Note that an affine automorphism $\varphi$ of $\tilde M$ always   
induces a real projective automorphism of $Q$ and $\varphi$ acts   
trivially on $Q$ if and only if $\varphi$ is a homothety.  
  
Let $\vth$ be an element of $G$. Then   
$\vth$ is not a homothety.  
If not, then the radial flow is periodic, and $M$ is easily shown  
to be a Benz\'ecri suspension by Proposition 3.3   
of \cite{barbot1}. (Note that for these arguments, there   
is no difference when $M$ is closed or has nonempty boundary.)  
We are done in this case.  
  
As $\tilde K$ is totally geodesic, $c = f(\tilde K)$ is   
a geodesic boundary component of $Q$.  
Suppose that $K$ is compressible in $M$. Then $\tilde K$   
is compressible in $\tilde M$. Let $D$ be an imbedded   
compressing disk with boundary in $\tilde M$.   
Then the radial projection $f|D: D \ra Q$ maps   
a disk $D$ onto $Q$ with boundary $\partial D$ onto  
a boundary component $c = f(\tilde K)$ of $Q$.  
This means that $Q$ is homeomorphic to a compact surface, and   
hence $Q$ is a compact disk. As $Q$ is bounded   
by a geodesic $c$, $Q$ must be projectively   
diffeomorphic to a $2$-hemisphere, and   
$\dev'$ is an embedding onto a $2$-hemisphere in $\SI^2$.  
(This can be seen by a doubling argument and the uniqueness of  
the projective structure on $\SI^2$.)  
This shows that $\tilde M$ is affinely diffeomorphic by $\dev$   
to an affine half-space with boundary containing $O$   
with $O$ removed.  
  
\begin{prop}\label{prop:hHopf}  
Suppose that $Q$ is real projectively diffeomorphic to   
a hemisphere in $\SI^2$, or equivalently $\tilde M$ is affinely  
diffeomorphic to an affine half-space $H_1$ with boundary  
containing $O$ with $O$ removed.  
Then $M$ is a half-Hopf manifold.  
\end{prop}   
\begin{proof}   
We will use the second hypothesis. The punctured 
half-plane $H_1 - \{O\}$   
includes a compact disk $D$ with boundary in $\partial H_1 -\{O\}$  
which is transversal to every ray.  
There exists a deck transformation $\vth$ of $H_1 - \{O\}$   
sending $D$ to $\vth(D)$ disjoint from $D$ as 
the deck-transformation groups are properly discontinuous.  
Clearly $H_1 -\{O\}$ quotient out by $\vth$ has   
a total cross-section corresponding to $D$, and is   
a half-Hopf manifold. Hence, $M$ is finitely covered by   
a generalized affine suspension, and so $M$ is   
a generalized affine suspension. The total cross-section   
has positive Euler characteristic, and hence is a compact   
disk transversal to radial flow. The lemma now follows.  
\end{proof}  
  
Assume that $K$ is incompressible from now on.  
First, suppose that $K$ is affinely homeomorphic   
to a quotient of $\bR^2 -\{O\}$ by an affine action.   
Now, we need to use the {\em holonomy cover}\/ $M_h$ of $M$,   
i.e., the cover of $M$ corresponding to the kernel of   
the holonomy homomorphism $h$. As we described   
in \cite{psconv}, the developing map  
$\dev$ induces an immersion $\dev_h:M_h \ra \bR^3$ and   
the holonomy homomorphism induces a homomorphism   
$h_h:\pi_1(M)/\pi_1(M_h) \ra \GL(3, \bR)$. Also   
$f$ induces a fibration $f_h:M_h \ra Q_h$ to a projective   
surface $Q_h$ covered by $Q$. The immersion   
$\dev_h$ induces an immersion $\dev'_h:Q_h \ra \SI^2$  
and $h_h$ induces a homomorphism   
$h'_h:\pi_1(M)/\pi_1(M_h) \ra \Aut(\SI^2)$.  
The surface $\tilde K$   
corresponds to a surface $K_h$ covering $K$.   
We see easily that $K_h$ is affinely diffeomorphic   
to $\bR^2 - \{O\}$. A deck-transformation group $G_h$ acts on   
$K_h$ so that $K$ is affinely homeomorphic to $K_h/G_h$.   
Hence $G_h$ is an infinite cyclic group.  
  
A {\em homothety} in $M_h$ is an affine automorphism $\varphi$ of  
$M_h$ acting on each ray and satisfying   
$\dev_h\circ \varphi = s\Idd \circ \dev_h$ for a positive scalar $s$.  
As above, if an element of $G_h$ is a homothety, then $M$   
is a generalized affine suspension.  
  
Assume now that no element of $G_h$ is a homothety.  
We claim that in this case $M$ is a half-Hopf manifold.  
Let $\vth$ be a generator of $G_h$. Then $h_h(\vth)$ acts   
on the totally geodesic plane $P$ including $\dev_h(K_h)$.   
We assume that $P$ is the $xy$-plane for simplicity.   
As no element of $G_h$ is a homothety, $G_h$ acts effectively  
on $Q_h$. We divide our cases according to the conjugacy classes of  
$h'_h(\vth)$ in $\Aut(\SI^2)$. Let $P'$ be the great   
circle in $\SI^2$ corresponding to $P$, and $H_1$ and   
$H_2$ the two hemispheres bounded by $P'$.  
We easily see that one of the following occurs:  
\begin{itemize}  
\item[(i)] The only $h'_h(\vth)$-invariant subset of $H_i$ including   
a neighborhood of $P'$ in $H_i$ is $H_i$ itself.  
\item[(ii)] The only subset of $H_i$ with this property   
equals $H_i - \{x\}$ for a point $x \in H_i^o$.   
\item[(iii)] $h'_h(\vth)$ under a suitable coordinates of   
the affine space $H_i^o$ is of form  
\[\begin{pmatrix} \cos 2\pi\theta & \sin 2\pi\theta \\  
-\sin 2\pi\theta & \cos 2\pi\theta \end{pmatrix}, \theta \ne 0. \]  
\end{itemize}  
  
Note that $\dev'_h|c$ is an imbedding onto $P'$,   
and hence a neighborhood $U$ of $c$ in $Q_h$ imbeds   
onto a neighborhood of $P'$ in say $H_1$.  
We see that $\dev'_h\circ \vth|U = h'_h(\vth) \circ \dev'_h| U$,   
and so $\dev'_h|\vth(U)$ is also an imbedding onto  
a neighborhood of $P'$. Therefore, by induction,   
we see that $\dev'_h| \bigcup_{i\in \bZ} \vth^i(U)$ is  
an imbedding onto an $h'_h(\vth)$-invariant subset of   
$H_1$ including a neighborhood of $P'$.  
  
In case (i), $Q_h$ maps homeomorphic onto $H_1$, and   
by Proposition \ref{prop:hHopf}, $M$ is a half-Hopf manifold.  
  
In case (ii), $Q_h$ maps homeomorphic to $H_1 - \{x\}$.  
Therefore $M_h$ under $\dev_h$ maps homeomorphic to $H'_1 - l$ for  
the upper half-space $H'_1$ corresponding to $H_1$   
and a line $l$ through $O$   
transversal to the $xy$-plane. Hence, we may   
identify $M_h$ with its image, and each deck transformation   
acts on $l$ and the upper-half space. We see that the group of  
deck transformations acting on $K_h$ equals the entire   
deck-transformation group of of $M_h$ and   
hence the group of deck transformations acting on $\tilde K$   
equals $\pi_1(M)$. This means that $K$ is homotopy   
equivalent to $M$ by the inclusion map where $K = \partial M$.  
We see easily by a topological doubling argument and   
the top dimensional $\bZ_2$-homology group computation that 
such a situation cannot   
happen.   
  
In case (iii), $\theta$ must be irrational since otherwise   
$h'_h(\vth)$ is periodic and hence $\vth$ must be periodic near   
$c_h$ and hence on $Q_h$; but this means that a power of   
$\vth$ is a homothety on $M_h$.   
  
Under the suitable coordinates,   
for some large $r$, there is an open neighborhood $U_r$ of $c$   
which under $\dev'_h$ maps homeomorphic to the set of form   
the union of $P'$ and the complement of a ball of radius $r$.   
  
Let $r_0$ be an infimum of possible values of $r$.  
Then we see easily that $U_{r_0}$ exists.   
We claim that (1) $Q_h$ equals this set $U_{r_0}$ or (2) $Q_h$ maps   
homeomorphic to $H_1$ under $\dev'_h$.   
Suppose that there   
exists a boundary point $x$ of $U_{r_0}$ in $Q_h$.  
As $\vth$ acts as an irrational rotation, we see   
that every point of the boundary of $\dev'(U_{r_0})$   
is an image of a boundary point of $U_{r_0}$ in $Q_h$.  
Assume that $r_0 > 0$.  
In this case, the boundary $\gamma$ of $U_{r_0}$ in $Q_h$   
maps homeomorphic to the closure $\gamma'$ of   
the orbit of $\vth$ in $H_1$, which is homeomorphic to   
a circle. If there exists a point of $\partial Q_h$  
in the boundary $\gamma$, then every point of   
$\gamma$ is in $\partial Q_h$ as $\partial Q_h$ is closed   
and $\vth$ acts on $\partial Q_h$. But this implies that   
a boundary component of $M_h$ is not totally geodesic.   
Therefore, there exists a regular neighborhood of $\gamma$   
mapping to a regular neighborhood of $\gamma'$.   
This contradicts the minimality of $r_0$ if $r_0 >0$.  
Hence, we obtain that $U_{r_0} = Q_h$.  
When $r_0 = 0$, the claim follows easily.   
  
In case (1), we obtain that $M_h$ maps homeomorphic   
under $\dev_h$ to the complement $L$ of a closed   
convex cone in $H'_1$ not meeting the boundary   
$P - \{O\}$. We identify $M_h$ with   
this set $L$. Then a contradiction as in (ii) occurs.   
  
In case (2), we see that $M$ is a half-Hopf manifold  
by Proposition \ref{prop:hHopf}.  
  
From now on, we will be working on $\tilde M$ (i.e., not on $M_h$)
and assume that $\tilde K$ is affinely homeomorphic   
to a convex cone in $\bR^2$.  
  
Since $h(G)$ acts on a convex cone $\dev(\tilde K)$, if   
$h(\varphi)$ for $\varphi \in G$ is a homothety, then   
$\varphi$ is a homothety   
near $\tilde K$ in $\tilde M$, and hence $\varphi$ is a homothety   
on $\tilde M$. Thus $q \circ h$ is injective.   
We identify $G$ with $h(G)$ from now on.  
  
Let $q:\GL(3, \bR) \ra \SL(3, \bR)$ be the homomorphism   
whose kernel consists of $s\Idd$ for $s \ne 0$.  
Since no element of $G$ is a homothety in $\tilde M$,   
Barbot \cite{barbot1} shows that the identity component $L$  
of the Zariski closure of $q(G)$ is conjugate to a subgroup of  
the following groups:   
\begin{description}  
\item[Case D] the group of all diagonal matrices with   
positive eigenvalues.  
\item[Case P] the group of matrices of form:   
\[\begin{pmatrix}  
e^u & t & 0 \\ 0 & e^u & 0 \\ 0 & 0 & e^v \end{pmatrix},   
u, v, t \in \bR, 2u + v = 0,\]  
\item[Case U] the group of matrices of form:   
\[\begin{pmatrix}   
1 & s & t \\ 0 & 1 & 0 \\ 0 & 0 & 1 \end{pmatrix}, s, t \in \bR, \]  
\item[Case C] the group of matrices of form:  
\[\begin{pmatrix}   
1 & s & t \\ 0 & 1 & s \\ 0 & 0 &1 \end{pmatrix}, s, t \in \bR,\]  
\item[Case S] the group of matrices of form:  
\[\begin{pmatrix}   
e^u\cos \theta & e^u \sin\theta & 0 \\   
-e^u \sin \theta  & e^u \cos\theta & 0 \\   
0 & 0 &e^v \end{pmatrix}, u, v, \theta \in \bR, 2u + v = 0,\] or  
\item[Case T] the group of matrices of form:  
\[\begin{pmatrix}   
1 & 0 & s \\ 0 & 1 & t \\ 0 & 0 &1 \end{pmatrix}, s, t \in \bR.\]  
\end{description}  
 
Since $q(G) \cap L$ must be a finite index subgroup of $q(G)$,   
we assume that $q(G)$ is a rank-two abelian subgroup of $L$   
by taking a finite cover of $M$ and choosing $(\dev, h)$   
carefully, i.e., conjugating $h$ by an element of $\GL(3, \bR)$. 

In cases S and T, if $L$ is one-dimensional, then 
$L$ can be considered a subgroup of $D$ and $C$ respectively 
by conjugations. 
In case S, since $\vth$ acts on $\dev(\tilde K)$ a convex 
two-dimensional subset, $\theta$ is always $0$ 
for each element of the group $G$ and we are reduced 
to the one-dimensional case.
In case T, assuming that $L$ is two-dimensional,
$G$ acts on the $xy$-plane as a homothety,  
and there are no other $G$-invariant subspaces of codimension-one.
But since $\dev(K)$ is a $G$-invariant subspace   
so that $\dev(K)/G$ is homeomorphic to the torus $K/G$,   
this is a contradiction by the following lemma. 
Hence the case T does not occur.  
  
\begin{lem}\label{lem:homoth}  
Let $D$ be a convex cone in $\bR^2$, and $G'$ an abelian   
group isomorphic to $\bZ + \bZ$. If $D/G$ is homeomorphic   
to a torus, then an element of $G'$ is not a homothety.  
\end{lem}  
\begin{proof}  
If each element of $G'$ is a homothety, then each   
element of $G'$ acts on each ray in $D$ ending   
at $O$. Since $\bZ + \bZ$ cannot act properly discontinuously  
and freely on a real line, this is absurd.  
\end{proof}  
  
We see that $G$ is a lattice  
in a connected two-dimensional subgroup   
$H$ of following group of matrices in $\GL(3, \bR)$:  
\begin{description}  
\item[Case D] the group of all diagonal matrices of   
positive eigenvalues.  
\item[Case P] the group of matrices of form   
\[\begin{pmatrix}   
e^u & t & 0 \\ 0 & e^u & 0 \\ 0 & 0 &e^v \end{pmatrix}, u, v, t \in \bR.\]  
\item[Case U] the group of matrices of form:   
\[\begin{pmatrix}   
e^u & s & t \\ 0 & e^u & 0 \\ 0 & 0 &e^u \end{pmatrix}, u, s, t \in \bR,\] or  
\item[Case C] the group of matrices of form:  
\[\begin{pmatrix}   
e^u & s & t \\ 0 & e^u & s \\ 0 & 0 &e^u \end{pmatrix}, u, s, t \in \bR.\]  
\end{description}  

The Lie group $q(H)$ is not zero dimensional, as $G$ does not contain
a homothety. $q(H)$ could be one-dimensional or two-dimensional.
If $q(H)$ is two-dimensional, then $q(G)$ is a lattice in $q(H)$; 
otherwise, $q(G)$ is dense in $q(H)$. 

The group $H$ acts on the projective sphere $\SI^2$ effectively  
as a subgroup of the group $\Aut(\SI^2)$ of projective   
automorphisms of $\SI^2$.  
The space $\SI^2$ decomposes into two-dimensional open   
orbits and one-dimensional orbits and zero-dimensional orbits   
under $H$. We have to divide our case to when $q(H)$ is 
one-dimensional and when two-dimensional. (1) indicates the one
dimensional cases which differ from (2) by having no 
two-dimensional orbits. The zero dimensional orbits are the same,
and there are additional one-dimensional orbits, which
foliate the regions corresponding to two-dimensional orbits.   
\begin{description} 
\item[Case D(2)] The zero-dimensional orbits are six  
points in $\SI^2$ comprising three pair of antipodal   
points, one-dimensional ones are twelve lines in $\SI^2$   
with endpoints in the three points, and two-dimensional   
orbits are eight open triangles bounded by the closures of the lines. 
(See Figure \ref{fig:dd}.) 
\begin{figure}[t] 
\centerline{\epsfxsize=12cm \epsfbox{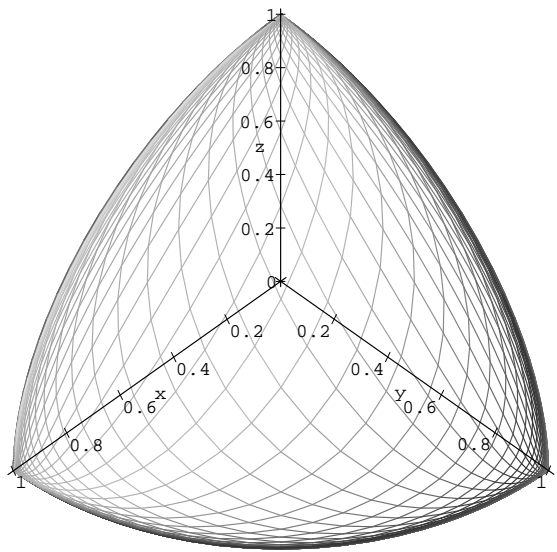}}  
\caption{\label{fig:dd} An example of an $H$-action on $\SI^2$ in case D(2).
The sphere itself is not drawn, and
the arcs correspond to orbits of one-parameter subgroups.} 
\end{figure} 
\item[Case D(1)] The additional one-dimensional orbits 
are curves curved in one direction connecting a fixed 
pair of zero-dimensional orbits, or lines that are components 
of great circles passing through a common antipodal
pair of zero-dimensional orbits with the union of 
the above twelve lines removed.
\item[Case P(2)] The zero-dimensional orbits comprise two   
pair of antipodal points $\{p, -p\}$ and $\{q, -q\}$,   
one-dimensional orbits are four lines with endpoints   
$\pm p$ and $\pm q$ and two lines that   
are components of a great circle containing $p, -p$ with $p$    
and $-p$ removed, and two-dimensional orbits are four open lunes   
bounded by the closure of the union of one-dimensional orbits.
(See Figure \ref{fig:pp}.)
\begin{figure}[t] 
\centerline{\epsfxsize=14cm \epsfbox{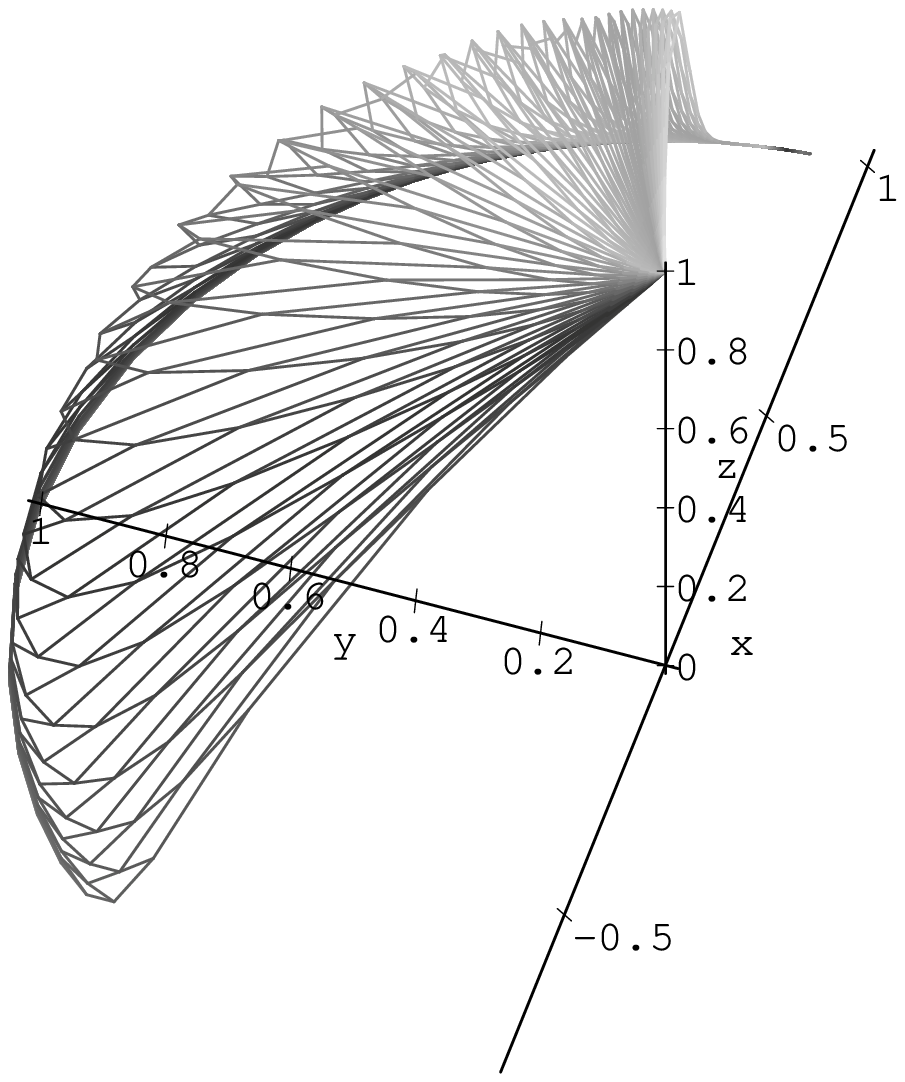}}  
\caption{\label{fig:pp}} The case P(2): the orbits are partially drawn here. 
\end{figure}   
\item[Case P(1)] The additional one-dimensional orbits are curves 
connecting a pair of zero-dimensional orbits curved in 
one-direction.
\item[Case U(2)] Zero-dimensional   
orbits are two points $p$ and $-p$. One-dimensional orbits   
are great circles containing $p$ and $-p$ with $p$ and $-p$ removed.   
Their union equals $\SI^2 - \{p, -p\}$ and there are no  
two-dimensional orbits.   
\item[Case U(1)] Zero-dimensional orbits are points of a great circle   
through two points $\{p, -p\}$. One dimensional orbits are   
the other great circles through $p$ and $-p$ with $p$   
and $-p$ removed. There are no two-dimensional orbits.    
\item[Case C(2)] Zero-dimensional orbits are two points   
$p$ and $-p$. One-dimensional orbits are two components of   
a great circle containing $p$ and $-p$ with $\{p, -p\}$ removed.   
Two-dimensional orbits are two open hemispheres bounded   
by the great circle. (See Figure \ref{fig:cc}.)
\begin{figure}[h] 
\centerline{\epsfxsize=12.5cm \epsfbox{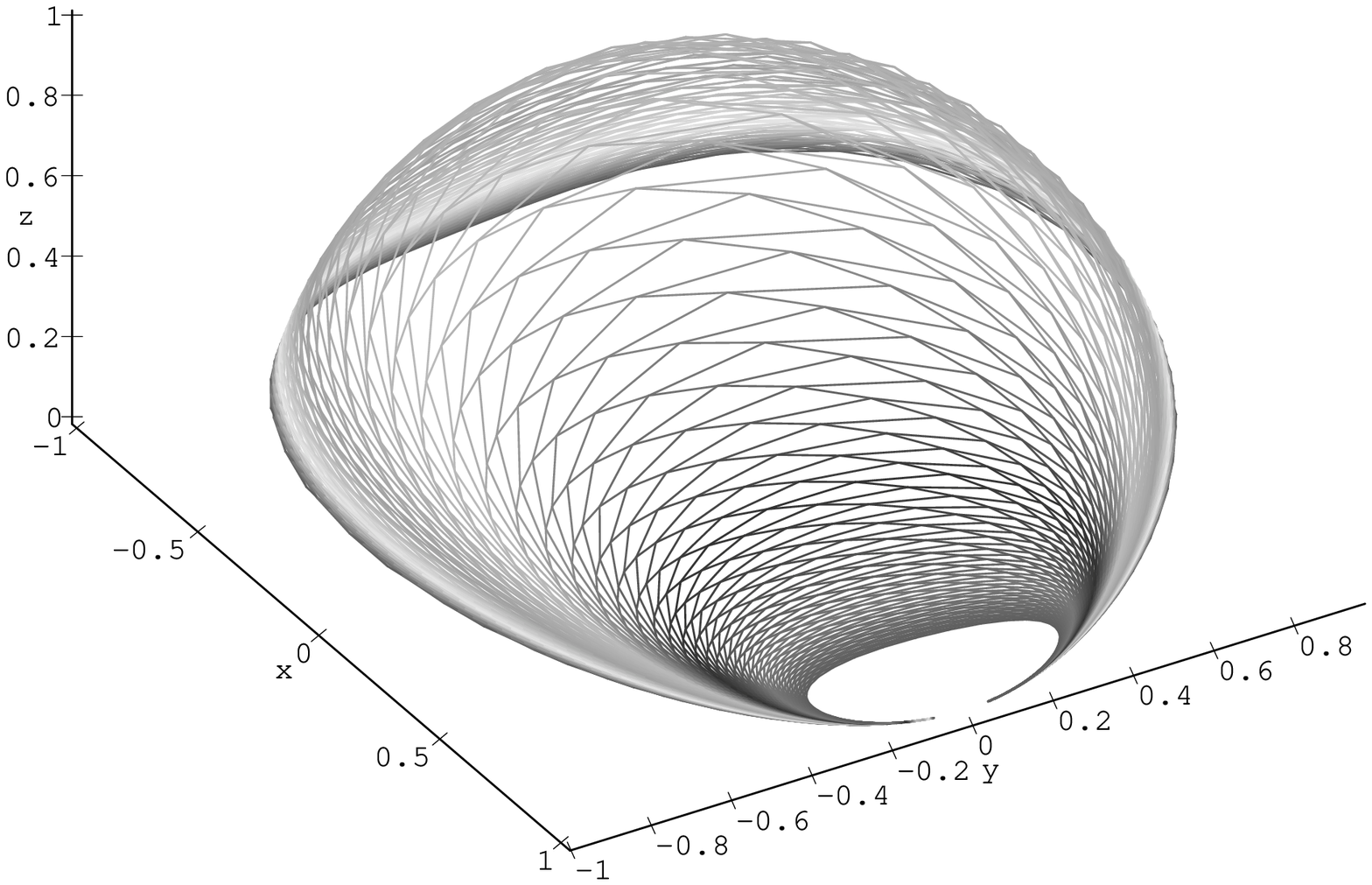}}  
\caption{\label{fig:cc} The case C(2).} 
\end{figure} 
\item[Case C(1)] The additional one-dimensional orbits are curves in 
the open hemisphere connecting $p$ to $p$ or $-p$ to $-p$,
they are curved in one directions.   
\end{description}  
We see that all one-dimensional orbits are simple geodesics   
in $\SI^2$ and are not closed ones when $q(H)$ is two-dimensional.
(To see some more figures of these kinds of actions, see \cite{cdcr3}.)
  
Recall that the developing map $\dev:\tilde M \ra \bR^3$   
induces an immersion $\dev':Q \ra \SI^2$ and 
the deck-transformation group acts on $Q$ as well so that   
$\dev' \circ \vth$ equals $h'(\vth) \circ \dev'$ where   
$h'(\vth) = q(h)(\vth)$ for each deck transformation   
$\vth$ of $\tilde M$.  
(The action is not necessarily proper.)  
  
\begin{lem}[Lemma 2.4.3 of Dupont \cite{dupont}]\label{lem:extend}  
Let $V$ be a $(U, X)$-manifold where $U$ is a Lie group  
acting on a space $X$. Let $L$ be a connected subgroup   
of $U$, and $\omega$ an $L$-orbit in $X$ and $\hat \omega$   
a connected component of $D^{-1}(\omega)$ where $(D, j)$ is   
a development pair of $V$. Suppose that the action of $L$   
on $\omega$ is covered by the action of $L$ on $\hat \omega$,   
and there exists a subgroup $\Gamma_0$ of the deck-transformation group   
so that $\hat \omega/\Gamma_0$ is compact and $j(\Gamma_0)$   
is in $L$. Then the action of $L$ on a neighborhood of   
$\omega$ is covered by an action of $L$ on a neighborhood of   
$\hat \omega$ such that, for any element $\gamma$  
of $\Gamma_{0}$, the action of $j(\gamma)$ coincides  
with the action of $\gamma$ as a deck transformation.  
\end{lem}  
  
Since $G$ acts on $f(\tilde K) \subset Q$ as in the premise of   
the above lemma, the extension argument (i.e., the proof of   
Theorem B in \cite{barbot1} using essentially Lemma \ref{lem:extend})   
shows that $H$ acts on $Q$, $Q$ is a union of   
orbits of dimension zero, one, or two, and each orbit maps   
homeomorphic to an orbit in $\SI^2$ under $\dev'$.   
The proof of this claim, similar to what is in \cite{barbot1},  
goes as follows:   
Let $\mathcal{O}$ be the maximal connected subset of   
$Q$ including $c$ where the $H$-action is defined  
and whose restriction to $G$  
coincide with the deck transformation action.  
Our claim follows from the following lemma:  
  
\begin{prop}\label{prop:o}  
The $H$-action is defined everywhere\/{\rm ;} 
i.e. $\mathcal{O}$ equals $Q$.  
\end{prop}

\begin{proof}  
We have to give a more precise definition of $\mathcal{O}$:  
this is the maximal connected subset including $c$ of   
the set of points $x$ in $Q$ such that:  
\begin{itemize}  
\item 
There is a continuous action of $H$ on $Q$,
denoted by $(k,x) \in H \times Q \mapsto k.x$ so that, for  
any $(k,x)$, we have $\dev'(k.x)=k\dev'(x)$.
\item for every element $g$ of $G$, we have the  
equality $g.x=gx$ (remember that $G$ is  
a group of deck transformations, and therefore acts  
on $Q$).  
\end{itemize}  
  
The fundamental group of an orbit of the radial flow   
is trivial or cyclic. It follows that  
there are no zero-dimensional $H$-orbits in $\mathcal{O}$.  
(An abelian group of rank $2$ such as $G$ cannot act on   
a connected one-dimensional space properly discontinuously   
and freely.)  
  
By Lemma \ref{lem:extend}, $\mathcal{O}$ is open.   
Since $Q$ is connected, the lemma will be proven  
if we show that $\mathcal{O}$ is closed.  
  
As $H$ acts on $\mathcal{O}$, $\mathcal{O}$ is a union of   
one- or two-dimensional $H$-orbits in $Q$.  
The restriction of $\dev'$ 
to each $H$-orbit maps homeomorphic   
to an $H$-orbit in $\SI^2$ by Lemma \ref{lem:inje}.  
In particular, two-dimensional $H$-orbits in $\mathcal{O}$   
are open surfaces, and $G$ acts properly   
discontinuously and freely on each of them.  
  
We note that each two-dimensional $H$-orbit in $\mathcal{O}$  
has to have at least one adjacent one-dimensional   
$H$-orbit. If not, the two-dimensional orbit   
is disconnected from $c$, a contradiction.  
  
We see that $\mathcal{O}$ is a union of 
two-dimensional orbits and adjacent one-dimensional   
orbits joined in a ``chain-like'' manner  
in cases D(2), P(2), and C(2) or is a union of
one-dimensional orbits in other cases.
  
Suppose that $x$ is a boundary point of $\mathcal{O}$ in $Q$.   
We aim to obtain a contradiction. If $\dev'(x)$ lies in   
a zero-dimensional orbit of $H$, then $x$ is a zero-dimensional  
orbit of $H$. This is a contradiction by above.   
The image $\dev'(x)$ does not lie in a two-dimensional orbit   
as each two-dimensional $H$-orbit of $Q$ is open in $Q$  
and $\mathcal{O}$ is a union of $H$-orbits of $Q$.  
Let $\dev'(x)$ be in a one-dimensional orbit $J$, 
and let $J'$ be a component of $\dev^{\prime -1}(J)$   
containing $x$. Then $J'$ is a closed subset of $Q$.  
The restriction of $\dev'$ to $J'$ is injective.  
  
Suppose that the intervals $J$ and $\dev'(J')$  
have a common endpoint. Then, for any element $g$  
of $G$, the intervals $\dev'(J')$ and $g \dev'(J')$  
are not disjoint. Actually,   
inverting $g$ if necessary, we can assume that  
$\dev'(J')$ contains $g \dev'(J')$.  
Since the restriction of $\dev'$ to every $H$-orbit in   
$\mathcal{O}$ is injective, a union $A$ of some $H$-orbits   
near $J'$ is an open set where $H$ acts, and   
the restriction of $\dev'$ to $A$ is injective.   
Since $\dev'$ restricted to $A \cup J'$ and $A \cup g(J')$   
are both injective, $J'$ contains $g(J')$.   
It follows that $J'$ is $G$-invariant.  
We deduce that $J'$ is contained in $\mathcal{O}$; a contradiction. (*)
   
Therefore, the endpoints of $\dev(J')$ both belong to $J$.  
  
Suppose that no element of $G$ fixes a point of $J$.   
Then, all the $G$-orbits in $J$ are dense.  
There is an element $g$ of $G$ such that  
$g(\dev'(J'))$ meets $\dev'(J')$. Since as above   
$\dev'$ restricted on $A \cup J'$ and   
$A \cup g(J')$ are both injective, $J'$ is $g$-invariant,  
meaning that $H$ acts on $J'$ again.  
  
This contradiction shows  
that some element $g$ of $G$ acts as an identity map on $J$.  
Therefore, the surface $D$ in $\tilde M$  
corresponding to $J'$ is filled with rays on which   
$g$ acts, and hence, maps to the union of   
closed orbits of radial flow in $M$.  
  
We claim that, for any deck transformation  
$\vth$,  $J'$ cannot meet $\vth(J')$ where $J'$ and $\vth(J')$ are 
distinct. Suppose not. Then $D$ and $\vth(D)$ meet at   
an isolated ray $l$ by the real analyticity of $J'$ and $\vth(J')$.
This ray is fixed by $g$; hence, $g$ is a power of 
a deck transformation $g'$ so that   
$l/\langle g'\rangle$ maps to the closed orbit in $M$.   
Similarly, $\vth \circ g \circ \vth^{-1}$ is a power of   
$g'$. Thus, a finite power $g''$ of $g$ acts trivially  
in $J'$ and $g(J')$ since a power of $\vth \circ g \circ \vth^{-1}$  
must equal a power of $g$. Since we can find four 
nearby fixed points in $Q$ of $g''$ contained in $J'$ and $\vth(J')$, 
$g''$ is a homothety, a contradiction. (To see this,
recall that orbits are real analytic submanifolds of $Q$.)
  
Also, the collection of sets of form $\vth(J')$ are   
locally finite. If not, then there exists a sequence of   
points $p_i \in \vpi_i(D)$ converging to $p \in \tilde M$  
for a sequence of deck transformations $\vpi_i$.   
As $\vpi_i \circ g \circ \vpi_i^{-1}$ acts on $\vpi_i(D)$   
as homotheties by a fixed factor $s, s > 0$, we see that   
$\vpi_i \circ g \circ \vpi_i^{-1}$ moves $p_i$ to a points   
in a compact subset of $\tilde M$. The properness of
the deck-transformation-group action shows that   
infinitely many of these deck transformations must be   
the same. We may choose $\vpi_i(D)$ and $\vpi_j(D)$   
sufficiently close so that $\vpi_i\circ g \circ \vpi_i^{-1}$   
equals $\vpi_j \circ g \circ \vpi_j^{-1}$.   
However this means that $\vpi_i\circ g \circ \vpi_i^{-1}$   
is a homothety since it acts trivially on $\vpi_i(J')$  
and another curve $\vpi_j(J')$, a contradiction.  
  
By above two conclusions, it follows that $D$ covers   
a closed surface in $M$ tangent to the radial flow,   
and a subgroup $G'$ of rank $2$ acts on $D$ and on $J'$.   
All we did above apply once more: there is   
a connected abelian group $H'$ including $G'$   
as a dense subgroup or a lattice, 
and $H'$ is again of type D, P, U, or C. 
(This is easy if $J'$ is geodesic. If not, we look at 
the line connecting the endpoints of $\dev(J')$.) 
We define similarly to $\mathcal{O}$ the locus $\mathcal{O}'$  
of definition of the $H'$-action on $Q$. It contains $J'$.  
Since $\dev'(J')$ has both endpoints inside $\dev(J)$,   
$H'$ cannot be in case U or C, since in these  
cases, $\dev'(J')$ should equal $J'$ by geometric reasons.  
  
Suppose that $H'$ is in case D(2) or P(2), and
$H$ is in case D(2), P(2), or C(2).  
Then, $J'$ is contained in the boundary of a two-dimensional   
$H$-orbit $A$.  
By looking at the image under $\dev'$ of $A$  
and the two-dimensional $H'$-orbit $B$ meeting $\dev'(A)$,   
we obtain that there exists a geodesic $c$ in $A$ mapping into   
a one-dimensional $H'$-orbit adjacent to $\dev'(B)$  
so that $\dev'(c)$ and $\dev'(J')$ share an endpoint.  
Since an endpoint of $\dev'(c)$ is an $H'$-orbit,  
$c$ must be included in an $H'$-orbit by the same   
reason as the paragraph marked by (*) with $G$ replaced by $G'$.  
By the following lemma \ref{lem:orbdim2}, this is   
a contradiction.   
  
If $H$ is in case D(1), P(1), C(1) or U(1)(2), then  
as in the above paragraph a one-dimensional $H'$-orbit meets   
a one-dimensional $H$-orbit transversally. (We can see this 
by the transversality theorem applied to the foliation by orbits,
i.e., a parameter of submanifolds \cite{GuiPol}.) 
This is a contradiction by Lemma \ref{lem:orbdim1}.  

Suppose that $H'$ is in case D(1) or P(1), and $H$ is
in case D(2), P(2), or C(2). Then clearly a one-dimensional
orbit of $H'$ near $J'$ meets a two-dimensional orbit of $H$ 
adjacent to $J'$, a contradiction by Lemma \ref{lem:orbdim2}. 
As above, $H$ cannot be in case D(1), P(1), C(1) or U(1)(2), 
This final contradiction achieves the proof of Proposition
\ref{prop:o}.  
\end{proof}  
    
\begin{lem}\label{lem:orbdim1}  
Let $H'$ be a connected two-dimensional abelian group 
including an abe\-lian group $G'$ of rank $2$ of 
deck transformations of $\tilde M$.
Then a one-dimensional orbit $K'$ of $H'$ in $\SI^2$   
does not meet a one-dimensional orbit $K$ of $H$ 
in $\SI^2$ transversally.  
\end{lem}  
\begin{proof}  
Suppose that $K$ and $K'$ meet at a point $x$ in $Q$.   
Then let $S$ and $S'$ be the corresponding
surfaces in $\tilde M$. They meet at a ray   
$l$ corresponding to $x$. As $S/G$ and $S'/G'$  
correspond to immersed tori, $l$ corresponds to   
a closed orbit of the radial flow in $M$.  
There is an element of $g$ in $G \cap G'$ corresponding   
to this orbit. Then $\langle g\rangle$ acts trivially on $K$ and $K'$  
as we can see from the matrices of form D, P, U, and C;  
i.e., $G$ acts on one-dimensional orbits   
in $\SI^2$ as translations with respect to certain coordinates.
Since we can find four nearby fixed points of $g$ in $K \cup K'$, 
the holonomy of $g$ must be a homothety. 
This contradicts our assumption.  
\end{proof}  
  
\begin{lem}\label{lem:orbdim2}  
Let $H'$ and $G'$ be as in the preceding lemma.   
Let $A$ be a two-dimensional $H$-orbit including   
in its boundary a one-dimensional $H$-orbit $C$  
included in $\mathcal{O}$. Then, no one-dimensional  
$H'$-orbit included in $\mathcal{O}'$ meets  
$A$ where $\mathcal{O}'$ is the region of $Q$ 
where $H'$ action is defined near $C$.
\end{lem}   
\begin{proof} 
Since $G$ is a lattice of $H$, for every element   
$x$ of $A$, there is a sequence of elements $h_{n}$  
of $G$ for which the sequence $h_{n}.x$ converges to some point   
of $C$ by Lemma \ref{lem:limit}. Assume that some   
one-dimensional $H'$-orbit $B$ meets $A$. Then, some   
iterates $h_{n}.b$, where $b$ belongs to $B \cap A$,  
accumulates to a point in $C$. But $B$ and $C$ correspond to some  
immersed tori in $M$; therefore, such an accumulation is impossible.  
\end{proof}

\begin{lem}\label{lem:limit}  
Let $H''$ be a connected two-dimensional abelian group  
of projective transformations of the projective sphere $\SI^2$.  
Let $G''$ be a lattice of $H''$, and let $A$ be an  
open orbit of $H''$ in the projective plane.  
Let $C$ be a one-dimensional orbit of $H''$  
contained in the boundary of $A$.  
Then, the closure of the $G''$-orbit of any element of 
$A$ contains a point of $C$.  
\end{lem}  
\begin{proof}  
Let $F$ be a compact connected fundamental domain for the action  
of $G''$ on $H''$ by translations.  
Let $a$ be any point of $A$, and $x$ a point of $C$.  
The orbit $F.x$ of $x$ by $F$ is a compact part  
of $C$.   
Let $U$ be an open neighborhood of $F.x$ in the projective plane.  
By continuity of the action of $H''$ on the  
real projective plane, and by compactness of $F$,  
there is an open neighborhood $V$ near $x$  
such that for every element $v$ of $V$, the  
orbit $F.v$ is contained in $U$. The neighborhood
$V$ certainly contains an element of the orbit of $a$ by $H''$.
Since $H''$ is abelian, and since $F$ is a fundamental  
domain, $F.V$, and thus $U$, must contain an element  
of the orbit of $a$ by $G$. Since this is true  
for any open neighborhood $U$ of $F.x$, it follows  
that the orbit $G''.a$ accumulates at least on some point   
of $F.x$.  
\end{proof}  
  
We claim that the decomposition of $Q$ into $H$-orbits   
are preserved under the action of the deck-transformation group.   
A one-dimensional $H$-orbit $l$ in $Q$ does not meet   
with an image $\vth(m)$ of a one dimensional orbit $m$ in $Q$   
transversally for a deck transformation $\vth$  
by Lemma \ref{lem:orbdim1}. A one-dimensional $H$-orbit   
does not meet an image of two-dimensional $H$-orbit   
under a deck transformation by Lemma \ref{lem:orbdim2}.   
Also, if a one-dimensional orbit $l$ meets an image of
$\vth(m)$ nontransversally with $\vth(m)$ not equal to $l$,
then either a transversal intersection 
of a one-dimensional orbit with $\vth(m)$
or the intersection with two-dimensional orbits must occur nearby.
Thus, $H$-orbits map to $H$-orbits under 
deck transformations as there are no zero dimensional   
$H$-orbits in $Q$.  

We can quickly achieve the proof of Theorem \ref{thm:barbot2II} 
when $q(H)$ is one-dimensional.
Each $q(H)$-orbit has a $G$-invariant one-dimensional volume 
form on it as $q(H)$ is a one-dimensional Lie group with 
a left-invariant volume form. Note that the volume form 
depends smoothly on orbits. As $G$ acts on each of 
the subsets of $\tilde M$ corresponding to one-dimensional
orbits in $Q$, these sets cover tori or Klein bottles under
the covering map $\tilde M \ra M$. 
These tori or Klein bottles fiber $M$.
By taking a finite cover of $M$ if necessary,
we can assume that $M$ is homeomorphic to a torus times 
an interval with fibers corresponding to the tori.
Hence, there exists a volume form transverse to the tori
which becomes a transverse volume form of the orbits.
Therefore, the volume form 
extends to a smooth two-dimensional volume form on $Q$. 
Hence, $M$ has a flow invariant three-dimensional volume form. 
By Proposition 4 of Carri\`ere \cite{Carlet}, 
we see that $M$ admits a total cross-section. 

From now on, we assume that $q(H)$ is two-dimensional.
We say that an affine $3$-manifold $M$ {\em decomposes} into   
affine $3$-manifolds $N_1, \dots, N_n$ if each $N_i$ is   
the closure of a component of $M$ with two-sided   
separating totally geodesic surfaces in $M$ removed.   

\begin{prop}\label{prop:decomps}   
Our manifold $M$ or a finite cover of $M$ decomposes into   
compact radiant affine $3$-manifolds $N_i$ each of which   
is affinely isomorphic to a quotient of a domain in   
$\bR^3 -\{O\}$ by an action of $G$. Every piece 
$N_i$ is homeomorphic   
to torus times an interval and has totally geodesic boundary   
which are tori isotopic to $K$ or a finite cover of $K$.  
\end{prop}  
\begin{proof}   
  
The $H$-orbits in $\SI^2$ correspond to $H$-invariant   
submanifolds in $\bR^3 - \{O\}$. Say these sets   
are {\em $H$-invariant sets}\/ in $\bR^3 - \{O\}$.  
We now choose domains in $\bR^3 - \{O\}$ consisting of   
adjacent three- and two-dimensional $H$-invariant sets.  
In cases D, P, C, an {\em $H$-domain} is the union of an $H$-invariant   
open set and $H$-invariant   
two-dimensional sets in the boundary of the set.  
In case U, we choose two antipodal two-dimensional   
$H$-invariant sets so that $\dev(\tilde K)$ is included in one of  
them.  
Their complement is   
the union of two convex $H$-invariant open sets.   
Call these sets {\em distinguished} two-dimensional   
$H$-invariant sets. An {\em $H$-domain} in case U is the union of   
one convex open set and the two two-dimensional $H$-invariant sets  
included in its boundary.  
    
We look at a component of the inverse image in $Q$ of   
one-dimensional orbits under $\dev'$ or distinguished orbits   
in case U. Then they are one-dimensional $H$-orbits mapping   
homeomorphic to the orbits below. A component of the complement   
of these orbits are either two-dimensional orbits or a union of   
one-dimensional orbits in case U.   
  
On $\tilde M$, these orbits correspond to a decomposition   
into $H$-invariant open sets and $H$-invariant totally
geodesic two-dimensional sets mapping homeomorphic 
to their images in $\bR^3$, which are also $H$-invariant.   
  
Recall that $H$-invariant two-dimensional sets map   
(distinguished ones in case U) to imbedded closed surfaces in $M$   
under the covering map. By taking a finite cover of $M$   
if necessary, we assume that these are two-sided tori in $M$.  
Take one, say $A$, of the open $H$-invariant sets in $M$. Then since   
the deck-transformation group acts on $\tilde M$ preserving   
the decomposition, it follows that $A$ union with   
adjacent two-dimensional orbits map to a closed   
submanifold in $M$ bounded by tori. Let us denote   
them by $N_1, \dots, N_h$. Then clearly, $M$ decomposes   
into $N_1, \dots N_n$.  
  
The manifolds $N_i$ are obviously affinely homeomorphic   
to the quotients of three-dimen\-sional closed domains which are   
unions of $H$-invariant sets.   
  
Take the $H$-invariant open set $A$ adjacent   
to $\tilde K$, a two-dimensional $H$-invariant set. Then   
the closure of $A$ is the union of $A$ and adjacent one, two, or   
three two-dimensional $H$-invariant sets. The closure covers   
a compact radiant affine manifold $N_1$,  
i.e, it may be identified with a universal cover $\tilde N_1$  
of $N_1$.  
The deck-transformation group of $M$ acting on $A$   
can be identified with the deck-transformation group of   
$\tilde N_1$. Then since the deck-transformation  
group acts on $A$ nontrivially, it acts on $\tilde K$ also.
(Remember that $G$ is a subgroup of $H$.)  
It follows that $\pi_1(K) \ra \pi_1(N_1)$ is   
an isomorphism. 
By three-manifold topology, we may assume that $N_1$ is   
homeomorphic to a torus times an interval.  
  
Removing $N_1$ from $M$ and taking the closure in $M$,  
we get a new radiant affine $3$-manifold which   
decomposes into $N_2, \dots, N_n$. By induction,   
we see that each $N_i$ is homeomorphic to  
a torus times an interval.
The proof of Proposition \ref{prop:decomps} is complete.  
\end{proof}  

\begin{cor}\label{cor:decomps}  
Let $p:\tilde M \ra M$ be the universal covering.  
Then $p^{-1}(N_i)$ is connected and under $p$ maps   
as a universal covering map onto $N_i$. Each $p^{-1}(N_i)$   
is a three-dimensional $H$-invariant domain, and maps   
homeomorphic to a three-dimensional $H$-invariant set   
union with two adjacent two-dimensional $H$-invariant sets  
under $\dev$.  
\end{cor}  

We denote by $\tilde N_i$ the set $p^{-1}(N_i)$  
for simplicity.  
  
From now on, we assume that $M$ satisfies the conclusion of   
Proposition \ref{prop:decomps}.  We end this section by
showing that $M$ must be a generalized affine suspension
in the cases D, P, U, and C. 
  
In case U, choose the open $H$-invariant set 
$\tilde A$ adjacent to $\tilde K$.   
Then $\tilde A$ is bounded by $\tilde K$ and another orbit  
$L$ so that the angle $\theta$ between $\tilde K$ and $L$ is   
less than or equal to $\pi$.   
  
If $\theta < \pi$, then $L$ is another boundary component  
of $\tilde{M}$, since if not we could have enlarged   
$\tilde{A}$.  
Thus, $N_1 = M$ and $M$ is a generalized affine suspension  
since we can use the $1$-form $dz/z$ for a linear function $z$   
to get a total cross-section where the plane $z = 0$ is disjoint   
from the image under $\dev$.  
  
Clearly, $G$ acts on $\dev(A) \cup \dev(\tilde K) \cup \dev(L)$   
properly discontinuously.  
If $\theta = \pi$, then this set is in the form of a half-space   
with a line in its boundary removed.  
By Lemma \ref{lem:pilune}, this is a contradiction.  
  
In case C, let $A$ be the open $H$-invariant set adjacent   
to $\tilde K$. The closure of $A$ in $\tilde M$ equals   
$A \cup \tilde K \cup L$ for a two-dimensional $H$-invariant   
set $L$. The developing map $\dev$ sends homeomorphic this set   
to a radiant half-space with a line in the boundary   
passing through $O$ removed. Again the following lemma gives   
us a contradiction.  

\begin{lem}\label{lem:pilune} 
Let $N$ be a radiant affine manifold homeomorphic to   
a torus times an interval with totally geodesic boundary.  
Let $(D, j)$ be the development pair of $N$ and   
$j(\pi_1(M))$ is in one of the above four cases D, P, U, or C.  
Let $x, y,$ and $z$ denote the standard coordinate functions   
of $\bR^3$.  
Then a developing map of $N$ can not be as follows\/{\rm :}  
it maps $\tilde M$ into {\rm (}\/but not necessarily onto\/{\rm )} 
a half-space given by $z \geq 0$, and   
the two boundary components of $\tilde M$ respectively   
homeomorphic onto two components of $z = 0$   
with the line given by $y = 0$ or $x = 0$ removed.   
\end{lem}   
\begin{proof}   
As $N$ is homeomorphic to a torus times an interval, 
an ordered choice of two generators of $G'$ induces orientations 
on each leaf torus. Given an orientation on $N$ and 
the boundary orientation on $\partial N$, the generator orientation 
agrees with the boundary orientation at one component of
$\partial N$ and disagrees at the other component. 

Look at $\tilde \partial N$ mapping homeomorphic to the two
components $A_1$ and $A_2$ of the $xy$-plane with a line removed, 
say $x=0$. Then since $G$ acts on $A_1$ and $A_2$, 
the ordered choice of generators induces an orientation 
on each component. As $-\Idd$ is orientation-preserving on
the $xy$-plane, and commutes with $G$, it follows that 
the generator orientation agrees with the boundary orientation 
of the half-space given by $z \geq 0$ or disagrees with it 
at both components $A_1$ and $A_2$. However, as $\dev(\tilde N)$ 
lies in the half-space, the boundary orientations of $A_1$ and 
$A_2$ are from the boundary orientation of $\partial N$, 
a contradiction. 
\end{proof}  
  
In case D, let $A$ and $L$ be as above. The open set $\dev(A)$  
is a cone with three adjacent two-dimensional $H$-invariant sets.  
In its boundary,  
$\dev(\tilde K \cup L)$ includes only two of   
the adjacent two-dimensional $H$-invariant sets.   
Let $m$ be the remaining $H$-invariant set. Then there exists  
a coordinate function $z$ on $\bR^3$ such that $z = 0$  
is a plane through $O$ including $m$.  
As above, $dz/z$ induces a closed $1$-form on $N_1$   
taking positive values under the radial vector field.   
Thus $N_1$ admits a total cross-section by Lemma \ref{lem:form}.  
  
Similarly, we can show that each $N_i$ admits a total   
cross-section to the radial flow.   
We can patch the cross-sections together at the tori, and   
obtain a total cross-section for $M$. The reason is   
that each $N_i$ is affinely homeomorphic to each other   
by maps induced by reflections along two-dimensional   
$H$-invariant sets. Hence, we may assume   
that the homotopy classes of the total cross-sections are   
the same. This shows that $M$ is a generalized affine suspension.   
There is a more detailed description of this case  
in section 3.2 of \cite{barbot1}.  
  
We now study the last but most complicated case P: Here, we are given   
standard coordinate functions $x, y,$ and $z$  
so that our group $H$ takes the form P.  
  
As the $xz$-plane and the $xy$-plane are $H$-invariant,   
the interior of $\dev(\tilde N_1)$ may be given by   
$y > 0$ and $z > 0$ (up to sign changes of coordinate functions).   
The boundary parts 
$\dev(\tilde K)$ and $\dev(L)$ are obtained from intersecting   
the planes given by $y=0$ and $z=0$ with $\dev(\tilde N_1)$.  
As $G$ acts on these sets so that the quotients are tori,   
we may assume without loss of generality   
that $\dev(\tilde K)$ is one of the following form:  
$y = 0, z > 0, x > 0$; $y=0, z > 0, x  <0$; or $y > 0, z = 0$;   
and similarly, $\dev(L)$ is of form: $y > 0, z = 0$;   
$y = 0, z > 0, x > 0$; or $y = 0, z > 0, x < 0$.  
  
As there are exactly two adjacent two-dimensional $H$-orbits,   
we have two cases:  
\begin{itemize}  
\item[(i)] $\dev(\tilde K)$ and $\dev(L)$ are both   
triangles given by $y = 0, z > 0, x > 0$ and  
$y = 0, z > 0, x < 0$ respectively.  
\item[(ii)] If $\dev(\tilde K)$ is the open lune given by   
$y> 0, z = 0$, then $\dev(L)$ must be a triangle   
given by $y =0, z > 0, x >0$ or $y=0, z > 0, x <0$ and vice versa.
\end{itemize}  
  
We define a Euclidean reflection $R$ about $xy$-plane.  
Then the group $F$ generated by $-\Idd$ and $R$ is of order four.  
As elements of $H$ commute with $F$, we see that   
the closure in $\tilde M$ of each three-dimensional   
$H$-invariant set in $\tilde M$ under $\dev$ maps homeomorphic to   
$\dev(A) \cup \dev(\tilde K) \cup \dev(L)$ as in (i) or (ii)   
or an image regions of type (i) or (ii) under an element of $F$.  
  
If every $N_i$ is of type (i), we easily see that   
$\dev(\tilde N_i)$ lies in the set $z > 0$ by induction.  
Thus $\dev(\tilde M)$ lies in the same set.  
Hence, the closed $1$-form $dz /z$ is $G$-invariant, and   
induces a closed $1$-form on $N_1$ which takes a non-zero   
value under radiant vectors. This means that $N_1$ has   
a total cross-section by Lemma \ref{lem:form}. Hence $M$   
is a generalized affine suspension.   
  
From now on we assume that each $N_i$ is of type (ii).  
Thus, $N_i$ and $N_{i+1}$ meet at a torus which   
is a quotient of a triangle or a lune alternatively according   
to $i$.  
  
Suppose that $N_j$ and $N_{j+1}$ meet at a torus which is   
a quotient of a triangle. Then we see that the remaining boundary   
components of $\tilde N_j \cup \tilde N_{j+1}$ map   
homeomorphic under $\dev$ to open lunes which are components   
of a plane with a line removed. This is a contradiction by   
Lemma \ref{lem:pilune}. Therefore $M$ equals $N_1$ or   
the union of only two $N_1$ and $N_2$ meeting a torus which is   
a quotient of a lune. (We remark that this corresponds 
to a generalized affine suspension of $\pi$-annulus of
type C.)

We will show that $N_1$ admits a total cross-section.   
Since $N_2$ can be obtained by the reflection $R$  
commuting with elements of $G$,   
this will show that $M$ has a total cross-section.  
  
We may assume without loss of generality by conjugations
that the connected two-dimensional   
subgroup $H$ of matrices of form $P$ is of form:  
\[\begin{pmatrix}  
e^a & b e^a & 0 \\ 0 & e^a & 0 \\ 0 & 0 & e^c   
\end{pmatrix}\]  
where $a, b, c$ lies in a two-dimensional subspace of   
$\bR^3$ with coordinate functions $a, b,$ and $c$.   
A lattice $L$ in $P$ determines a subspace of $\bR^3$,   
to be denoted by $P(L)$.   
Let the $ac$-plane have the orientation given by   
$\{e_a, e_c\}$ and the $ab$-plane have the orientation   
by $\{e_a, e_b\}$ where $e_a, e_b,$ and $e_c$ are unit   
vectors in the positive parts of the $a$-, $b$-, and $c$-axes respectively.  
  
\begin{lem}\label{lem:positivity}  
Let $L$ be a lattice in a connected   
two-dimensional subgroup of group of matrices of form P.   
Let $U$ be the domain given by $y > 0, z > 0$ union with   
the set $U_{xz}$ given by $x > 0, y = 0, z > 0$ and   
the set $U_{yz}$ given by $y > 0, z = 0$.   
Then $L$ acts on $U$ so that $U/L$   
is a manifold if and only if   
for the projections $p_{ac}:\bR^3 \ra \bR^2$ to   
the $ac$-plane and $p_{ab}:\bR^3 \ra \bR^2$ to  
the $ab$-plane, $g = p_{ac}\circ (p_{ab}|P(L))^{-1}$  
is orientation-reversing.  
\end{lem}  
\begin{proof}   
Suppose that $U/L$ is a manifold. Then $U/L$ is homeomorphic   
to a torus times an interval. Let $L'$ be a lattice in $P(L)$   
corresponding to $L$ by the above description. Then   
a connected two-dimensional group $\tilde H$ of elements of   
the above form with $a, b, c \in P(L)$ acts on $U$ properly and without   
fixed points. We see that $U$ is foliated by $\tilde H$-orbits.  
Thus, $U/L$ is foliated by leaves that are $\tilde H$-orbits   
quotient out by $L$, homeomorphic to tori.   
Using an orientation on $P(L)$,  
each leaf has an induced orientation. We see   
easily that the leaf-space is homeomorphic to   
an interval and the boundary of $U/L$ corresponds  
to the endpoints of the interval.   
  
Consider the map $\mathcal{F}: \bR^3 \ra U^o$ given by   
\[(a, b, c)\mapsto   
\begin{pmatrix} e^a & be^a & 0 \\   
0 & e^a & 0 \\ 0 & 0 & e^c \end{pmatrix}   
\begin{pmatrix} 1 \\ 1 \\ 1 \end{pmatrix}   
= ((b+1)e^a, e^a, e^c)  \]  
which is a homeomorphism.   
Also, define $\mathcal{F}_{ac}: \bR^2_{ac} \ra U_{xz}$, where   
$\bR^2_{ab}$ denotes the $ab$-plane, by   
\[(a,0,c) \mapsto   
\begin{pmatrix} e^a & be^a & 0 \\   
0 & e^a & 0 \\ 0 & 0 & e^c \end{pmatrix}   
\begin{pmatrix} 1 \\ 0 \\ 1 \end{pmatrix}  
 = (e^a, 0 , e^c) \]  
Define $\mathcal{F}_{ab}:\bR^2_{ab} \ra U_{xy}$, where   
$\bR^2_{ab}$ denotes the $xy$-plane, by   
\[(a,b,0) \mapsto   
\begin{pmatrix} e^a & be^a & 0 \\   
0 & e^a & 0 \\ 0 & 0 & e^c \end{pmatrix}   
\begin{pmatrix} 1 \\ 1 \\ 0 \end{pmatrix}  
 = ((b+1)e^a, e^a , 0). \]  
  
Under $\mathcal{F}$, each $\tilde H$-orbit corresponds to  
a translate of $P(L)$.  
We may add $\bR^2_{ac}$ and $\bR^2_{ab}$ to $\bR^3$  
by considering $(a + t, e^{-t}-1, c)$, $t < 0$,   
to converge to $(a, 0, c)$   
as $t \ra -\infty$ and considering $(a, b, t)$, $t < 0$, to  
converge to $(a, b, 0)$ as $t \ra -\infty$.   
Then $\bR^3$ is an open subspace of the completed space $C$.  
By adding $\mathcal{F}_{ac}$ and $\mathcal{F}_{ab}$ to   
$\mathcal{F}$, we obtain a homeomorphism $C \ra U$ with   
appropriate topology on $C$.  
  
In $U^o$, a $\tilde H$-orbit separates $U_{xz}$ and $U_{xy}$   
since they correspond to two boundary components of   
a torus times an interval. In $C$, $P(L)$ must separate $\bR^{ac}$   
and $\bR^{ab}$. Considering $g$ to be a map $\bR^2 \ra \bR^2$,   
as $g(a,b) = (a, b')$ for each $(a, b)$ and some $b'$,   
it follows that $g$ is represented by   
a matrix \[\begin{pmatrix} 1 & l \\ 0 & m \end{pmatrix}\]   
where the plane $P(L)$ is given by $c = la + mb$.   
  
Consider $\bR^3$ as an open upper-hemisphere, i.e.,   
as an affine patch, in the projective sphere $\SI^3$ with   
boundary $\SI^2$. 
Since the point that $(x+t, e^{-t}-1, z)$ converges   
on $\SI^2$ as $t \ra -\infty$ equals the ray through   
$(0,1,0,0)$ and $(x, y, t)$ converges to the ray through   
$(0,0,-1,0)$, for separation to hold, we must have that   
the function $z - lx - my$ takes different signs at these   
two points. This means that $m < 0$. Therefore   
the determinant of the matrix equals $m$ which is negative.  
  
We now prove the converse. Using the notation above,   
we see that $m < 0$. Let $l_1$ be the arc given by   
$\{(t,e^{-t}-1,0)| t \leq 0\}$ and $l_2$ one by $\{(0, 0, s)| s \leq 0\}$.   
Then $l_1 \cup l_2$ is an arc with well-defined endpoints   
at $\bR^2_{ac}$ and $\bR^2_{ab}$ in $C$. The arc $\alpha$ meets   
each of the translates of $P(L)$ at a unique point eventually  
far away. We can easily choose an arc $\alpha'$ which eventually agrees   
with $\alpha$ far away and meets each translate of $P(L)$   
at a unique point. The image of 
$\alpha'$ by $\mathcal{F}$ is an arc in $U$ with endpoints   
in $U_{xz}$ and $U_{xy}$ which meets each $\tilde H$-orbit   
at exactly one point. As $L$ acts on each $\tilde H$-orbit   
to produced a torus, it follows that this condition   
is enough to give us a compact fundamental domain in  
$U$ of the $L$-action. Hence, $U/L$ is a compact manifold  
homeomorphic to a torus times an interval.  
\end{proof}  
  
\begin{rem}   
There are analogous statements in case D also. But   
we omit them here.  
\end{rem}  
  
Let $L'$ be the lattice in $\bR^3$ corresponding to our   
group $G$ by above correspondence. Let $\vth$ be an element of $L'$   
so that $a-c$ and $b$ values are positive on it. We can always choose   
such $\vth$ since $m < 0$ for our group $L'$ by the above   
lemma \ref{lem:positivity}. We may further assume that $\vth$   
is not a power of an element of $L'$. Let $\vth'$ be the   
corresponding element of $G$. The condition implies that   
for the projective automorphism $\vth''$ corresponding to $\vth'$   
acting on $\SI^2$, $\langle \vth''\rangle$ acts properly and freely on   
the subset $U'$ of $\SI^2$ corresponding to $U$ under   
the radial projection as $[0,0,1]$ corresponds to 
a repelling fixed point of $\vth''$ (see Section 1.4 of \cite{cdcr2}).
Hence, $U'/\langle\vth''\rangle$ is a compact real projective annulus
with geodesic boundary (of type $IIb$ in \cite{cdcr2}).
Hence, $U/\langle\vth'\rangle$ is an $\bR$-bundle over this quotient space   
$U'/\langle\vth''\rangle$. We see easily that 
$U/\langle\vth'\rangle$ has a total   
cross-section $A$ with real projective structure isomorphic   
to $U'/\langle \vth''\rangle$. As the deck-transformation group is abelian 
of rank $2$, the group also acts on $U/\langle \vth'\rangle$ properly
discontinuously.
There exists an element $\vpi$ of $G$ inducing an automorphism
$\vpi'$ on $U/\langle\vth'\rangle$ so that $(\vpi^{\prime})^i(A)$ is 
disjoint from $A$ for $i \ne 0$. 
Since $U'/\langle \vth, \vpi\rangle$ is a finite cover of $N_1$ with
a total cross section coming from $A$,
$N_1$ admits a total cross-section.   

\begin{exmp}  
Let $\vth$ be equal to \[\begin{pmatrix} e^2 & e^2 & 0 \\ 0 & e^2 & 0 \\  
0 & 0 & e^{-4} \end{pmatrix} \] and $\vpi$ equal to   
\[\begin{pmatrix} e^3 & -e^3 & 0 \\ 0 & e^3 & 0   
\\ 0 & 0 & e^3 \end{pmatrix}.\] Then the group   
$\langle \vth, \vpi\rangle$ acts on $U$ freely and properly  
discontinuously. The above coefficients are given by   
$l = -1/5$ and $m = -18/5$. This has a total cross-section.  
\end{exmp}  

Finally, suppose that both types (i) or (ii) occur for 
our given manifold $M$, which decomposes into 
$N_1, \dots, N_n$ in a sequential manner.
We note that type (ii) ones occur in adjacent pairs or adjacent to
boundary components of $M$. Type (i) 
submanifolds can occur in sequence or alone adjacent to boundary 
components. Pick $N_j$ which is of type (ii), and choose a generator
$\vth$ as above. We can suppose that $\dev(\tilde N_j)$ is given 
by the region $U$ in Lemma \ref{lem:positivity}.
Then $[0,0,1]$ is a repelling fixed point of 
the projective automorphism $\vth'$ on $\SI^2$ corresponding to $\vth$. 

Recall that $A_i$ is the $H$-invariant domain given by 
$\dev(\tilde N_i)$ for $i=1, \dots, n$, and
the group $F$ generated by a reflection $R$ and $-\Idd$.
Recall also that up to an action of $F$, $A_i$ can be of the form
the open set $y>0$ and $z >0$ union with
\begin{itemize}
\item the two sets in $z>0$ and $y =0$; i.e.,
the open disk $D_1$ given by $x > 0$ and $D_2$ given by $x <0$
respectively.
\item the open disk $D_3$ given by $y > 0$ and $z=0$, and $D_1$.
\item $D_3$, and $D_2$.
\end{itemize}

Clearly $\langle \vth' \rangle$ acts properly on the domains 
in $\SI^2$ corresponding to first and second cases as it 
acts properly on $U'$. We claim that the third case is actually
impossible to occur as $\dev(\tilde N_i)$ for some $i$.
Since $H$ acts properly and freely on $\tilde N_i$, $\tilde N_i$ is 
foliated by $H$-orbits which map to tori foliating $N_i$
under the covering map.
Let $L$ be an $H$-orbit in $A_i^o$ with $A_i$ in the above forms. 
$L/G$ is a torus in $N_i$ separating two boundary components 
of $N_i$ regardless of what form $A_i$ is.
As the first two cases are possible, $L$ separates
$D_1$ and $D_2$, and $L$ separates $D_1$ and $D_3$. 
We can see that $L$ cannot separate $D_2$ and $D_3$ as well
using the $\bZ_2$-intersection theory.
But $L$ corresponds to an $H$-orbit in each $\dev(\tilde N_i)$. 
(As a dual reasoning, one can consider $H$-orbits of a point 
mapped by $F$ lifted to $A_i$s.)

Since every $\dev(\tilde N_i)$ is of first two forms up
to an action of $F$, it follows that $\langle \vth' \rangle$
acts properly on $Q$ which is a union of domains corresponding
to regions of the first two types up to elements of $F$. 

This means as before that $\tilde M/\langle \vth \rangle$ is 
an $\bR$-bundle over a compact surface $Q/\langle \vth' \rangle$ with
radial lines forming the fibers. 
Choosing a section and finding another deck transformation
$\vpi$ which induces an automorphism $\vpi'$ on 
$\tilde M/\langle \vth \rangle$ so that the images of the section 
under $\vpi^{\prime i}$ are all disjoint from each other;
$\tilde M/\langle \vth, \vpi \rangle$ is a radiant affine $3$-manifold
with a total cross-section covering $M$ finitely.
We obtain that $M$ admits a total cross-section.

\end{appendix} 
 
\backmatter 

\bibliographystyle{amsplain} 
\bibliography{rdsv.bib}

\end{document}